\documentclass[10pt,a4paper]{article}
\usepackage[english]{babel}
\usepackage{commath}
\usepackage{graphicx}
\usepackage{tikz}
\usetikzlibrary{shapes.geometric}

\usepackage{xcolor}
\newcommand\bldxi{\boldsymbol\xi}
\newcommand\bldomega{\boldsymbol\omega}
\newcommand\bldx{\boldsymbol x}
\newcommand\bldy{\boldsymbol y}
\newcommand\blda{\boldsymbol a}
\newcommand\bldb{\boldsymbol b}
\newcommand\bldu{\boldsymbol u}
\newcommand\bldv{\boldsymbol v}
\newcommand\bldt{\boldsymbol t}
\newcommand\bldh{\boldsymbol h}
\newcommand\A{a}
\newcommand\B{b}
\newcommand\uu{u}
\newcommand\vv{v}
\newcommand\xx{x}
\newcommand\yy{y}
\newcommand\zz{z}
\newcommand\tee{t}
\newcommand\xxi{\xi}
\newcommand\oomega{\omega}
\newcommand\Ma{A}
\newcommand\Mb{B}
\newcommand\Mc{C}
\newcommand\Md{D}

\DeclareMathOperator{\sgn}{sgn}
\newcommand{\LmA}{A}
\newcommand{\LmB}{B}
\newcommand{\LmC}{C}
\newcommand{\LmD}{D}
\newcommand{\LmE}{E}

\usepackage{enumitem}
\usepackage{float}
\usepackage{amsthm,amssymb,mathrsfs,amsmath}
\usepackage[a4paper,bindingoffset=0.2in,left=0.5in,right=0.5in,top=1in,bottom=1in,footskip=.25in]{geometry}
\usepackage{cite}
\usepackage[colorlinks,citecolor=blue,urlcolor=blue,bookmarks=false,hypertexnames=true]{hyperref}
\hypersetup{citecolor=blue}

\usepackage{array}
\usepackage{MnSymbol,bbding,pifont}

\newtheorem{definition}{Definition}[section]
\newtheorem{theorem}{Theorem}[section]

\title{A mathematical survey on Fourier type integral transform and their offshoots: windowed Fourier transform, wavelet transform and Stockwell transform}
\author{Bivek Gupta$^a$\thanks{bivekgupta040792@gmail.com}, Amit K. Verma$^b$\thanks{Corresponding Author Email : $^b$akverma@iitp.ac.in}\footnote{Amit K Verma dedicates this paper to his supervisor Professor Rajni Kant Pandey}\\\small{\it{$^{a,b}$ Department of Mathematics,}} \\\small{\it{Indian Institute of Technology Patna,}}\\\small{\it{ Bihta, Patna 801103, (BR) India.}}}
\date{\today}

\begin{document}
\maketitle
\begin{abstract}
This comprehensive review paper delves into the intricacies of advanced Fourier type integral transforms and their mathematical properties, with a particular focus on fractional Fourier transform (FrFT), linear canonical transform (LCT), quadratic phase Fourier transform (QPFT), and their associated offshoots: windowed Fourier transform, wavelet transform, and Stockwell transform. In the pursuit of a deeper understanding of these transformations, we explore their convolution properties, shedding light on their capacity to define windowed, wavelet and Stockwell transforms in the realm of Fourier, fractional Fourier and quadratic phase Fourier transforms. This review also expands its purview to the realm of uncertainty principles. Several uncertainty principles, like Heisenberg, logarithmic, local, Rényi uncertainty principles, etc., within the context of fractional Fourier, linear canonical, and quadratic phase Fourier transforms, as well as their derivative offshoots are presented in the paper both for the functions of complex as well as quatenrion valued. In particular, the counterpart of several important inequalities of classical Fourier transform are also presented in details for the quaternion case. This article also reviews that multiresolution analysis that has been developed in the literature so far.
\end{abstract}

We first mention below some important abbreviations that will be used throughout this paper.\\~\\
\textit{List of Abbreviations}\\
FT - Fourier transform\\
WFT - Windowed Fourier transform\\
WT - Wavelet transform\\
ST - Stockwell transform\\
FrFT - Fractional Fourier transform\\
WFrFT - Windowed fractional Fourier transform\\
FrWT - Fractional wavelet transform\\
FrST - Fractional Stockwell transform\\
LCT - Linear canonical transform\\
WLCT - Windowed linear canonical transform\\
LCWT - Linear canonical wavelet transform\\
LCST - Linear canonical Stockwell transform\\
OLCT - Offset linear canonical transform\\
WOLCT - Windowed offset linear canonical transform\\
OLCWT - Offset linear canonical wavelet transform\\
OLCST - Offset linear canonical Stockwell transform\\
QFT - Quaternion Fourier transform\\
QWFT - Quaternion window Fourier transform\\
QWT - Quaternion wavelet transform\\
QFrST - Quaternion fractional Stockwell transform\\
QFrFT - Quaternion fractional Fourier transform\\
QWFrT - Quaternion windowed fractional Fourier transform\\
QFrWT - Quaternion fractional wavelet transform\\
QLCT - Quaternion linear canonical transform\\
QWLCT - Quaternion windowed linear canonical transform\\
QLCWT - Quaternion linear canonical wavelet transform\\
QLCST - Quaternion linear canonical Stockwell transform\\
QQPFT - Quaternion offset linear canonical transform\\
WQQPFT - Windowed quaternion quadratic phase Fourier transform\\
QQPWT - Quaternion quadratic phase wavelet transform\\
QQPST - Quaternion quadratic phase Stockwell transform\\
%QOLCT - Quaternion offset linear canonical transform\\
%QWOLCT - Quaternion windowed offset linear canonical transform\\
%QOLCWT - Quaternion offset linear canonical wavelet transform\\
%QOLCST - Quaternion offset linear canonical Stockwell transform\\
RKHS - Reproducing kernel Hilbert space\\
UP - Uncertainty principle\\
\begin{center}
\includegraphics[scale=0.64]{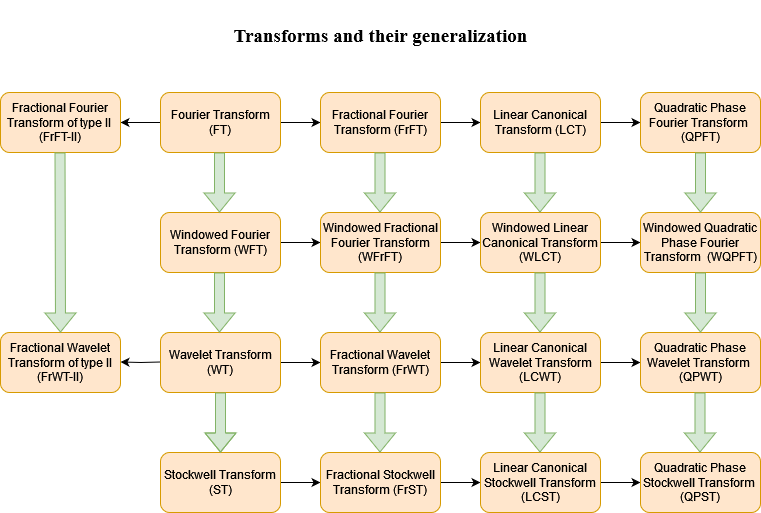}
\end{center}
The Fourier transform (FT) is a spectacular mathematical analysis tool that has long been acknowledged as having the greatest significance to contemporary science and technology developments. There have been numerous mathematical advances in the field of Fourier analysis over time, which have had major effects on the traditional FT. In the ever-expanding domain of signal processing and analysis, the quest for deeper insights, finer resolutions, and more comprehensive interpretations remains unceasing. This review paper embarks on an exploration of advanced signal transformation techniques that transcend the classical boundaries of the Fourier domain. We delve into the profound world of signal transformation, examining the Fractional Fourier Transform (FrFT), Linear Canonical Transform (LCT), Quadratic Phase Fourier Transform (QPFT), and their derived offshoots: the Windowed Fourier Transform, Wavelet Transform, and Stockwell Transform. These transforms offer a sophisticated toolkit that transcends the conventional transforms. The Fractional Fourier Transform, for instance, introduces fractional orders, allowing us to analyze signals with dynamic and complex frequency content, such as chirp signals and those influenced by dispersion. The Linear Canonical Transform unifies diverse mathematical transforms into a unified framework, providing a comprehensive approach to handling signals with intricate spatiotemporal characteristics. It transcends the limitations of individual transforms, making it a valuable asset in fields as diverse as optics, signal processing, and quantum mechanics. The most recent edition, and the generalization of all, i.e., Fourier transform, fractional Fourier transform and the linear canonical transform, the Quadratic Phase Fourier Transform refines our understanding of signals by focusing on the often-overlooked quadratic phase components. This transformation is indispensable for applications like radar, sonar, and image processing, where precise characterization of signals with varying phase characteristics is paramount. 

The introduction of the windowed Fourier transform, wavelet transform, and Stockwell transform arises from the inherent properties of the classical Fourier transform, specifically its use of a global kernel. Because of the global kernel, Fourier transform fails in giving local properties of signal. To overcome these limitations, the Windowed Fourier Transform (WFT) was introduced. It incorporates a window function that acts as a localized filter to select a specific segment of the signal. By applying the Fourier transform to this windowed segment, the WFT allows for the analysis of signal characteristics within this localized region. This approach provides more detailed and localized information about the signal's frequency content. The Wavelet Transform further extends the concept of localized analysis. Wavelets are a family of functions that are adjustable in scale and position. By using wavelets as the analysis function, the Wavelet Transform offers both time and frequency localization. It can capture high and low-frequency components with varying time resolutions, which is particularly useful for understanding transient and non-stationary features in signals. Meanwhile, the Stockwell transform is designed to address signals with rapidly changing frequency components. It employs a windowed approach to capture local characteristics of the signal in the time-frequency domain. The Stockwell Transform, like the WFT and Wavelet Transform, enables a localized view of the signal's frequency content, facilitating the analysis of non-stationary signals. The introduction of the windowed Fourier transform, wavelet transform, and Stockwell transform is driven by the need to go beyond the global view of the classical Fourier transform. These techniques provide localized analysis, offering a more detailed and nuanced understanding of signals, particularly those with complex, non-stationary, or transient characteristics. They are indispensable tools in signal processing, offering the ability to extract specific features from signals in diverse applications, from image analysis to audio processing. Similar, limitation of the fractional Fourier transform, linear canonical transform and the quadratic phase Fourier transform as the Fourier transform in the fractional Fourier domain, linear canonical domain and the quadratic phase Fourier domain, results in the development of their offshoots as windowed Fourier transform, wavelet transform and the Stockwell transform, in the realm of the fractional Fourier transform, linear canonical transform and the quadratic phase Fourier transform. 

Apart from the function of complex-valued ($\mathbb{C}$-valued), researchers are also actively involved in studying these transforms for quaternion-valued ($\mathbb{H}$-valued) functions. Indeed due to non-commutativity of the quaternion multiplication, extension of the complex Fourier transform to quaternion valued functions can be at least of three different types, left sided, right sides and the two sided. The non-commutativity of quaternion multiplication has given rise to distinct variants of quaternion Fourier transform including windowed Fourier transform, wavelet transform, and the Stockwell transform for the quaternion valued functions. These quaternion-based transforms consider the non-commutative properties of quaternion algebra and provide unique insights into signals characterized by both magnitude and phase components. The quaternion Fourier transform extends traditional Fourier analysis to quaternion-valued signals, offering a unified framework for time-frequency analysis, while windowed quaternion Fourier transforms focus on localized time-frequency features within quaternion signals. QWTs and the QSTs enable multi-resolution analysis of quaternion signals, capturing both amplitude and phase information facilitate the study of non-stationary signals in the quaternion domain. These specialized approaches have found applications in diverse fields, particularly in multidimensional and multichannel signal analysis, where the traditional Fourier transforms fail in doing so.

This review paper is intended to discuss the important mathematical properties of the fractional Fourier transform, linear canonical transform, quadratic phase Fourier transform and their associated offshoots: windowed Fourier transform, wavelet transform and Stockwell transform for the functions of complex-valued and and the quaternion-valued. In what follows, we mention the theme of this review paper point wise.
\begin{itemize}
\item This article is devoted to the survey of the convolutions and the corresponding convolution theorem associated with the fractional Fourier transform, linear canonical transform and the quadratic phase Fourier transform both for complex and the quaternion-valued functions.
\item  To survey the convolution based joint time-frequency analysis tools like wavelet transform and the Stockwell  transform in the Fractional Fourier domain, linear canonical domain and the quadratic phase Fourier domain, both for the complex and the quaternion valued functions.
\item Multiresolution analysis (MRA) are the fundamental tools in the construction of a wavelet. Several MRA are developed in the literature in the last few year, associated with the FrWT and LCWT. This paper collects  the MRA developed so far in the literature.
\item Important inequalities like Hausdorff-Young inequality, Pitts inequality are of fundamental importance in Fourier analysis. This review paper is also devoted to survey the counterpart of these inequalities for the FrFT, LCT and the QPFT, both for the function of complex as well as quaternion valued functions. Also, we survey the Lieb's inequality for the v in the content of the FrFT, LCT and the QPFT.
\item Several qualitative and quantitative uncertainty principle for the Fourier transform has been explored in the literature for the functions of complex valued. We survey the Heisenberg, local, logarithmic and R\'{e}nyi (in particular Shannon) uncertainty principles and their variants for the FrFT, LCT and the QPFT for the complex and quaternion valued functions.
\item We also see the boundedness of the WT, FrWT and the LCWT in the various function spaces like Sobolev space, Morrey space, BMO space, Hardy space, generalized Sobolev space,etc.
\end{itemize}
Through this exploration, we aim to provide a comprehensive resource for researchers, engineers, and mathematicians interested in advanced signal processing and analysis. By elucidating the mathematical foundations and properties of these transforms and the associated uncertainty principles, we endeavor to advance the understanding and application of these tools in a wide range of scientific and technological domains. 

Several sections comprise our survey. Each section is devoted to a single integral transform and its corresponding essential results. Before we move to the section-wise survey, we recall some basic definitions and results.
\\
\section{Preliminaries}
Let $\mathbb{R}^n$  be the $n$-dimensional Euclidean space with the euclidean norm $|\cdot|,$ i.e., for $\bldx=(\xx_{1},\xx_{2},\ldots,\xx_{n})\in\mathbb{R}^n$, $$|\bldx|=\sqrt{\sum_{i=1}^n \xx^2_{i}}.$$
We define $|\bldx|_{p}=|\xx_{1}\xx_{2}\ldots \xx_{n}|$ and $\mathbb{R}^n_{0}=\{\bldx\in\mathbb{R}^n:|\bldx|_{p}\neq 0\}.$\\ 
For $\bldx=(\xx_{1},\xx_{2},\ldots,\xx_{n}),$ $\bldy=(\yy_{1},\yy_{2},\ldots,\yy_{n})\in\mathbb{R}^n,$ $\bldx+\bldy=(\xx_{1}+\yy_{1},\xx_{2}+\yy_{2},\ldots,\xx_{n}+\yy_{n}),$ $\bldx\bldy=(\xx_{1}\yy_{1},\xx_{2}\yy_{2},\ldots,\xx_{n}\yy_{n}).$ If further $\bldy\in\mathbb{R}^n_{0}$, then $\frac{\bldx}{\bldy}=(\frac{\xx_{1}}{\yy_{1}},\frac{\xx_{2}}{\yy_{2}},\ldots,\frac{\xx_{n}}{\yy_{n}}).$
\begin{definition}
For $1\leq p<\infty,$ $L^p(\mathbb{R}^n)$ is a Banach space of all complex valued measurable function defined on $\mathbb{R}^n$ such that
$$\int_{\mathbb{R}^n}|f(\bldt)|^pd\bldt<\infty,$$ with norm defined by
$$\|f\|_{L^p(\mathbb{R}^n)}=\left(\int_{\mathbb{R}^n}|f(\bldt)|^p d\bldt \right)^{\frac{1}{p}}.$$
In particular, $L^2(\mathbb{R}^n)$ is a Hilbert space, where the inner product inducing the norm is given by 
$$\langle f,g\rangle_{L^2(\mathbb{R}^n)}=\int_{\mathbb{R}^n}f(\bldt)\overline{g(\bldt)}d\bldt,$$
where $\overline{g(\bldt)}$ is the conjugate of $g(\bldt).$
\end{definition}
\begin{definition}
Let $A$ be an arbitrary set and $X$ be a Hilbert space of complex valued functions defined on $A$ with the inner product $\langle\cdot,\cdot\rangle_{X}.$  Then a complex valued function $K$ defined on $A\times A$ is called a reproducing kernel of X if it satisfies the following condition:\\
\quad For any fixed $q\in A,$ we have $K(\cdot,q)$ is in $X$ and $f(q)=\langle f(\cdot),K(\cdot,q)\rangle_{X}$ for all $f\in X$. 
\end{definition}
\section{Quaternion Algebra and some important notions}
The field of real and complex numbers are respectively denoted by $\mathbb{R}$ and $\mathbb{C}.$ Let 
$$\mathbb{H}=\{r=r_0+ir_1+jr_2+kr_3 :r_0,r_1,r_2,r_3\in\mathbb{R}\},$$ where $i,j$ and $k$ are the imaginary units such that they satisfy the following Hamilton's multiplication rule
$$ij=k=-ji,~jk=i=-kj,~ki=j=-ik,~i^2=j^2=k^2=1.$$

For a quaternion $r=r_0+ir_1+jr_2+kr_3,$ we call $r_0$ the real scalar part of $r,$ and denote it by $Sc(r).$ The scalar part satisfies the following cyclic multiplication symmetry
\begin{align}\label{RQuaternionCyclicMultiplicationSymmetry}
Sc(pqr)=Sc(qrp)=Sc(rpq),~\forall~p,q,r\in\mathbb{H}.
\end{align}

We denote the quaternion conjugate of $r$ as $\bar{r}$ and is defined as
$$\bar{r}=r_0-ir_1-jr_2-kr_3.$$
The quaternion conjugate satisfy the following
\begin{align}\label{RPropertyOfQuaternionConjugate}
\overline{qr}=\bar{r}\bar{q},~\overline{q+r}=\bar{q}+\bar{r},~\bar{\bar{q}}=q,~\forall~q,r\in\mathbb{H}.
\end{align}
The modulus of $r\in\mathbb{H}$ is defined as 
\begin{align}\label{RQuaternionModules}
|r|=\sqrt{r\bar{r}}=\left(\sum_{l=0}^3r_l^2\right)^{\frac{1}{2}},
\end{align}
and it satisfies $|qr|=|q||r|,~\forall~q,r\in\mathbb{H}.$

Each quaternion can be split by
$$q=q_{+}+q_{-},~q_{\pm}=\frac{1}{2}(q\pm iqj).$$
By the real components $r_{0},r_{1},r_{2},r_{3}\in\mathbb{R},$ we have
\begin{align}
r_{\pm}=\left\{r_{0}\pm r_{3}+i(r_{1}\mp r_{2})\right\}\frac{1\pm k}{2}=\frac{1\pm k}{2}\left\{r_{0}\pm r_{3}+j(r_{1}\mp r_{2})\right\},
\end{align}
where $q_{\pm}$ is orthogonal in in the sense that $[r_{+}\overline{r_{-}}]_{0}=0.$ It is to be noted that $r_{+}$ and $r_{-}$ lies in the planes spanned by $\{i-j,1+ij\}$ and $\{i+j,1-ij\}$ respectively. These planes are orthogonal which spans the whole quaternion. This leads to the following modular identity \cite{hitzer2010directional}
\begin{align}\label{RModularIdentityQuaternion}
|r|^2=|r_{+}|^2+|r_{-}|^2.
\end{align}
A quaternion valued function $h$ defined on $\mathbb{R}^n$  can be written as
\begin{align*}
h(\bldx)=h_0(\bldx)+ih_1(\bldx)+jh_2(\bldx)+kh_3(\bldx),~\bldx\in\mathbb{R}^n,
\end{align*}
where $h_0,h_1,h_2$ and $h_3$ are real valued function on $\mathbb{R}^n.$ If $1\leq q<\infty,$ then the $L^q-$norm of $h$ is defined by
\begin{align}\label{RQuaternionLpNorm}
\|h\|_{L^q_\mathbb{H}(\mathbb{R}^n)}
&=\left(\int_{\mathbb{R}^n}|h(\bldx)|^qd\bldx\right)^\frac{1}{q}\notag\\
%&=\left\{\int_{\mathbb{R}^n}\left(|h_0(\bldx)|^2+|%h_1(\bldx)|^2+|h_2(\bldx)|^2+|h_3(\bldx)|%^2\right)^\frac{q}{2}d\bldx\right\}^\frac{1}{q}\notag\\
&=\left\{\int_{\mathbb{R}^n}\left(\sum_{l=0}^3|h_l(\bldx)|^2\right)^\frac{q}{2}d\bldx\right\}^\frac{1}{q}
\end{align}
and $L^q_\mathbb{H}(\mathbb{R}^n)$ is a Banach space of all measurable quaternion valued functions $f$ having finite $L^q-$norm.
$L^\infty_\mathbb{H}(\mathbb{R}^n)$ is the set of all essentially bounded quaternion valued measurable functions with norm 
\begin{align}\label{RQuaternionLinfinityNorm}
\|f\|_{L^\infty_\mathbb{H}(\mathbb{R}^n)}=\mbox{ess~sup}_{\bldx\in\mathbb{R}^n}|f(\bldx)|.
\end{align}
Moreover, the quaternion valued inner product 
\begin{align}\label{RQuaternionInnerProduct}
(f,g)=\int_{\mathbb{R}^n}f(\bldx)\overline{g(\bldx)}d\bldx,
\end{align}
with symmetric real scalar part
\begin{align}\label{RQuaternionScalerProduct}
\langle f,g\rangle
&=\frac{1}{2}[(f,g)+(g,f)]\notag\\
&=\int_{\mathbb{R}^n}Sc\left[f(\bldx)\overline{g(\bldx)}\right]d\bldx\notag\\
&=Sc\left(\int_{\mathbb{R}^n}f(\bldx)\overline{g(\bldx)}d\bldx\right)
\end{align}
turns $L^2_\mathbb{H}(\mathbb{R}^n)$ to a Hilbert space, where the norm in equation \eqref{RQuaternionLpNorm} can be expressed as 
\begin{align}
\|f\|_{L^2_\mathbb{H}(\mathbb{R}^n)}=\sqrt{\langle f,f \rangle}=\sqrt{(f,f)}=\left(\int_{\mathbb{R}^n}|f(\bldx)|^2d\bldx\right)^\frac{1}{2}.
\end{align} 
\begin{definition}
A function $f\in L^2_{\mathbb{H}}(\mathbb{R}^{2d})$  is said to be $\epsilon-$ concentrated on a measurable set $T\subset \mathbb{R}^{2d},$ if 
$$\left(\int_{T^c}|f(\bldx_1,\bldx_2)|^2d\bldx_1d\bldx_2\right)^{\frac{1}{2}}\leq\epsilon\|f\|_{L^2_{\mathbb{H}}(\mathbb{R}^{2d})}$$
\end{definition}
If $0\leq\epsilon\leq\frac{1}{2},$ then most of the energy is concentrated on $T$ and $T$ is indeed and essential support of $f.$ If $\epsilon=0,$ then $T$ is the exact support of $f.$
\\

%%%%%%%%%%%%%%%%%%%%%%%%%%%%%%%%%%%%%%%%%%%%%%%%%%%%%%%%%%%%%%%%%%%%%%%%%%%%%%%%%%%%%%%%%%%%%%%%%%%%%%
%%%%%%%%%%%%%%%%%%%%%%%%%%%%%%%%%%%%%%%%%%%%%%%%%%%%%%%%%%%%%%%%%%%%%%%%%%%%%%%%%%%%%%%%%%%%%%%%%%%%%%
%%%%%%%%%%%%%%%%%%%%%%%%%%%%%%%%%%%%%%%%%%%%%%%%%%%%%%%%%%%%%%%%%%%%%%%%%%%%%%%%%%%%%%%%%%%%%%%%%%%%%%
%%%%%%%%%%%%%%%%%%%%%%%%%%%%%%%%%%%%%%%%%%%%%%%%%%%%%%%%%%%%%%%%%%%%%%%%%%%%%%%%%%%%%%%%%%%%%%%%%%%%%%
%%%%%%%%%%%%%%%%%%%%%%%%%%%%%%%%%%%%%%%%%%%%%%%%%%%%%%%%%%%%%%%%%%%%%%%%%%%%%%%%%%%%%%%%%%%%%%%%%%%%%% 
%%%%%%%%%%%%%%%%%%%%%%%%%%%%%%%%%%%%%%%%%%%%%%%%%%%%%%%%%%%%%%%%%%%%%%%%%%%%%%%%%%%%%%%%%%%%%%%%%%%%%%
\section{Fractional Fourier Transform (FrFT)}
For a suitable complex valued function $f$ of real variable we define the FT by
\begin{align}\label{RFTDefn}
(\mathcal{F}f)(\xxi)=\frac{1}{\sqrt{2\pi}}\int_{\mathbb{R}}f(\tee)e^{-i\tee\xxi}d\tee.
\end{align}
It is known that the operator $\mathcal{F}$ is a homeomorphism on $L^2(\mathbb{R}),$ and has $\{e^{-in\frac{\pi}{2}}\}_{n=0}^{\infty}$ as the eigen values with corresponding eigen functions $\{e^{-\frac{1}{2}\tee^2}H_n(\tee)\}_{n=0}^{\infty},$ where $H_n$ is the Hermite polynomial of degree $n,$ i.e., 
$$H_n(\tee)=(-1)^ne^{\tee^2}\frac{d^n}{d\tee^n}e^{-\tee^2},~n\in\mathbb{N}.$$
Renaming the operator $\mathcal{F}$ in equation \eqref{RFTDefn} as $\mathcal{F}_{\frac{\pi}{2}},$ we see that
\begin{align}\label{REgnvaleqn}
\mathcal{F}_{\frac{\pi}{2}}\left[e^{-\frac{1}{2}\tee^2}H_n(\tee)\right](\xxi)=e^{-in\frac{\pi}{2}}e^{-\frac{1}{2}\xxi^2}H_n(\xxi).
\end{align} 
Namias (\cite{namias1980fractional}), considered the operator $\mathcal{F}_{\alpha},$ given by
\begin{align}\label{REqn3}
\mathcal{F}_{\alpha}\left[e^{-\frac{1}{2}\tee^2}H_n(\tee)\right](\xxi)=e^{-in\alpha}e^{-\frac{1}{2}\xxi^2}H_n(\xxi),~\alpha\in\mathbb{C}
\end{align} 
obtained by replacing $\frac{\pi}{2}$ by $\alpha$ in the eigen value equation \eqref{REgnvaleqn}.
Any function $f\in L^2(\mathbb{R}),$ can be written as 
\begin{align}\label{REqn4}
f(\tee)=\sum_{n=0}^\infty a_ne^{-\frac{\tee^2}{2}}H_n(\tee),
\end{align}
where $a_n$ is obtained, using the orthogonal property of the Hermite polynomials, as 
$$a_n=\frac{1}{2^n n!\sqrt{\pi}}\int_{\mathbb{R}}e^{-\frac{\tee^2}{2}}H_n(\tee)f(\tee)d\tee.$$
Taking $\mathcal{F}_{\alpha}$ on both sided of equation \eqref{REqn4} and using equation \eqref{REqn3}, results in 
\begin{align}\label{REqn5}
(\mathcal{F}_{\alpha}f)(\xxi)=\sum_{n=0}^\infty a_ne^{-in\frac{\pi}{2}}e^{-\frac{\xxi^2}{2}}H_n(\xxi)
\end{align}
The calculation of $\mathcal{F}_{\alpha}$ using \eqref{REqn5} is usually not practical. Namias \cite{namias1980fractional} obtained the following integral representation of the $\mathcal{F}_{\alpha}$ using the Mehler formula of the  integral representation of Hermite polynomials
\begin{align}\label{REqn6}
(\mathcal{F}_{\alpha}f)(\xxi)=\frac{e^{i\left(\frac{\pi}{4}-\frac{\alpha}{2}\right)}}{\sqrt{2\pi\sin{\alpha}}}e^{-\frac{1}{2}i\xxi^2\cot{\alpha}}\int_{\mathbb{R}}e^{\frac{i\tee\xxi}{\sin{\alpha}}-\frac{1}{2}i\tee^2\cot{\alpha}}f(\tee)d\tee
\end{align}
and called it the fractional Fourier transform (FrFT). McBride et al. \cite{mcbride1987namias} pointed out that the $\mathcal{F}_{\alpha}$ in \eqref{REqn5} and $\mathcal{F}_{\alpha}$ in \eqref{REqn6} are not the same for all $\alpha\in\mathbb{C},$ because the former has period $2\pi$ while the later has period $4\pi.$ Thus equation \eqref{REqn6} is not, after all, an integral representation of the original $\mathcal{F}_\alpha$ for all $\alpha\in\mathbb{C}.$ Also problem arises when $\sin\alpha$ becomes zero or negative. The right hand side of \eqref{REqn6} makes sense only when $\alpha\neq n\pi,~n\in\mathbb{Z},$ and we need to specify which branch of $(\sin\alpha)^{\frac{1}{2}}$ is used when $\sin\alpha<0.$ So they considered a new family of FrFT $\mathcal{F}_{\alpha},~\alpha\in\mathbb{R},$ given by
\begin{align}
(\mathcal{F}_{\alpha}f)(\xxi)=\frac{e^{i(\frac{\pi}{4}\tilde{\alpha}-\frac{\alpha}{2})}}{\sqrt{2\pi|\sin\alpha|}}e^{-\frac{i\xxi^2}{2}\cot\alpha}\int_{\mathbb{R}}e^{\frac{i\tee\xxi}{\sin{\alpha}}-\frac{i\tee^2}{2}\cot{\alpha}}f(\tee)d\tee,~\tilde{\alpha}=\sgn \alpha, 0<\alpha<\pi,
\end{align}
\begin{align*}
(\mathcal{F}_{0}f)(\xxi)=f(\xxi)\\
(\mathcal{F}_{\pi}f)(\xxi)=f(-\xxi)\\
(\mathcal{F}_{\alpha+2\pi}f)(\xxi)=(\mathcal{F}_\alpha f)(\xxi),~\forall~\alpha\in\mathbb{R}
\end{align*}
and proved that the transforms $\mathcal{F}_\alpha,~\alpha\in\mathbb{R}$ are homeomorphisms on the Schwartz space $\mathcal{S}(\mathbb{R}), $ along with the index law
$$\mathcal{F}_{\alpha}\mathcal{F}_{\beta}f=\mathcal{F}_{\alpha+\beta}f,~\forall \alpha,\beta\in\mathbb{R}, f\in\mathcal{S}(\mathbb{R}).$$
Kerr \cite{kerr1988namias} study the $L^2(\mathbb{R})$ theory of $\mathcal{F}_\alpha$ by defining 
\begin{align*}
(\mathcal{F}_{\alpha}f)(\xxi)=\frac{e^{i(\frac{\pi}{4}\tilde{\alpha}-\frac{\alpha}{2})}}{\sqrt{2\pi|\sin\alpha|}}e^{-\frac{i\xxi^2}{2}\cot\alpha}\lim_{R\rightarrow \infty}\int_{-R}^{R}e^{\frac{i\tee\xxi}{\sin{\alpha}}-\frac{i\tee^2}{2}\cot{\alpha}}f(\tee)d\tee,~\tilde{\alpha}=sgn \alpha, 0<|\alpha|<\pi
\end{align*} 
and proved that for $\alpha\in\mathbb{R},$ $\mathcal{F}_{\alpha}$ is a homeomorphism on $L^2(\mathbb{R})$ and satisfies
\begin{align*}
\mathcal{F}_{\alpha}\mathcal{F}_{\beta}f=\mathcal{F}_{\alpha+\beta}f,~\forall~\alpha,\beta\in\mathbb{R},~f\in L^2(\mathbb{R})\\
\mathcal{F}_0f=f,~\forall f\in L^2(\mathbb{R})\\
\|\mathcal{F}_{\alpha}f-f\|_{L^2(\mathbb{R})}\rightarrow 0~\mbox{as}~ \alpha\rightarrow 0,~f\in L^2(\mathbb{R}),
\end{align*}
i.e., $\{\mathcal{F}_\alpha\}$ is a strongly continuous unitary group of operators on $L^2(\mathbb{R}).$ Almeida in \cite{almeida1994fractional} gave the kernel of the FrFT in the following form 
\begin{equation}\label{R1994_Almeida-IEEETSP_FrFT_Kernel}
K_{\alpha}(\tee,\xxi)=
\begin{cases}
\sqrt{\frac{1-i\cot}{2\pi}}e^{\frac{i}{2}(\tee^2+\xxi^2)\cot\alpha-i\tee\xxi\csc\alpha},&\mbox{if}~\alpha\neq m\pi\\
\delta(\tee-\xxi),&\mbox{if}~\alpha=2m\pi\\
\delta(\tee+\xxi),&\mbox{if}~\alpha=(2m+1)\pi,m\in\mathbb{Z}.
\end{cases}
\end{equation}
He obtained some important properties of the kernel like $K_{\alpha}(\tee,\xxi)=K_{\alpha}(\xxi,\tee),$ $K_{-\alpha}(\tee,\xxi)=\overline{K_{\alpha}(\tee,\xxi)},$ $\int_{\mathbb{R}}K_{\alpha}(\tee,\xxi')K_{\beta}(\xxi',\xxi)d\xxi'=K_{\alpha+\beta}(\tee,\xxi)$ and $\int_{\mathbb{R}}K_{\alpha}(\tee,\xxi)K_{\alpha}(\tee,\xxi')d\tee=\delta(\xxi-\xxi')$ and also gave the formula for the FrFT of translation and modulation i.e.,  
\begin{align*}
\left(\mathcal{F}_\alpha (f(\tee-\tau))\right)(\xxi)=(\mathcal{F}_\alpha f)(\xxi-\tau\cos\alpha)e^{i\frac{\tau^2}{2}\sin\alpha\cos\alpha-i\xxi\tau\sin\alpha}
\end{align*}
and 
\begin{align*}
\left(\mathcal{F}_\alpha (f(\tee)e^{i\nu\tee})\right)(\xxi)=(\mathcal{F}_\alpha f)(\xxi-\tau\sin\alpha)e^{-i\frac{\tau^2}{2}\sin\alpha\cos\alpha+i\xxi\nu\cos\alpha}.
\end{align*}
Almeida \cite{almeida1997product} obtained the expression of the FrFT of the product and the convolution of two functions in the space $L^1(\mathbb{R})\cap\mathcal{W},$ where $\mathcal{W}$ is the Wigner algebra consisting of functions which are the FT of functions in $L^1(\mathbb{R}).$ The expression of the FrFT of the product and convolution of $f,g\in L^1(\mathbb{R})\cap\mathcal{W},$ are given respectively by the following :
\begin{align}\label{RAlmeidaProductFrFT}
\left(\mathcal{F}_\alpha(fg)\right)(\xxi)=\frac{|\csc\alpha|}{\sqrt{2\pi}}e^{i\frac{\xxi^2}{2}\cot\alpha}\int_{\mathbb{R}}\left(\mathcal{F}_{\alpha}f\right)(\nu)\left(\mathcal{F}_{\frac{\pi}{2}}g\right)[(\xxi-\nu)\csc\alpha]e^{-i\frac{\nu^2}{2}\cot\alpha}d\nu
\end{align}
and 
\begin{align}\label{RAlmeidaConvolutionFrFT}
\left(\mathcal{F}_\alpha(f\ast g)\right)(\xxi)=|\sec\alpha|e^{-i\frac{\xxi^2}{2}\tan\alpha}\int_{\mathbb{R}}\left(\mathcal{F}_{\alpha}f\right)(\nu)g[(\xxi-\nu)\sec\alpha]e^{-i\frac{\nu^2}{2}\tan\alpha}d\nu,
\end{align}
where the convolution $\ast$ is defined by 
\begin{align}\label{RConvolutionFT}
(f\ast g)(\xx)=\int_{\mathbb{R}}f(\tee)g(\xx-\tee)d\tee.
\end{align}
This can equivalently written as 
\begin{align}\label{RConvolutionFTEquivalent}
(f\ast g)(\xx)=\int_{\mathbb{R}}f(\tee)(\tau_\tee g)(\xx)d\tee,
\end{align}
where $\tau_\tee,$ is the time shifted operator, given as $(\tau_\tee g)(\xx)=g(\xx-\tee).$

Several methods can be adopted to extend the integral transform to a class of generalized functions. Important among them is 
\begin{itemize}
\item Embedding or kernel method: The kernel of an integral transform is embedded in a suitable space of test functions. The dual of a test function space is a space of generalized functions. The integral transform of the member of this dual space is defined as the function obtained by taking the action of the member on the kernel of the transform. Precisely,  If $\mathcal{A}$ is a space of test function with $\mathcal{A}^\star$ as its dual, then with the assumption that the kernel of the integral transforms $k(\tee,\xx)\in\mathcal{A},$ the integral transform of $f\in\mathcal{A}^\star$ is given by $F(\xx)=\langle f(\tee),k(\tee,\xx) \rangle.$

\item Adjoint Operator method: This method involves the construction of two test function spaces $\mathcal{A}$ and $\mathcal{B}$ such that the integral transform $\mathcal{I}$ maps $\mathcal{A}$ onto $\mathcal{B}.$ The integral transform $\mathcal{I}f$ of any generalized function $f\in\mathcal{B}^\star$ is defined as an element of $\mathcal{A}^\star,$ given by $\langle \mathcal{I}f,\phi\rangle=\langle f,\mathcal{I}\phi\rangle,~\phi\in\mathcal{A}.$

L. Schwartz first applied this method to extend the FT to tempered distribution, dual of the test function space consisting of $C^\infty$ function of rapid decay. This method can mainly be used when the kernel of the transform is not necessarily an infinitely differential function, as in the case of the Hilbert transform.
\item Algebraic method: This method is different from the one discussed above; it uses the algebraic structure of the integral transform. This method results in the space of generalized functions, the elements of which are called Boehmians. The theory of Boehmians, introduced by J. Mikusin\`ski in \cite{mikusinski1981quotients}, as a space of generalized functions, was motivated by the concept of Boehme's regular operators \cite{boehme1973support}. The algebraic construction of Boehmians is similar to the construction of the field of quotients. The main difference is that Boehmians can be constructed even if the ring has zero divisors.   If the construction is applied to a function space and the multiplication is interpreted as convolution, the construction yields a space of generalized functions. These spaces provide a natural setting for the extension of the transform like FT, Radon transform, Zak transform, etc. \cite{mikusinski1995transforms}.
\end{itemize} 
For other methods like sequence or limit method and indirect method (where the considered transform is expressed in terms of another integral transform which is defined in some space of generalized functions), we refer the reader to \cite{bryvckov1989integral},\cite{zayed2019handbook}. Zayed \cite{i1998fractional} extended the FrFT to different spaces of generalized function. Similar to the method adopted by Zemanian to extend the Laplace transform to generalized functions, Zayed used the embedding method to extend the FrFT to the generalized function with compact support. He also used the algebraic method to construct suitable space of Boehmians. For that, he introduced a new convolution structure \cite{i1998fractional},\cite{zayed1998convolution} 
\begin{align}\label{RConvolutionZayedFrFT}
\left(f\star g\right)(\xx)=\frac{c(\alpha)}{\sqrt{2\pi}} e^{-ia(\alpha)\xx^2}\left(\tilde{f}\ast\tilde{g} \right)(\xx),
\end{align}
where $\ast$ is the classical convolution \eqref{RConvolutionFT}, $\tilde{h}(\xx)=h(\xx)e^{ia(\alpha)\xx^2},$ $a(\alpha)=\frac{\cot\alpha}{2}$ and $c(\alpha)=\sqrt{1-i\cot\alpha}$ 
for the FrFT, as the FrFT does not behave very nicely under the convolution \eqref{RConvolutionFT}. The FrFT of the convolution $\star$ is given by
\begin{align}\label{RConvolutionTheoremZayedFrFT}
\left(\mathcal{F}_\alpha(f\star g)\right)(\xxi)=\left(\mathcal{F}_\alpha f\right)(\xxi)\left(\mathcal{F}_\alpha g\right)(\xxi)e^{-ia(\alpha)\xxi^2}.
\end{align}
This equation \eqref{RConvolutionTheoremZayedFrFT} finds its application in filter design. Like, if we are interested only in the frequency spectrum of the FrFT of $f$ in the region $[\xxi_1,\xxi_2],$ we choose g (the filter impulse response) in such a way that $\mathcal{F}_\alpha g$ is constant over $[\xxi_1,\xxi_2],$ and zero or rapidly decays outside that reason. Thus we have the spectrum of $f$ over $[\xxi_1,\xxi_2],$ passing output of the filter through the chirp multiplier, $e^{ia(\alpha)\xxi^2}.$ This is easier as compared to the one suggested in \cite{almeida1997product}. Taking advantage of the newly defined convolution \eqref{RConvolutionZayedFrFT} and the associated convolution formula \eqref{RConvolutionTheoremZayedFrFT}, he extended the definition of the FrFT to that of the Boehmians. Zayed in \cite{zayed1996relationship} obtained the relation between the FrFT and the FT, which is given by
\begin{align}\label{RRelationb/wFrFT&FT}
e^{-ia(\alpha)\xxi^2}\left(\mathcal{F}_\alpha f\right)(\xxi)=c(\alpha)\left(\mathcal{F}_{\frac{\pi}{2}}\tilde{f}\right)(2a(\alpha)b(\alpha)\xxi),~\tilde{f}(\tee)=e^{ia\tee^2}f(\tee).
\end{align}
Based on relation \eqref{RRelationb/wFrFT&FT}, Zayed obtained the sampling theorem that can be used to reconstruct the signal or its FrFT from their samples at a discrete set of points satisfying the Nyquist rate and the inversion formula for the FrFT from the corresponding results in Fourier analysis \cite{zayed2018advances} by suitable change of variable. The proof obtained is easier as compared to the one obtained by Xia in \cite{xia1996bandlimited}. Xia \cite{xia1996bandlimited} studied the properties of bandlimited and time-limited signals and generalized the fact that a non-zero function cannot be both time-limited as well as bandlimited by showing that a non-zero function $f$ cannot be bandlimited simultaneously with FrFT $\mathcal{F}_\alpha$ for a certain value of $\alpha,$ and FrFT $\mathcal{F}_\beta,$ with $\beta\neq \pm\alpha+m\pi,~\forall~n\in\mathbb{Z}.$ Mendlovic et al. in \cite{mendlovic1993fractional} gave a definition of fractional convolution of functions $f$ and $g$ in terms of functions in fractional Fourier domain as 
\begin{align}\label{RMendlovicFrFTConvolution}
(f\mbox{\EightStarTaper}_\alpha g)(\tee)=\mathcal{F}_{-\alpha}\left(\mathcal{F}_\alpha f(\xxi)\mathcal{F}_\alpha g(\xxi)\right)(\tee),
\end{align}
followed by the definition of the fractional correlation defined as 
\begin{align}\label{RMendlovicFrFTCorrelation}
CORR^\alpha(f,g)=(f\mbox{\EightStarTaper}_\alpha\bar{g}).
\end{align}
Equation \eqref{RMendlovicFrFTConvolution} can be expressed in terms of integral as 
\begin{align*}
(f\mbox{\EightStarTaper}_\alpha g)(\tee)=\frac{|\csc \alpha|}{2\pi}\int_{\mathbb{R}}\int_{\mathbb{R}}\int_{\mathbb{R}}e^{i\frac{\tee'2+\uu^2+\tau^2-\tee^2}{2}\cot\alpha}e^{-i\uu(\tee-\tau-\tee')\csc\alpha}f(\tee')h(\tau)d\tee'd\uu d\tau,
\end{align*}
which is, unfortunately, complicated to reduce into single integral as in the case of ordinary FT. Ozaktas et al. \cite{ozaktas1994convolution} introduced an alternative fractional convolution as 
\begin{align}\label{ROzaktasFrFTConvolution}(f\mbox{\EightStar}_\alpha g)(\tee)&=\left(\mathcal{F}_{-\alpha}\left[\left(\mathcal{F}_\alpha f\right)(\xxi)\ast\left(\mathcal{F}_\alpha g\right)(\xxi)\right]\right)(\tee)\notag\\
&=\int_{\mathbb{R}}f(\tau)g(\tee\cos\alpha-\tau)e^{i\pi \tee^2\cos\alpha\sin\alpha-i2\pi \tee\tau\sin\alpha} d\tau.
\end{align}
This definition of fractional convolution in equation \eqref{ROzaktasFrFTConvolution} lacks the classical convolution theorem, i.e., the FT of convolution of the functions is the product of their corresponding FT. Shi et al. \cite{shi2010multichannel} defined a simple fractional convolution structure given by
\begin{align}\label{RShiFrFTConvolution}
(f \star_\alpha h)(\tee)&=e^{-\frac{i}{2}\tee^2\cot\alpha}\left[\left(f(\cdot)e^{\frac{i}{2}(\cdot)^2\cot\alpha}\right)\ast h(\cdot)\right]\notag\\
&=\int_{\mathbb{R}} h(\tau) f(\tee-\tau) \mathrm{e}^{-i \tau(\tee-\frac{\tau}{2}) \cot \alpha} d\tau,
\end{align}
which satisfies
\begin{align}\label{RShiFrFTConvolutionTheorem}
\mathcal{F}_{\alpha}[(f \star_\alpha h)](\xxi)=\sqrt{2 \pi} (\mathcal{F}_{\alpha}f)(\xxi) (\mathcal{F}h)(\xxi \csc \alpha).
\end{align}
Based on the fractional time shift operator given by $(T_\tau^\alpha f)(\tee)=f(\tee-\tau)e^{-i\tau\left(\tee-\frac{\tau}{2}\right)\cot\alpha},$ which generalizes the classical translation operator $(T_\tau f)(\tee)=f(\tee-\tau),$ the fractional convolution in \eqref{RShiFrFTConvolution} can be written in the form
\begin{align}\label{RShiFrFTConvolutionI}
(f \star_\alpha h)(\tee)=\int_{\mathbb{R}} h(\tau) \left(T_{\tau}^{\alpha} f\right)(\tee)d\tau.
\end{align}
Shi et al. \cite{shi2014generalized} gave a generalized definition of convolution associated with the FrFT `$\star_{\alpha,\beta,\gamma}$' as 
\begin{align}\label{RShiFrFTConvolutionII}
\left(f \star_{\alpha, \beta, \gamma} h\right)(\tee)&=\int_{\mathbb{R}} h(\tau)\left(T_{\tau}^{\beta} f\right)(\tee) \phi_{\alpha, \beta, \gamma}(t, \tau) d \tau\notag\\
&=e^{-\frac{i}{2}\tee^2\cot\alpha}\left[g(\cdot)e^{\frac{i}{2}(\cdot)^2\cot\beta}\star h(\cdot)e^{\frac{i}{2}(\cdot)^2\cot\gamma}\right](\tee),
\end{align}
where $\star$ denotes the classical convolution operator and
\begin{align}\label{RShiFrFTConvolutionTheoremII}
\phi_{\alpha, \beta, \gamma}(\uu,\vv)= \mathrm{e}^{i \frac{\vv^{2}}{2} \cot \gamma+i \frac{\uu^{2}}{2}(\cot \beta-\cot \alpha)}.
\end{align}
The generalized convolution theorem associated to \eqref{RShiFrFTConvolutionII} is given by
\begin{align}
\mathcal{F}_{\alpha}\left[\left(f \star_{\alpha, \beta, \gamma} h\right)(\tee)\right](\uu)=\epsilon_{\alpha, \beta, \gamma}(\uu) (\mathcal{F}_{\beta}f)\left(\frac{\uu \csc \alpha}{\csc \beta}\right) (\mathcal{F}_{\gamma}h)\left(\frac{\uu \csc \alpha}{\csc \gamma}\right),
\end{align}
where
$\epsilon_{\alpha, \beta, \gamma}(\uu)= \frac{A_{\alpha}}{A_{\beta} A \gamma} \mathrm{e}^{i \frac{(\uu)^{2}}{2}}\left(\cot \alpha-\left(\frac{\csc \alpha}{\csc \beta}\right)^{2} \cot \beta-\left(\frac{\csc \alpha}{\csc \gamma}\right)^{2} \cot \gamma\right).$
The generalized convolution theorem given in \eqref{RShiFrFTConvolutionTheoremII} reduces to the one studied in \cite{kraniauskas1998method},\cite{shi2010multichannel},\cite{torres2010fractional},
\cite{zayed1998convolution} for some constraints on $\alpha,\beta$ and $\gamma.$
%%%%%%%%%%%%%%%%%%%%%%%%%%%%%%%%%
%%%%%%%%%%%%%%%%%%%%%%%%%%%%%%%%%%
%%%%%%%%%%%%%%%%%%%%%%%%%%%%%%%%%%%%
%%%%%%%%%%%%%%%%%%%%%%%%%%%%%%%%%
%%%%%%%%%%%%%%%%%%%%%%%%%%%%%%%%%%
%%%%%%%%%%%%%%%%%%%%%%%%%%%%%%%%%%%%
%%%%%%%%%%%%%%%%%%%%%%%%%%%%%%%%%
%%%%%%%%%%%%%%%%%%%%%%%%%%%%%%%%%%
%%%%%%%%%%%%%%%%%%%%%%%%%%%%%%%%%%%%
%%%%%%%%%%%%%%%%%%%%%%%%%%%%%%%%%
%%%%%%%%%%%%%%%%%%%%%%%%%%%%%%%%%%
%%%%%%%%%%%%%%%%%%%%%%%%%%%%%%%%%%%%
The fractional translation operator
$$
\left(\mathrm{T}_{\tau}^{\alpha} f\right)(\tee)=f(\tee-\tau) \mathrm{e}^{-i \tau\left(\tee-\frac{\tau}{2}\right) \cot \alpha}
$$
satisfies $ \left(\mathrm{T}_{\tau}^{\alpha} f\right)(\tee) \stackrel{\mathcal{F}^{\alpha}}{\leftrightarrow} F_{\alpha}(\uu) \mathrm{e}^{-i \uu \tau \csc \alpha}.$
Similarly, we can formulate a fractional frequency-shift operator, denoted by $\mathrm{U}_{\vv}^{\alpha}$, as follows:
$$
\left(\mathrm{U}_{\vv}^{\alpha} F_{\alpha}\right)(\uu)=F_{\alpha}(\uu-\vv) \mathrm{e}^{i \vv\left(\uu-\frac{\vv}{2}\right) \cot \alpha}
$$
where it is easy to verify that $
f(\tee) \mathrm{e}^{i \tee \vv \csc \alpha} \stackrel{\mathcal{F}^{\alpha}}{\leftrightarrow}\left(\mathrm{U}_{\vv}^{\alpha} F_{\alpha}\right)(\uu).$

For $\alpha,\beta\in \mathbb{R}$ such that $\alpha +\beta \nin 2\pi\mathbb{Z},$ let 
$\gamma=\frac{\alpha+\beta}{2},~\gamma=\frac{\alpha-\beta}{2},~\tilde{a}=\frac{\cot\gamma}{2},~\tilde{b}=\frac{\cos\delta}{sin\gamma},~d=\frac{ie^{-i\gamma}}{2\pi\sin\gamma},~\boldsymbol\zeta=
\begin{pmatrix}
\tilde{b} & \tilde{c} \\
-\tilde{c} & \tilde{b}
\end{pmatrix},$
then the FrFT $\mathcal{F}_{\alpha,\beta}$ of $f\in L^1(\mathbb{R}^2)$ is defined by Zayed \cite{zayed2018two} as
\begin{align}\label{R2022_Kamalakkannan-Roopkumar-Zayed_Coupled-FrFT-Definition}
(\mathcal{F}_{\alpha,\beta}f)(\bldu)=\int_{\mathbb{R}^2}f(\bldx)K_{\alpha,\beta}(\bldx,\bldu)d\bldx,
\end{align}
where $K_{\alpha,\beta}(\bldx,\bldu)=de^{-i[\tilde{a}(|\bldx|^2+|\bldu|^2)-\bldx\cdot\boldsymbol\zeta\bldu]}.$ This definition is not the tensor product of two one-dimensional transforms and is based on using the Hermite functions of complex variables \cite{ismail2016analytic} as eigen functions of the transform. Kamalakkannan et al. \cite{kamalakkannan2022extension} gave the rigorous proof to extend the definition to the functions in $f\in L^2(\mathbb{R}^2)$ in the natural way i.e., $\mathcal{F}_{\alpha,\beta}f$ is the $L^2(\mathbb{R}^2)$ limit of $\mathcal{F}_{\alpha,\beta}f_n,$ where $\{f_n\}$ is a sequence in $L^1(\mathbb{R}^2)\cap L^2(\mathbb{R}^2)$ such that $f_n\rightarrow f$ in $L^2(\mathbb{R}^2)$ as $n\rightarrow \infty.$ Moreover, $\mathcal{F}_{\alpha,\beta}$ is a unitary operator on $L^2(\mathbb{R}^2).$ It is easy to observe that when $\alpha=\beta$ and $\alpha=\beta=\frac{\pi}{2},$ then $\mathcal{F}_{\alpha,\beta}$ becomes the tensor product of two one-dimensional FrFT and the classical two-dimensional FT respectively. The coupled FrFT \eqref{R2022_Kamalakkannan-Roopkumar-Zayed_Coupled-FrFT-Definition} can be written in terms of the FT as 
\begin{align}\label{R2022_Kamalakkannan-Roopkumar-Zayed_Coupled-FrFT-FT_Relation}
(\mathcal{F}_{\alpha,\beta}f)(\bldu)=2\pi dE_\gamma^{-1}(\bldu)\widehat{fE_\gamma^{-1}}(-\boldsymbol\zeta\bldu),
\end{align}
where $E_\gamma(\bldx)=e^{i\tilde{a}|\bldx|^2}$ and $E^{-1}_\gamma(\bldx)=e^{i\tilde{a}|\bldx|^2},~\bldx\in \mathbb{R}^2.$ Furthermore, $\mathcal{F}_{\alpha,\beta}$ satisfies  $\mathcal{F}_{-\alpha,-\beta}\mathcal{F}_{\alpha,\beta}f=f,$ for all $f\in L^2(\mathbb{R}^2).$ In particular, if $\mathcal{F}_{\alpha,\beta}f\in L^1(\mathbb{R}^2),$ then $f(\bldx)=\int_{\mathbb{R}^2}(\mathcal{F}_{\alpha,\beta}f)(\bldu)K_{-\alpha,-\beta}(\bldu,\bldx)d\bldu,$ a.e. (\cite{zayed2018two}). Authors in \cite{kamalakkannan2022extension} also obtained the additive property $\mathcal{F}_{\alpha_2,\beta_2}(\mathcal{F}_{\alpha_1,\beta_1}f)=\mathcal{F}_{\alpha_1+\alpha_2,\beta_1+\beta_2}(f),~f\in L^2(\mathbb{R}^2),$ where the parameters satisfy the conditions $\alpha_1+\beta_1,\alpha_2+\beta_2,\alpha_1+\beta_1+\alpha_2+\beta_2\nin2\pi\mathbb{Z}.$
For $f\in L^p(\mathbb{R}^2),~p=1,2$ and $g\in L^1(\mathbb{R}),$ the convolution associated with the coupled FrFT  is given in \cite{zayed2018two} and \cite{kamalakkannan2022extension} as 
\begin{align}\label{R2022_Kamalakkannan-Roopkumar-Zayed_Coupled-FrFT-Convolution}
\left(f\bigotimes g\right)(\bldx)=2\pi dE_\gamma (\bldx)[f(\cdot)E^{-1}_\gamma(\cdot)\star g(\cdot)E^{-1}_\gamma(\cdot)](\bldx),
\end{align}
where $\star$ is the convolution given by $(f\star g)(\bldx)=\frac{1}{2\pi}\int_{\mathbb{R}^2}f(\bldx-\bldy)g(\bldy)d\bldy$ and the associated convolution theorem is given by $\mathcal{F}_{\alpha,\beta}(f\bigotimes g)=E_\gamma\mathcal{F}_{\alpha,\beta}f\mathcal{F}_{\alpha,\beta}g.$ Recently, Topan et al.\cite{topan2023two} obtained the sharp Hausdorff-Young inequality for the coupled FrFT \eqref{R2022_Kamalakkannan-Roopkumar-Zayed_Coupled-FrFT-FT_Relation} which is given for $f\in L^s(\mathbb{R}^2),$ where $\frac{1}{s}+\frac{1}{r}=1$ and $1\leq r\leq 2$ as $\left(\int_{\mathbb{R}^2}\left|\left(\mathcal{F}_{\alpha,\beta}f\right)(\bldu)\right|^s d\bldu \right)^{\frac{1}{s}}\leq \frac{\sin^{\frac{2}{s}}\gamma}{2\pi\sin\gamma}(c(r))^2\|f\|_{L^r(\mathbb{R}^2)},$ where $c(r)=\left(\frac{r^{\frac{1}{r}}}{s^{\frac{1}{s}}}\right)^{\frac{1}{2}}.$ Moving forward to the higher dimensional FrFT, Verma et al.\cite{verma2022note} introduced the FrFT for the function $f\in L^2(\mathbb{R})$ of real order $\alpha,$ using the $n-$copies of the kernel of one dimensional fractional Fourier kernel and is defined as 
\begin{align}\label{R2022_Verma-Gupta_IJWMIP_nD-FrFT-Definition}
\left(\mathfrak{F}_\alpha f\right)(\bldxi)=\int_{\mathbb{R}^n}f(\bldt)K_\alpha(\bldt,\bldxi)d\bldxi,~\bldxi\in\mathbb{R}^n,
\end{align}
where 
\begin{equation}
K_{\alpha}(\bldt,\bldxi)=
\begin{cases}
C_{\alpha}e^{\frac{i}{2}(|\bldt|^2+|\bldxi|^2)\cot\alpha-i\langle\bldt,\bldxi\rangle\csc\alpha},&\mbox{if}~\alpha\neq m\pi\\
\delta(\bldt-\bldxi),&\mbox{if}~\alpha=2m\pi\\
\delta(\bldt+\bldxi),&\mbox{if}~\alpha=(2m+1)\pi,m\in\mathbb{Z},
\end{cases}
\end{equation}
$\delta$ is a Dirac delta function in $n$-dimensions and $C_{\alpha}=\left(\frac{1-i\cot\alpha}{2\pi}\right)^{\frac{n}{2}}.$ It is to be noted that the definition \eqref{R2022_Verma-Gupta_IJWMIP_nD-FrFT-Definition} reduces to the classical FT in $n-$dimension if $\alpha=\frac{\pi}{2}.$ Analogous to the classical definition of the FT the definition of the FrFT \eqref{R2022_Verma-Gupta_IJWMIP_nD-FrFT-Definition} can be interpreted as $\mathfrak{F}_\alpha f=\lim_{n\rightarrow\infty}\mathfrak{F}_\alpha f_n$ in $L^2(\mathbb{R}^n),$ where $\{f_n\}_n$ is a sequence in $L^1(\mathbb{R}^n)\cap L^2(\mathbb{R}^n)$ such that $f_n\rightarrow f$ in $L^2(\mathbb{R}^n).$ The existence and uniqueness of the $\lim_{n\rightarrow\infty}\mathfrak{F}_\alpha f_n$ follows from the identity $\|f\|_{L^2(\mathbb{R}^n)}=\|\mathfrak{F}_\alpha f\|_{L^2(\mathbb{R}^n)},~\forall~f\in L^1(\mathbb{R}^n)\cap L^2(\mathbb{R}^n).$ Authors in \cite{verma2022note} proved the Parseval's formula $\langle f,g\rangle_{L^2(\mathbb{R}^n)}=\langle \mathfrak{F}_\alpha f,\mathfrak{F}_\alpha g\rangle_{L^2(\mathbb{R}^n)}$ followed by the additive property $\mathfrak{F}_\alpha(\mathfrak{F}_\beta f)=\mathfrak{F}_{\alpha+\beta}f,~f\in L^2(\mathbb{R}^n).$ Based on the convolution $\star_\alpha$ defined by $(f\star_\alpha g)(\bldt)=\int_{\mathbb{R}^n}e^{-\frac{i}{2}(|\bldt|^2-|\bldxi|^2)\cot\alpha}f(\bldy)g(\bldt-\bldy)d\bldy,$ provided the integral is defined, author obtained the convolution theorem for the transform \eqref{R2022_Verma-Gupta_IJWMIP_nD-FrFT-Definition} given as $\left(\mathfrak{F}_\alpha(f\star_\alpha g)\right)(\bldxi)=\frac{1}{C_\alpha}e^{-\frac{i}{2}|\bldxi|^2\cot\alpha}(\mathfrak{F}_\alpha f)(\bldxi)\left(\mathfrak{F}_\alpha\{e^{-\frac{i}{2}|\cdot|^2\cot\alpha}g(\cdot)\}\right)(\bldxi),$ for $f\in L^p(\mathbb{R}^n),~p=1,2$ and $g\in L^1(\mathbb{R}^n).$ Kamalakkannan et al. \cite{kamalakkannan2020multidimensional} studied the properties of the FrFT defined for the function $f\in L^1(\mathbb{R}^n)$ with respect to the parameter $\boldsymbol \alpha=(\alpha_1,\alpha_2,\cdots,\alpha_n)\in (-\pi,\pi)^n$ as 
\begin{align}\label{R2020-Kamalakkanan-Roopkumar_ITSF_nD-FrFT_nalpha-Definition}
\left(\mathfrak{F}_{\boldsymbol\alpha} f\right)(\bldxi)=\int_{\mathbb{R}^n}f(\bldt)K_{\boldsymbol\alpha}(\bldt,\bldxi)d\bldxi,~\bldxi\in\mathbb{R}^n,
\end{align}
where $K_{\boldsymbol\alpha}(\bldt,\bldxi)=\Pi_{l=1}^n K_{\alpha_l}(\tee_l,\xxi_l),$ such that 
\begin{equation}
K_{\alpha_l}(\tee_l,\xxi_l)=
\begin{cases}
\frac{c_{\alpha_l}}{\sqrt{2\pi}}e^{i\{a(\alpha_l)[\tee_l^2+\xxi_l^2-2b(\alpha_l)\tee_k\xxi_k]\}},&\mbox{if}~\alpha\neq m\pi\\
\delta(\tee_l-\xxi_l),&\mbox{if}~\alpha=2m\pi\\
\delta(\tee_l+\xxi_l),&\mbox{if}~\alpha=(2m+1)\pi,m\in\mathbb{Z},
\end{cases}
\end{equation}
where $\bldt=(\tee_1,\tee_2,\cdot,\tee_n),$ $a(\alpha_l)=\cot\alpha_l/2,$ $b(\alpha_l)=\sec\alpha_l,$ $c(\alpha_l)=\sqrt{1-i\cot\alpha_l},$ $l=1,2,\cdots,n.$ With the assumption that $\lambda\in \mathbb{R}_{0},$ authors defined $e_{\boldsymbol\alpha,\lambda}(\bldt)=e^{i\lambda\sum_{l=1}^na(\alpha_l)\tee_l^2},~\bldt\in\mathbb{R}^n,$ based on which the definition can be rewritten as
$$\left(\mathfrak{F}_{\boldsymbol\alpha} f\right)(\bldxi)=c(\boldsymbol\alpha)e_{\boldsymbol\alpha,1}(\bldxi)(e_{\boldsymbol\alpha,1}f\hat{)}(\xxi_1\csc\alpha_1,\xxi_2\csc\alpha_2,\cdot,\xxi_n\csc\alpha_n),$$
where $c(\boldsymbol\alpha)=c(\alpha_1)c(\alpha_2)\cdots c(\alpha_n)$ and $(e_{\boldsymbol\alpha,1}f\hat{)}$ denotes the FT of $e_{\boldsymbol\alpha,1}f.$ Also the kernel $K_{\boldsymbol\alpha}(\bldt,\bldxi)$ can be written in the form $K_{\boldsymbol\alpha}(\bldt,\bldxi)=\frac{c(\boldsymbol\alpha)}{(\sqrt{2\pi})^n}e_{\boldsymbol\alpha,1}(\bldt)e_{\boldsymbol\alpha,1}(\bldxi)e^{i\sum_{l=1}^na(\alpha_l)\tee_l^2}.$ It is to be noted that the definition \eqref{R2020-Kamalakkanan-Roopkumar_ITSF_nD-FrFT_nalpha-Definition} reduces to the classical FT in $n-$dimension if $\alpha_l=\frac{\pi}{2},~l=1,2,\cdots,n.$ Analogous to the classical definition of the FT the definition of the FrFT \eqref{R2020-Kamalakkanan-Roopkumar_ITSF_nD-FrFT_nalpha-Definition} can be interpreted as $\mathfrak{F}_{\boldsymbol\alpha} f=\lim_{n\rightarrow\infty}\mathfrak{F}_{\boldsymbol\alpha} f_n$ in $L^2(\mathbb{R}^n),$ where $\{f_n\}_n$ is a sequence in $L^1(\mathbb{R}^n)\cap L^2(\mathbb{R}^n)$ such that $f_n\rightarrow f$ in $L^2(\mathbb{R}^n).$ The existence and uniqueness of the $\lim_{n\rightarrow\infty}\mathfrak{F}_{\boldsymbol\alpha} f_n$ follows from the identity $\|f\|_{L^2(\mathbb{R}^n)}=\|\mathfrak{F}_{\boldsymbol\alpha} f\|_{L^2(\mathbb{R}^n)},~\forall~f\in L^1(\mathbb{R}^n)\cap L^2(\mathbb{R}^n).$ Authors in \cite{kamalakkannan2020multidimensional} proved the Parseval's formula $\langle f,g\rangle_{L^2(\mathbb{R}^n)}=\langle \mathfrak{F}_{\boldsymbol\alpha} f,\mathfrak{F}_{\boldsymbol\alpha} g\rangle_{L^2(\mathbb{R}^n)}$ followed by the additive property $\mathfrak{F}_{\boldsymbol\alpha}(\mathfrak{F}_{\boldsymbol\beta} f)=\mathfrak{F}_{\boldsymbol\alpha+\boldsymbol\beta}f,~f\in L^2(\mathbb{R}^n).$ They also introduced a generalized fractional convolution for $f\in L^2(\mathbb{R}^n),~p=1,2,$ $g\in L^1(\mathbb{R}^n)$ given by
\begin{align}\label{R2020_Kamalakkannan-Roopkumar_ITSF_nD-FrFT_Convolution}
f\star_{\boldsymbol\alpha,\lambda}g=c(\boldsymbol\alpha_\lambda)e_{\boldsymbol\alpha,-\lambda^2}\left[(e_{\boldsymbol\alpha,\lambda^2})\star(e_{\boldsymbol\alpha,\lambda^2}g)\right],
\end{align}
where $\lambda\in\mathbb{R}_{0},$ $\boldsymbol\alpha_\lambda=\left(\cot^{-1}(\lambda^2\cot\alpha_1),\cot^{-1}(\lambda^2\cot\alpha_2),\cdots,\cot^{-1}(\lambda^2\cot\alpha_n)\right)$ and $\star$ is the convolution given by $(f\star g)(\bldx)=\frac{1}{(\sqrt{2\pi})^n}\int_{\mathbb{R}^n}f(\bldx-\bldy)g(\bldy)d\bldy$ and obtained the associated convolution theorem for the FrFT \eqref{R2020-Kamalakkanan-Roopkumar_ITSF_nD-FrFT_nalpha-Definition} as  
\begin{align}\label{R2020-Kamalakkanan-Roopkumar_ITSF_nD-FrFT_nalpha-Convolution-Theorem}
\mathfrak{F}_{\boldsymbol\alpha_\lambda}\left(f\star_{\boldsymbol\alpha,\lambda}g\right)=e_{\boldsymbol\alpha,-\lambda^2}\left(\mathfrak{F}_{\boldsymbol\alpha_\lambda}f\right)\left(\mathfrak{F}_{\boldsymbol\alpha_\lambda}g\right).
\end{align}
It is to be noted that for $\lambda=1,~\alpha_l=\frac{\pi}{2},~l=1,2,\cdots,n$ the convolution $\star_{\boldsymbol\alpha,\lambda}$ reduces to the classical convolution and when $n=1,~\lambda=1,$ the convolution \eqref{R2020_Kamalakkannan-Roopkumar_ITSF_nD-FrFT_Convolution} reduces to the fractional convolution given in \eqref{RConvolutionZayedFrFT}. Moreover, for $f,f_1,f_2\in L^p(\mathbb{R}^n),~p=1,2$ $g,g_1,g_2\in L^1(\mathbb{R}),$ $\gamma\in\mathbb{C}$ and $s\in \mathbb{R}_{0},$ the convolution $\star_{\boldsymbol\alpha,\lambda}$ enjoys the properties $(f_1+f_2)\star_{\boldsymbol\alpha,\lambda}g=(f_1\star_{\boldsymbol\alpha,\lambda}g)+(f_2\star_{\boldsymbol\alpha,\lambda}g),$
$f\star_{\boldsymbol\alpha,\lambda}(g_1\star_{\boldsymbol\alpha,\lambda}g_2)=(f\star_{\boldsymbol\alpha,\lambda}g_1)\star_{\boldsymbol\alpha,\lambda}g_2,$ $\gamma(f\star_{\boldsymbol\alpha,\lambda}g)=(\gamma f)\star_{\boldsymbol\alpha,\lambda}g=f\star_{\boldsymbol\alpha,\lambda}(\gamma g),$ $f\star_{\boldsymbol\alpha,\lambda}g=g\star_{\boldsymbol\alpha,\lambda}f$ and also $\tau_{\bldt}(f\star_{\boldsymbol\alpha,\lambda}g)=(\tau_{\bldt}f)\star_{\boldsymbol\alpha,\lambda}(E_{\boldsymbol\alpha,\lambda^2\bldt}g)=(E_{\boldsymbol\alpha,\lambda^2\bldt}f)\star_{\boldsymbol\alpha,\lambda}(\tau_{\bldt}g),$ where $E_{\boldsymbol\alpha,\lambda^2\bldt}(\boldsymbol x)=e^{2i\lambda^2\sum_{l=1}^na(\alpha_l)\tee_l x_l}$ and $M_s(f\star_{\boldsymbol\alpha,\lambda s}g)=|s|^n\frac{c(\boldsymbol\alpha_\lambda)}{c(\boldsymbol\alpha_\lambda s)}(M_sf)\star_{\boldsymbol\alpha,\lambda s}(M_sg),$ where $(M_sf)(\bldx)=f(s\bldx).$ As a consequence of the convolution result \eqref{R2020-Kamalakkanan-Roopkumar_ITSF_nD-FrFT_nalpha-Convolution-Theorem}, author also obtained the product theorem $\mathfrak{F}_{\boldsymbol\alpha_\lambda}f\star_{-\boldsymbol\alpha,\lambda}\mathfrak{F}_{\boldsymbol\alpha_\lambda}g=\mathfrak{F}_{\boldsymbol\alpha_\lambda}(e_{\boldsymbol\alpha,\lambda^2}fg)$ for $f,g\in L^2(\mathbb{R}^n)$ such that $\mathfrak{F}_{\boldsymbol\alpha_\lambda}g\in L^1(\mathbb{R}^n).$\\\\
%%%%%%%%%%%%%%%%%%%%%%%%%%%%%%%%%
%%%%%%%%%%%%%%%%%%%%%%%%%%%%%%%%%%
%%%%%%%%%%%%%%%%%%%%%%%%%%%%%%%%%%%%
%%%%%%%%%%%%%%%%%%%%%%%%%%%%%%%%%
%%%%%%%%%%%%%%%%%%%%%%%%%%%%%%%%%%
%%%%%%%%%%%%%%%%%%%%%%%%%%%%%%%%%%%%
%%%%%%%%%%%%%%%%%%%%%%%%%%%%%%%%%
%%%%%%%%%%%%%%%%%%%%%%%%%%%%%%%%%%
%%%%%%%%%%%%%%%%%%%%%%%%%%%%%%%%%%%%
%%%%%%%%%%%%%%%%%%%%%%%%%%%%%%%%%
%%%%%%%%%%%%%%%%%%%%%%%%%%%%%%%%%%
%%%%%%%%%%%%%%%%%%%%%%%%%%%%%%%%%%%%
\noindent\textbf{Uncertainty Principles:} The UP in the time-frequency plane plays a vital role in signal processing. It states that for a given function $f\in L^2(\mathbb{R})$ such that $\|f\|_{L^2(\mathbb{R})}=1,$ the product of its spread in the time domain and frequency domain is bounded below and is given by
\begin{align}\label{RUPClassicalFT}
\left(\int_{\mathbb{R}}(\tee-\tee_0)^2|f(\tee)|^2d\tee\right)\left(\int_{\mathbb{R}}(\xxi-\xxi_0)^2|(\mathcal{F}f)(\xxi)|^2d\xxi\right)\geq \frac{1}{4},
\end{align}
where $\tee_0=\int_{\mathbb{R}}\tee|f(\tee)|^2d\tee$ and $\xxi_0=\int_{\mathbb{R}}\xxi|(\mathcal{F}f)(\xxi)|^2d\xxi.$ It has been seen in the \cite{aytur1995non}, \cite{ozaktas1995fractional} that the lower bound of the product of the spread in the time domain and the FrFT domain is given by $\frac{\sin^2{\alpha}}{4}.$ This can also be obtained using the UP for the FT and the relation between FrFT and the FT given by \eqref{RRelationb/wFrFT&FT}. Shinde et al. \cite{shinde2001uncertainty} gave the lower bound of the product of the spread of real unit energy signal $f$ in two FrFT domains as
\begin{align}\label{RShindeUP2FrFTDomain}
\Delta \xxi_{\alpha,f}^{2} \Delta \xxi_ {\beta,f}^{2} \geq\left[\Delta \tee_f^{2}\cos\alpha \cos \beta\right.&\left.+\frac{\sin \alpha \sin \beta}{4 \Delta \tee_f^{2}}\right]^{2}+\frac{\sin^{2}(\alpha-\beta)}{4}\left(=L(\Delta \tee_f,\alpha,\beta)~say\right),
\end{align}
where  $\Delta \xxi_{\alpha,f}^{2}=\int_{\mathbb{R}}(\xxi-\xxi_{f,\alpha 0})^2|(\mathcal{F}^\alpha f)(\xxi)|^2d\xxi,$ $\Delta \xxi_{\beta,f}^{2}=\int_{\mathbb{R}}(\xxi-\xxi_{f,\beta 0})^2|(\mathcal{F}^\beta f)(\xxi)|^2d\xxi$ and $\Delta \tee_f^2=\int_{\mathbb{R}}(\tee-\tee_{f,0})^2|f(\tee)|^2d\tee$ with $\xxi_{f,\alpha 0}=\int_{\mathbb{R}}\xxi|(\mathcal{F}^\alpha f)(\tee)|^2d\tee,$ $\xxi_{f,\beta 0}=\int_{\mathbb{R}}\xxi|(\mathcal{F}^\beta f)(\tee)|^2d\tee$ and $\tee_{f,0}=\int_{\mathbb{R}}\tee|f(\tee)|^2d\tee$ respectively. The equality in \eqref{RShindeUP2FrFTDomain} is achieved when
$$
f(\tee)=\left(\frac{1}{\pi \sigma^{2}}\right)^{\frac{1}{2}} e^{-\frac{\tee^{2}}{2 \sigma^{2}}},
$$
where $\sigma$ is an arbitrary real constant. The lower bound given in \eqref{RShindeUP2FrFTDomain} is tighter than the one given in \cite{ozaktas1995fractional}. As it is known that the shifting and scaling, given by $(\tau_\tau f)(\tee)=f(\tee-\tau)$ and $f_c(\tee)=c^{\frac{1}{2}}f(c\tee)$ respectively does not change the inequality \eqref{RUPClassicalFT}. Shinde et al.\cite{shinde2001uncertainty} also studied the effect of the shifting and scaling of the signal in the UP \eqref{RShindeUP2FrFTDomain} and obtained that the shifting has no effect on the inequality, whereas the scaling of the signal has the following form.
\begin{align}\label{RShindeUPDilationof Signal}
\Delta \xxi_{\alpha,f_c}^{2} \Delta \xxi_ {\beta,f_c}^{2}\geq \frac{c^4\sin^2\alpha\sin^2\beta}{\sin^2\theta\sin^2\gamma} L(\Delta \tee_{f_c},\theta,\gamma),
\end{align}
where $\tan\theta=c^2\tan\alpha$ and $\tan\gamma=c^2\tan\beta.$
%\cite{guanlei2009logarithmic}
Guanlei et al.\cite{guanlei2009generalized} obtained the Shannon UP in two FrFT domains.
\begin{align}\label{RGuanleiShannonEntropyUP}
E\left\{\left|(\mathcal F ^{\alpha}f)(\xxi)\right|^{2}\right\} -E\left\{\left|(\mathcal F ^{\beta}f)(\xxi)\right|^{2}\right\} \geq \ln (\pi e|\sin (\alpha-\beta)|),
\end{align}
where $E(\rho(\xx))=-\int_{\mathbb{R}}\rho(\xx)\ln\rho(\xx)d\xx$ is the Shannon Entropy for the probability density function $\rho$. If $\alpha\rightarrow 0$ and $\beta\rightarrow \frac{\pi}{2},$ \eqref{RGuanleiShannonEntropyUP} reduces to the traditional case. They also obtained the Renyi UP in two FrFT domains.
\begin{align}\label{RGuanleiRenyiEntropyUP}
H_\theta\{|(\mathcal{F}^\alpha f)(\xxi)|^2\}+H_\lambda\{|(\mathcal{F}^\alpha f)(\xxi)|^2\}\geq \frac{\ln(\frac{\theta}{\pi})}{2(\theta-1)}+\frac{\ln (\frac{\lambda}{ \pi})}{2(\lambda-1)}+\ln |\sin (\alpha-\beta)|,
\end{align}
where $H_\mu(\rho(\xx))=\frac{1}{1-\mu}\ln\int_{\mathbb{R}}|\rho(\xx)|^{\mu}d\xx$ is the Renyi Entropy \cite{dembo1991information} for the probability density function $\rho.$ Letting $\mu\rightarrow 1,$ $H_{\mu}(\rho(\xx))$ reduces to $E(\rho(\xx)).$ If $\alpha\rightarrow 0$ and $\beta\rightarrow \frac{\pi}{2},$ \eqref{RGuanleiRenyiEntropyUP} reduces to the traditional case. If $\sin(\alpha-\beta)=\pm(\frac{\theta}{\pi})^{-\frac{1}{2(\theta-1)}}(\frac{\lambda}{\pi})^{-\frac{1}{2(\lambda-1)}},$ then the bound will be zero, i.e., the bound of the Renyi UP in two FrFT domains can be zero, which also provides the possible better selection in FrFT domain for the information processing.
\\
\\
For the function $f\in L^2(\mathbb{R}^2)$ and $2\leq k<3,$ Shah et al. \cite{shah2022class} obtained the following Heisenberg-type UP for the coupled FrFT
\begin{align}\label{R2022_Shah-Lone-Nisar-Abdeljawad_JoIA_CFrFT_Heisenberg-UP}\left(\int_{\mathbb{R}^2}|\bldx|^k\left|f(\bldx)\right|^2d\bldx\right)^{\frac{1}{2}}\left(\int_{\mathbb{R}^2}|\bldu|^k\left|(\mathcal{F}_{\alpha,\beta}f)(\bldu)\right|^2d\bldu\right)^{\frac{1}{2}}\geq \frac{\sin^{2k}\gamma}{8\pi^k}\left[\frac{\Gamma\left(\frac{k}{4}\right)}{\Gamma\left(\frac{4-k}{4}\right)}\right]^2\|f\|^2_{L^2(\mathbb{R}^2)}.
\end{align}
For $k=2,$ equation \eqref{R2022_Shah-Lone-Nisar-Abdeljawad_JoIA_CFrFT_Heisenberg-UP} reduces to the classical Heisenberg UP for the coupled FrFT 
%as
%\textcolor{purple}{\begin{align}\label{R2022_Kamalakkannan-Roopkumar_CFrFT-Heisenberg-UP}
%\Delta\bldx^2\Delta\bldu^2_{\alpha,\beta}\geq \sin^2\gamma\|f\|^4_2,
%\end{align}
%the inequality is not reducing as a particular case of the above equation} where $\Delta\bldx^2=\frac{1}{2\pi}\int_{\mathbb{R}^2}|\bldx|^2|f(\bldx)|^2d\bldx$ and $\Delta\bldu^2_{\alpha,\beta}=\frac{1}{2\pi}\int_{\mathbb{R}^2}|\bldu|^2|(\mathcal{F}_{\alpha,\beta}f)(\bldu)|^2d\bldu,$ 
%which was 
obtained by Kamalakkannan et al. in \cite{kamalakkannan2022two} using the relation \eqref{R2022_Kamalakkannan-Roopkumar-Zayed_Coupled-FrFT-FT_Relation} and the Heisenberg UP for the FT. Moreover, for $k=2$ and $\alpha=\beta$ equation \eqref{R2022_Shah-Lone-Nisar-Abdeljawad_JoIA_CFrFT_Heisenberg-UP} reduces to the Heisenberg UP for the FrFT. Using the logarithmic UP for the FT for the function $f\in \mathcal{S}(\mathbb{R}^2)$ and the relation \eqref{R2022_Kamalakkannan-Roopkumar-Zayed_Coupled-FrFT-FT_Relation} authors obtained
\begin{align}
4\pi^2\int_{\mathbb{R}^2}\ln|\bldu|\left|\left(\mathcal{F}_{\alpha,\beta}f\right)(\bldu)\right|^2d\bldu+\int_{\mathbb{R}^2}\ln|\bldx|\left|f(\bldx)\right|^2d\bldx\geq \left(\frac{\Gamma'(\left(\frac{1}{2}\right))}{\Gamma\left(\frac{1}{2}\right)}-\ln\pi+8\pi^2\ln|\sin\gamma|\right)\|f\|^2_{L^2(\mathbb{R}^2)}
\end{align}
called the logarithmic UP for the coupled FrFT, which reduces to the corresponding inequality for the FrFT when $\alpha=\beta.$ Authors in \cite{shah2022class} also obtained the local uncertainty inequality for the coupled FrFT as 
\begin{align}
\int_{\mathbb{R}^2}|\bldx|^{2k}\left|f(\bldx)\right|^2d\bldx\geq \frac{4\pi^2}{C_k|E|^{2k}\sin^2\gamma}\int_{E}\left|\left(\mathcal{F}_{\alpha,\beta}f\right)(\bldu)\right|^2d\bldu,
\end{align}
where $E\subset \mathbb{R}^2$ is a set of finite measure, $0<k<1$ and $C_k$ is some constant. Using the classical Amrein-Berthier-Benedicks UP \cite{amrein1977support} for the classical FT and the relation \eqref{R2022_Kamalakkannan-Roopkumar-Zayed_Coupled-FrFT-FT_Relation}, author in \cite{shah2022class} proved that if $f\in L^1(\mathbb{R}^2)$ and $E_1,E_2\subset\mathbb{R}^2$ with the finite measure satisfying $supp(f)\subset E_1$ and $supp\left(\mathcal{F}_{\alpha,\beta}f\right)\subset E_2.$ If $|E_1||E_2|<\infty$ then $f=0.$ Apart from these UPs authors in \cite{shah2022class} proved the Donoho-Stark's UP for the transform \eqref{R2022_Kamalakkannan-Roopkumar-Zayed_Coupled-FrFT-Definition}, using the result in the classical setting in \cite{donoho1989uncertainty}, which says that if $E_1,E_2\subset \mathbb{R}^2$ be any two measurable set such that a non-zero $f\in L^2(\mathbb{R}^2)$ is $\epsilon_{E_1}-$concentrated on $E_1$ and the corresponding coupled FrFT $\mathcal{F}_{\alpha,\beta}f$ is $\epsilon_{E_2}-$concentrated on $E_2,$ then $|E_1||E_2|\geq \sin^2\gamma\left(1-\epsilon_{E_1}-\epsilon_{E_2}\right)^2.$ They also obtained the entropy based UP for the coupled FrFT \eqref{R2022_Kamalakkannan-Roopkumar-Zayed_Coupled-FrFT-Definition} for the function $f\in L^2(\mathbb{R}^2)$ satisfying $\|f\|_{L^2(\mathbb{R}^2)}=1$ as $E\left(|f|^2\right)+4\pi^2 E\left(|\mathcal{F}_{\alpha,\beta}f|^2\right)\geq \ln(\pi e)+4\pi^2\ln \left(4\pi^2\sin^2\gamma\right).$ For a real valued function $f\in L^2(\mathbb{R}^2),$ satisfying $P_lf,~P_l\hat{f}\in L^2(\mathbb{R}^2),~l=1,2,$ where $P_l(\xx_1,\xx_2)=\xx_l,$ Kamalakkannan et al. \cite{kamalakkannan2021short} obtained the Heisenberg UP for the the coupled FrFT \eqref{R2022_Kamalakkannan-Roopkumar-Zayed_Coupled-FrFT-Definition}, in two coupled fractional Fourier domain as 
\begin{align*}
&\Delta\bldu^2_{\alpha_1,\beta_1}\Delta\bldu^2_{\alpha_2,\beta_2}\geq \bigg(\|f\|^2_{L^2(\mathbb{R}^2)}+|\bldxi^0_f|\Delta\bldx^2_f\bigg)\|f\|^2_{L^2(\mathbb{R}^2)}\sin^2(\gamma_1-\gamma_2)\notag\\
&\hspace{6.5cm}+\bigg[(\sin\gamma_1 \sin\gamma_2)\bigg(\Delta\bldxi^2_f+|\bldxi^0_f|^2\|f\|^2_{L^2(\mathbb{R}^2)}\bigg)+(\cos\gamma_1 \cos\gamma_2)\Delta\bldx^2_f\bigg]^2,
\end{align*}
where 
$$\Delta\bldu^2_{\alpha,\beta}=\int_{\mathbb{R}^2}|\bldu-\bldu^{\alpha,\beta,0}|^2\left|\bigg(\mathcal{F}_{\alpha,\beta}f\bigg)(\bldu)\right|^2d\bldu,$$
with $\bldu^{\alpha,\beta,0}=(\uu_1^{\alpha,\beta,0},\uu_2^{\alpha,\beta,0}),$ $\uu_l^{\alpha,\beta,0}=\frac{1}{\|f\|^2_{L^2(\mathbb{R}^2)}}\int_{\mathbb{R}^2}\uu_l\left|\bigg(\mathcal{F}_{\alpha,\beta}f\bigg)(\bldu)\right|^2d\bldu,$ $\gamma_l=\frac{\alpha_l+\beta_l}{2}~l=1,2$ and 
\begin{align}\label{R2021_Kamalakkannan-Roopkumar-Zayed_Coupled_FrFT_Spread_Time_Domain}
\Delta \bldx^2_f=\int_{\mathbb{R}^2}|\bldx-\bldx^0_f|^2|f(\bldx)|^2d\bldx,~\mbox{with}~\bldx^0_f=(\xx^0_1,\xx^0_2),~\xx^0_l=\frac{1}{\|f\|^2_2}\int_{\mathbb{R}^2}\xx_k|f(\bldx)|^2d\bldx.
\end{align}
and 
\begin{align}\label{R2021_Kamalakkannan-Roopkumar-Zayed_Coupled_FrFT_Spread_Frequency_Domain}
\Delta \bldxi^2_f=\int_{\mathbb{R}^2}|\bldxi-\bldxi^0_f|^2|\hat{f}(\bldxi)|^2d\bldxi,~\mbox{with}~\bldxi^0_f=(\xxi^0_1,\xxi^0_2),~\xxi^0_l=\frac{1}{\|f\|^2_{L^2(\mathbb{R}^2)}}\int_{\mathbb{R}^2}\xxi_k|\hat{f}(\bldxi)|^2d\bldxi.
\end{align}
In particular, $\Delta\bldu^2_{\alpha_1,\beta_1}\Delta\bldu^2_{\alpha_2,\beta_2}\geq\|f\|^4_{L^2(\mathbb{R}^2)}\sin^2(\gamma_1-\gamma_2).$ Using the relation between the FrFT \eqref{R2022_Verma-Gupta_IJWMIP_nD-FrFT-Definition} and the classical FT in $n-$dimension and adopting the method used by Guanlei \cite{guanlei2009logarithmic}, authors in \cite{verma2022note} obtained the Heisenberg's UP in fractional Fourier domain as 
\begin{align}\label{R2022_Verma-Gupta_IJWMIP_nD-FrFT-Heisenberg_UP_2FrFDomain}
\left(\int_{\mathbb{R}^n}|\bldt|^2|(\mathfrak{F}_\alpha f)(\bldt)|^2d\bldt\right)\left(\int_{\mathbb{R}^n}|\bldxi|^2|(\mathfrak{F}_\beta f)(\bldxi)|^2d\bldxi\right)\geq \frac{n^2}{4}|\sin(\alpha-\beta)|^2\|f\|^4_{L^2(\mathbb{R}^n)}
\end{align} 
where $\alpha$ and $\beta$ are such that none of $\alpha,~\beta$ and $\alpha-\beta$ are integral multiple of $\pi.$ Moreover, with the similar assumption on $\alpha$ and $\beta$ they obtained the local uncertainty inequality for the transform \eqref{R2022_Verma-Gupta_IJWMIP_nD-FrFT-Definition} in two FrFT domain which says that for a given $\theta>0$ there exists a constant $A_\theta$ such that for any measurable subset $E$ of $\mathbb{R}^n$ and all $f\in L^2(\mathbb{R}^n)$
\begin{align}\label{R2022_Verma-Gupta_IJWMIP_nD-FrFT-Local_UP_2FrFDomain}
\int_E|(\mathfrak{F}_{\alpha}f)(\bldxi)|^2d\bldxi
\leq
\begin{cases}
\frac{A_{\theta}}{|\sin(\alpha-\beta)|^{2\theta}}(\lambda(E))^\frac{2\theta}{n}\left\||\cdot|^\theta(\mathfrak{F}_{\beta}f)(\cdot)\right\|^2_{L^2(\mathbb{R}^n)} &\mbox{if}~ 0<\theta<\frac{n}{2},\\
\frac{A_{\theta}}{|\sin(\alpha-\beta)|^{n}}\lambda(E)\|f\|^{2-\frac{n}{\theta}}_{L^2(\mathbb{R}^n)}\left\||\cdot|^\theta(\mathfrak{F}_{\beta}f)(\cdot)\right\|^{\frac{n}{\theta}}_{L^2(\mathbb{R}^n)} &\mbox{if}~\theta>\frac{n}{2}.
\end{cases}
\end{align}
It is to be noted that for $n=1,$ the Heisenberg's UP \eqref{R2022_Verma-Gupta_IJWMIP_nD-FrFT-Heisenberg_UP_2FrFDomain} reduces to the one obtained by Guanlei \cite{guanlei2009logarithmic} for one-dimensional FrFT. Apart from this, \cite{guanlei2009logarithmic} also contains the logarithmic UP in two FrFT domains for FrFT in one dimension
\begin{align}\label{R2008_Guanlei-Xiaotong-Xiaogang_SP_FrFT_Logarithmic-UP}
\int_{\mathbb{R}}\ln|\tee||(\mathfrak{F}_\alpha f)(\tee)|^2d\tee+\int_{\mathbb{R}}\ln|\xxi||(\mathfrak{F}_\beta f)(\xxi)|^2d\xxi\geq \ln|\sin(\alpha-\beta)|+\frac{\Gamma\left(\frac{1}{4}\right)}{\Gamma\left(\frac{1}{4}\right)},
\end{align}
where $f\in L^2(\mathbb{R})$ is such that $\|f\|_{L^2(\mathbb{R})}=1.$ It follows that if $\alpha\rightarrow 0$ and $\beta\rightarrow \frac{\pi}{2},$ then \eqref{R2008_Guanlei-Xiaotong-Xiaogang_SP_FrFT_Logarithmic-UP} reduces to the logarithmic UP for the classical FT in one dimension. Using the relation between the FrFT $\mathfrak{F}_{\boldsymbol\alpha_\lambda}$ given in equation \eqref{R2020-Kamalakkanan-Roopkumar_ITSF_nD-FrFT_nalpha-Convolution-Theorem} and the classical FT in $n-$dimension, Kaur et al.\cite{kaur2022multidimensional} obtained the Heisenberg's UP for the FrFT $\mathfrak{F}_{\boldsymbol\alpha_\lambda},$ defined in \eqref{R2020-Kamalakkanan-Roopkumar_ITSF_nD-FrFT_nalpha-Definition}, as
\begin{align}\label{R2022-Kaur-Gupta-Verma_nD-FrFT_nalpha-Heisenberg-UP}
\left(\int_{\mathbb{R}^n}|\bldt|^2|f(\bldt)|^2d\bldt\right)\left(\int_{\mathbb{R}^n}|\bldxi|^2|(\mathfrak{F}_{\boldsymbol\alpha_\lambda}f)(\bldxi)|^2 d\bldxi\right)\geq P_{\boldsymbol\alpha,\lambda}\frac{n^2}{4}\|f\|^4_{L^2(\mathbb{R}^n)},
\end{align}
for $f\in L^2(\mathbb{R}^n)$ where $P_{\boldsymbol\alpha,\lambda}=\frac{|c(\boldsymbol\alpha_\lambda)|^2|\sin\boldsymbol\alpha|_m}{M^2(\lambda)^{n+2}},$ $M=\mbox{max}\left\{\frac{1}{|\sin\alpha_l|},l=1,2,\cdots,n\right\}$ and with $P_{\boldsymbol\alpha,\lambda}'=\frac{\lambda^{2n}}{|c(\boldsymbol\alpha_\lambda)|^2|\sin\boldsymbol\alpha|_m},$ $M=\mbox{max}\left\{\frac{1}{|\sin\alpha_l|},l=1,2,\cdots,n\right\}$ they also obtained the logarithmic UP 
\begin{align}\label{R2022-Kaur-Gupta-Verma_nD-FrFT_nalpha-Logarithmic-UP}
\int_{\mathbb{R}^n}\ln|\bldx||f(\bldx)|^2d\bldx+P_{\boldsymbol\alpha,\lambda}'\int_{\mathbb{R}^n}\int_{\mathbb{R}^n}\ln|\bldxi||(\mathfrak{F}_{\boldsymbol\alpha_\lambda}f)(\bldxi)|^2d\bldxi\geq \left(D-P_{\boldsymbol\alpha,\lambda}'\ln(\lambda^2M)\right)\int_{\mathbb{R}^n}|f(\bldx)|^2d\bldx
\end{align}
followed by the local uncertainty inequality, which says that for a given $\theta>0$, there exists a constant $A_\theta$ such that for any measurable subset $E$ of $\mathbb{R}^n$ and all $f\in L^2(\mathbb{R}^n)$
\begin{align}\label{R2022-Kaur-Gupta-Verma_nD-FrFT_nalpha-Local-UP}
\int_E|(\mathcal{F}_{\boldsymbol\alpha_\lambda}f)(\bldxi)|^2d\bldxi
\leq
\begin{cases}
|c(\boldsymbol\alpha_\lambda)|^2A_{\theta}\frac{|\sin\boldsymbol\alpha|_m}{(\lambda^2)^n}(\lambda(E))^\frac{2\theta}{n}\left\||\cdot|^\theta f(\cdot)\right\|^2_{L^2(\mathbb{R}^n)} &\mbox{if}~ 0<\theta<\frac{n}{2},\\
|c(\boldsymbol\alpha_\lambda)|^2A_{\theta}\frac{|\sin\boldsymbol\alpha|_m}{(\lambda^2)^n}\lambda(E)\|f\|^{2-\frac{n}{\theta}}_{L^2(\mathbb{R}^n)}\left\||\cdot|^\theta f(\cdot)\right\|^{\frac{n}{\theta}}_{L^2(\mathbb{R}^n)} &\mbox{if}~\theta>\frac{n}{2}.
\end{cases}
\end{align} 

\subsection{Fractional Fourier transform of type II}

If $\mathcal{S}(\mathbb{R})$ is the space of rapidly decreasing test functions, then the Lizorkin space $\Phi(\mathbb{R})$ is the Fourier pre-image of the space $V(\mathbb{R})$ in the space $\mathcal{S}(\mathbb{R}),$ i.e., 
$$\Phi(\mathbb{R})=\left\{\phi\in\mathcal{S}(\mathbb{R}):\hat{\phi}\in V(\mathbb{R})\right\},$$
where $V(\mathbb{R})$ is the set of all $v\in \mathcal{S}(\mathbb{R})$ satisfying the conditions $\frac{d^nv(x)}{dx^n}|_{x=0}=0,~n=0,1,2,\cdots.$ It is to be noted that the space $\Phi(\mathbb{R})$ is invariant with respect to the fractional integration and differentiation operators, which is not the case for the whole space $\mathcal{S}(\mathbb{R}).$ Luchko et al. \cite{luchko2008fractional} introduced a new definition of the FrFT of real order $\alpha (0<\alpha\leq 1)$ for a function $u\in\Phi(\mathbb{R})$ as 
\begin{align}\label{R2008_Luchko-Martinez-Trujillo_FCAA_NewFrFT_Definition}
(\mathit{F}_\alpha u)(\oomega)=\int_{\mathbb{R}}u(\tee)e_{\alpha}(\oomega,\tee)d\tee,~\oomega\in\mathbb{R},
\end{align}
where 
\begin{align*}
e_{\alpha}(\oomega,\tee)=e^{i\sgn(\oomega)|\oomega|^{\frac{1}{\alpha}}\tee}=
\begin{cases}
e^{-i|\oomega|\tee},&\oomega\leq 0,\\
e^{i|\oomega|\tee},&\oomega\geq 0.
\end{cases}
\end{align*}
It is to be noted that for $\alpha=1,$ the kernel $e_{\alpha}(\oomega,\tee)$ coincides with that of the kernel of the conventional FT, i.e., the FrFT \eqref{R2008_Luchko-Martinez-Trujillo_FCAA_NewFrFT_Definition} reduces to the classical FT when $\alpha=1.$ Moreover the relation between the FrFT \eqref{R2008_Luchko-Martinez-Trujillo_FCAA_NewFrFT_Definition} and the classical FT is given by the formula
\begin{align}\label{R2008_Luchko-Martinez-Trujillo_FCAA_Relation_NewFrFT&FT}
(\mathit{F}_\alpha u)(\oomega)=\hat{u}(x)
\end{align}
where $
x=\sgn(\oomega)|\oomega|^{\frac{1}{\alpha}}=
\begin{cases}
-|\oomega|^{\frac{1}{\alpha}},&\oomega\leq 0,\\
|\oomega|^{\frac{1}{\alpha}},&\oomega\geq 0.
\end{cases}
$
Using this formula and the known properties of the FT, one can calculate the FrFT of the concrete functions and determine the inversion of the FrFT. Kilbas et al. \cite{kilbas2010fractional} obtained the inversion formula of \eqref{R2008_Luchko-Martinez-Trujillo_FCAA_NewFrFT_Definition} as 
\begin{align}
(\mathit{F}^{-1}_\alpha u)(\xx)=\frac{1}{2\pi \alpha}\int_{\mathbb{R}}e^{-i\sgn(\oomega)|\oomega|^{\frac{1}{\alpha}}\xx}|\oomega|^{\frac{1}{\alpha}-1}u(\oomega)d\oomega,~\xx\in\mathbb{R},\alpha>0,
\end{align}
for $u\in \Phi(\mathbb{R})$ and also obtain the  FrFT of the translation and dilation of $u\in \Phi(\mathbb{R})$ respectively as 
\begin{align}\label{R2010_Kilabs-Luchko-Martinez-Trujillo_ITSF_NewFrFT_Translation-Formula}
\left(\mathit{F}_\alpha (\tau_hu)\right)(\oomega)=e^{i\sgn(\oomega)|\oomega|^{\frac{1}{\alpha}}h}(\mathit{F}_\alpha u)(\oomega)
\end{align}
and
\begin{align}\label{R2010_Kilabs-Luchko-Martinez-Trujillo_ITSF_NewFrFT_Dilation-Formula}
\left(\mathit{F}_\alpha (u(\lambda\tee))\right)(\oomega)=\frac{1}{\lambda}(\mathit{F}_\alpha u)\left(\frac{\oomega}{\lambda^\alpha}\right),~\lambda>0.
\end{align}
Apart from this authors also obtained the convolution theorem for the FrFT as 
\begin{align}\label{R2010_Kilabs-Luchko-Martinez-Trujillo_ITSF_NewFrFT_Convolution-Theorem}
\left(\mathit{F}_\alpha (k\star u)\right)(\oomega)=(\mathit{F}_\alpha k)(\oomega)(\mathit{F}_\alpha u)(\oomega)
\end{align}
followed by the formula 
\begin{align}\label{R2010_Kilabs-Luchko-Martinez-Trujillo_ITSF_NewFrFT_Differentiation-Formula}
\left(\mathit{F}_\alpha u^{m}\right)(\oomega)=\left(-i\sgn(\oomega)|\oomega|^{\frac{1}{\alpha}}\right)^m(\mathit{F}_\alpha u)(\oomega),
\end{align}
where $\star$ denotes the classical convolution and $k,u\in \Phi(\mathbb{R}).$ Observe that the formulas given by equations \eqref{R2010_Kilabs-Luchko-Martinez-Trujillo_ITSF_NewFrFT_Translation-Formula}, \eqref{R2010_Kilabs-Luchko-Martinez-Trujillo_ITSF_NewFrFT_Dilation-Formula}, \eqref{R2010_Kilabs-Luchko-Martinez-Trujillo_ITSF_NewFrFT_Convolution-Theorem} and \eqref{R2010_Kilabs-Luchko-Martinez-Trujillo_ITSF_NewFrFT_Differentiation-Formula} reduces to the corresponding formulas in the case of classical FT when $\alpha=1.$ We know that the conditions of their validity can be improved, for example the convolution theorem for the classical FT is valid for either $k\in L^1(\mathbb{R})$ and $u\in L^1(\mathbb{R}),$ or $k\in L^2(\mathbb{R})$ and $u\in L^1(\mathbb{R})$ or $k\in L^2(\mathbb{R})$ and $u\in L^2(\mathbb{R}).$ Thus, using the relation \eqref{R2008_Luchko-Martinez-Trujillo_FCAA_Relation_NewFrFT&FT}, authors in \cite{kilbas2010fractional} mentioned that the condition of validity of the formulas given in equations \eqref{R2010_Kilabs-Luchko-Martinez-Trujillo_ITSF_NewFrFT_Translation-Formula}, \eqref{R2010_Kilabs-Luchko-Martinez-Trujillo_ITSF_NewFrFT_Dilation-Formula}, \eqref{R2010_Kilabs-Luchko-Martinez-Trujillo_ITSF_NewFrFT_Convolution-Theorem} and \eqref{R2010_Kilabs-Luchko-Martinez-Trujillo_ITSF_NewFrFT_Differentiation-Formula} can be improved as in the case of FT. Srivastava et al. \cite{srivastava2017family} proved that the transform \eqref{R2008_Luchko-Martinez-Trujillo_FCAA_NewFrFT_Definition} and its inverse is an isomorphism from $S(\mathbb{R})$ onto itself. Srivastava et al.\cite{srivastava2019certain} studied the FrFT given by \eqref{R2008_Luchko-Martinez-Trujillo_FCAA_NewFrFT_Definition} for the functions in $L^2(\mathbb{R})$ and obtained the Parseval formula for $f,g\in L^2(\mathbb{R})$ as  
\begin{align*}
\int_{\mathbb{R}}f(\tee)\overline{g(\tee)}d\tee=\frac{1}{2\pi\alpha}\int_{\mathbb{R}}|\oomega|^{\frac{1}{\alpha}-1}(\mathit{F}_\alpha f)(\oomega)\overline{(\mathit{F}_\alpha g)(\oomega)}d\oomega.
\end{align*}

\section{Linear Canonical Transform (LCT)}
Moshinsky and Quesne proposed linear canonical transform (LCT) \cite{moshinsky1971linear} in the early 1990s, which is a promising linear integral transform with three free parameters. It is a generalization of several integral transforms like FT, FrFT, Fresnel transform, and the scaling operations. For a real, unimodular matrix $M=(\Ma,\Mb,\Mc,\Md),$ the LCT of a function $f\in L^2(\mathbb{R}^2)$ is defined by
\begin{align}\label{RDefinitionLCT}
(\mathcal{L}^Mf)(\xxi)=\int_{\mathbb{R}}f(\tee)\mathcal{K}_M(\tee,\xxi)d\tee,
\end{align}
where the kernel $\mathcal{K}_M(\tee,\xxi)$ is given by
\begin{align}\label{RDefinitionLCTKernel}
\mathcal{K}_M(\tee,\xxi)=\frac{1}{\sqrt{2\pi i\Mb}}e^{\frac{i(\Ma\tee^2-2\tee\xxi+\Md\xxi^2)}{2\Mb}},~\Mb\neq 0.
\end{align}
It can be seen from \eqref{RDefinitionLCT} and \eqref{RDefinitionLCTKernel}, that for special choices of $M$, \eqref{RDefinitionLCT} reduces to several known integral transforms like
\begin{enumerate}
\item for $M=(0,1,-1,0),$ definition \eqref{RDefinitionLCT} reduces to classical FT
\item for $M=(\cos\theta,\sin\theta,-\sin\theta,\cos\theta),\theta\neq n\pi,~n\in\mathbb{Z}$ definition \eqref{RDefinitionLCT} reduces to FrFT
\item for $M=(1,\Mb,0,1),~\Mb\neq 0$ definition \eqref{RDefinitionLCT} reduces to classical Fresnel transform.
\end{enumerate}
Determining a convolution theorem by defining a suitable convolution associated with a given transform has become an active area of research. Different transform requires different definitions of convolutions to obtain the associated convolution theorem. Like in the case of FT, it gives the powerful result that the convolution $``\star"$, in \eqref{RConvolutionFT}, of two signals in the time domain results in simple multiplication of their FTs in the FT domain. However, the LCT of $(f\star g)(\tee)$ is given by \cite{deng2006convolution}
\begin{align}
\left(\mathcal{L}^M(f*g)\right)(\uu)=\left|\frac{1}{\Ma}\right| e^{-i \frac{\Mc}{2 \Ma} \uu^2} \int_{-\infty}^{\infty} \left(\mathcal{L}^Mf\right)(\vv) e^{-\i \frac{\Mc}{2 \Ma} \vv^2} g\left(\frac{\uu-\vv}{\Ma}\right) d\vv.
\end{align}
This formula lacks the requirement that the convolution in the time domain is the product of LCTs in the LCT domain. Moreover, it is inconvenient for the analysis of multiplicative filters in the LCT domain. Deng et al. \cite{deng2006convolution} also gave another convolution theorem for the LCT
\begin{align}\label{RCovolutionTheoremLCT2}
&\left(\mathcal{L}^M(f\star_M g)\right)(\uu)=\sqrt{i 2 \pi \Mb} \left(\mathcal{L}^Mf\right)(\uu) \left(\mathcal{L}^Mg\right)(\uu) e^{-i \frac{\Md}{2 \Mb} \uu^2}
\end{align}
based on the modified convolution operation given by
\begin{align}\label{RCovolutionLCT2}
\left(f\star_M g\right)(\tee)&=e^{-i \frac{\Ma}{2 \Mb} \tee^2}\left(\left(f(\cdot) e^{i\frac{\Ma}{2 \Mb} (\cdot)^2}\right) *\left(g(\cdot) e^{i\frac{\Ma}{2 \Mb} (\cdot)^2}\right)\right)(\tee)\notag\\
&=\int_{\mathbb{R}} f(\yy)g(\tee-\yy)e^{i\yy(\yy-\tee)\frac{\Ma}{\Mb}}d\yy
\end{align}
Even though the convolution theorem for the LCT given in \eqref{RCovolutionTheoremLCT2} satisfies our need for getting the product of LCTs in the LCT domain, during the process of evaluation of the convolution \eqref{RCovolutionLCT2}, the signal are multiplied three times by the different chirp signals. From realization point of view, this may create difficulty as in communication system, it is nearly impossible to generate a chirp signal accurately. Wei et al. \cite{wei2009convolution} generalized the convolution theorem for the LCT as 
\begin{align}\label{RCovolutionTheoremLCT3}
\left(\mathcal{L}^M\left(f\stackrel{M}{\Theta} g\right)\right)(\uu) = \left(\mathcal{L}^Mf\right)(\uu) \left(\mathcal{L}^Mg\right)(\uu)
\end{align}
where $\stackrel{M}{\Theta}$ denotes generalized convolution $\left(f\stackrel{M}{\Theta} g\right)(\tee)=\int_{\mathbb{R}}f(\tau)g(\tee\theta\tau)d\tau$ in the LCT domain, where $g(\tee\theta\tau)$ is the $\tau-$generalized translation of signal $g(\tee)$ by $\tau$. This convolution theorem has the needed form that the convolution in the time domain is equal to the product of the LCTs in the LCT domain, as in the case of the classical result for the FT. However, it is complicated to reduce the expression of the generalized convolution to a single integral form as in ordinary convolution \eqref{RConvolutionFT}. So, it is difficult to implement it compared to the FT. For the function $f$ with $\tilde{f}(\tee)=e^{i \frac{\Ma\tee^2}{2\Mb}},$ Wei et al.\cite{wei2011multichannel} proposed a new convolution for the LCT  of $f$ and $g$ is defined as
\begin{align}\label{R2011_Wei-Ran-Li-LCT-Convolution}
(f\star'_M g)(\tee)=e^{-i\frac{\Ma \tee^2}{2\Mb}}\left(\tilde{f}\star g\right)(\tee)
\end{align} 
and obtain the following convolution theorem  \begin{align}\label{R2011_Wei-Ran-Li-LCT-Convolution-Theorem}
\left(L_M(f\star'_M g)\right)(\uu)=\left(L_M f\right)(\uu)(\mathcal{F}g)\left(\frac{\uu}{\Mb}\right),
\end{align} 
where $(\mathcal{F}g)(\uu)=\int_{\mathbb{R}}g(t)e^{-itu}dt.$ Thus, the LCT of the convolution $\star'$ given by \eqref{R2011_Wei-Ran-Li-LCT-Convolution} is the multiplication of the LCT of the one with the FT of the other. Moreover, the convolution can be expressed by one-dimensional integrals, as in the ordinary convolution expression, which is useful and easy to implement, particularly in filter designing.

The classical convolution has translational invariance property, i.e., 
$$\tau_\zz(f\star g)=(\tau_\zz f)\star g=g\star (\tau_\zz g),$$
where $\star$ denotes the convolution given by \eqref{RConvolutionFT} and $(\tau_\zz f)(\tee)=f(\tee-\zz).$ So many convolution structures are defined for LCT in the literature, as discussed. Wei et al. \cite{wei2012new} pointed out that the convolution structure $\star_M$ given by \eqref{RCovolutionLCT2} is not translational invariance with respect to the translation operator $\tau_\zz,$ so they generalized the notion of translation and called it as linear canonical translation which is as follows:

For the LCT parameter $M=(\Ma,\Mb,\Mc,\Md)$ and value $\zz$, 
$T_{M;\zz} f:=(f\star_{M}\delta) e^{-i\frac{\Ma\zz^2}{2\Mb}},$ which can be equivalently written as 
$$(T_{M;\zz} f)(\tee)=f(\tee-\zz)e^{-i\frac{z(t-\frac{z}{2})\Ma}{\Mb}}.$$
For $M=(0,1,-1,0)$ $T_{M,\zz}=\tau_\zz$(which is the usual translation by $\zz$).
With this linear canonical translation operator, they obtained the translation invariant property for LCT i.e.,
$$T_{M;\zz}(f\star_{M}g)=(T_{M;\zz}f)\star_{M}g=f\star_{M}(T_{M;\zz} g)$$
Motivated by the definition \eqref{RConvolutionFTEquivalent} of classical convolution, Shi et al.\cite{shi2014generalizedconvolution} introduced a generalized canonical convolution operator, $\star_{M_1, M_2, M_3},$ for $f,g\in L^2(\mathbb{R})$ in the time domain as 
\begin{align}\label{RConvolutionGeneralizedLCT}
\left(f \star_{M_1, M_2, M_3} h\right)(\tee)&=\int_{-\infty}^{+\infty} h(\tau)\left(\mathrm{T}_\tau^{M_1} f\right)(\tee) \phi_{M_1, M_2, M_3}(\tee, \tau) \mathrm{d} \tau\notag\\
&=e^{-i\frac{\Ma_3}{2\Mb_3}\tee^2}\left[f(\cdot)e^{i\frac{\Ma_1}{2\Mb_1}(\cdot)^2}\star g(\cdot)e^{i\frac{\Ma_2}{2\Mb_2}(\cdot)^2}\right](\tee)
\end{align}
where $\left(\mathrm{T}_\tau^{M_1} f\right)(\tee)=f(\tee-\tau)e^{-\frac{i\Ma_1}{\Mb_1}\tau(\tee-\frac{\tau}{2})}$ and
$\phi_{M_1,M_2, M_3}(\uu, \vv) = \mathrm{e}^{i \frac{\Ma_2}{2 \Mb_2}\vv^2+i\left(\frac{\Ma_1}{2 \Mb_1}-\frac{\Ma_3}{2 \Mb_3}\right)\uu^2}$
is a domain-independent weighted function. It is to be noted that the convolution $\star_{M_1,M_2,M_3}$ given by equation \eqref{RConvolutionGeneralizedLCT} reduces to the convolution operator $\star_{\alpha,\beta,\gamma}$ \eqref{RShiFrFTConvolutionTheoremII}, for the choice of the matrices $M_1=(\cos\beta,\sin\beta;-\sin\beta,\cos\beta),$ $M_2=(\cos\gamma,\sin\gamma;-\sin\gamma,\cos\gamma)$ and $M_3=(\cos\alpha,\sin\alpha;-\sin\alpha,\cos\alpha)$. They also obtained the convolution theorem associated with the generalized convolution \eqref{RConvolutionGeneralizedLCT} as 
\begin{align}\label{RConvolutionTheoremGeneralizedLCT}
\mathcal{L}^{M_3}\left[f \star_{M_1, M_2, M_3} h\right](\uu)=\epsilon_{M_1, M_2, M_3}(\uu) \left(\mathcal{L}^{M_1}f\right)\left(\frac{\uu \Mb_1}{\Mb_3}\right) \left(\mathcal{L}^{M_1}h\right)\left(\frac{\uu \Mb_2}{\Mb_3}\right)
\end{align}
where $\left(\mathcal{L}^{M_1}f\right)(\uu)$ and $\left(\mathcal{L}^{M_2}h\right)(\uu)$ respectively denote the LCT with parameter $M_1$ of $f(\tee)$ and the LCT with parameter $M_2$ of $h(\tee)$ and  
$$
\epsilon_{M_1, M_2, M_3}(\uu) = \frac{M_{\Mb_3}}{M_{\Mb_1} M_{\Mb_2}} \mathrm{e}^{i \uu^2\left(\frac{\Md_3}{2 \Mb_3}-\frac{\Md_1}{2 \Mb_1}\left(\frac{\Mb_1}{\Mb_3}\right)^2-\frac{\Md_2}{2 \Mb_2}\left(\frac{\Mb_2}{\Mb_3}\right)^2\right)},~M_{\Mb_l}=\frac{1}{\sqrt{i2\pi \Mb_l}},~l=1,2,3,
$$
is a domain-independent weighted function. The convolution result for the LCT given in \eqref{RConvolutionTheoremGeneralizedLCT} generalized several convolution theorems for particular choices of the matrix parameter $M_1,~M_2$ and $M_3$ \cite{deng2006convolution},\cite{kraniauskas1998method},\cite{shi2010multichannel},\cite{torres2010fractional},\cite{wei2012new},\cite{zayed1998convolution}.\\

The sampling process is central in almost all domains of mathematics. One of the most fundamental concepts in digital signal processing is the sampling process, which serves as a link between a continuous and discrete signal. Before we proceed with our discussion on the sampling theorem, we need to know that a signal is said to be time-limited if it is compactly supported in the time domain, whereas it is band-limited if it is compactly supported in the transformed domain. We have interesting results concerning the support of the signal in the LCT domain. It says that if a nonzero signal $f(t)$ is band-limited in the LCT domain with parameter 
$M_1=\left(\begin{array}{ll}\Ma_1 & \Mb_1 \\ \Mc_1 & \Md_1\end{array}\right)$, 
then $f(t)$ cannot be band-limited with another parameter 
$M_2=\left(\begin{array}{ll}\Ma_2 & \Mb_2 \\ \Mc_2 & \Md_2\end{array}\right)$, 
satisfying the relationship $\Ma_1 / \Mb_1 \neq \Ma_2 / \Mb_2$. By virtue of this, it follows that if a signal is non-bandlimited in the Fourier domain, then it can be band-limited in the LCT domain. Healy et al. \cite{healy2008cases} presented non-band-limited signals in the Fourier domain that are band-limited in the LCT domain. So, it is important to establish the sampling theorem for the band-limited function in the LCT domain. Stern \cite{stern2006sampling} obtained the sampling theorem for the LCT of the compactly supported functions in the time domain or the LCT domain. Even though the formulae given in \cite{stern2006sampling} can be used to reconstruct the original signal or its LCT from their samples at a discrete set of points that satisfies the Nyquist rate of that particular domain, it is not strictly in an interpolation form. Keeping this in mind Li et al.\cite{li2007new} developed a sampling theorem for $f(\tee)$  band-limited to $\Omega_{M}$ in the LCT domain with parameter 
$M=\left(\begin{array}{ll}\Ma & \Mb \\ \Mc & \Md\end{array}\right)$ and $\Mb>0$ as
\begin{align}\label{RLiSamplingTheoremLCT}
f(\tee)=e^{-i \frac{\Ma}{2 \Mb} \tee^2} \sum_{n=-\infty}^{+\infty} f(n T) e^{i \frac{\Ma}{2 \Mb}(n T)^2} \frac{\sin \left[\frac{\Omega_M(\tee-n T)}{\Mb\pi}\right]}{\frac{\Omega_M(\tee-n T)}{\Mb \pi}}
\end{align}
where $T$ is the sampling period and satisfies $T=\pi \Mb / \Omega_{M}$; and the Nyquist rate of sampling theorem associated with the LCT is $f_s=\Omega_{M} / \pi \Mb$. The result \eqref{RLiSamplingTheoremLCT}, which is derived using the relation between the LCT and the FT, is more explicit than that given in \cite{stern2006sampling}, and its derivation is easier than that in \cite{deng2006convolution}. It also generalizes the sampling theorem for the FrFT studied in \cite{xia1996bandlimited},\cite{zayed1996relationship}. From an application point of view, when sampling non-stationary signals or broadband signals, the sampling condition in \eqref{RLiSamplingTheoremLCT} may not always hold because of the requirement of the low sampling rate. So, Li et al.\cite{li2007new} also developed the following sampling formulae that are used to reconstruct the signal from the samples of the signal and its first derivatives or its generalized Hilbert transform \cite{li2006hilbert} associated with the LCT each at half of the Nyquist rate.
\begin{enumerate}
\item If $f(\tee)$ is band-limited to $\Omega_{M}$ in LCT sense and has a continuous first derivative. Then the following sampling formula for $f(\tee)$ holds
\begin{align}\label{R}
f(\tee)=\mathrm{e}^{-i\left(\frac{\Ma}{2 \Mb}\right) \tee^2} \sum_{k=-\infty}^{\infty} \mathrm{e}^{i\left(\frac{\Ma}{2 \Mb}\right)\left(k T_1\right)^2}\left\{\left[1+i\left(k T_1\right)\left(\tee-k T_1\right) \frac{\Ma}{\Mb}\right]f\left(k T_1\right)+\left(\tee-k T_1\right) f^{\prime}\left(k T_1\right)\right\}\left[\operatorname{sinc}\left(\frac{\tee}{T_1}-k\right)\right]^2
\end{align}
where $T_1=2 T=2 \pi \Mb / \Omega_{M}, f^{\prime}\left(k T_1\right)=\left.f^{\prime}(\tee)\right|_{\tee=k T_1}.$

\item If a signal $f(\tee)$ is band-limited to $\Omega_{M}$ in the LCT sense. Then $f(\tee)$ can be represented by the samples of $f(\tee)$ and its generalized Hilbert transform $H_{M}[f(\tee)]$ at half Nyquist rate as follows:
\begin{align}\label{RLiSamplingTheoremLCTHT}
f(\tee)=&\mathrm{e}^{-i\left(\frac{\Ma}{2 \Mb}\right) \tee^2} \sum_{k=-\infty}^{+\infty}\left\{\mathrm{e}^{i\left(\frac{\Ma}{2 \Mb}\right)\left(k T_2\right)^2} f\left(k T_2\right)\cos \left[\frac{\pi}{T_2}\left(\tee-k T_2\right)\right]-\mathrm{e}^{i\left(\frac{\Md}{2 \Mb}\right)\left(k T_2\right)^2} H_{\mathrm{A}}[f]\left(k T_2\right)\sin \left[\frac{\pi}{T_2}\left(\tee-k T_2\right)\right]\right\}\notag\\
&\hspace{10cm}\times\operatorname{sinc}\left[\frac{\pi}{T_2}\left(\tee-k T_2\right)\right],
\end{align}
where $T_2=2 T=2 \pi \Mb / \Omega_{M}$.
Equation \eqref{RLiSamplingTheoremLCTHT} is an extension of the result in \cite{zayed1999new} associated with the FrFT.
\end{enumerate}
Equation \eqref{RLiSamplingTheoremLCT} establishes the relationship between the signal samples and the original signal; that is to say, we can reconstruct the original signal or the LCT of a band-limited signal in the LCT domain from the uniform sampling points of the signal, provided the sampling interval satisfies the uniform sampling conditions. However, when applying this formula in practical situations, two problems are encountered. The first is concerned with nonuniform sampling signal processing, and the other is related to the fact that only a finite number of samples are practically available. This problem is addressed in \cite{tao2008sampling}.
\\
\\
\noindent\textbf{Uncertainty Principles:} Guanlei et al.\cite{guanlei2009new} extended the traditional Pitt's inequality in two LCT domain as 
\begin{align}\label{R2009Guanlie-Xiaotong-XiaogangPitt'sInequalityin2LCTDomain}
\int_{\mathbb{R}}|\uu|^{-\lambda}|\left(\mathcal{L}^{M_1}f\right)(\uu)|^2d\uu\leq\frac{R_\lambda}{|\Ma_1\Mb_2-\Ma_2\Mb_1|^\lambda}\int_{\mathbb{R}}|\vv|^{\lambda}|\left(\mathcal{L}^{M_2}f\right)(\vv)|^2d\vv,
\end{align}
where $R_\lambda=\left[\frac{\Gamma(\frac{1-\lambda}{4})}{\Gamma\left(\frac{1+\lambda}{4}\right)}\right]^2,$ $0\leq \lambda<1.$ Setting $S(\lambda)=|\Ma_1\Mb_2-\Ma_2\Mb_1|^\lambda\int_{\mathbb{R}}|\uu|^{-\lambda}|\left(\mathcal{L}^{M_1}f\right)(\uu)|^2d\uu-R_\lambda\int_{\mathbb{R}}|\vv|^{\lambda}|\left(\mathcal{L}^{M_2}f\right)(\vv)|^2d\vv$ one can obtain $S'(0+)\leq 0,$ which on simplification gives
\begin{align}\label{R2009Guanlie-Xiaotong-XiaogangLogatithmicUPin2LCTDomain}
\int_{\mathbb{R}}\ln|\uu|\left|\left(\mathcal{L}^{M_1}f\right)(\uu)\right|^2d\uu+\int_{\mathbb{R}}\ln|\vv|\left|\left(\mathcal{L}^{M_2}f\right)(\vv)\right|^2d\vv\geq \ln\left(|\Ma_1\Mb_2-\Ma_2\Mb_1|\right)+\frac{\Gamma'\left(\frac{1}{4}\right)}{\Gamma\left(\frac{1}{4}\right)}
\end{align}
which is called the logarithmic UP in LCT domain. If $M_1=(0,-1,1,1)$ and $M_2=(\nu,\nu-1,1,1),$ where $\nu=\sqrt{-\frac{2\Gamma'\left(\frac{1}{4}\right)}{\Gamma\left(\frac{1}{4}\right)}},$ then $\ln\left(|\Ma_1\Mb_2-\Ma_2\Mb_1|^2\right)+\frac{2\Gamma'\left(\frac{1}{4}\right)}{\Gamma\left(\frac{1}{4}\right)}=0,$ i.e., the bound of \eqref{R2009Guanlie-Xiaotong-XiaogangLogatithmicUPin2LCTDomain} may be zero. Authors in \cite{guanlei2009new} also obtained the entropy UP for the LCT as 
\begin{align}\label{R2009Guanlie-Xiaotong-XiaogangGEntropyUPLCT}
E(|\left(\mathcal{L}^{M_1}f\right)(\uu)|^2)+E(|\left(\mathcal{L}^{M_2}f\right)(\vv)|^2)\geq \ln(\pi e|\Ma_1\Mb_2-\Ma_2\Mb_1|),
\end{align}
where $E(\rho(\xx))=-\int_{\mathbb{R}}\rho(\xx)\ln\rho(\xx)d\xx,$ $\rho(\xx)$ is the probability density function and the Heisenberg's UP in two LCT domain as 
\begin{align}\label{R2009Guanlie-Xiaotong-XiaogangHeisenbergUPin2LCTDomain}
\left(\int_{\mathbb{R}}|\uu|^2|\left(\mathcal{L}^{M_1}f\right)(\uu)|^2d\uu\right)\left(\int_{\mathbb{R}}|\vv|^2|\left(\mathcal{L}^{M_2}f\right)(\vv)|^2d\vv\right)\geq\frac{|\Ma_1\Mb_2-\Ma_2\Mb_1|^2}{4},
\end{align}
for function $f$ satisfying $\|f\|_{L^2(\mathbb{R})}=1.$ Equation \eqref{R2009Guanlie-Xiaotong-XiaogangHeisenbergUPin2LCTDomain}, for more general centres can be found in \cite{guanlei2009uncertainty}. It is to be noted that \eqref{R2009Guanlie-Xiaotong-XiaogangPitt'sInequalityin2LCTDomain}, \eqref{R2009Guanlie-Xiaotong-XiaogangGEntropyUPLCT} and \eqref{R2009Guanlie-Xiaotong-XiaogangHeisenbergUPin2LCTDomain} depends only on the parameters $\Ma$ and $\Mb$ but not on $\Mc$ and $\Md,$ because from the definition of the LCT it is seen that the parameters $\Mc$ and $\Md$ only play the role of scaling and modulation and that they have no effect on the inequality \eqref{R2009Guanlie-Xiaotong-XiaogangPitt'sInequalityin2LCTDomain}, \eqref{R2009Guanlie-Xiaotong-XiaogangGEntropyUPLCT} and \eqref{R2009Guanlie-Xiaotong-XiaogangHeisenbergUPin2LCTDomain}. If $M_1=(0,-1,1,1)$ and $M_2=(\frac{1}{\pi e},\frac{1}{\pi e}-1,1,1)$ then we have $\ln(\pi e|\Ma_1\Mb_2-\Ma_2\Mb_1|)=0,$ which means that the bound of the entropic UP may be zero. Moreover, if $M_1=(\Ma_1,\Mb_1,\Mc_1,\Md_1)=(\cos\alpha,\sin\alpha,-\sin\alpha,\cos\alpha)$ and $M_2=(\Ma_2,\Mb_2,\Mc_2,\Md_2)=(\cos\beta,\sin\beta,-\sin\beta,\cos\beta)$ then \eqref{R2009Guanlie-Xiaotong-XiaogangLogatithmicUPin2LCTDomain},\eqref{R2009Guanlie-Xiaotong-XiaogangGEntropyUPLCT} and \eqref{R2009Guanlie-Xiaotong-XiaogangHeisenbergUPin2LCTDomain} reduces respectively to the logarithmic UP, entropy UP and Heisenberg's UP \cite{guanlei2009generalized},\cite{guanlei2009logarithmic} for the FrFT and to the traditional case \cite{folland1997uncertainty} when $M_1=(1,0,0,1)$ and $M_2=(0,1,-1,0),$ i.e., in the later case \eqref{R2009Guanlie-Xiaotong-XiaogangPitt'sInequalityin2LCTDomain}, \eqref{R2009Guanlie-Xiaotong-XiaogangLogatithmicUPin2LCTDomain}, \eqref{R2009Guanlie-Xiaotong-XiaogangGEntropyUPLCT} and \eqref{R2009Guanlie-Xiaotong-XiaogangHeisenbergUPin2LCTDomain} reduces, respectively, to the traditional Pitt's inequality  \cite{beckner1995pitt}, logarithmic UP by Beckner \cite{beckner2012pitt}, Shannon entropy \cite{hirschman1957note} and Heisenberg UP. Guanlei et al.\cite{guanlei2009uncertainty} obtained the R\'enyi entropy-based UP in two LCT domains as 
\begin{align}\label{Rguanlei2009uncertainty_Renyi_UP_2LCTD}
H_\theta\left\{\left|(\mathcal{L}^{M_2}f)(\uu)\right|^2\right\}+H_\gamma\left\{\left|(\mathcal{L}^{M_1}f)(\uu)\right|^2\right\}\geq \frac{\ln\left(\frac{\theta}{\pi}\right)}{2(\theta-1)}+\frac{\ln\left(\frac{\gamma}{\pi}\right)}{2(\gamma-1)}+\ln|\Ma_1\Mb_2-\Ma_2\Mb_1|,
\end{align}
using the Hausdorff-Young inequality in two LCT domain, for $1<q\leq 2$ and $\frac{1}{p}+\frac{1}{q}=1,$ 
\begin{align}\label{Rguanlei2009uncertainty_H-Y_inequality_2LCTD}
\|(\mathcal{L}^{M_2}f)(\uu)\|_{L^p(\mathbb{R})}\leq \left(\frac{\left(\frac{q}{2\pi}\right)^{\frac{1}{q}}}{\left(\frac{p}{2\pi}\right)^{\frac{1}{p}}}\right)^{\frac{1}{2}}\left(\frac{1}{|\Ma_1\Mb_2-\Ma_2\Mb_1|}\right)^{\frac{1}{q}}|\Ma_1\Mb_2-\Ma_2\Mb_1|^{\frac{1}{2}}\|(\mathcal{L}^{M_2}f)(\uu)\|_{L^q(\mathbb{R})}.
\end{align}
For $\theta\rightarrow 1$ and $\gamma\rightarrow 1$ equation \eqref{Rguanlei2009uncertainty_Renyi_UP_2LCTD} reduces to \eqref{R2009Guanlie-Xiaotong-XiaogangGEntropyUPLCT}. Guanlei et al. \cite{guanlei2010uncertainty} obtained the UPs for the LCT with a bound sharper than that appearing in \eqref{R2009Guanlie-Xiaotong-XiaogangHeisenbergUPin2LCTDomain} for a complex-valued function of the form $f(\tee)=s(\tee)e^{i\phi(\tee)},$ where $s(\tee)$ and $\phi(\tee)$ are real-valued and differentiable functions and $s(\tee)\in L^1_{\mathbb{R}}(\mathbb{R})\cap L^2_{\mathbb{R}}(\mathbb{R})$ and with the assumption that $\frac{\Ma_1}{\Mb_1}\neq \frac{\Ma_2}{\Mb_2}$ as
\begin{align}\label{R2010_guanlei-Xiaotong-Xiaogang_LCT_UP2LCT-Domain-SpecialComplex-Functions}
&\left(\int_{\mathbb{R}}|\uu|^2\left|\left(\mathcal{L}^{M_1}f\right)(\uu)\right|^2d\uu\right)\left(\int_{\mathbb{R}}|\vv|^2\left|\left(\mathcal{L}^{M_2}f\right)(\vv)\right|^2d\vv\right)\notag\\
&\hspace{3cm}\geq\frac{|\Ma_1\Mb_2-\Ma_2\Mb_1|^2}{4}+\left(\Ma_1\Ma_2\Delta\tee^2+\Mb_1\Mb_2\Delta\oomega_s^2+\Mb_1\Mb_2K_1+(\Ma_1\Mb_2+\Ma_2\Mb_1)K_2\right)^2
\end{align}
where $K_1=\int_{\mathbb{R}}\left[\phi'(\tee)s(\tee)\right]^2d\tee,$ $K_2=\int_{\mathbb{R}}\tee(s(\tee))^2\phi'(\tee)d\tee$ and with $\tee_0=\int_{\mathbb{R}}\tee|s(\tee)|^2d\tee$ and $\oomega_{s,0}=\int_{\mathbb{R}}\oomega|\hat{s}(\oomega)|^2d\oomega,$ $\Delta\tee^2=\int_{\mathbb{R}}(\tee-\tee_0)^2|s(\tee)|^2d\tee$ and $\Delta\oomega_s^2=\int_{\mathbb{R}}(\oomega-\oomega_0^2)^2|\hat{s}(\oomega)|^2d\oomega.$ For the function $f(\tee)=\left(\frac{2}{\sqrt{\pi}}\right)^{\frac{1}{2}}\tee e^{-\frac{\tee^2}{2}}e^{i\frac{\tee^2}{2}},$ it can be shown that $\Delta\tee^2=\Delta_s^2=K_1=K_2=\frac{3}{2}.$
The equality in \eqref{R2010_guanlei-Xiaotong-Xiaogang_LCT_UP2LCT-Domain-SpecialComplex-Functions} is obtained for the function $f(\tee)=k_2e^{-i\left[\frac{k_1 a_2-a_1}{k_1 b_2-b_1}\right]\frac{\tee^2}{2}},$ where $k_1,k_2$ are arbitrary complex constants. 
%\textcolor{violet}{This is important: Motivated from \cite{guanlei2009new}(\eqref{R2009Guanlie-Xiaotong-XiaogangPitt'sInequalityin2LCTDomain}, \eqref{R2009Guanlie-Xiaotong-XiaogangLogatithmicUPin2LCTDomain}, \eqref{R2009Guanlie-Xiaotong-XiaogangGEntropyUPLCT} and \eqref{R2009Guanlie-Xiaotong-XiaogangHeisenbergUPin2LCTDomain}, this paper is mentioned above),  \cite{kundu2022uncertainty} In this paper we find some UPs like Pitt's inequality,Logarithmic Uncertainty inequality, Entropy UP using  Hausdorff Young inequality, that will be used in the explanation of the paper of Kundu et al.\cite{kundu2022uncertainty} in discussing the QLCT.} \\

\section{Quadratic Phase Fourier Transform (QPFT)}

Recently,  Castro et al.\cite{castro2018new} developed a new integral transform called the QPFT (QPFT), which generalizes several transforms like FT, FrFT, and LCT. For the function $f \in L^1(\mathbb{R})$ or $f \in L^2(\mathbb{R})$ and given parameter $\Lambda=(\LmA,\LmB,\LmC,\LmD,\LmE),~\LmB\neq 0,$ the QPFT is defined by 
\begin{align}\label{RDefinitionQPFT}
(Q_{\Lambda} f)(\oomega):=\int_{\mathbb{R}} \mathcal{K}^i_{\Lambda}(\xx,\oomega)f(\xx)d\xx,
\end{align}
where the kernel is given by 
\begin{align}\label{RDefinitionofKernel_QPFT}
\mathcal{K}^i_{\Lambda}(\xx,\oomega)=\frac{1}{\sqrt{2\pi}}e^{-i(\LmA\xx^2+\LmB \xx \oomega+\LmC \oomega^2+\LmD \xx+\LmE \oomega)},
\end{align}
they proved the Riemann-Lebesgue lemma, inversion formula,
Plancherel's extension theorem and also the Parseval identity given by 
$$
\langle Q_{\Lambda} f, Q_{\Lambda} g\rangle_{L^2(\mathbb{R})}=\frac{1}{|\LmB|}\langle f, g\rangle_{L^2(\mathbb{R})},~f, g \in L^2(\mathbb{R})$$
denotes the usual inner product in $L^2(\mathbb{R})$ given by $\langle f, g\rangle_{L^2(\mathbb{R})}:=$ $\int_{\mathbb{R}} f(\xx) \overline{g(\xx)} d\xx$. In the special case of $f=g$, we then have
$$
\|Q_{\Lambda} f\|_{L^2(\mathbb{R})}^2=\frac{1}{|\LmB|}\|f\|_{L^2(\mathbb{R})}^2 .
$$
They also gave several new convolutions followed by the corresponding convolution theorem, which, in particular cases, was the convolution theorem for the FT, FrFT, and LCT. They also obtained the Young's convolution inequalities associated with these convolutions. Shah et al.\cite{shah2021uncertainty} proved that the transform \eqref{RDefinitionQPFT} is linear and has the translation formula $Q_{\Lambda}\left[f(\tee-k)\right](\oomega)=e^{-i(\LmA k^2+\LmB k\oomega+\LmD k)}Q_{\Lambda}\left[e^{-2i \LmA k\tee}f(\tee)\right](\oomega),$  modulation formula $Q_{\Lambda}\left[e^{i\mu\tee}f(\tee)\right](\oomega)=e^{\frac{i}{\LmB}(\mu^2+2\mu \LmB \oomega+\LmE\mu \LmB)}\left(Q_{\Lambda}f\right)\left(\oomega+\frac{\mu}{\LmB}\right)$ and scaling formula 
$$Q_{\Lambda}\left[f(\mu\tee)\right](\oomega)=\frac{1}{\mu}e^{-i \LmC \oomega^2\left(1-\frac{1}{\mu^2}\right)-i \LmE \oomega\left(1-\frac{1}{\mu}\right)}\left(Q_{\Lambda}F\right)\left(\frac{\oomega}{\mu}\right),$$ 
where $F(\tee)=e^{i\left(\LmA \tee^2\left(1-\frac{1}{\mu^2}\right)+\LmD \tee\left(1-\frac{1}{\mu}\right)\right)}f(\tee).$ They also gave the following convolution theorem
\begin{align}\label{RShahConvolutionTheoremQPFT}
\left(Q_{\Lambda}\left[f \otimes_{\Lambda} g\right]\right)(\oomega)=e^{i\left(\LmC \oomega^2+\LmE \oomega\right)} \left(Q_{\Lambda}f\right)(\oomega)\left( Q_{\Lambda}g\right)(\oomega),
\end{align} 
with the convolution $\otimes_{\Lambda}$ given by
\begin{align}\label{RShahConvolutionQPFT}
\left(f \otimes_{\Lambda} g\right)(\tee)=\frac{1}{\sqrt{2 x}} \int_{\mathbb{R}} f(\xx) g(\tee-\xx) e^{2 i \LmA \xx(\tee-\xx)} d\xx.
\end{align} 
Kumar et al.\cite{kumar2022quadratic} proved the continuity of the linear map $Q_{\LmD, \LmE}^{\LmA, \LmB, \LmC}: \mathcal{S}(\mathbb{R}) \longmapsto \mathcal{S}_{\Omega}(\mathbb{R}),$
where $Q_{\LmD, \LmE}^{\LmA, \LmB, \LmC}$ is defined as 
\begin{align}\label{R2022_Kumar-Pradhan_QPFT-of-TemperedDistribution-DefinitionQPFT}
\left(Q_{\LmD, \LmE}^{\LmA, \LmB, \LmC}f\right)(\yy)=\sqrt{\frac{\LmB}{2\pi i}}\int_{\mathbb{R}}e^{-i(\LmA \xx^2+\LmB \xx \yy+\LmC \yy^2+\LmD \xx+\LmE \yy)}f(\xx)d\xx
\end{align}
where $\mathcal{S}(\mathbb{R})$ is a Schwartz space and $\mathcal{S}_{\Omega}(\mathbb{R})$ is the Schwartz-type space  which contains a complex-valued $C^{\infty}$-function $\varphi$ on $\mathbb{R}$, such that for every pair $(\alpha, \beta)$ of non-negative integers, it satisfies
$$
\Gamma_{\Omega}^{\alpha, \beta}(\varphi)=\sup _{\xx \in \mathbb{R}}\left|\xx^\alpha \Delta_{\xx, a}^\beta \varphi(\xx)\right|<\infty,
$$
where
$$
\Delta_{\xx, \LmA}=\left(\frac{\mathrm{d}}{\mathrm{d} x}-2 i \LmA \xx\right)
$$
and $\Delta_{\xx, \LmA}^\beta=\sum_{j=0}^\beta p_j^\beta(\xx) D_\xx^{\beta-j} ; p_j^\beta$ denotes the polynomial of order $j$ with $p_0^\beta(\xx)=1$. They also extended the definition of the QPFT to the space of all tempered distribution $\mathcal{S}^{\prime}(\mathbb{R}),$ given by 
$$
\left\langle Q_{\LmD, \LmE}^{\LmA, \LmB, \LmC} f, \varphi\right\rangle:=\left\langle f, Q_{\LmD, \LmE}^{\LmA, \LmB, \LmC} \varphi\right\rangle
$$
for all $\varphi \in \mathcal{S}(\mathbb{R}),$ where $\mathcal{S}^{\prime}(\mathbb{R})$ consists of all continuous linear functions $f: \mathcal{S}(\mathbb{R}) \rightarrow \mathbb{C}$ and proved that QPFT is a continuous linear map from  $\mathcal{S}^{\prime}(\mathbb{R})$ onto itself. Prasad et al. \cite{prasad2020quadratic} proved the convolution theorem for the QPFT as
\begin{align}\label{R2018_Prasad-Sharma_MMAS_QPWT_QP-ConvolutionTheorem}
\left(Q_{\LmD, \LmE}^{\LmA, \LmB, \LmC}\left(f\star_{\Omega_{\LmD}^{\LmA}}\psi\right)\right)(\oomega)=\sqrt{\frac{2\pi i}{\LmB}}e^{-i \LmC \oomega^2-i \LmE \oomega}\left(Q_{\LmD, \LmE}^{\LmA, \LmB, \LmC}f\right)(\oomega)\left(Q_{\LmD, \LmE}^{\LmA, \LmB, \LmC}\left(e^{-i \LmA (\cdot)^2-i \LmD (\cdot)}\psi(\cdot)\right)\right)(\oomega),
\end{align}
for $f,\psi\in \mathcal{S}_\Omega(\mathbb{R}),$ where the convolution $\star_{\Omega_{\LmD}^{\LmA}}$ is defined as 
\begin{align}\label{R2018_Prasad-Sharma_MMAS_QPWT_QP-Convolution}
(f\star_{\Omega_{\LmD}^{\LmA}}\psi)(\xx)=\int_{\mathbb{R}}f(\yy)\psi(\xx-\yy)e^{-i \LmA (\xx^2-\yy^2)-i \LmD (\xx-\yy)}d\yy.
\end{align} 
Shah et al. \cite{shah2022sampling} defined a quadratic-phase convolution of $f,g\in L^2(\mathbb{R})$ as 
\begin{align}\label{R2023_SIVP_Shah-Tantary_QPFT_Convolution}
\left(f\star'_{\Lambda}g\right)(\zz)=\frac{1}{\sqrt{2\pi}}\int_{\mathbb{R}}f(\tee)g(\zz-\tee)e^{-2i\LmA\tee(\tee-\zz)}d\tee
\end{align}
and proved that the QPFT \eqref{RDefinitionQPFT} satisfies the formula $\left(Q_{\Lambda}(f\star'_{\Lambda}g)\right)(\oomega)=e^{i\left(\LmC\oomega^2+\LmE \oomega\right)}\left(Q_{\Lambda} f\right)(\oomega)\left(Q_{\Lambda} g\right)(\oomega).$ They also defined the correlation of $f,g\in L^2(\mathbb{R})$ as $\left(f\circ'_{\Lambda}g\right)(\zz)=\frac{1}{\sqrt{2\pi}}\int_{\mathbb{R}}f(\tee)\overline{g(\tee-\zz)}e^{2i\LmA \zz(\zz-\tee)}d\tee$ and the corresponding correlation formula as $\left(Q_{\Lambda}(f \circ'_{\Lambda} g)\right)(\oomega)=e^{i\left(\LmC\oomega^2+\LmE\oomega\right)}\left(Q_{\Lambda} f\right)(\oomega)\overline{\left(Q_{\Lambda} g\right)(\oomega)}.$ Lone et al. \cite{lone2022weighted} proposed the convolutions associated with the QPFT followed by the corresponding convolution theorem. These convolution operations enjoy both the commutative and the associative properties. We mention below such convolutions along with the associated convolution theorems.
\begin{enumerate}
\item For any pair of square integrable functions $f_1$ and $f_2$, the quadratic-phase convolution is denoted by $\oplus_\Lambda$ and is defined as
\begin{align}\label{RLoneConvolutionQPFT3}
\left(f_1 \oplus_\Lambda f_2\right)(\zz)=|\LmB| e^{i \LmA\zz^2}\left(\check{f}_1(\LmB~\cdot) \star \check{f}_2(\cdot)\right)(\zz) e^{i \LmD\zz},
\end{align}
where $\breve{f}(\tee)=f(\tee) e^{-i \LmA \tee^2}$ and $\star$ represents the classical convolution operator.
The corresponding convolution theorem is given by
\begin{align}\label{RLoneConvolutionTheoremQPFT3}
\left(Q_{\Lambda}[f_1 \oplus_\Lambda f_2]\right)(\oomega)=e^{i\left\{\LmC\left(\frac{\oomega}{\LmB}\right\}^2+\LmE\left(\frac{\oomega}{\LmB}\right)\right)}\left(Q_{\Lambda}F_1\right)\left(\frac{\oomega}{\LmB}\right) \left(Q_{\Lambda}F_2\right)\left(\frac{\oomega}{\LmB}\right),
\end{align} 
where $F_1(\tee)=e^{i\LmD\tee}f_1(\tee)$ and $F_2(\tee)=e^{i\LmD\tee}f_2(\tee).$
\item For a set of parameters $\Lambda=(\LmA, \LmB, \LmC, \LmD, \LmE),~\LmB\neq 0$ and any pair of functions $f_1, f_2 \in L^2(\mathbb{R})$, the weighted quadratic-phase convolution is denoted by $\ast_\Lambda$ and is defined as
\begin{align}\label{RLoneConvolutionQPFT4}
\left(f_1 \ast_\Lambda f_2\right)(\zz)=\frac{1}{\sqrt{2 \pi}} \int_{\mathbb{R}} f_1(\tee) f_2(\zz-\tee) e^{i\LmA\left(\zz^2-\tee^2\right)} \mathrm{d\tee},
\end{align}
with the corresponding convolution theorem as 
\begin{align}\label{RLoneConvolutionTheoremQPFT4}
\left(Q_{\Lambda}\left[f_1 \ast_\Lambda f_2\right]\right)(\oomega)=e^{i\left(\LmC\oomega^2+\LmE\oomega\right)}\left(Q_\Lambda f_1\right)(\oomega) \left(Q_\Lambda f_2\right)(\oomega).
\end{align}
\end{enumerate}
The convolution theorem in equations \eqref{RLoneConvolutionTheoremQPFT3} and \eqref{RLoneConvolutionTheoremQPFT4} are such that they contain an extra chirp. So, Lone et al. \cite{lone2022weighted},\cite{lone2023non} also gave some convolution in the context of QPFT such that the corresponding convolution theorem has no chirp multiplier. In what follows, we mention such convolution followed by the associated convolution theorem 
\begin{enumerate}
\item For a set of parameters $\Lambda=(\LmA,\LmB, \LmC, \LmD, \LmE)$ and any pair of functions $f_1, f_2 \in L^2(\mathbb{R})$, the chirp-free weighted convolution associated with the QPFT is denoted by $\star_\Lambda$ and is defined as
\begin{align}\label{RLoneConvolutionQPFT1}
\left(f_1 \star_\Lambda f_2\right)(\zz)=\frac{1}{\sqrt{2 \pi}} \int_{\mathbb{R}} f_1(\tee) f_2(\zz-\tee) e^{i\LmA\left(\zz^2-\tee^2\right)} \mathrm{d\tee},
\end{align}
with the corresponding convolution theorem as 
\begin{align}\label{RLoneConvolutionTheoremQPFT1}
\left(Q_{\Lambda}\left[f_1 \star_\Lambda f_2\right]\right)(\oomega)=\left(Q_\Lambda f_1\right)(\oomega) \left(\mathscr{F} F_2\right)(\LmB \oomega),
\end{align}
where $\mathscr{F}\left[F_2\right]$ represents the FT of $F_2(\tee)=e^{-i \LmD \tee} f_2(\tee)$. The convolution operation given in \eqref{RLoneConvolutionQPFT1} do not satisfy commutative nor associative properties. The convolution given below is a chirp free convolution operation which satisfy commutative as well as associative properties.
\item For a parametric set $\Lambda=(\LmA, \LmB, \LmC, \LmD, \LmE), \LmB \neq 0$, the chirp-free quadratic-phase convolution of functions $f_1, f_2 \in L^2(\mathbb{R})$ is denoted by $\odot_\Lambda$ and is given by
\begin{align}\label{RLoneConvolutionQPFT2}
\left(f_1 \odot_\Lambda f_2\right)(\zz)=\frac{1}{\sqrt{2\pi}} \int_{\mathbb{R}} f_1(\tee) f_2(\sqrt{2} \zz-\tee) e^{-2i\LmA\left(\frac{\zz}{\sqrt{2}}-\tee\right)^2} \mathrm{d\tee},
\end{align}
with the corresponding convolution theorem as 
\begin{align}\label{RLoneConvolutionTheoremQPFT2}
\left(Q_{\Lambda}\left[f_1 \odot_\Lambda f_2\right]\right)(\oomega)=\frac{1}{\sqrt{2}} \left(Q_{{\Lambda^{\prime}}}f_1\right)\left(\frac{\oomega}{\sqrt{2}}\right) \left(Q_{\Lambda^{\prime}}f_2\right)\left(\frac{\oomega}{\sqrt{2}}\right), \quad \Lambda^{\prime}=\left(\LmA, \LmB, \LmC, \frac{\LmD}{\sqrt{2}},\frac{\LmE}{\sqrt{2}}\right)
\end{align}
\end{enumerate}
Lone et al. \cite{lone2022shift} defined the discrete-time QPFT of any sequence $u(n) \in \ell^2(\mathbb{Z}),$ denoted by $D_{\Omega}[u(n)],$ as
\begin{align}\label{RLoneDefinitionDTQPFT}
\mathcal{D}_{\Omega}[u(n)](\oomega)=\frac{1}{\sqrt{2 \pi}} \sum_{n \in \mathbb{Z}} u(n) e^{-i\left(\LmA n^2+\LmB n \oomega+\LmC \oomega^2+\LmD n+\LmE \oomega\right)},
\end{align}
which generalized the discrete-time FrFT \cite{bhandari2011shift}, discrete-time FrFT \cite{bhandari2011shift} and discrete-time FT as \cite{debnath2017lecture}. The inversion formula corresponding to \eqref{RLoneDefinitionDTQPFT} is given by
\begin{align}\label{RLoneDefinitionInverseofDTQPFT}
u(n)=\int_0^{2 \pi / \LmB} D_{\Omega}[u(n)](\oomega) \mathcal{K}_{\Omega}(n, \oomega) d \oomega,
\end{align}
where $\mathcal{K}_{\Omega}(n, \oomega)=e^{-i(\LmA n^2+\LmB n \oomega+\LmC \oomega^2+\LmD n+\LmE \oomega)}.$
%Definition $2.3$ allows us to make the following comments:\\\\(i). For $\Omega=(-\LmA / 2 \LmB, 1 / \LmB,-\LmD / 2 \LmB, 0,0)$ and multiplying the right side of (2.8) by $1 / \sqrt{i \LmB}$, Definition $2.3$ yields the discrete-time LCT \cite{huo2020nonuniform}:
%$$
%\mathcal{D}_{\Omega}[u(n)](\omega)=\frac{1}{\sqrt{2 \pi i \LmB}} \sum_{n \in \mathbb{Z}} u(n) e^{i\left(\LmA n^2-2 n \omega+\LmD \omega^2\right) / 2 \LmB} .
%$$
%(ii). For $\Omega=(-\cot \theta / 2, \csc \theta,-\cot \theta / 2,0,0), \theta \neq n \pi$ and multiplying the right side of (2.8) by $\sqrt{1-i \cot \theta}$, Definition $2.3$ boils down to the discrete-time FrFT \cite{bhandari2011shift}:
%$$
%D_{\Omega}[u(n)](\omega)=\sqrt{\frac{1-i \cot \theta}{2 \pi}} \sum_{n \in \mathbb{Z}} u(n) e^{l\left(n^2+\omega^2\right) \cot \theta / 2-i \omega n \csc \theta} .
%$$
%(iii). For $\Omega=(0,1,0,0,0)$, Definition $2.3$ reduces to the conventional discrete-time FT as \cite{debnath2017lecture}:
%$$
%D_{\Omega}[u(n)](\omega)=\frac{1}{\sqrt{2 \pi}} \sum_{n \in \mathbb{Z}} u(n) e^{-i n \omega} .
%$$
They also introduced the notion of discrete convolution and semi-discrete convolution in the context of QPFTs for $u(n),v(n)\in\ell^2(\mathbb{Z})$ and $\phi\in L^2{(\mathbb{R})}$, denoted by $\oplus_D$ and $\circledast_S$ respectively, as
\begin{align}\label{RLoneDefinitionDiscreteQPConvolution}
u(n) \oplus_D v(n)=\frac{1}{\sqrt{2 \pi}} \sum_{n \in \mathbb{Z}} u(n) v(k-n) e^{-2 \pi i \LmA n(n-k)}
\end{align}
and 
\begin{align}\label{RDefinitionSemiDiscreteTQPConvolution}
u(n) \circledast_S \phi(\tee)=\frac{1}{\sqrt{2 \pi}} \sum_{n \in \mathbb{Z}} s(n) \phi(\tee-n) e^{-2 \pi i \LmA n(n-\tee)} .
\end{align}
The operation of convolution $\oplus_D$ defined in \eqref{RLoneDefinitionDiscreteQPConvolution} is commutative and associative. Moreover, the convolution theorem for \eqref{RLoneDefinitionDTQPFT} associated with \eqref{RLoneDefinitionDiscreteQPConvolution} and \eqref{RDefinitionSemiDiscreteTQPConvolution} are given respectively by
\begin{align}\label{RDefinitionDistreteQPConvolutionTheorem}
\mathcal{D}_{\Omega}\left[u(n) \oplus_D v(n)\right](\omega)=e^{i\left(\LmC \omega^2+\LmE \omega\right)} \mathcal{D}_{\Omega}[u(n)](\omega) \mathcal{D}_{\Omega}[v(n)](\omega)
\end{align}
and 
\begin{align}\label{RDefinitionSemiDiscreteTQPConvolutionTheorem}
\mathcal{L}_{\Omega}[u(n) \circledast_S \phi](\oomega)=e^{i\left(\LmC \oomega^2+\LmE \oomega\right)} \mathcal{D}_{\Omega}[u(n)](\oomega) \mathcal{L}_{\Omega}[\phi](\oomega).
\end{align}
Moreover, $\left|\mathcal{D}_{\Omega}\left[ u(n)\right](\oomega)\right|$ is periodic with periodicity $\frac{2\pi}{\LmB}$. Authors in \cite{lone2023non} have introduced a pair of novel discrete and semi-discrete convolution structures in the framework of QPFT for $x(n),y(n)\in\ell^2(\mathbb{Z})$ and $\phi\in L^2(\mathbb{R}),$ denoted  respectively by $\smallstar_D$ and $\smallstar_S$ and are defined by
\begin{align}\label{RLoneNewDiscreteCOnvolutionQQPT}
x(n)\smallstar_D y(n)=\frac{1}{\sqrt{2 \pi}} \sum_{n \in \mathrm{Z}} x(n) y(k-n) e^{i\LmA\left(k^2-n^2\right)}
\end{align}
and 
\begin{align}\label{RLoneNewSemiDiscreteCOnvolutionQQPT}
x(n)\smallstar_S\phi(\tee)=\frac{1}{\sqrt{2\pi}}\sum_{n\in\mathbb{Z}}x(n)\phi(\tee-n)e^{i\LmA(\tee^2-n^2)}.
\end{align}
The operation of convolution $\smallstar_D$ defined in \eqref{RLoneNewDiscreteCOnvolutionQQPT} is neither commutative nor associative. Moreover, the convolution theorem associated with \eqref{RLoneNewDiscreteCOnvolutionQQPT} and \eqref{RLoneNewSemiDiscreteCOnvolutionQQPT} are given respectively by 
\begin{align}\label{RLoneNewDiscreteCOnvolutionTheoremQQPT}
\mathbf{D}_{\Lambda}\left[x(n) \smallstar_D y(n)\right](\oomega)=\mathbf{D}_{\Lambda}[x(n)](\oomega) \mathscr{F}\left[y_{\Lambda}(n)\right](\LmB \oomega),
\end{align}
and 
\begin{align}\label{RLoneNewSemiDiscreteCOnvolutionTheoremQQPT}
\mathscr{I}_{\mathrm{N}}\left[x(n) \smallstar_S \phi\right](\oomega)=\mathbf{D}_{\Lambda}[x(n)](\oomega) \mathscr{F}[\Phi](\LmB \oomega),
\end{align}
where $\mathscr{F}[y_{\Lambda}(n)]$ and $\mathscr{F}[\Phi]$ represent respectively the discrete-time FT  of $y_{\Lambda}(n)=e^{-i \LmD n} y(n)$ and  FT of $\Phi(\tee)=e^{-i\LmD\tee}\phi(\tee).$ Observe that the convolution results given in \eqref{RLoneNewDiscreteCOnvolutionTheoremQQPT} and \eqref{RLoneNewSemiDiscreteCOnvolutionTheoremQQPT} does not contain any chirp multipliers as was not the case for the convolution result given in \eqref{RDefinitionDistreteQPConvolutionTheorem} and \eqref{RDefinitionSemiDiscreteTQPConvolutionTheorem}.
\\

%Theorem 2.8. For any sequence $u(n) \in \ell^2(\mathbb{Z})$ and a function $\phi \in L^2(\mathbb{R})$, we have
%\begin{align}\label{RDefinitionSemiDiscreteTQPConvolutionTheorem}
%\mathcal{L}_{\Omega}[u(n) \circledast s \phi](\omega)=e^{i\left(\LmC \omega^2+\LmE \omega\right)} \mathcal{D}_{\Omega}[u(n)](\omega) \mathcal{L}_{\Omega}[\phi](\omega) .
%\end{align}
%Moreover, $\left|\mathcal{D}_{\Omega}\right| u(n)|(\omega)|$ is periodic with periodicity $2 \pi /$ \LmB.

Shanon's sampling theorem is one of the remarkable, profound, and elegant concepts of digital signal processing, which serves as a bridge between analog and digital signals. The theorem asserts that if a bandlimited signal has a range of $n$ frequencies, it can be perfectly reconstructed by taking $2 n$ evenly spaced samples. If we take any additional samples, it will turn out to be redundant, while as if we take fewer samples, we will lose some quality. This is a fundamental result in the field of information theory, in particular, telecommunications and signal processing \cite{debnath2017lecture}. The Shanon's sampling theorem for the QPFT serves as a unified sampling theorem for the wide class of integral transforms. Authors in \cite{lone2022shift}, \cite{gupta2023quadratic} and \cite{shah2022sampling} proved that a bandlimited function $f(\tee) \in L^2(\mathbb{R})$ with bandwidth $\sigma$ in the quadratic-phase Fourier domains can be reconstructed from its equally spaced samples $f(n \tau), n \in \mathbb{Z}$ by the formula
\begin{align}\label{RLeonSamplingTheoremQPFT}
f(t)=e^{i\left(\LmA \tee^2+\LmD \tee\right)} \sum_{n \in \mathbb{Z}} e^{-i\left(\LmA(n \tau)^2+\LmD(n \tau)\right)} \frac{\sin (\LmB \sigma(\tee-n \tau))}{\LmB \sigma(\tee-n \tau)} f(n \tau),
\end{align}
where $\tau=\pi / \LmB \sigma$ is the sampling period.
\\ 
A closed subspace $\mathcal{V}$ of $L^2(\mathbb{R})$ is said to be shift-invariant if $f(\tee-n) \in \mathcal{V}$, whenever $f \in \mathcal{V}$. A closed shift-invariant subspace $\mathcal{V} \subset L^2(\mathbb{R})$ is said to be generated by $\phi \in L^2(\mathbb{R})$ if $\mathcal{V}=\overline{\operatorname{span}}\{\phi(\tee-n): n \in \mathbb{Z}\}$. Lone et al. \cite{lone2022shift} introduced the chirp-modulated shift-invariant subspace of $L^2(\mathbb{R})$ in the quadratic-phase Fourier domains $\mathcal{V}_{\Omega}(\phi)$ for any sequence $u(n) \in \ell^2(\mathbb{Z})$ and $\phi(\tee) \in L^2(\mathbb{R})$ as 
\begin{align}\label{RLoneShiftModulatedInvariantSpaceQPFT}
\mathcal{V}_{\Omega}(\phi)=\text { closure }\left\{f \in L^2(\mathbb{R}): f(\zz)=\left(u(n) \circledast_S \phi(\tee)\right)(\zz)\right\},
\end{align}
where $\circledast_S$ is the semi-discrete convolution given by \eqref{RDefinitionSemiDiscreteTQPConvolution} and established that the necessary and sufficient condition for $\left\{e^{-2 i \LmA n(n-\tee)} \phi(\tee-n)\right\}$ to be a Riesz basis for $\mathcal{V}_{\Omega}(\phi)$ is that there exists two positive constants $\Gamma_1$ and $\Gamma_2$ such that
\begin{align}
\Gamma_1 \leq \sum_{n \in \mathbb{Z}}\left|\mathcal{L}_{\Omega}[\phi](\oomega+n)\right|^2 \leq \Gamma_2, \quad \forall \oomega \in[0,2 \pi / \LmB]
\end{align}
For any sequence $u(n) \in \ell^2(\mathbb{Z})$ and a collection of functions $\phi_{\ell}(\tee) \in L^2(\mathbb{R}), 1 \leqslant \ell \leqslant M$, the generalized SISs of $L^2(\mathbb{R})$ with respect to a parametric set $\Lambda=(\LmA, \LmB,\LmC, \LmD, \LmE), \LmB>0$ is denoted as $\gamma_{\Lambda}\left(\phi_1, \ldots, \phi_M\right)$ and is defined by
\begin{align}
\gamma_{\Lambda}\left(\phi_1, \ldots, \phi_M\right)=\text {closure}\left\{f \in L^2(\mathbb{R}): f(\tee)=\sum_{\ell=1}^M \sum_{n\in\mathbb{Z}}u_l(n) \smallstar_S \phi_l(\tee)\right\},
\end{align}\label{RLoneNewShiftModulatedInvariantSpaceQPFT}
where $\smallstar_S$ denoted the semi-discrete convolution as defined by \eqref{RLoneNewSemiDiscreteCOnvolutionQQPT}. Let $\gamma_{\Lambda}\left(\phi_1, \ldots, \phi_M\right) \subset L^2(\mathbb{R})$ be the generalized SISs associated with the QPFT. Then, the family $\quad\left\{e^{\LmA n\left(\tee^2-n^2\right)} \phi_\ell(\tee-n): 1 \leqslant \ell \leqslant M\right\}$ forms a Riesz basis for $\gamma_{\Lambda}\left(\phi_1, \ldots, \phi_M\right)$ if and only if there exists a pair of constants $\Gamma_1, \Gamma_2>0$ such that 
\begin{align}
\Gamma_1 \leqslant \sum_{\ell=1}^M \sum_{n\in\mathbb{Z}} \left|\mathscr{F}\left[\Phi_{\ell}\right](\LmB \oomega+n)\right|^2 \leqslant \Gamma_2, \quad 1 \leqslant \ell \leqslant M, \forall \oomega \in[0,2 \pi / \LmB] \text {, }
\end{align}
where $\mathscr{F}[\Phi_{\ell}]$ represents the FT of $\Phi(t)=e^{-i\LmD\tee}\phi_{\ell}(\tee).$
%\cite{shah2022sampling}
\\\\
\noindent\textbf{Uncertainty Principles:} Using the relation between the QPFT and the classical FT and the Heisenberg's and logarithmic UP for the classical FT, Shah et al.\cite{shah2021uncertainty} obtained the Heisenberg's  UP for the QPFT as
\begin{align}\label{R2021_Shah-Nisar-Lone-Tantary_MMAS_QPFT_Heisenberg-UP}
\left(\int_{\mathbb{R}}|\tee-\tee_0|^2|f(\tee)|^2d\tee\right) \left(\int_{\mathbb{R}}|\oomega-\oomega_0|^2|(Q_{\Lambda}f)(\oomega)|^2d\tee\right)\geq \frac{1}{4\LmB^2}\|f\|^4_{L^2(\mathbb{R})}
\end{align}
for $f\in L^2(\mathbb{R})$ and logarithmic UP as
\begin{align}\label{R2021_Shah-Nisar-Lone-Tantary_MMAS_QPFT_Logarithmic-UP}
\int_{\mathbb{R}}\ln|\tee||f(\tee)|^2d\tee+|\LmB|\int_{\mathbb{R}}\ln|\oomega|\left|(Q_{\Lambda}f)(\oomega)\right|^2d\oomega\geq \left(\frac{\Gamma'\left(\frac{1}{4}\right)}{\Gamma\left(\frac{1}{4}\right)}-\ln\pi-\ln|\LmB|\right)\int_{\mathbb{R}}|f(\tee)|^2d\tee
\end{align}
for $f\in \mathcal{S}(\mathbb{R}).$ Adopting the same method, authors obtained the local UP for $f\in L^2(\mathbb{R})$ and $0<\alpha<\frac{1}{2}$ as
\begin{align}\label{R2021_Shah-Nisar-Lone-Tantary_MMAS_QPFT_Local-UP}
\int_{\mathbb{R}}|\tee|^{2\alpha}|f(\tee)|^2d\tee\geq \frac{|\LmB|}{C_\alpha|E_0|^{2\alpha}}\int_{E_0}\left|(Q_{\Lambda}f)(\oomega)\right|^2d\oomega.
\end{align}
Furthermore, they also obtained the Donoho-Stark's UP for the QPFT, which says that if $E_1,E_2$ be any two measurable subsets of $\mathbb{R}$ such that $f$ is $\epsilon_{E_1}-$concentrated on $E_1$ and $Q_{\Lambda}f$ is $\epsilon_{E_2}-$concentrated on $E_2,$ then $|E_1||E_2|\geq \frac{1}{|\LmB|}(1-\epsilon_{E_1}-\epsilon_{E_2})^2,$ where $|E_l|,$ denote the Lebesgue measure of $E_l,~l=1,2.$ Several other concentration-based UPs like Nazarov's and  Amrein-Berthier-Benedick's  UPs and Hardy's and Beurling's UPs are also established in \cite{shah2021uncertainty}.
 
\section{Quaternion Fourier Transform (QFT)}
In classical theory of Fourier analysis, it is known that if $f\in L^1_{\mathbb{C}}(\mathbb{R}),$ then the FT $\hat{f}$ given by $\hat{f}(\xi)=\int_{\mathbb{R}}f(t)e^{-it\xi}dt$ is well defined. Moreover, if $\hat{f}$ is integrable then $f$ is equivalent to a continuous function and is such that $f(\tee)=\frac{1}{2\pi}\int_{\mathbb{R}}\hat{f}(\xxi)e^{i\tee\xxi}d\xxi.$ There is a natural way of defining the FT of function in $L^2_{\mathbb{C}}(\mathbb{R})$ using the density of $L^1_{\mathbb{C}}(\mathbb{R})\cap L^2_{\mathbb{C}}(\mathbb{R}).$ In fact, the FT is not just isometric but also onto from $L^2_{\mathbb{C}}(\mathbb{R})$ to $L^2_{\mathbb{C}}(\mathbb{R}).$ So, there is a natural question: whether there is an analogous theory of FT for the quaternion-valued function. Because of the non-commutativity  of quaternion multiplication, there are at least three different types of quaternion FT for the quaternion-valued function defined on $\mathbb{R}^2,$ namely  right-sided, two-sided, and left-sided QFT, which are respectively given as 
\begin{align}\label{R2019Chen-KouDefinitionRsQFT}
\left(\mathcal{F}_rf\right)(\oomega_1,\oomega_2)=\frac{1}{2\pi}\int_{\mathbb{R}^2}f(\xx_1,\xx_2)e^{-i\oomega_1\xx_1}e^{-j\oomega_2\xx_2}d\xx_1d\xx_2,
\end{align} 
\begin{align}\label{R2019Chen-KouDefinitionTsQFT}
\left(\mathcal{F}_sf\right)(\oomega_1,\oomega_2)=\frac{1}{2\pi}\int_{\mathbb{R}^2}e^{-i\oomega_1\xx_1}f(\xx_1,\xx_2)e^{-j\oomega_2\xx_2}d\xx_1d\xx_2
\end{align} 
and
\begin{align}\label{R2019Chen-KouDefinitionLsQFT}
\left(\mathcal{F}_lf\right)(\oomega_1,\oomega_2)=\frac{1}{2\pi}\int_{\mathbb{R}^2}e^{-i\oomega_1\xx_1}e^{-j\oomega_2\xx_2}f(\xx_1,\xx_2)d\xx_1d\xx_2,
\end{align} 
where $f\in L^1_{\mathbb{H}}(\mathbb{R}^2).$ Ernst\cite{ernst1987principles} and Delucs\cite{delsuc1988spectral} first applied quaternions to FT in the late 1980s. In 1992, Ell\cite{ell1992hypercomplex} analyzed two-sided QFT and applied it to color image processing. In \cite{ell2006hypercomplex}, Ell et al. gave a detailed interpretation of the Fourier coefficient obtained from a QFT. Following this work, Hitzer \cite{hitzer2007quaternion} studied the right-sided and two-sided QFTs. Chen et al.\cite{cheng2019plancherel} systematically developed the QFT theory for the square-integrable quaternion-valued function. One of the important results in Fourier analysis is the multiplication formula $\int_{\mathbb{R}}\hat{f}(\xx)g(\xx)d\xx=\int_{\mathbb{R}}f(\oomega)\hat{g}(\oomega)d\oomega,~f,g\in L^2_{\mathbb{C}}(\mathbb{R}),$ which is not true for right-sided QFT. So, they introduced an auxiliary transform $\boldsymbol\alpha$ given by $\boldsymbol\alpha f(\xx_1,\xx_2)=f_0(\xx_1,\xx_2)+if_1(\xx_1,-\xx_2)+jf_2(-\xx_1,\xx_2)+kf_3(-\xx_1,-\xx_2),$ for  $f(\xx_1,\xx_2)=f_0(\xx_1,\xx_2)+if_1(\xx_1,\xx_2)+jf_2(\xx_1,\xx_2)+kf_k(\xx_1,\xx_2),$ and obtained the multiplication formula as
\begin{align}\label{R2019Chen-KouModifiedMultipolicationFormulaRsQFT}
\int_{\mathbb{R}^2}F_r(\xx_1,\xx_2)g(\xx_1,\xx_2)d\xx_1d\xx_2=\int_{\mathbb{R}^2}f(\xx_1,\xx_2)H_r(\xx_1,\xx_2)d\xx_1d\xx_2,
\end{align}
where $f,g\in L^1_{\mathbb{H}}(\mathbb{R}^2),$ $h=\boldsymbol\alpha g,~H_r=\mathcal{F}_rh$ and $F_r=\mathcal{F}_rf.$ Moreover, if $g\in L^2_{\mathbb{H}}(\mathbb{R}^2)$ then $\|g\|_{L^2_{\mathbb{H}}(\mathbb{R}^2)}=\|h\|_{L^2_{\mathbb{H}}(\mathbb{R}^2)}.$ This multiplication formula has been applied to prove that the right-sided QFT is unitary on $L^2_{\mathbb{H}}(\mathbb{R}^2).$ To study the theory of the two-sided QFT, authors in \cite{cheng2019plancherel} have obtained the relation of $\mathcal{F}_s$ with that of $\mathcal{F}_r$ as $\mathcal{F}_sf=\mathcal{F}_r(\boldsymbol\beta f),$ where $\boldsymbol\beta f(\xx_1,\xx_2)=f_0(\xx_1,\xx_2)+if_1(\xx_1,\xx_2)+jf_2(-\xx_1,\xx_2)+kf_3(-\xx_1,\xx_2).$

The two-sided QFT $\mathcal{F}_Q(f):\mathbb{R}^{2d}\mapsto\mathbb{H}$ of a quaternion function $f\in L^1_{\mathbb{H}}(\mathbb{R}^{2d})$ is defined by \cite{hitzer2007quaternion}
\begin{align}\label{RChenDefinitionQFT}
\mathcal{F}_Q(f)(\bldu,\bldv)=\hat{f}(\bldu,\bldv)=\int_{\mathbb{R}^{2d}}e^{-2\pi i \bldx\cdot \bldu}f(\bldx,\bldy)e^{-2\pi j \bldy\cdot \bldv}d\bldx d\bldy.
\end{align}
Definition \eqref{RChenDefinitionQFT} has a symmetrical representation 
$$\hat{f}(\bldu,\bldv)=\hat{f}_1(\bldu,\bldv)+i\hat{f}_2(\bldu,\bldv)+j\hat{f}_3(\bldu,\bldv)+i\hat{f}_4(\bldu,\bldv)j,$$
where $\hat{f}_l(\bldu,\bldv)=\mathcal{F}_Q(f_l)(\bldu,\bldv),~l=1,2,3,4.$ Chen et al.\cite{chen2015pitt} developed the new module and $L^p-$norm of $\hat{f}$ given respectively as
\begin{align}\label{R2015ChenNewModuleQFT}
|\hat{f}(\bldu,\bldv)|_Q=\left(|\hat{f_1}(\bldu,\bldv)|^2+|\hat{f_2}(\bldu,\bldv)|^2+|\hat{f_3}(\bldu,\bldv)|^2+|\hat{f_4}(\bldu,\bldv)|^2\right)^{\frac{1}{2}},
\end{align}
and 
\begin{align}\label{R2015ChenNewNormQFT}
\|\hat{f}\|_{Q,p}=\left(\int_{\mathbb{R}^{2d}}|\hat{f}(\bldu,\bldv)|^p_Qd\bldu d\bldv\right)^{\frac{1}{p}}
\end{align}
which are not equivalent to $|\hat{f}|_Q$ and $\|\hat{f}\|_{L^{p}_{\mathbb{H}}(\mathbb{R}^{2d})}$  respectively. Apart from several important properties like $\|\hat{f}\|_{Q,\infty}\leq\|f\|_{L^1_{\mathbb{H}}(\mathbb{R}^{2d})},$ continuity of $\hat{f}$ and the Riemann-Lebesgue lemma i.e., $\hat{f}(\bldu,\bldv)\rightarrow 0$ as $|(\bldu,\bldv)|\rightarrow\infty$ for $f\in L^1_{\mathbb{H}}(\mathbb{R}^{2d}),$ using the density of $L^1_{\mathbb{R}}(\mathbb{R}^{2d})\cap L^2_{\mathbb{R}}(\mathbb{R}^{2d})$ in $L^2_{\mathbb{R}}(\mathbb{R}^{2d})$ and the property that $\mathcal{F}_Q(h)\in L^2_{\mathbb{R}}(\mathbb{R}^{2d})$ whenever $h\in L^1_{\mathbb{R}}(\mathbb{R}^{2d})\cap L^2_{\mathbb{R}}(\mathbb{R}^{2d})$ they also extended the definition \eqref{RChenDefinitionQFT} of QFT to the functions in $L^2_{\mathbb{R}}(\mathbb{R}^{2d})$ and then to $f\in L^2_{\mathbb{H}}(\mathbb{R}^{2d}),$ followed by the equality $\|\hat{f}\|_{Q,2}=\|f\|_{L^2_{\mathbb{H}}(\mathbb{R}^{2d})},$ the inversion formula
$\left(\mathcal{F}^{-1}_Q\right)(\bldx)=\mathcal{F}_Q(-\bldx),$ for $f\in L^2_{\mathbb{H}}(\mathbb{R}^{2d}).$ 
%and the QFT \textcolor{red}{of $|\bldx|^{-\alpha}.$}
The generalization of two-sided quaternion FT \eqref{RChenDefinitionQFT} by replacing $i,j$ with a pair of unit quaternions $\mu,\nu,$ i.e.,
\begin{align}\label{RHitzerSteerable2sQFT}
\mathcal{F}_{\mu,\nu}(f)(\uu,\vv)=\int_{\mathbb{R}^{2d}}e^{-2\pi \mu \xx\cdot \uu}f(\xx,\yy)e^{-2\pi \nu \xx\cdot \uu}d\xx d\yy,
\end{align}
which is called steerable quaternion FT, can be found in \cite{hitzer2013orthogonal}. Using the splitting method of the quaternion signal, \eqref{RHitzerSteerable2sQFT} can be written in the following complex form 
$$\mathcal{F}_{\mu,\nu}(f)(\uu,\vv)=\mathcal{F}_{\mu,\nu}(f_{-})(\uu,\vv)+\mathcal{F}_{\mu,\nu}(f_{+})(\uu,\vv),$$
where 
$$\mathcal{F}_{\mu,\nu}(f_{\pm})=\int_{\mathbb{R}^2}e^{2\pi \mu(\xx\uu\mp \yy\vv)}f_{\pm}d\xx d\yy.$$
Recall that Bebenko\cite{babenko1961inequality} first proved the sharp Hausdorff Youngs inequality for the complex FT in $\mathbb{R}$ with even positive number $q$ and later for the general case was proved by Beckner \cite{beckner1975inequalities}. Using the Minkowski inequality, sharp Hausdorff-Young inequality\cite{beckner1975inequalities}  for the complex FT and the modulus identities 
$|f(\xx,\yy)|^2=|f_{-}(\xx,\yy)|^2+|f_{+}(\xx,\yy)|^2$ 
and 
$|\mathcal{F}_{\mu,\nu}(f)(\uu,\vv)|^2=|\mathcal{F}_{\mu,\nu}(f_{-})(\uu,\vv)|^2+|\mathcal{F}_{\mu,\nu}(f_{+})(\uu,\vv)|^2,$ which was first proved for $\mu=i,\nu=j$ in \cite{hitzer2010directional}, Lian\cite{lian2020sharp0}
proved the sharp Hausdorff Young inequality for the quaternion FT as
\begin{align}\label{RLianSharpHausdorffYoungInequality2sQFT}
\|\mathcal{F}_{\mu,\nu}f\|_{L^q_{\mathbb{H}}(\mathbb{R}^2)}\leq A_p^2\|f\|_{L^p_{\mathbb{H}}(\mathbb{R}^2)},
\end{align}
where $A_p=\left(\frac{p^{\frac{1}{p}}}{q^{\frac{1}{q}}}\right)^{\frac{1}{2}},$ $1\leq p\leq 2,~\frac{1}{p}+\frac{1}{q}=1$ and $f\in L^p_{\mathbb{H}}(\mathbb{R}^2).$ The constant $A_p^2$ appearing in the inequality \eqref{RLianSharpHausdorffYoungInequality2sQFT} is sharp as can be seen that choosing one of $f_{\pm}$ to be zero $\mathcal{F}_{\mu,\nu}$ reduces to the complex FT and $A_p^2$ serves as a sharp constant in the Hausdorff-Young inequality for the complex FT\cite{beckner1975inequalities}. Using the relation $\mathcal{F}_l(f)=\mathcal{F}_s\tilde{f}$ and using the sharp inequality \eqref{RLianSharpHausdorffYoungInequality2sQFT}, Lian\cite{lian2020sharp0} also obtained the sharp Hausdorff Young inequality, with the same assumption on $p$ and $f,$ for the left-sided quaternion FT as $\|\mathcal{F}_l(f)\|_{L^q_{\mathbb{H}}(\mathbb{R}^2)}\leq A_p^2\|\tilde{f}\|_{L^p_{\mathbb{H}}(\mathbb{R}^2)},$ where $\tilde{f}(\xx,\yy)=f_0(\xx,\yy)+if_1(\xx,-\yy)+jf_2(\xx,\yy)+kf_3(\xx,-\yy).$ The norms $\|\cdot\|_{L^p_{\mathbb{H}}(\mathbb{R}^2)}$ and $\|\tilde{\cdot}\|_{L^p_{\mathbb{H}}(\mathbb{R}^2)}$ are equivalent but they are not in general equal. However $\|f\|_{L^2_{\mathbb{H}}(\mathbb{R}^2)}=\|\tilde{f}\|_{L^2_{\mathbb{H}}(\mathbb{R}^2)}.$\\

\noindent\textbf{Uncertainty Principles:} The UP describes a feature of a function and its quaternion FT. Beckner \cite{beckner1995pitt} has shown that sharp Pitt's inequality yields a short proof of a logarithmic uncertainty estimate, and Heisenberg-Weyl's inequality follows by using logarithmic inequality. In \cite{beckner1995pitt}, with the aid of rearrangement and symmetrization, Beckner has proved a sharp Pitt's inequality by applying the sharp $L^1$ Young's inequality for convolution on $\mathbb{R}_{+}.$ As in \cite{beckner2008pitt}, based on Young's inequality for the convolution on $\mathbb{R}_+\times S^{2d-1},$ Chen et al.\cite{chen2015pitt} gave a different proof of the Pitt's inequality for QFT. They obtained the Pitt's  inequality 
\begin{align}\label{R2015ChenPitt'sInequalityReal}
\int_{\mathbb{R}^{2d}}|\bldxi|^{-\alpha}|\hat{f}(\bldxi)|^2d\bldxi\leq C_\alpha\int_{\mathbb{R}^{2d}}|\bldx|^{\alpha}|f(\bldx)|^2d\bldx
\end{align} 
for $f\in S(\mathbb{R}^{2d},\mathbb{R}),~0\leq \alpha\leq 2d,$ where $C_\alpha=\pi^\alpha\left[\Gamma\left(\frac{2d-\alpha}{4}\right)/\Gamma\left(\frac{2d+\alpha}{4}\right)\right]$. Using \eqref{R2015ChenPitt'sInequalityReal}, \eqref{R2015ChenNewModuleQFT} they established the Pitt's inequality for the quaternion valued function $f\in S(\mathbb{R}^{2d},\mathbb{H})$ as 
\begin{align}\label{R2015ChenPitt'sInequalityQuaternion}
\int_{\mathbb{R}^{2d}}|\bldxi|^{-\alpha}|\hat{f}(\bldxi)|_Q^2d\bldxi\leq C_\alpha\int_{\mathbb{R}^{2d}}|\bldx|^{\alpha}|f(\bldx)|^2d\bldx.
\end{align}
By differentiating the sharp Pitt's inequality at $\alpha=0,$ they also gave the logarithmic UP as 
\begin{align}\label{R2015ChenLogarithmicUPQFT}
\int_{\mathbb{R}^{2d}}\ln|\bldx||f(\bldx)|^2d\bldx +\int_{\mathbb{R}^{2d}}\ln|\bldxi||f(\bldxi)|^2_Qd\bldxi\geq D\int_{\mathbb{R}^{2d}}|f(\bldx)|^2d\bldx,
\end{align}
$D=\psi(\frac{2d}{4})-\ln(\pi)$ and $\psi=\frac{d}{dt}[\ln\Gamma(t)].$
The Heisenberg-Weyl's UP for $f\in L^2_{\mathbb{H}}(\mathbb{R}^{2d}),$ with $\|f\|^2_2=1$ is
\begin{align}\label{R2015ChenHeisenbergUPQFT}
\Delta_f\bldx\cdot \Delta_f\bldxi\geq \frac{d}{2\pi},
\end{align}
where the standard deviations $\Delta_f\bldx=\left(\int_{\mathbb{R}^{2d}}|\bldx-\bldx_0|^2|f(\bldx)|^2d\bldx\right)^{\frac{1}{2}}$ and $\Delta_f\bldxi=\left(\int_{\mathbb{R}^{2d}}|\bldxi-\bldxi_0|^2|f(\bldxi)|_Q^2d\bldxi\right)^{\frac{1}{2}}$ are smallest for the expectations $\bar{\bldx}=\left(\int_{\mathbb{R}^{2d}}\xx_l|f(\bldx)|^2d\bldx\right)_{l=1}^{2d},~\bldx=(\xx_1,\xx_2,\cdots,\xx_{2d})$ and $\bar{\bldxi}=\left(\int_{\mathbb{R}^{2d}}\xxi_l|f(\bldxi)|_Q^2d\bldxi\right)_{l=1}^{2d},$ $\bldxi=(\xxi_1,\xxi_2,\cdots,\xxi_{2d})$ respectively is obtained in \cite{chen2015pitt}. The inequality \eqref{R2015ChenHeisenbergUPQFT} is based on the interpretation of the standard deviation as the size of the essential support of $f.$ With the different notion of the support gives different version of the UP like, for a non-zero function $f\in L^2_{\mathbb{H}}(\mathbb{R}^{2d})$ which is $\epsilon_T-$concentrated on a measurable set $T\subset \mathbb{R}^{2d}$ and $\hat{f}$ is $\epsilon_{\Omega}-$concentrated on a measurable set $\Omega\subset \mathbb{R}^{2d},$ 
$$|T||\Omega|\geq (1-\epsilon_t-\epsilon_{\Omega})^2,$$
which is the Donoho-Stark's UP, for an alternative proof see \cite{kou2017uncertainty},\cite{lian2018uncertainty}. It happens that if $supp f\subset T$ and $supp \hat{f}\subset\Omega,$ then $|T||\Omega|\geq 1.$
%%%%%%%%%%%%%%%%%%%%%%%%%%%%%%%%%%%%%%%%%%%%%%%%%%%%%%%%%%%%%%%%%%%%%%%%%%%%%%%%%%%%%%%%%%%%%%%%%%%
%%%%%%%%%%%%%%%%%%%%%%%%%%%%%%%%%%%%%%%%%%%%%%%%%%%%%%%%%%%%%%%%%%%%%%%%%%%%%%%%%%%%%%%%%%%%%%%%%%%
%%%%%%%%%%%%%%%%%%%%%%%%%%%%%%%%%%%%%%%%%%%%%%%%%%%%%%%%%%%%%%%%%%%%%%%%%%%%%%%%%%%%%%%%%%%%%%%%%%%
%%%%%%%%%%%%%%%%%%%%%%%%%%%%%%%%%%%%%%%%%%%%%%%%%%%%%%%%%%%%%%%%%%%%%%%%%%%%%%%%%%%%%%%%%%%%%%%%%%%
%{Check:}%%%%%%%%%%%%%%%%%%%%%%%%%%%%%%%%%%%%%%%%%%%%%%%%%%%%%%%%%%%%%%%%%%%%%%%%%%%%%%%%%%%%%%%%%%
In \cite{kou2017uncertainty}, authors discussed the recovery of a quaternion-valued signal in $\mathbb{R}^2$ despite a significant amount of missing information. Like, say, if a signal $f\in L^2_{\mathbb{H}}(\mathbb{R}^2)$ is transmitted to a receiver who knows that it is bandlimited, meaning that $f$ was synthesized using only frequencies in a set $\Omega\subset \mathbb{R}^2.$ Now, suppose that the receiver is unable to observe all the data of $f,$ a certain subset $T$ of $\bldx-$values is unobserved. Moreover, the observed signal is contaminated by observational noise $n(\bldx)\in L^2_{\mathbb{H}}(\mathbb{R}^2).$ Thus the received signal $r(\bldx)\in L^2_{\mathbb{H}}(\mathbb{R}^2)$ satisfies
$$
r(\bldx)=
\begin{cases}
f(\bldx)+n(\bldx),&\bldx\nin T,\\
0,&\bldx\in T.
\end{cases}
$$
The receiver's aim is to reconstruct the transmitted signal $f$ from the noisy received signal $r(\bldx).$ Although it may seem that information  about $f$ for $\bldx\in T$ is unavailable,  the UPs say that the recovery is possible provided that $|T||\Omega|<1.$
%%%%%%%%%%%%%%%%%%%%%%%%%%%%%%%%%%%%%%%%%%%%%%%%%%%%%%%%%%%%%%%%%%%%%%%%%%%%%%%%%%%%%%%%%%%%%%%%%%%
%%%%%%%%%%%%%%%%%%%%%%%%%%%%%%%%%%%%%%%%%%%%%%%%%%%%%%%%%%%%%%%%%%%%%%%%%%%%%%%%%%%%%%%%%%%%%%%%%%%
%%%%%%%%%%%%%%%%%%%%%%%%%%%%%%%%%%%%%%%%%%%%%%%%%%%%%%%%%%%%%%%%%%%%%%%%%%%%%%%%%%%%%%%%%%%%%%%%%%%
%%%%%%%%%%%%%%%%%%%%%%%%%%%%%%%%%%%%%%%%%%%%%%%%%%%%%%%%%%%%%%%%%%%%%%%%%%%%%%%%%%%%%%%%%%%%%%%%%%%
Donoho and Stark first proved this result \cite{donoho1989uncertainty} in the classical Fourier domain for the function on $\mathbb{R}.$ Some generalizations of the Donoho-Stark UP for the right-sided QFT have been studied in \cite{abouelaz2020donoho}.

The following result, in \cite{chen2015pitt}, is a generalization in the quaternion setting of the result due to Benedicks \cite{benedicks1985fourier} (also see \cite{folland1997uncertainty}) for ordinary FT whose elegant proof was first circulated as a preprint in 1974 but not formally published for another decade.
\begin{align}\label{R2015ChenBenedicksUPQFT}
\mbox{If}~f\in L^1_{\mathbb{H}}(\mathbb{R}^{2d})~ \mbox{with}~supp f\subset T~\mbox{and}~supp \hat{f}\subset \Omega~\mbox{such that}~|T||\Omega|<\infty, ~\mbox{then}~f=0.
\end{align}
Note that if $f\in L^p_{\mathbb{H}}(\mathbb{R}^{2d})(p>1)$ and $|T|<\infty,$ then $f\in L^1_{\mathbb{H}}(\mathbb{R}^{2d})$ by Holder's inequality and if $f\in L^1_{\mathbb{H}}(\mathbb{R}^{2d})$ and $|\Omega|<\infty$ then $f\in L^p_{\mathbb{H}}(\mathbb{R}^{2d})$ for all $p>1$ hence the result \eqref{R2015ChenBenedicksUPQFT} applies equally to $L^p_{\mathbb{H}}(\mathbb{R}^{2d})$ functions. The proof of the result \eqref{R2015ChenBenedicksUPQFT} relies on the Poisson summation formula for the quaternion valued function, i.e., if $f\in L^1_{\mathbb{H}}(\mathbb{R}^{2d}),$ the series $\displaystyle \phi(\xx_1,\xx_2)=\sum_{(k_1,k_2)\in\mathbb{Z}^{2d}}f(\xx_1+k_1,\xx_2+k_2)$ converges in $L^1_{\mathbb{H}}([0,1]^{2d}),$ then the quaternion Fourier series of $\phi$ is $\displaystyle\sum_{(k_1,k_2)\in\mathbb{Z}^{2d}}e^{2\pi ik_1\cdot \xx_1}\hat{f}(k_1,k_2)e^{2\pi jk_2\cdot \xx_2}.$ The local uncertainty inequalities for the classical FT, which states that if $f\in L^2(\mathbb{R}^d,\mathbb{C})$ is highly localized then its FT can be concentrated in a small neighborhood of two or more separated points, was first obtained by Faris \cite{faris1978inequalities} and sharper by Price \cite{price1983inequalities},\cite{price1987sharp}. Lian \cite{lian2018uncertainty} obtained the local UP for the two-sided QFT as
%If $0<\alpha<1,$ there is a constant $K_\alpha$ such that for all $f\in L^2_{\mathbb{H}}(\mathbb{R}^2)$ and all measurable $E\subset\mathbb{R}^2,$
\begin{align}\label{R2018LianLocalUP0<alpha<1}
\int_{E}|(\mathcal{F}f)(\bldx)|^2d\bldx\leq K_\lambda|E|^{\lambda}\||\bldx|^{\lambda}f\|^2_2,~0<\lambda<1
\end{align}
and
\begin{align}\label{R2018LianLocalUPalpha>1}
\int_{E}|(\mathcal{F}f)(\bldx)|^2d\bldx\leq K_\lambda|E|\|f\|^{2-\frac{2}{\lambda}}_2\||\bldx|^{\lambda}f\|^{\frac{2}{\lambda}}_2,~\lambda>1
\end{align}
for all $f\in L^2_{\mathbb{H}}(\mathbb{R}^2)$ and all measurable $E\subset\mathbb{R}^2,$ where
\begin{align}\label{R2018LianLocalUPConstant}
K_\lambda=
\begin{cases}
\frac{(1+\lambda)^2}{\lambda^{2\lambda}}(2-2\lambda)^{\lambda-2},&0<\lambda<1\\
\frac{\pi}{\lambda\Gamma(\frac{1}{\lambda})}\Gamma(\frac{1}{\lambda})\Gamma(1-\frac{1}{\lambda})(\lambda-1)^\lambda(1-\frac{1}{\lambda})^{-1},&\lambda>1.
\end{cases}
\end{align}
Using this local UP and following the lines of  \cite{ghobber2012uncertainty}, one can get the global UP, which includes different powers of norm as
%For $s,\beta>0,$ there exists a constant $C_{s,\beta}$ such that for all $f\in L^2_{\mathbb{H}}(\mathbb{R}^2),$ we have
\begin{align}\label{R2018LianGlobalUP}
\||\bldx|^sf\|^{\frac{2\beta}{s+\beta}}_{L^2_{\mathbb{H}}(\mathbb{R}^2)}\cdot \||\bldy|^{\beta}\mathcal{F}f\|^{\frac{2s}{s+\beta}}_{L^2_{\mathbb{H}}(\mathbb{R}^2)}\geq C_{s,\beta}\|f\|^2_{L^2_{\mathbb{H}}(\mathbb{R}^2)},
\end{align}
where $s,\beta>0,$ $f\in L^2_{\mathbb{H}}(\mathbb{R}^2)$ and $C_{s,\beta}$ is some constant.

The Pitt's inequality and the logarithmic UPs given in \eqref{R2015ChenPitt'sInequalityQuaternion} and \eqref{R2015ChenLogarithmicUPQFT} are expressed in terms of the module given by \eqref{R2015ChenNewModuleQFT}. Moreover, it was first proved for real-valued function under two-sided QFT following Beckner's lines with quite a lot of effort. Thus, a genuine analogues of Pitt's inequality and the logarithmic UP for quaternion-valued function on $\mathbb{R}^2$ is obtained in \cite{lian2018uncertainty}, using the two-dimensional plane split of a quaternion and the results for the complex FT, which are same as \eqref{R2015ChenPitt'sInequalityQuaternion} and \eqref{R2015ChenLogarithmicUPQFT} respectively with the norm $|\cdot|_Q$ replaced as $|\cdot|,$ furthermore the constant $D$ in equation \eqref{R2015ChenLogarithmicUPQFT} is $\psi(1)-\ln(\pi).$ The method adopted by Lian\cite{lian2018uncertainty} is easier and also applicable to the one-sided QFT, general two-sided QFT or QFT on higher dimension $\mathbb{R}^d.$ Without using the sharp Hausdorff-Young's inequality of the two-sided QFT, Lian also obtained the Entropy UP
\begin{align}\label{R2018LianEntropyUP2sQFT}
H(|f|^2)+H(|\mathcal{F}(f)|^2)\geq 2(1-\ln 2),
\end{align}
for $f\in S(\mathbb{R}^2,\mathbb{H})$ and $\|f\|_{L^2_{\mathbb{H}}(\mathbb{R}^2)}=1,$ where $H(f)=-\int_{\mathbb{R}^2}f(\bldx)\ln f(\bldx)d\bldx.$ 
%followed by one of a qualitative UP called the Benedicks-Amrein-Berther principle which says that 
%If $S,V$ be a measurable subset of $\mathbb{R}^2$ with \textcolor{purple}{$\mu(S),\mu(V)<\infty.$} There exists a constant $C=C(S,V)$ such that for all $f\in L^2_{\mathbb{H}}(\mathbb{R}^2),$
%\begin{align}\label{R2018LianBenedicksAmreinBertherUP2sQFT}
%\|f\|^2_2\leq C\left(\|f\|^2_{L^2(\mathbb{R}^2\setminus S)}+\|\mathcal{F}f\|^2_{L^2(\mathbb{R}^2\setminus V)}\right).
%\end{align}

\subsection{One dimensional quaternion Fourier transform}
So far, we have seen that the quaternion FT was studied for the quaternion-valued function defined on $\mathbb{R}$ or $\mathbb{R}^{2d}.$ The quaternion FT for the quaternion-valued function defined on $\mathbb{R}$ was first defined by He et al.\cite{he2004continuous} as 
\begin{align}\label{R2002HeDefinitionQuaternionFT}
\hat{f}(\xxi)=\hat{f_1}(\xxi)+j \hat{f_2}(\xxi)
\end{align}
for $f\in L^2_{\mathbb{H}}(\mathbb{R}).$ They obtained the Parseval's formula $\langle f,g\rangle_{L^2_{\mathbb{H}}(\mathbb{R})}=\langle f,g\rangle_{L^2_{\mathbb{H}}(\mathbb{R})},$ for $f,g\in L^2_{\mathbb{H}}(\mathbb{R}),$ and proved that $\widehat{(f\star g)}(\xxi)\neq \hat{f}(\xxi)\hat{g}(\xxi)$ in general, for the convolution $\star$ defined by 
\begin{align}\label{R2002HeConvolutionQuaternion1DFT}
(f\star g)(\xx)=\int_{\mathbb{R}}f(\xx-\xx')g(\xx')d\xx'.
\end{align}
Also if we take $f=f_1+jf_2$ and $g=g_1+jg_2,$ then we have
$$\widehat{(f\star g)}(\xxi)=\left[\hat{f_1}(\xxi)\hat{g_1}(\xxi)-\hat{\bar{f_2}}(\xxi)\hat{g_2}(\xxi)\right]+j \left[\hat{f_2}(\xxi)\hat{g_1}(\xxi)-\hat{\bar{f_1}}(\xxi)\hat{g_2}(\xxi)\right]$$
$$\hat{f}(\xxi)\hat{g}(\xxi)=\left[\hat{f_1}(\xxi)\hat{g_1}(\xxi)-\bar{\hat{f_2}}(\xxi)\hat{g_2}(\xxi)\right]+j \left[\hat{f_2}(\xxi)\hat{g_1}(\xxi)-\bar{\hat{f_1}}(\xxi)\hat{g_2}(\xxi)\right].$$ 
Thus $\widehat{(f\star g)}(\xxi)=\hat{f}(\xxi)\hat{g}(\xxi)$ if $\bar{\hat{f_1}}(\xxi)=\hat{\bar{f_1}}(\xxi)$ and $\bar{\hat{f_2}}(\xxi)=\hat{\bar{f_2}}(\xxi).$ Bahri et al.\cite{bahri2019one} also studied one-dimensional quaternion FT with some convolution $\star $ for which $f\star g\neq g\star f$ in general and also faces the same problem as that in \cite{he2004continuous}. Akila et al.\cite{akila2014natural} introduced a novel natural convolution given as 
\begin{align}\label{R2014AkilaConvolution1DQFT}
(f\ostar g)(\xx)=(f_1\star g_1)(\xx)-(\check{\overline{f_2}}\star g_2)(\xx)+j\left((\check{\overline{f_1}}\star g_2)(\xx)+(f_2\star g_1)(\xx)\right)
\end{align}
and proved that $f\ostar\delta\in L^p_{\mathbb{H}}(\mathbb{R}),$ $\alpha(f\ostar g)=(\alpha f)\ostar g,$ $(f+g)\ostar\delta=f\ostar\delta+g\ostar\delta$ and $(f\ostar g)\ostar\epsilon=f\ostar(\delta\ostar\epsilon)$ for all $\alpha\in\mathbb{H},$ $\delta,\epsilon\in L^1_{\mathbb{H}}(\mathbb{R}),$ and $f\in L^p_{\mathbb{H}}(\mathbb{R}),~p=1,2.$ They also proved the convolution theorem $\mathcal{F}(f\ostar g)=\mathcal{F}(f)\mathcal{F}(\delta),~\delta\in L^1_{\mathbb{H}}(\mathbb{R})$ and $f\in L^2_{\mathbb{H}}(\mathbb{R}),$ and hence $\ostar$ defined in \eqref{R2014AkilaConvolution1DQFT} serves as a better convolution structure for the transform defined in \eqref{R2002HeDefinitionQuaternionFT} than the one given in \eqref{R2002HeConvolutionQuaternion1DFT} and also in \cite{bahri2019one}. In general, $\ostar$ is non-commutative, i.e., $f\ostar g$ need not be equal to $g\ostar f.$ If both $f$ and $g$ are complex-valued or if one of them is a real-valued even function, then we get $f\ostar g=g\ostar f.$ Finding a suitable non-trivial necessary and sufficient condition on $f$ and $g$ to get $f\star g=g\star f$ is still and interesting open problem. The definition \eqref{R2002HeDefinitionQuaternionFT} together with the convolution structure \eqref{R2014AkilaConvolution1DQFT} and its properties for the arbitrary dimension $d,$ i.e., for the functions defined on $\mathbb{R}^d$ can be found in \cite{akila2016multidimensional}. 
\section{Quaternion Fractional Fourier Transform (QFrFT)}
Guanlei et al. \cite{guanlei2008fractional} proposed a definition of two-sided QFrFT. It extends the definition of classical FrFT for the complex-valued function to that of the quaternion-valued functions. This transform is a natural extension of the two-sided QFT \eqref{R2019Chen-KouDefinitionTsQFT}. The idea of defining the transform is by replacing the Fourier kernel in \eqref{R2019Chen-KouDefinitionTsQFT} with that of the fractional Fourier kernel. Thus the definition of the two-sided QFrFT for the function $f\in L^2_{\mathbb{H}}(\mathbb{R})$ is defined by
\begin{align}\label{R2008_Guanlei-Xiaotong-Xiaogang_SP_Definition_2sQFrFT}
\left(\mathcal{F}^{ij}_{\alpha,\beta}f\right)(\uu_1,\uu_2)=\int_{\mathbb{R}^2}K^i_{\alpha}(\xx_1,\uu_1)f(\xx_1,\xx_2)K^j_{\beta}(\xx_2,\uu_2)d\xx_1d\xx_2,
\end{align} 
where $K^i_{\alpha}(\xx_1,\uu_1)=\sqrt{\frac{1-i\cot\alpha}{2\pi}}e^{i\left(\frac{\xx_1^2+\uu_1^2}{2\tan\alpha}-\frac{\xx_1\uu_1}{\sin\alpha}\right)}$ and $K^j_{\beta}(\xx_2,\uu_2)=\sqrt{\frac{1-j\cot\alpha}{2\pi}}e^{j\left(\frac{\xx_2^2+\uu_2^2}{2\tan\alpha}-\frac{\xx_2\uu_2}{\sin\alpha}\right)}$ and $\alpha,\beta\neq n\pi,~n\in \mathbb{Z}.$ In particular, the QFrFT \eqref{R2008_Guanlei-Xiaotong-Xiaogang_SP_Definition_2sQFrFT}, reduces to QFT for $\alpha=\beta=\frac{\pi}{2},$ reduces to the classical $1-D$ FT of $f(\xx_1,\xx_2)$ for variable $\xx_1,$ when $\alpha=\frac{\pi}{2}$ and $\beta=0,$ reduces to the $1-D$ FT of $f(\xx_1,\xx_2)$ for variable $\xx_2.$ Moreover, the transform \eqref{R2008_Guanlei-Xiaotong-Xiaogang_SP_Definition_2sQFrFT} is such that it has the property that $\left(\mathcal{F}^{ij}_{\alpha,\beta}f\right)(\uu_1,\uu_2)=f(\xx_1,\xx_2)$ for $(\alpha,\beta)=(0,0)~\mbox{or}~(2\pi,2\pi)$ and $\left(\mathcal{F}^{ij}_{\pi,\pi}f\right)(\uu_1,\uu_2)=f(-\xx_1,-\xx_2),$ $\left(\mathcal{F}^{ij}_{3\frac{\pi}{2},3\frac{\pi}{2}}f\right)(\uu_1,\uu_2)=(\mathcal{F}_sf)(-\uu_1,-\uu_2),$ where $\mathcal{F}_s$ is a two sided QFrFT given by \eqref{R2019Chen-KouDefinitionTsQFT}. Authors in \cite{guanlei2008fractional}, obtained the $\mathbb{R}-$linear property, odd-even invariant property, i.e., $\left(\mathcal{F}^{ij}_{\alpha,\beta}f(\epsilon_1\xx_1,\epsilon_2\xx_2)\right)(\uu_1,\uu_2)=\left(\mathcal{F}^{ij}_{\alpha,\beta}f\right)(\epsilon_1\uu_1,\epsilon_2\uu_2),~\epsilon_l=\pm 1,~l=1,2$ and also the additivity property, i.e., $\mathcal{F}^{ij}_{\alpha_1,\beta_1}\mathcal{F}^{ij}_{\alpha_2,\beta_2}=\mathcal{F}^{ij}_{\alpha_1+\alpha_2,\beta_1+\beta_2}$ of the transform \eqref{R2008_Guanlei-Xiaotong-Xiaogang_SP_Definition_2sQFrFT}. Wei et al. \cite{wei2013different}, obtained the inversion formula for the transform \eqref{R2008_Guanlei-Xiaotong-Xiaogang_SP_Definition_2sQFrFT} as
$$f(\xx_1,\xx_2)=\int_{\mathbb{R}^2}K^i_{-\alpha}(\xx_1,\uu_1)\left(\mathcal{F}^{ij}_{\alpha,\beta}f\right)(\uu_1,\uu_2)K^j_{-\beta}(\xx_2,\uu_2)d\uu_1d\uu_2,$$
and proved the Parseval's formula $\langle \mathcal{F}^{ij}_{\alpha,\beta}f,\mathcal{F}^{ij}_{\alpha,\beta}g\rangle =\langle f,g\rangle,~f,g\in L^2_{\mathbb{H}}(\mathbb{R})$ with respect to the scalar product \eqref{RQuaternionScalerProduct}. In particular, when $f=g,$ it follows that 
$ \|\mathcal{F}^{ij}_{\alpha,\beta}f\|_{ L^2_{\mathbb{H}}(\mathbb{R})}= \|f\|_{ L^2_{\mathbb{H}}(\mathbb{R})}.$ From the signal processing point of view, the Parseval's formula states that the signal energy is preserved by the transform \eqref{R2008_Guanlei-Xiaotong-Xiaogang_SP_Definition_2sQFrFT}. It can be seen that it is not possible to establish a general Parseval's formula for the inner product given by \eqref{RQuaternionInnerProduct} because it lacks the cyclic symmetry \eqref{RQuaternionCyclicMultiplicationSymmetry}. Thus, to obtain such Parseval's formula, it is therefore either necessary to modify the symmetric property of the inner product as in \eqref{RQuaternionScalerProduct} or to modify the QFrFT itself. Thus, author in \cite{wei2013different} gave the definition of right-sided QFrFT 
\begin{align}\label{R2008Wei-Li_Optik_Definition_RQFrFT}
\left(\mathcal{Q}^{ij}_{\alpha,\beta}f\right)(\uu_1,\uu_2)=\int_{\mathbb{R}^2}f(\xx_1,\xx_2)K^i_{\alpha}(\xx_1,\uu_1)K^j_{\beta}(\xx_2,\uu_2)d\xx_1d\xx_2.
\end{align}
The transform enjoys the property $\left( \mathcal{Q}^{ij}_{\alpha,\beta}f,\mathcal{Q}^{ij}_{\alpha,\beta}g\right) =\left( f,g\right),~f,g\in L^2_{\mathbb{H}}(\mathbb{R})$ according to the inner product \eqref{RQuaternionInnerProduct}, as a consequence of which it follows that $\|\mathcal{Q}^{ij}_{\alpha,\beta}f\|_{L^2_{\mathbb{H}}(\mathbb{R})}=\|f\|_{L^2_{\mathbb{H}}(\mathbb{R})}.$ In what follows, we see the progress of the FrFT for the quaternion-valued function defined on $\mathbb{R}.$

Using the idea of He et al.\cite{he2004continuous}, Roopkumar \cite{roopkumar2016quaternionic} gave the definition of FrFT of type II of quaternion-valued function $f=f_1+jf_2\in L^2_{\mathbb{H}}(\mathbb{R})$ as
\begin{align}\label{R2016Roopkumar_1DFrFT_Definition}
(\mathcal{F}_\alpha f)(\xxi)=(\mathcal{F}_\alpha f_1)(\xxi)+j(\mathcal{F}_\alpha f_2)(\xxi),
\end{align}
where $\mathcal{F}_\alpha f_l,~l=1,2$ is the  classical FrFT. The transform \eqref{R2016Roopkumar_1DFrFT_Definition} is $\mathbb{H}-$linear satisfying the property $\mathcal{F}_\alpha\mathcal{F}_\beta f=\mathcal{F}_{\alpha+\beta}f,$ for $f\in L^2_{\mathbb{H}}(\mathbb{R}).$ Apart from this, the transform also satisfies the following Parseval's identity
\begin{align}\label{R2016Roopkumar_1DFrFT_Parseval's-Formula}
(\mathcal{F}_\alpha f,\mathcal{F}_\alpha g)=(f,g),~f,g\in L^2_{\mathbb{H}}(\mathbb{R}),
\end{align}
where $(\cdot,\cdot)$ denotes the inner product given by \eqref{RQuaternionInnerProduct}. Motivated by the idea of defining the convolution of one-dimensional quaternion-valued function for the QFT given in \eqref{R2014AkilaConvolution1DQFT}, author in\cite{roopkumar2016quaternionic} defined the  convolution associated with the QFrFT as 
\begin{align}\label{R2016RoopkumarConvolution1DQFrFT}
(f\ostar_\alpha g)(\xx)=(f_1\star_\alpha g_1)(\xx)-(\mathcal{F}_{-2\alpha}\overline{f_2}\star_\alpha g_2)(\xx)+j\left((\mathcal{F}_{-2\alpha}\overline{f_1}\star_\alpha g_2)(\xx)+(f_2\star_\alpha g_1)(\xx)\right),
\end{align}
where $\star_\alpha$ is the convolution for the complex valued function given by $(f_1\star_\alpha g_1)(\xx)=\frac{c_\alpha}{\sqrt{2\pi}}e^{-ia_\alpha\xx^2}(\tilde{f_1}\star\tilde{g_1})(\xx),$ where $\tilde{f_1}(\xx)=f_1(\xx)e^{ia_\alpha\xx^2},$ $c_\alpha=\sqrt{1-i\cot\alpha},$ $a_\alpha=\frac{\cot\alpha}{2}$ and $\star$ is the classical convolution given as $(f_1\star g_1)(\xx)=\int_{\mathbb{R}}f_1(\xx-\yy)g(\yy)d\yy.$ Also, they considered the following operation between two quaternion-valued functions as
 \begin{align}\label{R2016RoopkumarProduct1DQFrFT}
(f\otimes_\alpha g)(\xx)=(f_1\star'_\alpha g_1)(\xx)-(\overline{\mathcal{F}_{-2\alpha}f_2}\star'_\alpha g_2)(\xx)+j\left((\overline{\mathcal{F}_{-2\alpha}f_1}\star'_\alpha g_2)(\xx)+(f_2\star'_\alpha g_1)(\xx)\right),
\end{align}
where $\star'_\alpha$ is the operation defined on complex valued functions as $(f_1\star'_\alpha g_1)(\xx)=\frac{1}{\sqrt{2\pi}}e^{ia_\alpha\xx^2}(\tilde{\tilde{f_1}}\star\tilde{\tilde{g_1}})(\xx),$ where $\tilde{\tilde{f_1}}(\xx)=f_1(\xx)e^{-ia_\alpha\xx^2}.$ For the convolution \eqref{R2016RoopkumarConvolution1DQFrFT} and the operation \eqref{R2016RoopkumarProduct1DQFrFT}, author in \cite{roopkumar2016quaternionic} obtained the convolution and product theorem given respectively as $(\mathcal{F}_\alpha(f\ostar_\alpha g))(\xi)=(\mathcal{F}_\alpha f)(\xi) (\mathcal{F}_\alpha g)(\xi)e^{-ia_\alpha\xi^2}$ and $\mathcal{F}_\alpha\left(\tilde{\tilde{(fg)}}\right)=c_{-\alpha}[\mathcal{F}_\alpha f\otimes_\alpha\mathcal{F}_\alpha g]$ for $f\in L^2_{\mathbb{H}}(\mathbb{R}),~g\in L^1_{\mathbb{H}}(\mathbb{R}).$ For the convolution \eqref{R2016RoopkumarConvolution1DQFrFT}, with $\star_\alpha$ for the complex valued functions given by $(f_1\star_\alpha g_1)(\xx)=e^{-ia_\alpha\xx^2}(\tilde{f_1}\star g_1)(\xx),$ author in \cite{roopkumar2018quaternionic} proved that  for $f\in L^2_{\mathbb{H}}(\mathbb{R})$ and $g\in L^2_{\mathbb{H}}(\mathbb{R})$ $f\ostar_\alpha g\in L^2_{\mathbb{H}}(\mathbb{R})$ and obtained the convolution theorem for the transform \eqref{R2016Roopkumar_1DFrFT_Definition} as 
\begin{align}\label{R2017_Roopkumar_QFrWT_Convolution-Theorem}
(\mathcal{F}_\alpha(f\ostar_\alpha g))(\xxi)=\sqrt{2\pi}(\mathcal{F}_\alpha f)(\xxi)(\mathcal{F}_\alpha g)(\xxi\csc\alpha).
\end{align}
They proved that the convolution $\ostar_\alpha$ satisfies the properties analogue to those for the convolution $\ostar$ given in  \eqref{R2014AkilaConvolution1DQFT}, like $(f+h)\ostar_\alpha g=f\ostar_\alpha g+h\ostar_\alpha g$ and $(\gamma f)\ostar_\alpha g=\gamma(f\ostar_\alpha g),$ for $f,h\in L^2_{\mathbb{H}}(\mathbb{R}),$ $g\in L^1_{\mathbb{H}}(\mathbb{R})$ and $\gamma\in \mathbb{H}.$ Motivated by He et al.\cite{he2004continuous} and Roopkumar \cite{roopkumar2016quaternionic}, Gupta et al.\cite{gupta2023quaternionic} extended the FrFT \cite{srivastava2019certain} for the complex-valued functions to the quaternion-valued function. They called the transform as the QFrFT and is defined for $f=f_1+jf_2\in L^2_{\mathbb{H}}(\mathbb{R})$ and $0<\theta\leq 1$ as 
\begin{align}\label{R2023_BookChapter_Gupta-Verma_Cattani_1DQFrFT_Definition}
\mathfrak{F}_\theta f =\mathfrak{F}_\theta f_1+j\mathfrak{F}_\theta f_2,
\end{align}
where $(\mathfrak{F}_\theta g)(\xxi)=\int_{\mathbb{R}}g(t)e^{-i(\sgn\xxi)|\xxi|^\frac{1}{\theta}\tee}d\tee$ for $g\in L^2_{\mathbb{C}}(\mathbb{R})$ is the FrFT given in \cite{srivastava2019certain}. They proved that the transform \eqref{R2023_BookChapter_Gupta-Verma_Cattani_1DQFrFT_Definition} is $\mathbb{H}-$linear and using the Riemann-Lebesgue lemma and the Parseval's formula for the classical case author proved the same for the quaternion case i.e., $\lim_{|\xxi|\rightarrow \infty}|(\mathfrak{F}_\theta f)(\xxi)|= 0,$ for $f\in L^1_{\mathbb{H}}(\mathbb{R})$ and $\int_{\mathbb{R}}f(\tee)(g(\tee))^cd\tee=\frac{1}{2\pi\theta}\int_{\mathbb{R}}|\xxi|^{\frac{1}{\theta}-1}(\mathfrak{F}_\theta f)(\xxi)((\mathfrak{F}_\theta g)(\xxi))^cd\xxi,$ for $f,g\in L^2_{\mathbb{H}}(\mathbb{R}).$ With the convolution $\otimes$ given by \eqref{R2014AkilaConvolution1DQFT}, the convolution theorem for the QFrFT \eqref{R2023_BookChapter_Gupta-Verma_Cattani_1DQFrFT_Definition} is given by $\mathfrak{F}_\theta (f\otimes g)=(\mathfrak{F}_\theta f)(\mathfrak{F}_\theta g),$ for $f\in L^2_{\mathbb{H}}(\mathbb{R}),~g\in L^1_{\mathbb{H}}(\mathbb{R}).$ Authors in \cite{gupta2023quaternionic} also obtained the Heisenberg type UP for the QFrFT \eqref{R2023_BookChapter_Gupta-Verma_Cattani_1DQFrFT_Definition} as
\begin{align}\label{R2023_Gupta-Verma-Cattani_BookChapter_QFrFT_UP}
\left(\int_{\mathbb{R}}|x|^2|f(x)|^2dx\right)\left(\int_{\mathbb{R}}|\xi|^{\frac{3}{\theta}-1}|(\mathfrak{F}_{\theta} f)(\xi)|^2\right)\geq \frac{\pi\theta}{4}\|f\|^4_{L^2_{\mathbb{H}}(\mathbb{R})},
\end{align}
for the function $f=f_1+jf_2$ such the $f_l,~l=1,2$ are in the Schwartz class.

\section{Quaternion Linear Canonical Transform (QLCT)}
The LCT has a degree of freedom and is more flexible than the FT and the FrFT, with similar computational costs. Thus, because of these advantages, the LCT has been extended to the quaternion algebra framework. This leads to the notion of the QLCT; furthermore, because of the non-commutativity of the quaternion multiplication, there are three different types of QLCTs, namely right-sided, two-sided, and left-sided QLCT. The following are their definitions of right-sided, left-sided, and two-sided QLCT, respectively, which are obtained by replacing the kernels of the FT with that of the kernels of LCT 
\begin{align}\label{R2016Kou-Ou-MoraisDefinitionLsQLCT}
\left(\mathcal{H}^{ij}_{M_1,M_2}f\right)(\uu_1,\uu_2)=\int_{\mathbb{R}^2}K^i_{M_1}(\xx_1,\uu_1)K^j_{M_2}(\xx_2,\uu_2)f(\xx_1,\xx_2)d\xx_1d\xx_2,
\end{align} 
\begin{align}\label{R2016Kou-Ou-MoraisDefinitionRsQLCT}
\left(\mathcal{Q}^{ij}_{M_1,M_2}f\right)(\uu_1,\uu_2)=\int_{\mathbb{R}^2}f(\xx_1,\xx_2)K^i_{M_1}(\xx_1,\uu_1)K^j_{M_2}(\xx_2,\uu_2)d\xx_1d\xx_2
\end{align} 
and
\begin{align}\label{R2016Kou-Ou-MoraisDefinitionTsQLCT}
\left(\mathcal{L}^{ij}_{M_1,M_2}f\right)(\uu_1,\uu_2)=\int_{\mathbb{R}^2}K^i_{M_1}(\xx_1,\uu_1)f(\xx_1,\xx_2)K^j_{M_2}(\xx_2,\uu_2)d\xx_1d\xx_2,
\end{align} 
where 
\begin{align}\label{R2016Kou-Ou-MoraisDefinitionTsQLCT-KernelI}
K^i_{M_1}(\xx_1,\uu_1)=
\begin{cases}
\frac{1}{\sqrt{i2\pi \Mb_1}}e^{i\left(\frac{\Ma_1}{2\Mb_1}\xx_1^2-\frac{1}{\Mb_1}\xx_1\uu_1+\frac{d_1}{2\Mb_1}\uu_1^2\right)},&~\mbox{for}~\Mb_1\neq 0 \\
\sqrt{\Md_1}e^{i\left(\frac{\Mc_1\Md_1}{2}\right)\uu_1^2}\delta(\xx_1-\Md_1\uu_1),&~\mbox{for}~\Mb_1=0,
\end{cases}
\end{align}
and
\begin{align}\label{R2016Kou-Ou-MoraisDefinitionTsQLCT-KernelII}
K^j_{M_2}(\xx_2,\uu_2)=
\begin{cases}
\frac{1}{\sqrt{i2\pi \Mb_2}}e^{j\left(\frac{\Ma_2}{2\Mb_2}\xx_2^2-\frac{1}{\Mb_2}\xx_2\uu_2+\frac{\Md_2}{2\Mb_2}\uu_2^2\right)},&~\mbox{for}~\Mb_2\neq 0 \\
\sqrt{\Md_2}e^{j\left(\frac{\Mc_2\Md_2}{2}\right)\uu_2^2}\delta(\xx_2-\Md_2\uu_2),&~\mbox{for}~\Mb_2=0,
\end{cases}
\end{align}
$M_l=(\Ma_l,\Mb_l,\Mc_l,\Md_l)\in\mathbb{R}^{2\times 2}$ are matrix parameters such that $det(M_l)=1,$ for $l=1,2$ and $f\in L^1_{\mathbb{H}}(\mathbb{R}^2).$ To obtain the Parseval theorem for the QFT, Hitzer \cite{hitzer2007quaternion} modified the inner product to be the scalar inner product, while Chen \cite{chen2015pitt} developed a strategy applying the $Q-$modules, i.e., \eqref{R2015ChenNewModuleQFT}, of QFT to define the $L^p,$ i.e., \eqref{R2015ChenNewNormQFT} to prove the Parseval formula. Bahri \cite{bahri2019relationships} obtained the relation between any two transforms given by \eqref{R2016Kou-Ou-MoraisDefinitionLsQLCT},\eqref{R2016Kou-Ou-MoraisDefinitionRsQLCT} and \eqref{R2016Kou-Ou-MoraisDefinitionTsQLCT}. Adopting the method of Chen\cite{chen2015pitt},  Kou et al. \cite{kou2016uncertainty} considered the following symmetrical representation
$$\mathcal{L}^{ij}_{M_1,M_2}f=\mathcal{L}^{ij}_{M_1,M_2}f_0+\mathcal{L}^{ij}_{M_1,M_2}f_1i+\mathcal{L}^{ij}_{M_1,M_2}f_2j+i\mathcal{L}^{ij}_{M_1,M_2}f_3j$$ 
and introduced the $Q-$modules and $L^2-$norm of QLCT $\mathcal{L}^{ij}_{A_1,A_2}f$ respectively as 
\begin{align}\label{R2016_Kou-Ou-Morais_Qmodulusof2sQLCT}
|\mathcal{L}^{ij}_{M_1,M_2}f|_Q=\left(|\mathcal{L}^{ij}_{M_1,M_2}f_0|^2+|\mathcal{L}^{ij}_{M_1,M_2}f_1|^2+|\mathcal{L}^{ij}_{M_1,M_2}f_2|^2+|\mathcal{L}^{ij}_{M_1,M_2}f_3|^2\right)^{\frac{1}{2}}
\end{align}
and 
\begin{align}\label{R2016_Kou-Ou-Morais_QNormof2sQLCT}
\|\mathcal{L}^{ij}_{M_1,M_2}f\|_{Q,p}=\left(\int_{\mathbb{R}^2}\left|\left(\mathcal{L}^{ij}_{M_1,M_2}f\right)(\uu_1,\uu_2)\right|^p_Qd\uu_1d\uu_2\right)^{\frac{1}{p}}
\end{align}
to prove the Parseval formula
$$\|f\|^2_{L^2_{\mathbb{H}}(\mathbb{R}^2)}=\|\mathcal{L}^{ij}_{M_1,M_2}f\|^2_{Q,2},~f\in L^2_{\mathbb{H}}(\mathbb{R}^2).$$
Bahri et al. \cite{bahri2016simplified} proved some important properties like linearity, shifting, and modulation followed by the continuity of the QLCT of the function in $L^1_{\mathbb{H}}(\mathbb{R}^2).$ Achak et al.\cite{achak2019uncertainty} proves that  Riemann-Lebesgue lemma i.e., 
$$\lim_{|\uu_1|\rightarrow\infty}|\left(\mathcal{L}^{ij}_{M_1,M_2}f\right)(\uu_1,\uu_2)|=0~\mbox{and}~\lim_{|\uu_2|\rightarrow\infty}|\left(\mathcal{L}^{ij}_{M_1,M_2}f\right)(\uu_1,\uu_2)|=0,~f\in L^1_{\mathbb{H}}(\mathbb{R}^2)$$ and the Hausdorff Young inequality
\begin{align}\label{R2019Achak-Abouelaz-SafouaneHYITsQLCT}
\|\mathcal{L}^{ij}_{M_1,M_2}f\|_{Q,r'}\leq\frac{|\Mb_1\Mb_2|^{\frac{1}{r'}-\frac{1}{2}}}{2\pi}\|f\|_{L^r_{\mathbb{H}}(\mathbb{R}^2)},
\end{align}
where $1\leq r\leq 2,$ $\frac{1}{r}+\frac{1}{r'}=1$ and $f\in L^r_{\mathbb{H}}(\mathbb{R}^2).$ In particular, $\|\mathcal{L}^{ij}_{M_1,M_2}f\|_{Q,\infty}\leq\frac{|\Mb_1\Mb_2|^{-\frac{1}{2}}}{2\pi}\|f\|_{L^1_{\mathbb{H}}(\mathbb{R}^2)}.$ Using the orthogonal plane split method Bahri \cite{bahri2019simplex} obtained the simplex complex form
\begin{align}\label{R2019BahriQuaterniontoComplexLCT-TsQLCT}
\left(\mathcal{L}^{ij}_{M_1,M_2}f\right)(\uu_1,\uu_2)=\left(\mathcal{L}^{ij}_{M_1,M_2}f_+\right)(\uu_1,\uu_2)+\left(\mathcal{L}^{ij}_{M_1,M_2}f_-\right)(\uu_1,\uu_2)
\end{align}
and, using \eqref{RModularIdentityQuaternion}, the corresponding modular form
\begin{align}\label{R2019BahriModularIdentityTsQLCT}
\left|\left(\mathcal{L}^{ij}_{M_1,M_2}f\right)(\uu_1,\uu_2)\right|^2=\left|\left(\mathcal{L}^{ij}_{M_1,M_2}f_+\right)(\uu_1,\uu_2)\right|^2+\left|\left(\mathcal{L}^{ij}_{M_1,M_2}f_-\right)(\uu_1,\uu_2)\right|^2,
\end{align}
where
$$\left(\mathcal{L}^{ij}_{M_1,M_2}f_{\pm}\right)(\uu_1,\uu_2)=\frac{1}{2\pi\sqrt{-\Mb_1(\mp \Mb_2)}}\int_{\mathbb{R}^2}e^{\frac{i}{2\Mb_1}(\Ma_1\xx_1^2+\Md_1\uu_1^2-2\xx_1\uu_1)\mp\frac{i}{2\Mb_2}(\Ma_2\xx_2^2+\Md_2\uu_2^2-2\xx_2\uu_2)}f_{\pm}(\xx_1,\xx_2)d\xx_1d\xx_2.$$
Kundu et al. \cite{kundu2022uncertainty} obtained the Plancherel theorem with respect the scalar product \eqref{RQuaternionScalerProduct} as $\langle \mathcal{L}^{ij}_{M_1,M_2}f ,\mathcal{L}^{ij}_{M_1,M_2}g\rangle=\langle f,g\rangle,~f,g\in L^2_{\mathbb{H}}(\mathbb{R}^2).$ Also taking the advantage of \eqref{R2019BahriModularIdentityTsQLCT} and the sharp Hausdorff-Young inequality for 2D LCT, they obtained the following sharp Hausdorff Young inequality 
\begin{align}\label{R2022Kundu-PrasadSHYITsQLCT}
\|\mathcal{L}^{ij}_{M_1,M_2}f\|_{L^q_{\mathbb{H}}(\mathbb{R}^2)}\leq|\Mb_1\Mb_2|^{\frac{1}{q}-\frac{1}{2}}\left(\frac{p^{\frac{1}{p}}}{q^{\frac{1}{q}}}\right)(2\pi)^{\frac{1}{q}-\frac{1}{p}}\|f\|_{L^p_{\mathbb{H}}(\mathbb{R}^2)}
\end{align}
where $1\leq p\leq 2,$ $\frac{1}{p}+\frac{1}{q}=1$ and $f\in L^p_{\mathbb{H}}(\mathbb{R}^2).$ Extending the idea of Guanlie et al.\cite{guanlei2009new} to prove \eqref{R2009Guanlie-Xiaotong-XiaogangPitt'sInequalityin2LCTDomain}, authors in \cite{kundu2022uncertainty} obtained the Pitt's inequality for  LCT of 2D function using which they obtained Pitt's inequality for the QLCT for $f\in  S(\mathbb{R}^2,\mathbb{H})$ and $0\leq\lambda<2$ as
\begin{align}\label{R2022Kundu-PrasadLocalPitt'sInequality}
\int_{\mathbb{R}^2}|\bldu|^{-\lambda}\left|\left(\mathcal{L}^{ij}_{M_1,M_2}f\right)(\bldu)\right|^2d\bldu\leq\frac{R_\lambda}{|\Mb_1\Mb_2|^\lambda}\int_{\mathbb{R}^2}|\bldt|^\lambda|f(\bldt)|^2d\bldt
\end{align}
where $R_\lambda=\left(\frac{\Gamma\left(\frac{1-\lambda}{4}\right)}{\Gamma\left(\frac{1+\lambda}{4}\right)}\right)^2.$\\

\noindent\textbf{Uncertainty Principles:} By differentiating the Pitt's inequality given by \eqref{R2022Kundu-PrasadLocalPitt'sInequality} for LCT of 2D function, Kundu et al. \cite{kundu2022uncertainty} obtained logarithmic UP for the LCT of two-dimensional complex valued function and using the orthogonal plane split method they obtained the following logarithmic UP for $f\in L^2_{\mathbb{H}}(\mathbb{R}^2)$ 
\begin{align}\label{R2022Kundu-PrasadLogarithmicUncertiantyPrinciple}
\int_{\mathbb{R}^2}\ln|\bldu|\left|\left(\mathcal{L}^{ij}_{M_1,M_2}f\right)(\bldu)\right|^2d\bldu+\int_{\mathbb{R}^2}\ln|\bldt||f(\bldt)|^2d\bldt\geq K_0'\int_{\mathbb{R}^2}|f(\bldt)|^2d\bldt,
\end{align}
where $K_0'=\frac{d}{d\lambda}(K_\lambda)$ at $\lambda=0$ and $K_\lambda=\frac{R_\lambda}{|b_1b_2|^\lambda}.$ Furthermore, they also obtained the local UP for the QLCT,
\begin{align}\label{R2022Kundu-PrasadLocalUPQLCT0<alpha<1}
\int_{E}\left|\left(\mathcal{L}^{ij}_{M_1,M_2}f\right)(\bldu)\right|^2d\bldu\leq \frac{K_\lambda}{|\Mb_1\Mb_2|}(2\pi)^{1-2\lambda}|E|^\lambda\||\bldt|^\lambda f\|_{L^2_{\mathbb{H}}(\mathbb{R}^2)}^2,~0<\lambda<1
\end{align}
and
\begin{align}\label{R2022Kundu-PrasadLocalUPQLCTalpha>1}
\int_{E}\left|\left(\mathcal{L}^{ij}_{M_1,M_2}f\right)(\bldu)\right|^2d\bldu\leq\frac{K_\lambda}{|\Mb_1\Mb_2|}|E|\|f\|_{L^2_{\mathbb{H}}(\mathbb{R}^2)}^{2-\frac{2}{\lambda}}\||\bldt|^\lambda f\|_{L^2_{\mathbb{H}}(\mathbb{R}^2)}^\frac{2}{\lambda},~\lambda>1,
\end{align}
where $K_\lambda$ is given equation \eqref{R2018LianLocalUPConstant},
%$$K_\lambda=
%\begin{cases}
%\frac{(1+\lambda)^2}{\lambda^{2\lambda}}(2-2\lambda)^{\lambda-2},&0<\lambda<1\\
%\frac{\pi}{\lambda\Gamma(\frac{1}{\lambda})}\Gamma(\frac{1}{\lambda})\Gamma(1-\frac{1}{\lambda})(\lambda-1)^\lambda(1-\frac{1}{\lambda})^{-1},&\lambda>1,
%\end{cases}
%$$
$\Mb_1\Mb_2\neq 0,$ $E\subset\mathbb{R}^2$ is measurable with measure $|E|$ and $f\in L^2_{\mathbb{H}}(\mathbb{R}^2).$ Without using the Hausdorff-Young inequality of the QLCT, authors in \cite{kundu2022uncertainty} obtained the entropy UP for $f\in S(\mathbb{R}^2,\mathbb{H})$ as
\begin{align}\label{R2022Kundu-PrasadEntropyUncertaintyPrincipleTsQFT}
\mathcal{H}(|f(\bldt)|^2)+\mathcal{H}\left(\left|\left(\mathcal{L}^{ij}_{M_1,M_2}f\right)(\bldu)\right|^2\right)\geq 1+\ln\left(\left|\frac{\Mb_1\Mb_2}{2}\right|\right),
\end{align}
where $\mathcal{H}(P(\bldt))=-\int_{\mathbb{R}^2}P(\bldt)\ln\left(P(\bldt)\right)d\bldt,$ where $P$ is a probability density function. Based on the relation between quaternion LCT and the  quaternion FT \cite{bahri2016simplified}, Zhang et al.\cite{zhang2018novel} obtained the directional UP for $f\in L^1_{\mathbb{H}}(\mathbb{R}^2)\cap L^2_{\mathbb{H}}(\mathbb{R}^2)$ as
\begin{align}\label{R2018_Zhang-LiUncertaintyPrinciple2sQLCT}
\int_{\mathbb{R}^2}(\boldsymbol\alpha\cdot\bldx)^2|f(\bldx)|^2d\bldx\cdot\int_{\mathbb{R}^2}(\boldsymbol\gamma\cdot\bldu)^2|\mathcal{L}_{M_1,M_2}^{\mathbb{H}}\{f\}(\bldu)|^2d\bldu\geq\frac{1}{4}[(\boldsymbol\alpha\cdot \boldsymbol b\boldsymbol\gamma)^2F_{-}^2+(\boldsymbol\alpha\cdot\boldsymbol b\boldsymbol\gamma')^2F_{+}^2],
\end{align}
where $\boldsymbol b\boldsymbol\gamma=(\Mb_1\gamma_1,\Mb_2\gamma_2),$ constant vectors $\boldsymbol\alpha,\boldsymbol\gamma\in\mathbb{R}^2,$
 $\boldsymbol\gamma'=-\gamma_1\boldsymbol e_1+\gamma \boldsymbol e_2,$ $F_{-}=\int_{\mathbb{R}^2}|f_{-}(\bldx)|^2d\bldx$ and $F_{+}=\int_{\mathbb{R}^2}|f_{+}(\bldx)|^2d\bldx.$ Choosing $\boldsymbol\alpha=\boldsymbol\gamma=\boldsymbol e_k,~k=1,2$ equation \eqref{R2018_Zhang-LiUncertaintyPrinciple2sQLCT} reduces to the component-wise UP 
\begin{align}\label{R2016Kou-Ou-MoraisComponentWiseUPTsQLCT}
\int_{\mathbb{R}^2}\xx_k^2|f(\bldx)|^2d\bldx\cdot\int_{\mathbb{R}^2}\uu_k^2|\mathcal{L}_{M_1,M_2}^{\mathbb{H}}\{f\}(\bldu)|^2d\bldu\geq\frac{\Mb_k^2}{4}F^2,~F^2=F_{-}^2+F_{+}^2
\end{align}
which was obtained in \cite{kou2016uncertainty}. Moreover, if $M_l=(\Ma_l,\Mb_l,\Mc_l,\Md_l)=(0,1,-1,0),~l=1,2$ then \eqref{R2018_Zhang-LiUncertaintyPrinciple2sQLCT} reduces to the directional UP proposed by Hitzer for two-sided QFT in \cite{hitzer2010directional}. Author in \cite{zhang2018novel} also obtained the modified UP for $f\in \mathcal{S}(\mathbb{R}^2,\mathbb{H})$ as
\begin{align}\label{R2018_Zhang-LiModifiedUP2sQLCT}
\int_{\mathbb{R}^2}\xx_1^2|f(\bldx)|^2\xx_2^2d\bldx\cdot\int_{\mathbb{R}^2}\uu_1^2|\mathcal{L}_{M_1,M_2}^{\mathbb{H}}\{f\}(\bldu)|^2\uu_2^2d\bldu\geq(\Mb_1\Mb_2)^2\left|\frac{1}{2}F-\int_{\mathbb{R}^2}|\xx_1J^{f}_{M_1,M_2}(\bldx)\xx_2|d\bldx\right|^2,
\end{align}
where $F=\int_{\mathbb{R}^2}|f(\bldx)|^2d\bldx$ and $$J^{f}_{M_1,M_2}(\bldx)=\frac{\Ma_1\xx_1}{\Mb_1}\overline{\frac{\partial f(\bldx)}{\partial\xx_2}}if(\bldx)+\overline{\frac{\partial f(\bldx)}{\partial\xx_2}}\frac{\partial f(\bldx)}{\partial\xx_1}-\frac{\Ma_2\xx_2}{\Mb_2}j\overline{f(\bldx)}if(\bldx)\frac{\Ma_1\xx_1}{\Mb_1}-\frac{\Ma_2\xx_2}{\Mb_2}j\overline{f(\bldx)}\frac{\partial f(\bldx)}{\partial\xx_1}.$$
which generalizes the modified UP for the two-sided QFT by Bahri \cite{bahri2016modified}. Equality in \eqref{R2018_Zhang-LiModifiedUP2sQLCT} holds for
$$f(\bldx)=\sqrt{\Mb_1}e^{i\frac{\pi}{4}}e^{-i\frac{\Ma_1}{2\Mb_1}\xx_1^2}K_0e^{-\left(\frac{\xx_1^2}{2\alpha_1}+\frac{\xx_2^2}{2\alpha_2}\right)}E^{-j\frac{\Ma_2}{2\Mb_2}\xx_2^2}e^{j\frac{\pi}{4}}\sqrt{\Mb_2},$$
where $\alpha_1>0,\alpha_2>0$ are real constants and $K_0$ is a quaternion constant.

$$\mathcal{S}(\mathbb{R}^2,\mathbb{H})=\left\{f\in\mathcal{C}^{\infty}(\mathbb{R}^2,\mathbb{H}):\sup_{x\in\mathbb{R}^2}(1+|\bldx|^k)|\partial^{\boldsymbol\alpha}f(\bldx)|<\infty\right\},$$
where $\bldx^{\boldsymbol\alpha}=\xx_1^{\alpha_1}\xx_2^{\alpha_2},$ $\partial^{\boldsymbol\alpha}=\frac{\partial^{|\boldsymbol\alpha|}}{\partial^{\alpha_1}\xx_1\partial^{\alpha_2}\xx_2},$ $\mathcal{C}^\infty(\mathbb{R}^2,\mathbb{H})$ is the set of smooth function from $\mathbb{R}^2$ to $\mathbb{H}.$
\\
\\
If $f\in L^2_{\mathbb{H}}(\mathbb{R}^2)$ such that
\begin{align}\label{R2018_Zhang-LiRelationb/w2sQLCTand2sFTin2LCDomain}
\mathcal{F}_Q\{g\}(\uu_1,\uu_2)=e^{-i\frac{d_5}{2b_5}\uu_1^2}\mathcal{L}^\mathbb{H}_{M_1,M_2}\{f\}(\uu_1,\uu_2)e^{-j\frac{d_6}{2b_6}\uu_2^2},
\end{align}
then 
\begin{align}\label{R2018_Zhang-LiRelationb/w2sQLCTand2sFTin2LCDomain&function_g}
g\left(\frac{\tee_1}{b_5},\frac{\tee_2}{b_6}\right)=\sqrt{b_5b_6}e^{i\frac{a_5}{2b_5}\tee_1^2-i\frac{\pi}{4}}\mathcal{L}^\mathbb{H}_{M_3,M_4}\{f\}(\tee_1,\tee_2)e^{j\frac{a_6}{2b_6}\tee_2^2-j\frac{\pi}{4}},
\end{align}
where $M_3=M_5^{-1}M_1,$ $M_4=M_6^{-1}M_2.$ Moreover, $\left|g\left(\frac{\tee_1}{b_5},\frac{\tee_2}{b_6}\right)\right|=\sqrt{b_5b_6}\left|\mathcal{L}^\mathbb{H}_{M_5^{-1}M_1,M_6^{-1}M_2}\{f\}(\tee_1,\tee_2)\right|,$ where $b_5=\Ma_3\Mb_1-\Ma_1\Mb_3$ and $b_6=\Ma_4\Mb_2-\Ma_2\Mb_4.$ This result was proved by Zhang et al. in \cite{zhang2018generalized}. Here, $\mathcal{L}^\mathbb{H}_{M_3,M_4}\{f\}(\tee_1,\tee_2)$ is the inverse two-sided quaternion LCT of $\mathcal{L}^\mathbb{H}_{M_1,M_2}\{f\}(\uu_1,\uu_2)$ with the parameters $M_5,M_6.$ The function $g(\tee_1,\tee_2)$ serves as a bridge to obtain the functional relation between $\mathcal{L}^\mathbb{H}_{M_1,M_2}\{f\}(\uu_1,\uu_2)$ and $\mathcal{L}^\mathbb{H}_{M_3,M_4}\{f\}(\tee_1,\tee_2).$ Note that this sort of relation is only available for the two-sided quaternion LCT. Relations \eqref{R2018_Zhang-LiRelationb/w2sQLCTand2sFTin2LCDomain} and \eqref{R2018_Zhang-LiRelationb/w2sQLCTand2sFTin2LCDomain&function_g} are important, as they help in obtaining several UPs using the same for the quaternion FT. For $f\in L^1_{\mathbb{H}}(\mathbb{R}^2)\cap L^2_{\mathbb{H}}(\mathbb{R}^2)$ with unit energy, they obtained the  UP, which presents lower bounds for the product of the spread of the quaternion-valued function in two different quaternion linear canonical domains, as
\begin{align}\label{R2008_Zhang-Li_GeneralizedUPs2sQLCT}
\int_{\mathbb{R}^2}\uu_l^2|\mathcal{L}^{\mathbb{H}}_{M_1,M_2}\{f\}(\bldu)|^2d\bldu\cdot\int_{\mathbb{R}^2}\vv_l^2|\mathcal{L}^{\mathbb{H}}_{M_3,M_4}\{f\}(\bldv)|^2d\bldv\geq P_l,~l=1,2
\end{align}
where $P_1=\frac{(\Ma_4\Mb_2-\Ma_2\Mb_4)^2}{4}$ and $P_2=\frac{(\Ma_3\Mb_1-\Ma_1\Mb_3)^2}{4}.$ In particular, if $M_3=M_4=(1,0,0,1),$ \eqref{R2008_Zhang-Li_GeneralizedUPs2sQLCT} reduces to \eqref{R2016Kou-Ou-MoraisComponentWiseUPTsQLCT} for $F=1.$ Using the relation \eqref{R2018_Zhang-LiRelationb/w2sQLCTand2sFTin2LCDomain} and the Hausdorff-Young inequality for the two-sided quaternion FT, Zhu et al. \cite{zhu2021uncertainty} introduced the Hausdorff-Young inequality for the quaternion LCT in two linear canonical domain as
\begin{align}\label{R2020_Zhu-Zheng_HausdorffYoungInequality2sQLCT}
\|\mathcal{L}^\mathbb{H}_{M_1,M_2}{f}(\bldu)\|_{L^q_{\mathbb{H}}(\mathbb{R}^2)}\leq \left(\frac{p}{2\pi}\right)^{\frac{1}{p}}\left(\frac{q}{2\pi}\right)^{-\frac{1}{q}}|b_5b_6|^{\frac{1}{2}-\frac{1}{p}}\|\mathcal{L}^\mathbb{H}_{A_3,A_4}\{f\}(\bldt)\|_{L^p_{\mathbb{H}}(\mathbb{R}^2)},
\end{align}
$b_5=\Ma_3\Mb_1-\Ma_1\Mb_3$ and $b_6=\Ma_4\Mb_2-\Ma_2\Mb_4,$ which is sharper than the one obtained in \cite{zhu2020uncertainty}. Zhu et al. also obtained the logarithmic UP in two LCT domains using the relation \eqref{R2018_Zhang-LiRelationb/w2sQLCTand2sFTin2LCDomain}, for $f\in \mathcal{S}(\mathbb{R}^2,\mathbb{H}),$ as 
\begin{align}\label{R2021_Zhu-Zheng_LogarithmicUP2sQLCTfor2LCDomain_Q}
\int_{\mathbb{R}^2}\ln\left|\frac{\bldt}{\boldsymbol b}\right||\mathcal{L}^\mathbb{H}_{M_3,M_4}\{f\}(\bldt)|^2d\bldt+\int_{\mathbb{R}^2}\ln|\bldu||\mathcal{L}^\mathbb{H}_{M_1,M_2}\{f\}(\bldu)|_Q^2d\bldu\geq D\int_{\mathbb{R}^2}|\mathcal{L}^\mathbb{H}_{M_3,M_4}\{f\}(\bldt)|^2d\bldt
\end{align}
%and \textcolor{red}{the following is same as \eqref{R2021_Zhu-Zheng_LogarithmicUP2sQLCTfor2LCDomain_Q} but norm $|\cdot|_Q$ \eqref{R2016_Kou-Ou-Morais_Qmodulusof2sQLCT} is replaced with $|\cdot|,$ so this can be written in language } 
%\begin{align}\label{R2021_Zhu-Zheng_LogarithmicUP2sQLCTfor2LCDomain}
%\int_{\mathbb{R}^2}\ln\left|\frac{\bldt}{\boldsymbol b}\right||\mathcal{L}^\mathbb{H}_{A_3,A_4}\{f\}(\bldt)|^2d\bldt+\int_{\mathbb{R}^2}\ln|\bldu||\mathcal{L}^\mathbb{H}_{A_1,A_2}\{f\}(\bldu)|^2d\bldu\geq D\int_{\mathbb{R}^2}|\mathcal{L}^\mathbb{H}_{A_3,A_4}\{f\}(\bldt)|^2d\bldt
%\end{align}
where $\frac{\bldt}{\boldsymbol b}=\left(\frac{\tee_1}{b_5},\frac{\tee_2}{b_6}\right)$ and $D=\psi(\frac{1}{2})+\ln2,$ $\psi=\frac{d}{d\tee}(\ln(\Gamma(\tee))),$ and the module $|\cdot|_Q$ on the second term of the left-hand side in equation \eqref{R2021_Zhu-Zheng_LogarithmicUP2sQLCTfor2LCDomain_Q}, given by \eqref{R2016_Kou-Ou-Morais_Qmodulusof2sQLCT} can also be replaced by $|\cdot|$ given by \eqref{RQuaternionModules}. Using the Heisenberg UP for the two-sided quaternion FT for the function $f\in L^2_{\mathbb{H}}(\mathbb{R}^2)$ and the relation \eqref{R2018_Zhang-LiRelationb/w2sQLCTand2sFTin2LCDomain}, authors in \cite{zhu2021uncertainty} obtained the Heisenberg-type UP
\begin{align}\label{R2021_Zhu-Zheng_HeisenbergUP2sQLCTfor2LCDomain}
\left(\int_{\mathbb{R}^2}(t_l^2)^\alpha|\mathcal{L}^{\mathbb{H}}_{M_3,M_4}\{f\}(\bldt)|^2d\bldt\right)^{\frac{\alpha}{\alpha+\beta}}\left(\int_{\mathbb{R}^2}(u_l^2)^\beta|\mathcal{L}^{\mathbb{H}}_{M_1,M_2}\{f\}(\bldt)|^2d\bldt\right)^{\frac{\alpha}{\alpha+\beta}}\geq P_l\left(\int_{\mathbb{R}^2}|\mathcal{L}^{\mathbb{H}}_{M_3,M_4}|\{f\}(\bldt)|^2d\bldt\right),~l=1,2
\end{align}
for $f\in L^1_{\mathbb{H}}(\mathbb{R}^2)\cap L^2_{\mathbb{H}}(\mathbb{R}^2)$ and $\alpha,\beta\geq 1,$ where $P_1=\left(\frac{(\Ma_3\Mb_1-\Ma_1\Mb_3)^2}{4}\right)^{\frac{\alpha\beta}{\alpha+\beta}}$ and $P_2=\left(\frac{(\Ma_4\Mb_2-\Ma_2\Mb_4)^2}{4}\right)^{\frac{\alpha\beta}{\alpha+\beta}},$ followed by the local UP associated with the two-sided quaternion LCT in two different linear canonical domain for all measurable subset $E$ of $\mathbb{R}^2$ as
\begin{align}\label{R2021_Zhu-Zheng_LocalUP2sQLCTfor2LCDomain0<Alpha<1}
\int_{E}|\mathcal{L}^{\mathbb{H}}_{M_1,M_2}\{f\}(\bldu)|^2d\bldu\leq K_\lambda|E|^\lambda\int_{\mathbb{R}^2}\left|\frac{\bldt}{\boldsymbol b}\right|^{2\lambda}|\mathcal{L}^{\mathbb{H}}_{M_3,M_4}\{f\}(\bldt)|^2d\bldt,~0<\lambda<1,
\end{align}
and
\begin{align}\label{R2021_Zhu-Zheng_LocalUP2sQLCTfor2LCDomainAlpha>1}
\int_{E}|\mathcal{L}^{\mathbb{H}}_{M_1,M_2}\{f\}(\bldu)|^2d\bldu\leq K_\lambda|E|^\lambda\|\mathcal{L}^{\mathbb{H}}_{M_3,M_4}\{f\}(\bldt)\|^{2-2\lambda}_2\left(\int_{\mathbb{R}^2}\left|\frac{\bldt}{\boldsymbol b}\right|^{2\lambda}|\mathcal{L}^{\mathbb{H}}_{M_3,M_4}\{f\}(\bldt)|^2d\bldt\right)^\frac{1}{\lambda},~\lambda>1
\end{align}
where $K_\lambda$ is given by \eqref{R2018LianLocalUPConstant} and $\boldsymbol b=(b_5,b_6).$ The inequalities given by equations \eqref{R2008_Zhang-Li_GeneralizedUPs2sQLCT},\eqref{R2008_Zhang-Li_GeneralizedUPs2sQLCT},\eqref{R2021_Zhu-Zheng_LogarithmicUP2sQLCTfor2LCDomain_Q},\eqref{R2021_Zhu-Zheng_HeisenbergUP2sQLCTfor2LCDomain},\eqref{R2021_Zhu-Zheng_LocalUP2sQLCTfor2LCDomain0<Alpha<1} and \eqref{R2021_Zhu-Zheng_LocalUP2sQLCTfor2LCDomainAlpha>1} reduces to the corresponding UP for the two-sided quaternion LCT in the spatial and frequency domain when the parameters $M_5=M_1,$ $M_6=M_2$ which yields $M_3=M_4=(1,0,0,1),$ if further  $M_1=(0,1,0,-1)=M_2=(0,1,0,-1)$ then these inequalities leads to the corresponding uncertainty inequalities for the two-sided quaternion FT. Even though the standard deviation is widely accepted as a measure of uncertainty, it is not an ideal tool to express UP in some situations; some examples are cited in \cite{bialynicki2011entropic} to express limitations on the use of standard deviation. Thus, entropy was proposed by Shannon and R\`enyi. In particular, R\`enyi entropy is one parameter extension of Shannon entropy.
\begin{definition}
(R\`enyi entropy) The entropy of a probability density function $P$ on $\mathbb{R}^2$ is defined by
$$H_\alpha(P)=\frac{1}{1-\alpha}\ln\left(\int_{\mathbb{R}^2}[P(\bldt)]^\alpha d\bldt\right),~\alpha\geq 1.$$
When $\alpha\to 1,$ R\`enyi entropy leads to Shannon entropy, given by $E(P)=-\int_{\mathbb{R}^2}\ln(P(\bldt))P(\bldt)d\bldt.$
\end{definition}
Using the Hausdorff-Young inequality \eqref{R2020_Zhu-Zheng_HausdorffYoungInequality2sQLCT}, the R\`enyi UPs for the two-sided quaternion LCT in different linear canonical domains are obtained in \cite{zhu2021uncertainty} as
\begin{align}\label{R2021_Zhu-Zheng_RenyiEntropyUP2sQLCTfor2LCDomain}
H_\alpha(|\mathcal{L}^\mathbb{H}_{M_3,M_4}{f}(\bldt)|^2)+H_\beta(|\mathcal{L}^\mathbb{H}_{M_1,M_2}{f}(\bldu)|^2)\geq\frac{\ln\alpha}{\alpha-1}+\frac{\ln\beta}{\beta-1}+2\ln\pi+\ln|b_5b_6|,
\end{align}
where $\alpha$ and $\beta$ are such that $\frac{1}{\alpha}+\frac{1}{\beta}=2$ and $f\in L^2_{\mathbb{H}}(\mathbb{R}^2)$ such that $\|f\|_{L^2_{\mathbb{H}}(\mathbb{R}^2)}=1.$ When $\alpha\to 1$ and $\beta\to 1,$ then \eqref{R2021_Zhu-Zheng_RenyiEntropyUP2sQLCTfor2LCDomain} becomes the Shannon UP
\begin{align}\label{R2021_Zhu-Zheng_ShannonEntropyUP2sQLCTfor2LCDomain}
E(|\mathcal{L}^\mathbb{H}_{M_3,M_4}{f}(\bldt)|^2)+E(|\mathcal{L}^\mathbb{H}_{M_1,M_2}{f}(\bldu)|^2)\geq 2+2\ln\pi+\ln|b_5b_6|,
\end{align}
in two different linear canonical domain. Equations \eqref{R2021_Zhu-Zheng_RenyiEntropyUP2sQLCTfor2LCDomain} and \eqref{R2021_Zhu-Zheng_ShannonEntropyUP2sQLCTfor2LCDomain} reduces to the corresponding UPs for the quaternion LCT in the spatial and frequency domain when $M_5=M_1,M_6=M_2$  and to that of the quaternion FT, if further $M_3=M_4=(1,0,0,1).$
%\cite{bahri2016simplified},\cite{zhu2020uncertainty}.

\section{Quaternion Quadratic Phase Fourier Transform (QQPFT)}
Recently, Gupta and Verma \cite{gupta2022short} extended QPFT \eqref{RDefinitionQPFT} to the quaternion-valued functions. They defined the two-sided quaternion quadratic phase FT by replacing the kernels $\mathcal{K}^i(\tee_1,\xxi_1)=\frac{1}{\sqrt{2\pi}} e^{-i\tee_1,\xxi_1}$ and $\mathcal{K}^j(\tee_2,\xxi_2)=\frac{1}{\sqrt{2\pi}}e^{-j\tee_2,\xxi_2},$ in the definition 
\begin{align}
(\mathcal{F}_{\mathbb{H}}f)(\bldxi)=\int_{\mathbb{R}^2}\mathcal{K}^i(\tee_1,\xxi_1)f(\bldt)\mathcal{K}^j(\tee_2,\xxi_2)d\bldt,~\bldxi=(\xxi_1,\xxi_2)\in\mathbb{R}^2,
\end{align}
of two-sided quaternion FT by the kernels in the quadratic phase FT given in \eqref{RDefinitionofKernel_QPFT} as
\begin{align}\label{RGupta-VermaDefinition2sQQPFT}
(\mathcal{Q}^{\wedge_1,\wedge_2}_{\mathbb{H}} f)(\bldxi)=\int_{\mathbb{R}^2}\mathcal{K}^i_{\wedge_1}(\tee_1,\xxi_1)f(\bldt)\mathcal{K}^j_{\wedge_2}(\tee_2,\xxi_2)d\bldt,~\bldxi=(\xxi_1,\xxi_2)\in\mathbb{R}^2
\end{align}
for $f(\bldt)\in L^2_{\mathbb{H}}(\mathbb{R}^2),~\bldt=(\tee_1,\tee_2),$ where $\wedge_l=(\LmA_l,\LmB_l,\LmC_l,\LmD_l,\LmE_l),\LmA_l,\LmB_l,\LmC_l,\LmD_l,\LmE_l\in\mathbb{R} ~\mbox{and}~\LmB_l\neq 0~\mbox{for}~l=1,2$
\begin{align}\label{RGupta_VermaKernelofQQPFTinI}
\mathcal{K}^i_{\wedge_1}(\tee_1,\xxi_1)=\frac{1}{\sqrt{2\pi}}e^{-i\left(\LmA_1\tee_1^2+\LmB_1\tee_1\xxi_1+\LmC_1\xxi_1^2+\LmD_1\tee_1+\LmE_1\xxi_1\right)}
\end{align}
and 
\begin{align}\label{RGupta_VermaKernelofQQPFTinJ}
\mathcal{K}^j_{\wedge_2}(\tee_2,\xxi_2)=\frac{1}{\sqrt{2\pi}}e^{-j\left(\LmA_2\tee_2^2+\LmB_2\tee_2\xxi_2+\LmC_2\xxi_2^2+\LmD_2\tee_2+\LmE_2\xxi_2\right)}
\end{align}
The corresponding inversion formula is given by
\begin{align}\label{RGupta_VermaInversionof2sQQPFT}
f(\bldt)=|\LmB_1\LmB_2|\int_{\mathbb{R}^2}\overline{\mathcal{K}^i_{\wedge_1}(\tee_1,\xxi_1)}(\mathcal{Q}^{\wedge_1,\wedge_2}_{\mathbb{H}} f)(\bldxi)\overline{\mathcal{K}^j_{\wedge_2}(\tee_2,\xxi_2)}d\bldxi.
\end{align}
Based on the relation between quaternion QPFT and the quaternion FT
\begin{align}\label{RGupta_VermaRelationb/w2sQQPFT&2sQFT}
(\mathcal{Q}^{\wedge_1,\wedge_2}_{\mathbb{H}} f)(\bldxi)=e^{-i\left(\LmC_1\xxi_1^2+\LmE_1\xxi_1\right)}\left(\mathcal{F}_{\mathbb{H}}\tilde{f}\right)(\LmB_1\xxi_1,\LmB_2\xxi_2)e^{-j\left(\LmC_2\xxi_2^2+\LmE_2\xxi_2\right)},
\end{align}
where $\left(\mathcal{F}_{\mathbb{H}}\tilde{f}\right)(\bldxi)=\int_{\mathbb{R}^2}\frac{1}{\sqrt{2\pi}} e^{-i\tee_1\xxi_1}\tilde{f}(\bldt)\frac{1}{\sqrt{2\pi}}e^{-j\tee_2\xxi_2}d\bldt$ with $\tilde{f}(\bldt)=e^{-i\left(\LmA_1\tee_1^2+\LmD_1\tee_1\right)}f(\bldt)e^{-j\left(\LmA_2\tee_2^2+\LmD_2\tee_2\right)}$ author in \cite{gupta2022short} obtained the Parseval's formula $\langle f,g\rangle=|\LmB_1\LmB_2|\langle \mathcal{Q}^{\wedge_1,\wedge_2}_{\mathbb{H}} f, \mathcal{Q}^{\wedge_1,\wedge_2}_{\mathbb{H}} g \rangle,$ for $f,g\in L^2_{\mathbb{H}}(\mathbb{R}^2)$ (in particular,  $\|f\|^2_{L^2_{\mathbb{H}}(\mathbb{R}^2)}=|\LmB_1\LmB_2|\|\mathcal{Q}^{\wedge_1,\wedge_2}_{\mathbb{H}} f\|^2_{L^2_{\mathbb{H}}(\mathbb{R}^2\times\mathbb{R}^2)}$), and the sharp Hausdorff-Young inequality for the two-sided quaternion QPFT as 
\begin{align}\label{RGupta_VermaSharpHUIfor2sQQPFT}
\|\mathcal{Q}^{\wedge_1,\wedge_2}_{\mathbb{H}} f\|_{L^q_{\mathbb{H}}(\mathbb{R}^2)}\leq \frac{(2\pi)^{\frac{1}{q}-\frac{1}{p}}A_p^2}{|\LmB_1\LmB_2|^\frac{1}{q}}\|f\|_{L^p_{\mathbb{H}}(\mathbb{R}^2)},
\end{align}
where $1\leq p\leq 2,$ $\frac{1}{p}+\frac{1}{q}=1$ and $f\in L^2_{\mathbb{H}}(\mathbb{R}^2)$ and $A_p=\left(\frac{p^{\frac{1}{p}}}{q^{\frac{1}{q}}}\right)^\frac{1}{2}.$ Recently, Bahri et al. \cite{bahri2023some} also obtained the same inversion formula, Parseval's formula and the sharp Hausdorff-Young inequality \eqref{RGupta_VermaSharpHUIfor2sQQPFT} for the QQPFT \eqref{RGupta-VermaDefinition2sQQPFT}. They proved the Riemann-Lebesgue lemma i.e., $\lim_{|\xxi_l|\rightarrow\infty}\left|(\mathcal{Q}^{\wedge_1,\wedge_2}_{\mathbb{H}} f)(\bldxi)\right|=0,~l=1,2$ followed by the Pitt's inequality $\int_{\mathbb{R}^2}|\bldxi|^{\alpha}\left|(\mathcal{Q}^{\wedge_1,\wedge_2}_{\mathbb{H}} f)(\bldxi)\right|^2d\bldxi\leq M(\alpha)|\LmB_1\LmB_2|^{1-\alpha}\int_{\mathbb{R}^2}|\bldt|^\alpha|f(\bldt)|^2d\bldt,~M(\alpha)=\pi^\alpha\left[\frac{\Gamma\left(\frac{2-\alpha}{4}\right)}{\Gamma\left(\frac{2+\alpha}{4}\right)}\right]^2,$ where $f\in \mathcal{S}(\mathbb{R}^2,\mathbb{H}),$ for the QQPFT using the same for the QFT. Bhat et al. \cite{younus2023towards} considered the two-sided quaternion QPFT \eqref{RGupta-VermaDefinition2sQQPFT} where the factor $\frac{1}{\sqrt{2\pi}}$ in the kernels $\mathcal{K}^i_{\wedge_1}(\tee_1,\xxi_1)$ and $\mathcal{K}^j_{\wedge_2}(\tee_2,\xxi_2)$ are replaced respectively with $\sqrt{\frac{\LmB_1i}{2\pi}}$ and $\sqrt{\frac{\LmB_2j}{2\pi}}.$ Let us denote the transform in this case by $^{(1)}\mathcal{Q}^{\wedge_1,\wedge_2}_{\mathbb{H}}.$ They also defined the right-sided and left-sided quaternion QPFT for the functions $L^2_{\mathbb{H}}(\mathbb{R}^2).$ Apart from the properties studied in \cite{gupta2022short}, like the inversion formula, Parseval's formula, and Hausdorff-Youngs inequality, authors also obtained the translation formula given by
\begin{align*}
&\left(^{(1)}\mathcal{Q}^{\wedge_1,\wedge_2}_{\mathbb{H}}f(\bldx-\boldsymbol k)\right)(\bldomega)=e^{-i\left[\LmA_1k_1^2+\LmD_1k_1+\LmB_1k_1\oomega_1-4\frac{\LmA_1^2\LmC_1}{\LmB_1^2}k_1^2-4\frac{\LmA_1}{\LmB_1}\LmC_1\oomega_1k_1-2\frac{\LmA_1}{\LmB_1}k_1\right]}\notag\\
&\hspace{3.5cm}\times\left(^{(1)}\mathcal{Q}^{\wedge_1,\wedge_2}_{\mathbb{H}}f\right)\left(\oomega_1+2\frac{\LmA_1}{\LmB_1}K_1,\oomega_2+2\frac{\LmA_2}{\LmB_2}K_2\right)e^{-j\left[\LmA_2k_2^2+\LmD_2k_2+\LmB_2k_2\oomega_2-4\frac{\LmA_2^2\LmC_2}{\LmB_2^2}k_2^2-4\frac{\LmA_2}{\LmB_2}\LmC_2\oomega_2k_2-2\frac{\LmA_2}{\LmB_2}k_2\right]}
\end{align*}
\begin{align*}
\left(^{(1)}\mathcal{Q}^{\wedge_1,\wedge_2}_{\mathbb{H}} (\mathbb{M}_{\bldomega_0}f)\right)(\bldomega)= e^{i\left[\frac{\LmC_1 u_0^2-2\LmB_1\LmC_1u_0\oomega_1-\LmB_1-\LmB_1 \LmE_1 u_0}{\LmB_1^2}\right]}\left(^{(1)}\mathcal{Q}^{\wedge_1,\wedge_2}_{\mathbb{H}}f\right)\left(\oomega_1-\frac{u_0}{\LmB_1},\oomega_2-\frac{v_0}{\LmB_2}\right)e^{j\left[\frac{\LmC_2 v_0^2-2\LmB_2\LmC_2v_0\oomega_2-\LmB_2-\LmB_2 \LmE_2 v_0}{\LmB_2^2}\right]}.
\end{align*}\\
\noindent\textbf{Uncertainty Principles:} Using the sharp Hausdorff Young inequality \eqref{RGupta_VermaSharpHUIfor2sQQPFT}, the following R\`eyni UP for the quaternion QPFT 
\begin{align}\label{RGupta_VermaRenyiUPfor2sQQPFT}
H_\alpha(|f|^2)+H_\beta\left(\left|\sqrt{|\LmB_1\LmB_2|}\left(\mathcal{Q}^{\wedge_1,\wedge_2}_\mathbb{H}f\right)(\bldxi)\right|^{2}\right)\geq-\log(|\LmB_1\LmB_2|)-2\log(2\pi)-\left(\frac{1}{1-\alpha}\log(2\alpha)+\frac{1}{1-\beta}\log(2\beta)\right),
\end{align}
where $\frac{1}{2}<\alpha<1$ and $\frac{1}{\alpha}+\frac{1}{\beta}=1$ was proved in \cite{gupta2022short}.  Equation \eqref{RGupta_VermaRenyiUPfor2sQQPFT} results in Shannon UP for the quaternion QPFT $E(|f|^2)+E\left(\left|\sqrt{|\LmB_1\LmB_2|}\left(\mathcal{Q}^{\wedge_1,\wedge_2}_\mathbb{H}f\right)(\bldxi)\right|^{2}\right)\geq\log \left(\frac{e^2}{16\pi^2|\LmB_1\LmB_2|}\right),$ when $\alpha\to 1$ and $\beta\to 1.$ It is to be noted that the Hausdorff-Youngs inequality \eqref{RGupta_VermaSharpHUIfor2sQQPFT}, obtained in the \cite{gupta2022short} is more sharper than the one obtained in \cite{younus2023towards}. Based on the relationship between QQPFT and the QFT and using the Heisenberg's and logarithmic UP for the QFT, authors in \cite{younus2023towards} (and also \cite{bahri2023some}) obtained the Heisenberg's UP for the transform $^{(1)}\mathcal{Q}^{\wedge_1,\wedge_2}_{\mathbb{H}}$ given as
\begin{align}\label{R2023_Bhat-Dar_MMAS_QQPFT_Heisenberg_UP}
\left(\int_{\mathbb{R}^2}\xx_l^2|f(\bldx)|^2d\bldx\right)\left(\int_{\mathbb{R}^2}\oomega_l^2\left|\left(^{(1)}\mathcal{Q}^{\wedge_1,\wedge_2}_{\mathbb{H}}f\right)(\bldomega)\right|^2d\bldomega\right)\geq \frac{1}{4\LmB_l^2}\left(\int_{\mathbb{R}^2}|f(\bldx)|^2d\bldx\right)^2,~l=1,2
\end{align}
for $f\in L^2_{\mathbb{H}}(\mathbb{R}^2)$ followed by the logarithmic UP given as 
\begin{align}\label{R2023_Bhat-Dar_MMAS_QQPFT_Logarithmic_UP}
\int_{\mathbb{R}^2}\ln|\bldx||f(\bldx)|^2d\bldx+\int_{\mathbb{R}^2}\ln|\bldomega|\left|\left(^{(1)}\mathcal{Q}^{\wedge_1,\wedge_2}_{\mathbb{H}}f\right)(\bldomega)\right|^2d\bldomega \geq (D-\ln|\boldsymbol{B}|)\int_{\mathbb{R}^2}|f(\bldx)|^2d\bldx,
\end{align}
for the functions $f\in \mathcal{S}(\mathbb{R}^2,\mathbb{H}),$ where $D=\ln(2\pi^2)-2\psi\left(\frac{1}{2}\right),~\psi=\frac{d}{dt}(\ln(\Gamma(t))).$ In fact, a more general version of the Heisenberg UP for $\mathcal{Q}^{\wedge_1,\wedge_2}_{\mathbb{H}}$ can be found in \cite{bahri2023some}. Apart from these UPs, authors in \cite{younus2023towards} also obtained the Hardy's UP and Beurling's UP for the transform $^{(1)}\mathcal{Q}^{\wedge_1,\wedge_2}_{\mathbb{H}},$ using the same for the QFT. Using the relation between $^{(1)}\mathcal{Q}^{\wedge_1,\wedge_2}_{\mathbb{H}}$ and the QFT and its associated Donoho-Stark's UP authors also proved that, if a non-zero $f\in L^2_{\mathbb{H}}(\mathbb{R}^2)$ is $\epsilon_{E_1}-$ concentrated on $E_1\subset \mathbb{R}^2$ and $^{(1)}\mathcal{Q}^{\wedge_1,\wedge_2}_{\mathbb{H}}f$ is $\epsilon_{E_2}-$concentrated on $E_2\subset\mathbb{R}^2,$ then $|E_1| |E_2|\geq \frac{2\pi}{|\boldsymbol{B}|}(1-\epsilon_{E_1}-\epsilon_{E_2})^2.$ This is called the Donoho-Stark's UP for the QQPFT $^{(1)}\mathcal{Q}^{\wedge_1,\wedge_2}_{\mathbb{H}}.$

\section{Windowed Fourier Transform (WFT)}
It is well known that the representation of a function $f$ of $\xx$ is usually called its time representation, while its FT $\hat{f}(\xxi)$ is called the frequency representation. From the application point of view, one often needs information regarding the frequency behavior of a signal at a certain time and the time behavior at certain frequency components of a signal. This leads to constructing several joint time-frequency representations like the Gabor transform, WT, ST, Wigner distribution, etc. Other names for the Gabor transform that are frequently used in the literature are the Weyl-Heisenberg transform, short-time FT, and WFT. The motivations for the development of the WT are of a similar nature. However, it is preferably called the time-scale representation since the parameter $a$ involved in its definition cannot completely be identified with an inverse frequency, as it is often done in the literature. 
\begin{definition}
The Gabor transform of $f\in L^2(\mathbb{R})$ with respect to the window function $\psi$ is defined by 
\begin{align}\label{R2000_Wilczok_DefinitionGaborTransform}
(\mathcal{G}_{\psi}f)(\oomega,\tee)=\frac{1}{\sqrt{2\pi}}\int_{\mathbb{R}}f(\xx)\overline{\psi(\xx-\tee)}e^{-i\oomega\xx}d\xx.
\end{align}
\end{definition}
Equivalently, \eqref{R2000_Wilczok_DefinitionGaborTransform} has several important forms. Before we write them, let us first mention some important operators that we will be frequently using throughout the paper.
\begin{enumerate}
\item \textbf{Translation operator:} \label{RTranslationOperatorFT}$(\tau_{\yy}\psi)(\xx)=\psi(\xx-\yy),~y\in\mathbb{R}.$
\item \textbf{Dilation Operator(s):}\label{RDilationOperatorFT} $(\delta_{\A}\psi)(x)=\frac{1}{\sqrt{|\A|}}\psi\left(\frac{x}{\A}\right)$ and $(\mathcal{D}_{\A}\psi)(x)=|\A|\psi(\A\xx),~\A\in\mathbb{R}_{0}.$
\item \textbf{Modulation Operator:}\label{RModulationOperatorFT} $(\mathcal{M}_{\A}\psi)(\xx)=e^{i\A\xx}\psi(\xx),~\A\in\mathbb{R}.$
\end{enumerate}
Based on these operators, equation \eqref{R2000_Wilczok_DefinitionGaborTransform} can be written as
\begin{align}\label{R2000_Wilczok_DefinitionGaborTransform1OperatorForm}
(\mathcal{G}_{\psi}f)(\oomega,\tee)=\frac{1}{\sqrt{2\pi}}\int_{\mathbb{R}}f(\xx)\overline{(\tau_\tee\psi)(\xx)}e^{-i\oomega\xx}d\xx
\end{align}
\begin{align}\label{R2000_Wilczok_DefinitionGaborTransform2OperatorForm}
(\mathcal{G}_{\psi}f)(\oomega,\tee)=\frac{1}{\sqrt{2\pi}}\int_{\mathbb{R}}f(\xx)\overline{(\mathcal{M}_a\tau_\tee\psi)(\xx)}d\xx
\end{align}
\begin{align}\label{R2000_Wilczok_DefinitionGaborTransformInnerProductForm}
(\mathcal{G}_{\psi}f)(\oomega,\tee)=\langle f,\psi_{\oomega,\tee} \rangle_{L^2(\mathbb{R})},
\end{align}
where $\psi_{\oomega,\tee}(x)=\frac{1}{\sqrt{2\pi}}(\mathcal{M}_a\tau_\tee\psi)(\xx)=\frac{1}{\sqrt{2\pi}}\psi(\xx-\tee)e^{i\oomega\xx}.$ It is  a matter of simple verification that $\mathcal{G}_{\psi}f$ is linear in $f$ and conjugate linear in $\psi.$ The following properties, called the covariance properties of the WFT, give respectively the formula for the WFT of the translation by $\xx_0$ and modulation by $\oomega_0$ of function $f$ with respect to the window function $\psi$
\begin{align}\label{R2000_Wilczok_GaborTransformforTranslateof_f}
[\mathcal{G}_{\psi}(\tau_{\xx_0}f)](\oomega,\tee)=e^{-i\oomega \xx_0}(\mathcal{G}_{\psi}f)(\oomega,\tee-\xx_0)
\end{align}
\begin{align}\label{R2000_Wilczok_GaborTransformforModulationof_f}
[\mathcal{G}_{\psi}(\mathcal{M}_{\oomega_0}f)](\oomega,\tee)=(\mathcal{G}_{\psi}f)(\oomega-\oomega_0,\tee).
\end{align}
Using \eqref{R2000_Wilczok_GaborTransformforTranslateof_f} and \eqref{R2000_Wilczok_GaborTransformforModulationof_f}, one can obtain $[\mathcal{G}_{\psi}(\mathcal{M}_{\oomega_0}\tau_{\xx_0}f)](\oomega,\tee)=e^{-i\oomega \xx_0}(\mathcal{G}_{\psi}f)(\oomega-\oomega_0,\tee-\xx_0),~(\oomega,\tee)\in \mathbb{R}^2.$
$$\int_{\mathbb{R}^2}(\mathcal{G}_{\psi}f)(\oomega,\tee)\overline{(\mathcal{G}_{\psi}g)(\oomega,\tee)}=\|\psi\|^2_{L^2(\mathbb{R})}\langle f,g\rangle_{L^2(\mathbb{R})}.$$
This in particular, for $f=g,$ implies that $\|\mathcal{G}_{\psi}f\|^2_{L^2(\mathbb{R}^2)}=\|\psi\|^2_{L^2(\mathbb{R})}\|f\|^2_{L^2(\mathbb{R})},$ which further shows that the normalized WFT $\frac{1}{\|\psi\|_{L^2(\mathbb{R})}}\mathcal{G}_{\psi}$ is an isometry from $L^2(\mathbb{R})$ into $L^2(\mathbb{R}^2).$ The WFT converts a function $f(\xx)$ of one variable into a function $(\mathcal{G}_{\psi}f)(\oomega,\tee)$ of two variables without changing the total energy.   Also, $\mathcal{G}_{\psi}f$ is in $L^2(\mathbb{R}^2)$ for $f\in L^2(\mathbb{R}),$ i.e., the range $\mathcal{G}_{\psi}\left(L^2(\mathbb{R})\right)$ of WFT of square integrable functions  is a subspace of $L^2(\mathbb{R}^2).$ In fact, it is a proper subspace of $L^2(\mathbb{R}^2),$ because $|(\mathcal{G}_{\psi}f)(\oomega,\tee)|\leq \frac{1}{\sqrt{2\pi}}\|f\|_{L^2(\mathbb{R}^2)}\|\psi\|_{L^2(\mathbb{R}^2)},~\oomega,\tee\in\mathbb{R}$ implies that $\mathcal{G}_{\psi}f$ is bounded, but $L^2(\mathbb{R}^2)$ contains functions that are unbounded. Thus, square integrability of $\mathcal{G}_{\psi}f$ is a necessary condition but not sufficient. The range $\mathcal{G}_{\psi}(L^2(\mathbb{R})),$  of bounded linear operator $\mathcal{G}_{\psi}f: L^2(\mathbb{R})\to L^2(\mathbb{R}^2)$ which is isometric upto some constant, is a closed subspace. Being closed subspace of the Hilbert space $L^2(\mathbb{R}^2),$ $\mathcal{G}_{\psi}(L^2(\mathbb{R}))$ is also a Hilbert space. Moreover, it is a RKHS, with the kernel given by   
\begin{align}\label{R2000_Wilczok_GaborTransformRKHS}
K_{\psi}(\oomega',\tee';\oomega,\tee)=\frac{1}{\|\psi\|^2_{L^2(\mathbb{R})}}\langle \psi_{\oomega,\tee},\psi_{\oomega',\tee'}\rangle_{L^2(\mathbb{R})}
\end{align}
which is pointwise bounded, i.e., $|K_{\psi}(\oomega',\tee';\oomega,\tee)|\leq 1,~\forall ~(\oomega',\tee'),(\oomega,\tee)\in \mathbb{R}^2.$ A strong inequality for the WFT, proved by Lieb \cite{lieb1990integral}, often called Lieb's inequality in the literature, for $f,\psi\in L^2(\mathbb{R}^d)$ and  $2\leq p<\infty$ as
\begin{align}\label{RDonoho-Stark_UP_Strong_Version_STFT}
\int_{\mathbb{R}^{2d}}|\mathcal{V}_{\psi}f(\xx,\oomega)|^pd\xx d\oomega\leq \left(\frac{2}{p}\right)^{d}\left(\|f\|_{L^2(\mathbb{R}^d)}\|\psi\|_{L^2(\mathbb{R}^d)}\right)^p.
\end{align}
For $1\leq p\leq 2,$
%for $f\in L^2(\mathbb{R}^d)$ and $2\leq p\leq\infty.$
\begin{align}
\int_{\mathbb{R}^{2d}}|\mathcal{V}_{\psi}f(\xx,\oomega)|^pd\xx d\oomega\geq \left(\frac{2}{p}\right)^{d}\left(\|f\|_{L^2(\mathbb{R}^d)}\|\psi\|_{L^2(\mathbb{R}^d)}\right)^p,
\end{align}
where $\mathcal{V}_{\psi}f$ for  the functions $f,\psi$ on $\mathbb{R}^d$ are defined similar to equation \eqref{R2000_Wilczok_DefinitionGaborTransform}.
\\\\
%\begin{theorem}
%Suppose that $\|f\|_2=\|g\|_2=1.$ If $U\subset \mathbb{R}^{2d}$ and $\epsilon\geq 0$ are such that 
%$$\int_{\mathbb{R}^{2d}}|\mathcal{V}_{\psi}f(\xx,\oomega)|^2d\xx d\oomega\geq 1-\epsilon,$$
%then $|U|\geq (1-\epsilon)^{\frac{p}{p-2}}\left(\frac{p}{2}\right)^{\frac{2d}{p-2}},$ 
%for all $p>2.$ In particular, $|U|\geq \sup_{p>0}(1-\epsilon)^{\frac{p}{p-2}}\left(\frac{p}{2}\right)^{\frac{2d}{p-2}}\geq (1-\epsilon)^22^d.$
%\end{theorem}
\noindent\textbf{Uncertainty Principles:} Wilczok \cite{wilczok2000new} proved that $\mathcal{G}_{\psi}f$ of a function $f\in L^2(\mathbb{R})\setminus\{0\}$ cannot have a support of finite Lebesgue measure. The author also proved that the portion of $\mathcal{G}_{\psi}f$ lying outside some set $M\subset \mathbb{R}^2$ cannot be arbitrarily small. For sufficiently small $M,$ this fact is justified by the following theorem, which can be proved by estimating the Hilbert-Schmidt norm of a suitably defined operator. 
\begin{theorem}
(Concentration of $\mathcal{G}_{\psi}f$ on a small set) Let $\psi$ be a window function and $M\subset \mathbb{R}^2$ with $\lambda^{(2)}(M)<1.$ Then, for $f\in L^2(\mathbb{R})$ arbitrary
\begin{align}\label{R2000_Wilczok_GaborTransformConcentration_on_small_set}
\|\mathcal{G}_{\psi}f-\chi_M\cdot \mathcal{G}_{\psi}f\|_{L^2(\mathbb{R}^2)}\geq \|\psi\|_{L^2(\mathbb{R})}(1-\lambda^{(2)}(M)^{\frac{1}{2}})\|f\|_{L^2(\mathbb{R})}.
\end{align}
\end{theorem}
It is to be noted that the lower bound of  $\|\mathcal{G}_{\psi}f-\chi_M\cdot \mathcal{G}_{\psi}f\|_{L^2(\mathbb{R}^2)}$ in \eqref{R2000_Wilczok_GaborTransformConcentration_on_small_set} is larger the smaller $\lambda^{(2)}(M)$ is. However, it can be seen that the above theorem requires the assumption that $\lambda^{(2)}(M)<1.$ So a more general result is established with the assumption that $M$ is an arbitrary subset of $\mathbb{R}^2$ with finite Lebesgue measure, i.e.,  
for a given window function $\psi$ and $M\subset \mathbb{R}^2$ with $\lambda^{(2)}(M)<\infty,$ there exists a constant $C^{\mathcal{G}}_{\psi, M}>0$ such that for $f\in L^2(\mathbb{R})$ arbitrary
\begin{align}\label{R2000_Wilczok_GaborTransformConcentration_on_small_set}
\|\mathcal{G}_{\psi}f-\chi_M\cdot \mathcal{G}_{\psi}f\|_{L^2(\mathbb{R}^2)}\geq C^{\mathcal{G}}_{\psi, M}\|f\|_{L^2(\mathbb{R})}.
\end{align}
This result was established by Wilczok \cite{wilczok2000new} using the abstract result of Havin and J\"oricke in \cite{havin2012uncertainty}.
%(mention which result, at least cite the proper page)
However, the author in \cite{wilczok2000new} did not give an estimate for the constant $C^{\mathcal{G}}_{\psi, M}$ involved in \eqref{R2000_Wilczok_GaborTransformConcentration_on_small_set} by the measure of $M.$ Motivated by the idea of Singer \cite{singer1999uncertainty} in the case of WT, Wilczok obtained a different class of UPs results for the comparing the localization of $f$ and $\hat{f}$ with the localization of the WFT regarded as function of one variable as
\begin{align}\label{R2000_Wilczok_UP_of_GaborTransformin_omega}
\left(\int_{\mathbb{R}^2}\oomega^2|(\mathcal{G}_{\psi}f)(\oomega,\tee)|^2d\oomega d\xx\right)^{\frac{1}{2}}\left(\xx^2|f(\xx)|^2d\xx\right)^{\frac{1}{2}}\geq\frac{1}{2}\|\psi\|_{L^2(\mathbb{R})}\|f\|^2_{L^2(\mathbb{R})}
\end{align}
and
\begin{align}\label{R2000_Wilczok_UP_of_GaborTransformin_t}
\left(\int_{\mathbb{R}^2}\tee^2|(\mathcal{G}_{\psi}f)(\oomega,\tee)|^2d\oomega d\tee\right)^{\frac{1}{2}}\left(\xxi^2|\hat{f}(\xxi)|^2d\xxi\right)^{\frac{1}{2}}\geq\frac{1}{2}\|\psi\|_{L^2(\mathbb{R})}\|f\|^2_{L^2(\mathbb{R})}
\end{align}
which are called the UPs of the Heisenberg type for the WFT in $\oomega$ and $\tee$, respectively. Using \eqref{R2000_Wilczok_UP_of_GaborTransformin_omega} and \eqref{R2000_Wilczok_UP_of_GaborTransformin_t} one gets $\left(\int_{\mathbb{R}^2}\oomega^2|(\mathcal{G}_{\psi}f)(\oomega,\tee)|^2d\oomega d\xx\right)^{\frac{1}{2}}\left(\int_{\mathbb{R}^2}\tee^2|(\mathcal{G}_{\psi}f)(\oomega,\tee)|^2d\oomega d\tee\right)^{\frac{1}{2}}\left(\xx^2|f(\xx)|^2d\xx\right)^{\frac{1}{2}}\left(\xxi^2|\hat{f}(\xxi)|^2d\xxi\right)^{\frac{1}{2}}\geq\frac{1}{4}\|\psi\|^2_{L^2(\mathbb{R})}\|f\|^4_{L^2(\mathbb{R})},$ which is called the phase space uncertainty of WFT. It can be interpreted by saying that the better the phase space localization of the pair $(f,\hat{f}),$ the worse the phase space localization of the WFT $(\mathcal{G}_{\psi}f)(\oomega,\tee).$ The time-frequency resolution of $\mathcal{G}_{\psi}f$ is limited by the Heisenberg UP for $\psi,$ which is based on the interpretation of the standard deviations as the size of the essential supports of $\psi$ and its Fourier transform. Donoho and Stark \cite{grochenig2001foundations},\cite{donoho1989uncertainty} used the concept “$\epsilon$-concentrated” to replace standard deviation and derived the UP that reads as 
\begin{align*}
|\Omega||\Gamma|\geq 2\pi(1-\epsilon_\Omega-\epsilon_\Gamma)^2,
\end{align*}
whenever a non-zero function $f\in L^2(\mathbb{R})$ is $\epsilon_\Omega-$concentrated on $\Omega\subset \mathbb{R}$ and its Fourier transform is $\epsilon_\Gamma-$concentrated on $\Gamma\subset\mathbb{R}.$ An analogous to which is given as
$$|U|\geq (1-\epsilon),$$
for $f,\psi\in L^2(\mathbb{R}^{d}),$ with $\|f\|_{L^2(\mathbb{R}^d)}=\|\psi\|_{L^2(\mathbb{R}^d)}=1,$ and $U\subset \mathbb{R}^{2d}$ such that  
\begin{align}\label{RDonoho-Stark_UP_WeakVersion_STFT}
\int\int_{U}\left|(\mathcal{V}_{\psi}f)(\xx,\oomega)\right|^2d\xx d\oomega\geq 1-\epsilon,
\end{align}
$\epsilon\geq 0.$ With the similar assumptions on $f,\psi,U$ and $\mathcal{V}_{\psi}f$ satisfying the assumption \eqref{RDonoho-Stark_UP_WeakVersion_STFT}, Lieb (see \cite{grochenig2001foundations}) obtain a stronger lower bound of essential support as $|U|\geq (1-\epsilon)^{\frac{p}{p-2}}\left(\frac{p}{2}\right)^{\frac{2d}{p-2}},$ for all $p>2.$ In particular, $|U|\geq \sup_{p>0}(1-\epsilon)^{\frac{p}{p-2}}\left(\frac{p}{2}\right)^{\frac{2d}{p-2}}\geq (1-\epsilon)^22^d.$ This is often called the Lieb's UP for short time Fourier transform in literature. With the assumption, that $f,\psi\in L^2(\mathbb{R}^{2d}),$ $U\subset \mathbb{R}^{2d}$ such that $\int\int_{U}\left|(\mathcal{V}_{\psi}f)(\xx,\oomega)\right|^2d\xx d\oomega\geq (1-\epsilon)\|f\|_{L^2(\mathbb{R}^{2d})}\|\psi\|_{L^2(\mathbb{R}^{2d})}$ and using the Lieb's inequality for $p\geq 2,$ Poria in \cite{poria2021uncertainty} estimated the size of the essential support of the STFT as $|U|\geq (1-\epsilon)^{\frac{p}{p-1}}\left(\frac{p}{2}\right)^{\frac{d}{p-1}},~p\geq 2.$ Furthermore, he also proved that if $\int\int_{U}\left|\mathcal{V}_{\psi}f(\xx,\oomega)\right|^pd\xx d\xxi\geq (1-\epsilon)\|\mathcal{V}_{\psi}f\|_{L^1(\mathbb{R}^{2d})}^p,$ then $|U|\geq 2^{\frac{2pd}{2-p}}(1-\epsilon)^{\frac{2}{2-p}}$ for all $1\leq p<2.$ 

\section{Wavelet Transform (WT)}
In signal analysis, $(G_\phi f)(\B,\xxi)$ gives the time-frequency content of a signal $f$ at time $\B$ and frequency $\xxi$ by placing the window $\phi$ at time $\B.$ The drawback here is that a window of fixed width is used for all time $\B.$ It is more accurate and desirable if we have an adaptive window that gives a wide window for low frequency and a narrow one for high frequency. This can be done using the WT.
\begin{definition}
A function $\psi\in L^2(\mathbb{R})\setminus\{0\}$ satisfying the admissibility condition 
$$c_{\psi}=2\pi\int_{\mathbb{R}}|\hat{\psi}(\xxi)|^2\frac{d\xxi}{|\xxi|}<\infty$$
is called a mother wavelet.
\end{definition}
For a given mother wavelet $\psi$ and $(\A,\B)\in\mathbb{R}^{\ast}\times\mathbb{R},$ the family 
$\psi_{\A,\B}(\xx)=\tau_{\B}\delta_{\A}\psi(\xx)=\frac{1}{\sqrt|\A|}\psi\left(\frac{\xx-\B}{\A}\right)$ are called the daughter wavelets.
\begin{definition}
The WT of $f\in L^2(\mathbb{R})$ with respect to the mother wavelet $\psi$ is defined by 
\begin{align}\label{R2000_Wilczok_DefinitionWaveletTransform}
(\mathcal{W}_{\psi}f)(\A,\B)=\int_{\mathbb{R}}f(\xx)\frac{1}{\sqrt{|\A|}}\overline{\psi\left(\frac{\xx-\B}{\A}\right)}d\xx.
\end{align}
\end{definition}
\noindent Equivalently, using the definition of translation and the dilation operator in \ref{RTranslationOperatorFT} and \ref{RDilationOperatorFT} respectively,  \eqref{R2000_Wilczok_DefinitionWaveletTransform} has several important forms. Like
\begin{align}\label{R2000_Wilczok_DefinitionWaveletTransform_using_Operators}
(\mathcal{W}_{\psi}f)(\A,\B)=\int_{\mathbb{R}}f(\xx)\overline{\tau_{\B}\delta_{\A}\psi(\xx)}d\xx.
\end{align}
and
\begin{align}\label{R2000_Wilczok_DefinitionWaveletTransform_using_Convolution}
(\mathcal{W}_{\psi}f)(\A,\B)=\langle f,\psi_{\A,\B}\rangle_{L^2(\mathbb{R})}.
\end{align}
Another important, equivalent form of the definition of the WT \eqref{R2000_Wilczok_DefinitionWaveletTransform} is in terms of the convolution given by
\begin{align}\label{R2000_Wilczok_DefinitionWaveletTransform_using_InnerProduct}
(\mathcal{W}_{\psi}f)(\A,\B)=\left(f\star\delta_{\A}\tilde{\psi}\right)(\B),
\end{align}
where $\tilde{\psi}(t)=\overline{\psi(-t)}.$ It is  a matter of simple verification that $\mathcal{W}_{\psi}f$ is linear in $f$ and conjugate linear in $\psi.$ The following properties, called the covariance properties of the WT, give the formula for the WT of the translation by $\xx_0$ and dilation by $c\in \mathbb{R}^{\ast}$ of the function $f$ with respect to $\psi$
$$[(\mathcal{W}_{\psi}(\tau_{\xx_0}f)](\A,\B)=(\mathcal{W}_{\psi}f)(\A,\B-\xx_0)$$
and
$$[(\mathcal{W}_{\psi}(\delta_cf)](\A,\B)=(\mathcal{W}_{\psi}f)\left(\frac{\A}{c},\frac{\B}{c}\right).$$
The transform  \eqref{R2000_Wilczok_DefinitionWaveletTransform}, satisfies the following inner product relations
$$\int_{\mathbb{R}}\int_{\mathbb{R}}(\mathcal{W}_{\psi}f)(\A,\B)\overline{(\mathcal{W}_{\psi}g)(\A,\B)}\frac{d\A d\B}{\A^2}=c_{\psi}\langle f,g\rangle_{L^2(\mathbb{R})}$$
This in particular, for $f=g,$ implies that $\|\mathcal{W}_{\psi}f\|^2_{L^2(\mathbb{R}^{\ast}\times\mathbb{R},\frac{d\A d\B}{\A^2})}=c_{\psi}\|f\|^2_{L^2(\mathbb{R})},$ which further shows that the operator $\frac{1}{\sqrt{c_{\psi}}}\mathcal{W}_{\psi}$ is an isometry from $L^2(\mathbb{R})$ into $L^2(\mathbb{R}^{\ast}\times\mathbb{R},d\mu_{aff}),$ where $d\mu_{aff}=\frac{d\A d\B}{\A^2}$ denotes the so called affine measure. The WT converts a function $f(\xx)$ of one variable into a function $(\mathcal{G}_{\psi}f)(\A,\B)$ of two variables without changing the total energy.   Also, the $\mathcal{W}_{\psi}f$ is in $L^2(\mathbb{R}^{\ast}\times\mathbb{R},\frac{d\A d\B}{\A^2})$ for $f\in L^2(\mathbb{R}),$ i.e., the space of all WFT of square-integrable functions $\mathcal{W}_{\psi}(L^2(\mathbb{R}))$ is a subspace of $L^2(\mathbb{R}^{\ast}\times\mathbb{R},\frac{d\A d\B}{\A^2}).$ In fact, it is a proper subspace because $|(\mathcal{W}_{\psi}f)(\A,\B)|\leq \|f\|_{L^2(\mathbb{R}^2)}\|\psi\|_{L^2(\mathbb{R}^2)},~(\A,\B)\in\mathbb{R}^{\ast}\times\mathbb{R}$ implies that $\mathcal{W}_{\psi}f$ is bounded, but $L^2(\mathbb{R}^{\ast}\times\mathbb{R},\frac{d\A d\B}{\A^2})$ contains functions that are unbounded. Thus, the belongingness of $\mathcal{W}_{\psi}f\in L^2(\mathbb{R}^{\ast}\times\mathbb{R},\frac{d\A d\B}{\A^2})$ is necessary condition but not sufficient.  The range $\mathcal{W}_{\psi}(L^2(\mathbb{R})),$  of  bounded linear operator $\mathcal{G}_{\psi}f:L^2(\mathbb{R})\to L^2(\mathbb{R}^{\ast}\times\mathbb{R},\frac{d\A d\B}{\A^2})$ which is isometric upto some constant, is a closed subspace. Being closed subspace of the Hilbert space $L^2(\mathbb{R}^{\ast}\times\mathbb{R},\frac{d\A d\B}{\A^2}),$ $\mathcal{G}_{\psi}(L^2(\mathbb{R}))$ is also a Hilbert space. Moreover, it is a RKHS, with the kernel given by

$$K_{\psi}(\A',\B';\A,\B)=\frac{1}{c_{\psi}}\langle\psi_{\A,\B},\psi_{\A',\B'}\rangle_{L^2(\mathbb{R})}.$$
The kernel is pointwise bounded, i.e., $|K_{\psi}(\A',\B';\A,\B)|\leq\frac{\|\psi\|^2_{L^2(\mathbb{R})}}{c_{\psi}},~\forall~(\A',\B'),(\A,\B)\in\mathbb{R}^{\ast}\times\mathbb{R}.$

%%%%%%%%%%%%%%%%%%%%%%%%%%%%%%%%%%%%%%%%%%%%%%%%%%
%%%%%%%%%%%%%%%%%%%%%%%%%%%%%%%%%%%%%%%%%%%%%%%%%%
%%%%%%%%%%%%%%%%%%%%%%%%%%%%%%%%%%%%%%%%%%%%%%%%%%
%%%%%%%%%%%%%%%%%%%%%%%%%%%%%%%%%%%%%%%%%%%%%%%%%%
%%%%%%%%%%%%%%%%%%%%%%%%%%%%%%%%%%%%%%%%%%%%%%%%%%
\begin{definition}
A function $f\in L^2(\mathbb{R}^n)$ is a window function if it satisfies the conditions $\xx_{l_1}\xx_{l_2}\cdots\xx_{l_r}f(\bldx)\in L^2(\mathbb{R}^n),$ for all possible choices of $l_1,l_2,\cdot,l_r\in \{1,2,\cdots,n\},$ with $r\in \{1,2,\cdots,n\}$ where $l_1,l_2,\cdot,l_r$ are all different. In particular, if $
n=2,$ then we say that $f\in L^2(\mathbb{R}^2)$ is a window function if $\xx_1f(\bldx),\xx_2f(\bldx)$ and $\xx_1\xx_2f(\bldx)\in L^2(\mathbb{R}^2).$ For example $\psi(\bldx)=\xx_1\xx_2\cdots\xx_n e^{-(\xx_1^2+\xx_2^2+\cdots+\xx_n^2)}$ is a window function.
\end{definition}
A function $f\in L^2(\mathbb{R}^n)$ is a basic wavelet if it satisfies  the following admissibility condition
\begin{align}\label{R2015_Pandey_Admissibility_Consition_WT_PoAMS}
\int_{\mathbb{R}^n}\frac{|\hat{\psi}(\bldxi)|^2}{|\bldxi|_p}d\bldxi<\infty.
\end{align}
If $f\in L^2(\mathbb{R}^n)$ be a window function then $f\in L^1(\mathbb{R}^n),$ this can be proved using the notion of pseudo-orthants \cite{pandey2014pseudo}. Pandey et al. \cite{pandey2015continuous} proved that a window function $\psi\in L^2(\mathbb{R}^n)$ satisfies the admissibility condition \eqref{R2015_Pandey_Admissibility_Consition_WT_PoAMS} if and only if $\int_{\mathbb{R}}\psi(\xx_1,\xx_2,\cdots,\xx_n)d\xx_l,~l=1,2,\cdots,n.$  Some examples of window functions that are wavelets are $\Pi_{i=1}^n\xx_i e^{-\xx_i^2}$ and $\Pi_{i=1}^n\frac{\xx_i}{1+\xx_i^4}.$ Pandey et al. \cite{pandey2016continuous} defined the WT of $f\in L^2(\mathbb{R}^n)$ with respect to the basic wavelet $\psi\in L^2(\mathbb{R}^n)$ by 
\begin{align}\label{R2015_Pandey_nDWT_PoAMS}
(W_\psi f)(\bldb,\blda)=\left\langle f(\bldx),\frac{1}{\sqrt{|\blda|_p}}\psi\left(\frac{\bldx-\bldb}{\blda}\right)\right\rangle
\end{align}
and proved that 
\begin{align}\label{R2015_Pandey_Inversion_nDWT_IJWMIP}
\int_{\mathbb{R}^n}\int_{\mathbb{R}^n}(W_\psi f)(\bldb,\blda)\overline{(W_\psi g)(\bldb,\blda)}\frac{d\blda d\bldb}{|\blda|_p^2}=C_\psi\langle f,g\rangle~\mbox{for all}~f,g\in L^2(\mathbb{R}^n)
\end{align}
followed by the inversion formula is 
\begin{align}\label{R2015_Pandey_Inversion_nDWT_PoAMS}
f=\frac{1}{C_\psi}\int_{\mathbb{R}^n}\int_{\mathbb{R}^n}(W_\psi f)(\blda,\bldb)|\blda|^{-\frac{1}{2}}\psi\left(\frac{\bldx-\bldb}{\blda}\right)\frac{d\blda d\bldb}{|\blda|_p^2},
\end{align}
where $C_\psi=(2\pi)^n\int_{\mathbb{R}^n}\frac{|\hat{\psi}(\bldxi)|^2}{|\bldxi|_p}d\bldxi$ and $|\blda|_p=|\A_1\A_2\cdots\A_n|.$ The convergence in \eqref{R2015_Pandey_Inversion_nDWT_PoAMS} is interpreted in the $L^2(\mathbb{R}^n)$ sense. Moreover, the convergence of the inversion formula \eqref{R2015_Pandey_Inversion_nDWT_PoAMS} is pointwise at all the points of continuity of $f.$ It is to be noted that the parameter $\blda$ was taken to be $(\A_1,\A_2,\cdots,\A_n)$ where  the wavelet $\psi$ is not spherically symmetric in \cite{pandey2014pseudo},\cite{pandey2015continuous} which was not the case in Daubechies \cite{daubechies1992ten}, they took $\blda=(\A,\A,\cdots,\A),$ with $\A>0$ and the wavelet $\psi$ as spherically symmetric to obtain 
\begin{align}
\int_0^\infty\frac{d\A}{\A^{n+1}}\int_{-\infty}^{\infty}d\bldb(T^{wav}f)(\A,\bldb)\overline{(T^{wav}g)(\A,\bldb)}=C_\psi\langle f,g\rangle
\end{align}
and 
\begin{align}
f=C_{\psi}^{-1}\int_{-\infty}^{\infty}d\bldb(T^{wav}f)(\A,\bldb)\psi^{\A,\bldb},
\end{align}
where $\psi^{\A,\bldb}(\bldx)=\A^{-\frac{n}{2}}\psi\left(\frac{\bldx-\bldb}{\A}\right)$ and $(T^{wav}f)(\A,\bldb)=\langle f,\psi^{\A,\bldb}\rangle$ is the WT of $f\in L^2(\mathbb{R}^n)$ with respect to the spherically symmetric function $\psi\in L^2(\mathbb{R}^n)$ satisfying the condition $C_{\psi}=(2\pi)^n\int_0^\infty\frac{dt}{t}|\eta(t)|^2<\infty$ and $\hat{\psi}(\bldxi)=\eta(|\bldxi|)$ for some $\eta,$ which is possible to choose because $\psi$ is spherically symmetric. Recently, Pandey et al. \cite{pandey2019continuous} defined the continuous WT of tempered distribution $f\in\mathcal{S}'(\mathbb{R}^n)$ with the wavelet kernel $\psi\in\mathcal{S}(\mathbb{R}^n).$ It is to be noted that a function in $\psi\in\mathcal{S}(\mathbb{R}^n)$ is naturally a window function, as per the definition given   in \cite{pandey2015continuous}. So it satisfied the admissible condition \eqref{R2015_Pandey_Admissibility_Consition_WT_PoAMS} if and only if its integral with respect to each variable $\xx_1,\xx_2,\cdots,\xx_n$ along the real line vanishes. Let $s(\mathbb{R}^n)$ contains all the wavelets in $\mathcal{S}(\mathbb{R}^n).$ For the distribution $f\in \mathcal{S}'(\mathbb{R}^n)$ and the wavelet $\psi\in s(\mathbb{R}^n),$ the WT is defined in  \cite{pandey2019continuous}(also see \cite{pandey2019continuouss}) with the corresponding inversion formula, where the convergence is interpreted in the weak topology of $\mathcal{S}'(\mathbb{R}^n).$ They found that the WT of a constant distribution is zero and that the inversion formula is not true for a constant distribution but is true for a non-constant distribution which is not the sum of a non-constant distribution and a non-zero constant distribution. Pandey \cite{pandey2021continuous} considered the distribution space $\mathcal{D}_{L^2}'(\mathbb{R}^n),$ where $\mathcal{D}_{L^2}(\mathbb{R}^n)$ consists of $\mathcal{C}^\infty(\mathbb{R}^n)$  function whose every derivative is in $L^2(\mathbb{R}^n)$ and is equipped with a topology generated by a separating family of semi-norms $\{\gamma_k\}_{|k|=0}^\infty$ given by $\gamma_k(\phi)=\|\phi^{(k)}\|_{L^2(\mathbb{R}^n)}$ for  $\phi\in \mathcal{D}_{L^2}(\mathbb{R}^n).$ They defined the WT of generalized functions in $\mathcal{D}_{L^2}'(\mathbb{R}^n)$ which respect to the wavelet $\psi(\bldx)$ which are functions $\mathcal{D}_{L^2}(\mathbb{R}^n)\cap \mathcal{D}_{L^1}(\mathbb{R}^n)$ with the property $\int_{\mathbb{R}}\psi(\bldx)d\xx_l=0,~l=1,2,\cdots,n$ and $\int_{\mathbb{R}}\psi(\bldx)d\bldx\neq 0,$ where $\bldx^{\boldsymbol m}=\xx_1^{m_1}\xx_2^{m_2}\cdots\xx_n^{m_n},$ $\boldsymbol m=(m_1,m_2,\cdots, m_n), m_1,m_2,\cdots,m_n\geq 1$  also obtained the corresponding inversion formula with the convergence interpreted  in the weak distributional sense. The space $\mathcal{D}'_{L^2}(\mathbb{R}^n)$ does not contain a non-zero constant, so it does not encounter a problem as seen with the WT of distribution in the space $\mathcal{D}'(\mathbb{R}^n)$ with respect to the wavelet $\phi$ in $\mathcal{D}$ satisfying the property that all its moments along the real lines are non-zero \cite{10.1216/RMJ-2019-49-6-2005}, i.e., two functions having the same WT may differ by the constant.

Rieder \cite{rieder1990wavelet} extended the WT $L_\psi f(\B,\A):=\frac{1}{C_{\psi}}\langle \tau_\B\delta_\A\psi,f\rangle_0$ to the Sobolev spaces $H^\alpha(\mathbb{R}),\alpha\in\mathbb{R},$ consisting of the tempered distribution having a regular and with respect to the weight $(1+\oomega^2)^\alpha$ square integrable FT, and interpreted its image as subset of the fiber space $\mathscr{F}^\alpha:=L_2\left(\left(\mathbb{R}_0,\frac{d\A}{\A^2}\right),H^\alpha(\mathbb{R})\right),$   which consists of functions $\phi\in H^\alpha(\mathbb{R})$ such that $\int_{\mathbb{R}_0}\|\phi\|_\alpha^2\frac{d\A}{\A^2}<\infty.$ This space is isomorphic to the tensor product $L_2\left(\mathbb{R}_0,\frac{d\A}{\A^2}\right)\otimes H^\alpha(\mathbb{R})$ as well as the Sobolev space of two variables $H^{0,\alpha}\left(\mathbb{R}^2,\frac{d\A d\B}{\A^2}\right).$ Using the density of the Schwartz function in $H^\alpha(\mathbb{R}),$ author proved that the for a fixed $\A,$ and an integrable and admissible $\psi,$ $L_\psi f(\cdot,\A)$ is a continuous mapping from $H^\alpha(\mathbb{R})$ to itself. Moreover, the operator $L_\psi$ is an isometry from $H^\alpha(\mathbb{R})$ to the space $\mathscr{F}^\alpha,$ i.e., 
\begin{align}\label{RRieder_WT_Isometry-Property_Sobolev_Space}
\|L_\psi f\|_{\mathscr{F}^\alpha}=\left(\int_{\mathbb{R}}\|L_\psi f(\cdot,\A)\|_{H^\alpha(\mathbb{R})}^2\frac{d\A}{\A^2}\right)^{\frac{1}{2}}=\|f\|_{H^\alpha(\mathbb{R})}
\end{align}
and Range $(L_\psi)\subset \mathscr{F}^\alpha$ is a Hilbert space with the reproducing kernel $P(\tilde{\B},\tilde{\A},\B,\A)=\frac{1}{\sqrt{C_\psi}}(L_\psi \psi)\left(\frac{\tilde{\B}-\B}{\A},\frac{\tilde{\A}}{\A}\right).$ The $H^s-$distance of two WT with different admissible and integrable $\psi$ and $\phi$ of different argument function $f$ and $g$ in $H^s(\mathbb{R})$ respectively were obtained to study the dependence of the transform on its wavelets and its argument as
\begin{align}\label{RRieder_WT_Argument-Dependence-Property}
\|(L_\psi f)(a,\cdot)-(L_\phi g)(a,\cdot)\|_{H^s(\mathbb{R})}\leq |a|^{\frac{1}{2}}\left(\left\|\frac{\psi}{\sqrt{C_\psi}}-\frac{\phi}{\sqrt{C_\phi}}\right\|_{L^1(\mathbb{R})}\|f\|_{H^\alpha(\mathbb{R})}+\left\|\frac{\phi}{\sqrt{C_\phi}}\right\|_{L^1(\mathbb{R})}\|f-g\|_{H^\alpha(\mathbb{R})}\right).
\end{align}
Perrier et al. \cite{perrier1996besov} extended the WT to the $L^p(\mathbb{R})$ spaces to characterize the Besov spaces $B_p^{\alpha,q}(\mathbb{R})$ by the behavior at the small scale of the wavelet coefficient.

\begin{definition}
The Hardy space $\mathcal{H}^1(\mathbb{R}^n)$ is defined as the space of all functions $f\in L^1(\mathbb{R}^n)$ such that $\|f\|_{\mathcal{H}^1(\mathbb{R}^n)}=\int_{\mathbb{R}^n}\sup_{t>0}|(f\star \phi_t)(\bldx)|d\bldx<\infty,$ where $\phi$ is a function in Schwartz space such that $\int_{\mathbb{R}^n}\phi(\bldx)d\bldx\neq 0$ and $\phi_t(\bldx)=t^{-n}\phi(\frac{\bldx}{t}),~t>0,~\bldx\in\mathbb{R}^n.$ 
\end{definition}
Chuong et al. \cite{chuong2013boundedness} proved that for any $a\neq 0$ and any basis wavelet $\psi,$ with the assumption that $\psi$ is integrable, the operator $W_\psi:\mathcal{H}^1(\mathbb{R}^n)\rightarrow \mathcal{H}^1(\mathbb{R}^n)$ given by $f\rightarrow (W_\psi f)(a,\cdot)$ are bounded such that 
$\|(W_\psi f)(a,\cdot)\|_{\mathcal{H}^1(\mathbb{R}^n)}\leq |a|^{\frac{n}{2}}\|\psi\|_{L^1(\mathbb{R}^n)}\|f\|_{\mathcal{H}^1(\mathbb{R}^n)}$ and $\|(W_\psi f)(a,\cdot)\|_{\mathcal{H}^1(\mathbb{R}^n)}=O(|a|^{\frac{n}{2}}).$  Moreover, 
\begin{align}\label{R_2013_Chuong_WT_Argument-Dependence-Property_Hardy_Space}
\|(W_\psi f)(a,\cdot)-(W_\phi g)(a,\cdot)\|_{\mathcal{H}^1(\mathbb{R}^n)}\leq |a|^{\frac{n}{2}}\left(\|\psi-\phi\|_{L^1(\mathbb{R}^n)}\|f\|_{\mathcal{H}^1(\mathbb{R}^n)}+\|\phi\|_{L^1(\mathbb{R}^n)}\|f-g\|_{\mathcal{H}^1(\mathbb{R}^n)}\right),~f,g\in \mathcal{H}^1(\mathbb{R}^n).
\end{align}
Similar result was also obtained with the same estimates of the bound in the case of the Besov space $B_p^{\alpha,q}(\mathbb{R}^n)$, bounded mean oscillation space $BMO(\mathbb{R}^n)$ and the $VMO(\mathbb{R}^n)$ space with the norm replaced with the corresponding norm and the basic wavelet compactly supported in the case of $BMO(\mathbb{R}^n)$ and $VMO(\mathbb{R}^n)$ spaces. We discuss these results for the weighted Besov spaces and the weighted BMO spaces associated with a tempered weight function. Let us recall the definitions of the weighted Besov spaces and the weighted BMO spaces associated with a tempered weight function, which was proposed by H\"ormander to study the theory of partial differential equations, i.e., a positive function $k$ defined in $\mathbb{R}^n$ for which there exist $M,N>0$ such that $k(\bldx+\bldy)\leq(1+M|\bldx|)^Nk(\bldy),~\forall~\bldx,\bldy\in\mathbb{R}^n.$ For example, $k(\bldx)=(1+|\bldx|^2)^{\frac{s}{2}},$ for every $s\in\mathbb{R}.$ We start with the following definitions.
\begin{definition}
The weighted Lebesgue space $L^p_k(\mathbb{R}^n),~1\leq p<\infty$ is the space of all measurable functions on $\mathbb{R}^n$ such that $\|f\|_{p,k}=\left(\int_{\mathbb{R}^n}|f(\bldx)|^pk(\bldx)d\bldx\right)^{\frac{1}{p}}<\infty.$
\end{definition}
\begin{definition}
Let $0<\alpha<1$ and $1\leq p,q<\infty,$ then the weighted Besov space $B_{p,k}^{\alpha,q}(\mathbb{R}^n)$ is defined as 
$$B_{p,k}^{\alpha,q}(\mathbb{R}^n)=\left\{f\in L^p_k(\mathbb{R}^n):\int_{\mathbb{R}^n}|w_{p,k}(f,\bldh)|^q\frac{d\bldh}{|\bldh|^{n+\alpha q}}<\infty\right\},~q<\infty$$ 
and 
$$B_{p,k}^{\alpha,\infty}(\mathbb{R}^n)=\left\{f\in L^p_k(\mathbb{R}^n):|\bldh|^{-\alpha} w_{p,k}(f,\bldh)\in L^\infty(\mathbb{R}^n\setminus\{0\})\right\}$$
with the respective norm as $\|f\|_{B_{p,k}^{\alpha,q}}=\|f\|_{p,k}+\left(\int_{\mathbb{R}^n}|w_{p,k}(f,\bldh)|^q\frac{d\bldh}{|\bldh|^{n+\alpha q}}\right)^{\frac{1}{q}}$ and $\|f\|_{B_{p,k}^{\alpha,\infty}}=\|f\|_{p,k}+\||\bldh|^{-\alpha} w_{p,k}(f,\bldh)\|_\infty,$ where $w_{p,k}(f,\bldh)=\|f(\cdot+\bldh)-f(\cdot)\|_{p,k}.$
\end{definition}
\begin{definition}
The weighted bounded mean oscillation space $BMO_k(\mathbb{R}^n)$ is defined as the space of all functions $f\in L_{loc}^{1,k}(\mathbb{R}^n)$ such that $\|f\|_{BMO_k}=\sup_{B\subset\mathbb{R}^n}=\frac{1}{|B|_k}\int_{B}|f(\bldx)-f_B|k(\bldx)d\bldx<\infty,$ where supremum is taken over all balls in $\mathbb{R}^n$ and $|B|_k=\int_{B}k(\bldx)d\bldx.$
\end{definition}
Authors in \cite{chuong2013boundedness} also proved that for any $\A\neq 0,$ and any compactly supported basis wavelet $\psi,$ with the support belonging to a ball of radius $r$ with center at the origin, the operator $W_\psi:B_{p,k}^{\alpha,q}(\mathbb{R}^n)\rightarrow B_{p,k}^{\alpha,q}(\mathbb{R}^n)$ and $W_\psi:BMO_k(\mathbb{R}^n)\rightarrow BMO_k(\mathbb{R}^n)$ given by $f\rightarrow (W_\psi f)(\A,\cdot)$ are bounded such that $\|(W_\psi f)(\A,\cdot)\|_{B_{p,k}^{\alpha,q}}\leq (1+Mr|\A|)^{\frac{N}{p}}|\A|^{\frac{n}{2}}\|\psi\|_{L^1(\mathbb{R}^n)}\|f\|_{B_{p,k}^{\alpha,q}}$ and $\|(W_\psi f)(\A,\cdot)\|_{BMO_k}\leq (1+Mr|\A|)^{2N}|\A|^{\frac{n}{2}}\|\psi\|_{L^1(\mathbb{R}^n)}\|f\|_{BMO_k}.$ Moreover, 
\begin{align*}
%\label{R_2013_Chuong_WT_Argument-Dependence-Property_Weighted_Besov_Space}
\|(W_\psi f)(\A,\cdot)-(W_\phi g)(\A,\cdot)\|_{B_{p,k}^{\alpha,q}}\leq (1+Mr|\A|)^{\frac{N}{p}}|\A|^{\frac{n}{2}}\left(\|\psi-\phi\|_{L^1(\mathbb{R}^n)}\|f\|_{B_{p,k}^{\alpha,q}}+\|\phi\|_{L^1(\mathbb{R}^n)}\|f-g\|_{B_{p,k}^{\alpha,q}}\right),~f,g\in B_{p,k}^{\alpha,q}(\mathbb{R}^n)
\end{align*}
and for $f,g\in BMO_k(\mathbb{R}^n)$
\begin{align}\label{R_2013_Chuong_WT_Argument-Dependence-Property_BMO_Space}
\|(W_\psi f)(\A,\cdot)-(W_\phi g)(\A,\cdot)\|_{BMO_k}\leq (1+Mr|\A|)^{2N}|\A|^{\frac{n}{2}}\left(\|\psi-\phi\|_{L^1(\mathbb{R}^n)}\|f\|_{BMO_k}+\|\phi\|_{L^1(\mathbb{R}^n)}\|f-g\|_{BMO_k}\right).
\end{align}
Chuong et. al \cite{chuong2002integral} studied the WT on the weighted Sobolev space $W_\omega^{m,p}(\mathbb{R}^n)$ which is defined as
 \begin{align}\label{R2001_Chuong_Tri_Weighted Sobolev Space_Definition}
 W_\omega^{m,p}(\mathbb{R}^n)=\left\{f\in L^p(\mathbb{R}^n)|\partial^kf\in L^p(\mathbb{R}^n),|k|\leq m\right\},
 \end{align}
which is equipped with the norm $\|f\|_{m,p,\omega}=\sum_{|\mu|\leq m}\left(\int_{\mathbb{R}^n}\omega_{\mu}(\bldx)|\partial^\mu f(\bldx)|^pd\bldx\right)^{\frac{1}{p}}<\infty,$ where $\mu=(\mu_1,\mu_2,\cdots,\mu_n),$ $|\mu|=\mu_1+\mu_2+\cdots+\mu_n,\mu_l\geq 0$ and $\omega_\mu(\bldx)\in L^\infty(\mathbb{R}^n)$ with $\omega_\mu (\bldx)>0,$ for almost all $\bldx\in \mathbb{R}^n$ and for each $\bldx\in\mathbb{R}^n,$ $\omega_\mu(\bldx+\bldy)\leq C_{1,\mu}\omega_\mu(\bldx)$ for almost all $\bldy\in\mathbb{R}^n.$ Authors in \cite{chuong2002integral} proved that the WT operator $(L_\psi f)(\bldb,\A)\left(:=\frac{1}{C_\psi^{\frac{1}{2n}}}\frac{1}{|\A|^{\frac{n}{2}}}\int_{\mathbb{R}^n}\bar{\psi}\left(\frac{\bldx-\bldb}{\A}\right)f(\bldx)d\bldx\right),$ for fixed $\A>0$ maps $W_\omega^{m,p}(\mathbb{R}^n)$ to itself and $\|(L_\psi f)(\cdot,\A))\|_{m,p,\omega}\leq C\|f\|_{m,p,\omega},$ where $C=\max_{|\mu|\leq m}C_\mu,$ $C_\mu=\frac{(C_{1,\mu})^{\frac{1}{p}}}{C_{\psi}^{\frac{1}{2n}}}|a|^{\frac{n}{2}}\|\psi\|_{L^1(\mathbb{R}^n)}$ and $C_\psi=(2\pi)^n\int_0^\infty\frac{|\hat{\psi}(a\bldxi)|^2}{\A}d\A$ is a non zero constant for every $\bldxi\neq 0.$ They also estimated the $\|\cdot\|_{m,p,\omega}-$distance of two WT with different basic wavelets $\psi$ and $\phi$ of different argument function $f$ and $g$ in $W_\omega^{m,p}(\mathbb{R}^n)$ respectively and studied the asymptotic behavior of the WT for small dilation parameter, analogous to the one studied in \cite{rieder1990wavelet}.\\

\noindent\textbf{Uncertainty Principles:} Wilczok \cite{wilczok2000new} proved that $\mathcal{W}_{\psi}f$ of a function $f\in L^2(\mathbb{R})\setminus\{0\}$ cannot have a support of finite affine measure. The author also proved that the portion of $\mathcal{W}_{\psi}f$ lying outside some set $M\subset \mathbb{R}^2$ cannot be arbitrarily small. For sufficiently small $M,$ this fact is justified by the following theorem, which can be proved by estimating the Hilbert-Schmidt norm of a suitably defined operator. 
\begin{theorem}
(Concentration of $\mathcal{W}_{\psi}f$ on a small set) Let $\psi$ be a window function and $M\subset \mathbb{R}^{\ast}\times\mathbb{R}$ with the assumption $\frac{\|\psi\|_{L^2(\mathbb{R})}}{\sqrt{c_{\psi}}}\mu_{aff}(M)^{\frac{1}{2}}<1.$ Then, for $f\in L^2(\mathbb{R})$ arbitrary
\begin{align}\label{R2000_Wilczok_WaveletTransformConcentration_on_small_set}
\|\mathcal{W}_{\psi}f-\chi_M\cdot \mathcal{W}_{\psi}f\|_{L^2\left(\mathbb{R}^{\ast}\times\mathbb{R},\frac{d\A d\B}{\A^2}\right)}\geq \sqrt{c_{\psi}}\left(1-\frac{\|\psi\|_{L^2(\mathbb{R})}}{\sqrt{c_{\psi}}}\mu_{aff}(M)^{\frac{1}{2}}\right)\|f\|_{L^2(\mathbb{R})}.
\end{align}
\end{theorem}
\noindent It is to be noted that the lower bound of  $\|\mathcal{W}_{\psi}f-\chi_M\cdot \mathcal{W}_{\psi}f\|_{L^2(\mathbb{R}^2)}$ in \eqref{R2000_Wilczok_WaveletTransformConcentration_on_small_set} is larger the smaller $\mu_{aff}(M)$ is. However, it can be seen that the above theorem requires the assumption that $\frac{\|\psi\|_{L^2(\mathbb{R})}}{\sqrt{c_{\psi}}}\mu_{aff}(M)^{\frac{1}{2}}<1.$ So a more general result is established with the assumption that $M$ is an arbitrary subset of $\mathbb{R}^{\ast}\times\mathbb{R}$ with finite affine measure $\mu_{aff}$, i.e.,  
for a given $\psi$ and $M\subset \mathbb{R}^{\ast}\times\mathbb{R}$ with $\mu_{aff}(M)<\infty,$ there exists a constant $C^{\mathcal{W}}_{\psi, M}>0$ such that for $f\in L^2(\mathbb{R})$ arbitrary
\begin{align}\label{R2000_Wilczok_WaveletTransformConcentration_on_small_set}
\|\mathcal{W}_{\psi}f-\chi_M\cdot \mathcal{W}_{\psi}f\|_{L^2(\mathbb{R}^2)}\geq C^{\mathcal{W}}_{\psi, M}\|f\|_{L^2(\mathbb{R})}.
\end{align}
This result was established by Wilczok \cite{wilczok2000new} using the abstract result of Havin and J\"oricke \cite{havin2012uncertainty}. However, the author in \cite{wilczok2000new} did not give an estimate for the constant $C^{\mathcal{W}}_{\psi, M}$ involved in \eqref{R2000_Wilczok_WaveletTransformConcentration_on_small_set} by the measure of $M.$ Motivated by the idea of Singer \cite{singer1999uncertainty} in the case of WT, Wilczok obtained a different class of UPs results for the comparing the localization of $f$ and $\hat{f}$ with the localization of the WT regarded as function of one variable as
\begin{align}\label{R2000_Wilczok_UP_of_WaveletTransformin_b}
\left(\int_{\mathbb{R}^2}\B^2|(\mathcal{W}_{\psi}f)(\A,\B)|^2\frac{d\A d\B}{|\A|^2}\right)^{\frac{1}{2}}\left(\xxi^2|\hat{f}(\xxi)|^2d\xxi\right)^{\frac{1}{2}}\geq\frac{1}{2}\sqrt{c_{\psi}}\|f\|^2_{L^2(\mathbb{R})}
\end{align}
and
\begin{align}
\int_0^\infty\int_{-\infty}^\infty\A^2|(\tilde{\mathcal{W}}_{\psi}f)(\A,\B)|^2d\A d\B\cdot\int_{-\infty}^\infty\xx^2|f(\xx)|^2d\xx\geq\pi\left(\mathcal{M}\left(|\hat{\psi}|^2\right)\right)(2)\|f\|^2_{L^2(\mathbb{R})},
\end{align}
where $\mathcal{M}$ denotes the classical Mellin transform given by $(\mathcal{M}f)(\sigma)=\int_0^\infty f(\xx)\xx^{-\sigma}\frac{d\xx}{\xx}.$

\section{Stockwell Transform (ST)}
The S transform, also known as the ST in the literature, was first introduced by Stockwell et al. \cite{stockwell1996localization} in 1996. This transform is basically an extension of the ideas of the continuous WT, which is based on a moving and scalable localization Gaussian window and is defined as
\begin{align}\label{R1996_Stockwell-Mansinha-Lowe_Definition_STransform}
(\mathcal{S}f)(\tau,\xxi)=\int_{-\infty}^\infty f(\tee)\frac{|\xxi|}{\sqrt{2\pi}}e^{-\frac{(\tau-\tee)^2 \xxi^2}{2}}e^{-i2\pi \xxi\tee}d\tee.
\end{align} 
%%%%%%%%%%%%%%%%%%%%%%%%%%%%%%%%%%%%%%%%%%%%%%%%%%%%%%%%%%%%%%%%%%%%%%%%%%%%%%%%%%%%%%%%%%%
%%%%%%%%%%%%%%%%%%%%%%%%%%%%%%%%%%%%%%%%%%%%%%%%%%%%%%%%%%%%%%%%%%%%%%%%%%%%%%%%%%%%%%%%%%%
%%%%%%%%%%%%%%%%%%%%%%%%%%%%%%%%%%%%%%%%%%%%%%%%%%%%%%%%%%%%%%%%%%%%%%%%%%%%%%%%%%%%%%%%%%%
%%%%%%check:%%%%%%%%%%%%%%%%%%%%%%%%%%%%%%%%%%%%%%%%%%%%%%%%%%%%%%%%%%%%%%%%%%%%%%%%%%%%%%%
%\textcolor{purple}{(Fourier kernel used here is different that that we are going to use below. So , I %think we must change it, for uniformity)}
%%%%%%%%%%%%%%%%%%%%%%%%%%%%%%%%%%%%%%%%%%%%%%%%%%%%%%%%%%%%%%%%%%%%%%%%%%%%%%%%%%%%%%%%%%%
%%%%%%%%%%%%%%%%%%%%%%%%%%%%%%%%%%%%%%%%%%%%%%%%%%%%%%%%%%%%%%%%%%%%%%%%%%%%%%%%%%%%%%%%%%%
%%%%%%%%%%%%%%%%%%%%%%%%%%%%%%%%%%%%%%%%%%%%%%%%%%%%%%%%%%%%%%%%%%%%%%%%%%%%%%%%%%%%%%%%%%%
%%%%check:%%%%%%%%%%%%%%%%%%%%%%%%%%%%%%%%%%%%%%%%%%%%%%%%%%%%%%%%%%%%%%%%%%%%%%%%%%%%%%%%%
One can infer from \eqref{R1996_Stockwell-Mansinha-Lowe_Definition_STransform} that the width of the window function $e^{-\frac{(\tau-\tee)^2 \xxi^2}{2}}$ shrinks as the analyzed frequency increases and provides a better time localization for high frequencies. This transform was introduced to provide frequency-dependent resolution while maintaining a direct relationship with the Fourier spectrum, which was absent in the continuous WT. The advantage of these transforms is due to the fact that the modulation sinusoidal is fixed with respect to the time axis, whereas the localizing scalable Gaussian window dilates and translates. If $f\in L^2(\mathbb{R}),$ then the ST given in equation \eqref{R1996_Stockwell-Mansinha-Lowe_Definition_STransform} can be given in the Fourier domain as 
\begin{align}\label{R1996_Stockwell-Mansinha-Lowe_Definition_STransform_Fourier-Domain}
(\mathcal{S}f)(\tau,\xxi)=e^{-2\pi i\xxi\tau}\int_{-\infty}^\infty \hat{f}(v)e^{-\frac{2\pi^2(v-\xxi)^2 }{\xxi^2}}e^{i2\pi \xxi\tau}dv,~\tau\in \mathbb{R},\xxi\in\mathbb{R}_{0}.
\end{align}
%%%%%%%%%%%%%%%%%%%%%%%%%%%%%%%%%%%%%%%%%%%%%%%%%%%%%%%%%%%%%%%%%%%%%%%%%%%%%%%%%%%%%%%%%%%
%%%%%%%%%%%%%%%%%%%%%%%%%%%%%%%%%%%%%%%%%%%%%%%%%%%%%%%%%%%%%%%%%%%%%%%%%%%%%%%%%%%%%%%%%%%
%%%%%%%%%%%%%%%%%%%%%%%%%%%%%%%%%%%%%%%%%%%%%%%%%%%%%%%%%%%%%%%%%%%%%%%%%%%%%%%%%%%%%%%%%%%
%%%%%%%%%%%%%%%%%%%%%%%%%%%%%%%%%%%%%%%%%%%%%%%%%%%%%%%%%%%%%%%%%%%%%%%%%%%%%%%%%%%%%%%%%%%
%%%\check:%%%%%%%%%%%%%%%%%%%%%%%%%%%%%%%%%%%%%%%%%%%%%%%%%%%%%%%%%%%%%%%%%%%%%%%%%%%%%%%%%
%\textcolor{purple}{(Fourier kernel used here is different that that we are going to use below. So , I think we must change it, for uniformity)}
%%%%%%%%%%%%%%%%%%%%%%%%%%%%%%%%%%%%%%%%%%%%%%%%%%%%%%%%%%%%%%%%%%%%%%%%%%%%%%%%%%%%%%%%%%%
%%%%%%%%%%%%%%%%%%%%%%%%%%%%%%%%%%%%%%%%%%%%%%%%%%%%%%%%%%%%%%%%%%%%%%%%%%%%%%%%%%%%%%%%%%%
%%%%%%%%%%%%%%%%%%%%%%%%%%%%%%%%%%%%%%%%%%%%%%%%%%%%%%%%%%%%%%%%%%%%%%%%%%%%%%%%%%%%%%%%%%%
%%%%%%%%%%%%%%%%%%%%%%%%%%%%%%%%%%%%%%%%%%%%%%%%%%%%%%%%%%%%%%%%%%%%%%%%%%%%%%%%%%%%%%%%%%%
Thus, the  ST can be obtained as the Fourier multiplier $\mathcal{F}^{-1}e^{-\frac{2\pi^2(\cdot -\xxi)^2}{\xxi^2}}\mathcal{F}$ followed by a modulation for non zero frequency $\xxi,$ here $(\mathcal{F}f)(\xxi)=\int_{\mathbb{R}}f(\tee)e^{-2\pi i\tee\xxi}d\tee.$ Moreover, the ST satisfy the property $A\mathcal{S}f=\hat{f},~f\in L^2(\mathbb{R}),$ where $A$ is called the time average operator defined by 
\begin{align}\label{R1996_Stockwell-Mansinha-Lowe_Definition_STransform_time_average_Operator}
(AF)(\xxi)=\int_{\mathbb{R}}F(\tau,\xxi)d\tau,~\xxi\in \mathbb{R}
\end{align}
for the measurable function $F:\mathbb{R}\times\mathbb{R}\rightarrow\mathbb{C},$ whenever the integral exists. As a result of which we have the inversion formula $f=\mathcal{F}^{-1}A\mathcal{S}f,$ for all $f\in L^2(\mathbb{R}),$ i.e.,
\begin{align}\label{R1996_Stockwell-Mansinha-Lowe_Definition_STransform_Inversion}
f(t)=\int_{\mathbb{R}}e^{2\pi it\xxi}\left(\int_{\mathbb{R}}(\mathcal{S}f)(\tau,\xxi)\right)d\xxi,~t\in\mathbb{R}
\end{align} 
%%%%%%%%%%%%%%%%%%%%%%%%%%%%%%%%%%%%%%%%%%%%%%%%%%%%%%%%%%%%%%%%%%%%%%%%%%%%%%%%%%%%%%%%%%%
%%%%%%%%%%%%%%%%%%%%%%%%%%%%%%%%%%%%%%%%%%%%%%%%%%%%%%%%%%%%%%%%%%%%%%%%%%%%%%%%%%%%%%%%%%%
%%%%%%%%%%%%%%%%%%%%%%%%%%%%%%%%%%%%%%%%%%%%%%%%%%%%%%%%%%%%%%%%%%%%%%%%%%%%%%%%%%%%%%%%%%%
%%%check:%%%%%%%%%%%%%%%%%%%%%%%%%%%%%%%%%%%%%%%%%%%%%%%%%%%%%%%%%%%%%%%%%%%%%%%%%%%%%%%%%%
%\textcolor{purple}{(Fourier kernel used here is different that that we are going to use below. So , I think we must change it, for uniformity)}
%%%%%%%%%%%%%%%%%%%%%%%%%%%%%%%%%%%%%%%%%%%%%%%%%%%%%%%%%%%%%%%%%%%%%%%%%%%%%%%%%%%%%%%%%%%
%%%%%%%%%%%%%%%%%%%%%%%%%%%%%%%%%%%%%%%%%%%%%%%%%%%%%%%%%%%%%%%%%%%%%%%%%%%%%%%%%%%%%%%%%%%
%%%%%%%%%%%%%%%%%%%%%%%%%%%%%%%%%%%%%%%%%%%%%%%%%%%%%%%%%%%%%%%%%%%%%%%%%%%%%%%%%%%%%%%%%%%
%%%%%%%%%%%%%%%%%%%%%%%%%%%%%%%%%%%%%%%%%%%%%%%%%%%%%%%%%%%%%%%%%%%%%%%%%%%%%%%%%%%%%%%%%%%
and also the Plancherel theorem for the ST as $\int_{\mathbb{R}}\left|\int_{\mathbb{R}}(\mathcal{S}f)(\tau,\xxi)d\tau\right|^2d\xxi=\|f\|^2_{L^2(\mathbb{R})}.$  Wong et al.\cite{wong2007characterization} considered the function space $M=\left\{F:\mathbb{R}\times\mathbb{R\rightarrow\mathbb{C}}:\int_{\mathbb{R}}\left|\int_{\mathbb{R}}F(\tau,\xxi)d\tau\right|^2d\xxi\right\},$ which is indefinite Hilbert space with the indefinite inner product $\langle,\rangle_M$ given by $\langle F,G\rangle_M=\langle AF,AG\rangle_{L^2(\mathbb{R})},~F,G\in M,$ where $AF$ is the average-time of $F$ given by \eqref{R1996_Stockwell-Mansinha-Lowe_Definition_STransform_time_average_Operator}. It is to be noted that $\langle ,\rangle_M$ is indefinite in the sense that there exist nonzero $F\in M$ such that $\langle F,F\rangle_M=0.$ They proved that the range $R(\mathcal{S})$ of the ST $\mathcal{S}:L^2(\mathbb{R}^2)\rightarrow M$ is a Hilbert space with respect to the inner product $\langle ,\rangle_M.$ In fact, $\mathcal{S}$ is an isometry from $L^2(\mathbb{R})$ onto $R(\mathcal{S}).$ As a consequence of $A\mathcal{S}f=\hat{f}$ and the Plancherel theorem for the ST, its range $R(\mathcal{S})$ is characterized in \cite{wong2007characterization} as $R(\mathcal{S})=M/Z,$ where $Z=\{F:\mathbb{R}\times\mathbb{R}\rightarrow C:AF=0\}.$ 
\\
%%%%%%%%%%%%%%%%%%%%%%%%%%%%%%%%%%%%%%%%%%%%%%%%%%%%%%%%%%%%%%%%%%%%%%%%%%%%%%%%%%%%%%%%%%%
%%%%%%%%%%%%%%%%%%%%%%%%%%%%%%%%%%%%%%%%%%%%%%%%%%%%%%%%%%%%%%%%%%%%%%%%%%%%%%%%%%%%%%%%%%%
%%%check:%%%%%%%%%%%%%%%%%%%%%%%%%%%%%%%%%%%%%%%%%%%%%%%%%%%%%%%%%%%%%%%%%%%%%%%%%%%%%%%%%%
%\textcolor{purple}{(What does this actually mean, need to check.)}
%%%%%%%%%%%%%%%%%%%%%%%%%%%%%%%%%%%%%%%%%%%%%%%%%%%%%%%%%%%%%%%%%%%%%%%%%%%%%%%%%%%%%%%%%%%
%%%%%%%%%%%%%%%%%%%%%%%%%%%%%%%%%%%%%%%%%%%%%%%%%%%%%%%%%%%%%%%%%%%%%%%%%%%%%%%%%%%%%%%%%%%
%%%%%%%%%%%%%%%%%%%%%%%%%%%%%%%%%%%%%%%%%%%%%%%%%%%%%%%%%%%%%%%%%%%%%%%%%%%%%%%%%%%%%%%%%%%
The ST with the general window in the place of the Gaussian function in the definition given in \cite{stockwell1996localization} can be defined below which is a hybrid of the Gabor transform and the WT
\begin{definition}
The ST of a function $f\in L^2(\mathbb{R})$ with respect to the window function $\varphi\in L^1(\mathbb{R})\cap L^2(\mathbb{R})$ is defined by 
\begin{align}\label{R2009_Molahajloo-Wong_DefnitionST_GeneralWindow}
(\mathcal{S}_\varphi f)(\B,\xxi)=\frac{|\xxi|}{\sqrt{2\pi}}\int_{\mathbb{R}}e^{-i\xx\xxi}f(\xx)\overline{\varphi(\xxi(\xx-\B))}d\xx,~\xx\in\mathbb{R},~\xxi\in\mathbb{R}_{0}.
\end{align}
\end{definition}
\noindent Equivalently, using the definition of translation, dilation  and the modulation operators in \ref{RTranslationOperatorFT}, \ref{RDilationOperatorFT} and \ref{RModulationOperatorFT} respectively,  \eqref{R2009_Molahajloo-Wong_DefnitionST_GeneralWindow} has several important forms.
\begin{align}\label{R2009_Molahajloo-Wong_DefnitionST_GeneralWindow_Operators}
(\mathcal{S}_{\varphi}f)(\B,\xxi)=\int_{\mathbb{R}}f(\xx)\overline{\mathcal{M}_{\xxi}\tau_{\B}\mathcal{D}_{\xxi}\varphi(\xx)}d\xx
\end{align}
which can also be written in the form
\begin{align}\label{R2009_Molahajloo-Wong_DefnitionST_GeneralWindow_InnerProduct}
(\mathcal{S}_{\varphi}f)(\B,\xxi)=\frac{1}{\sqrt{2\pi}}\langle f,\varphi^{\B,\xxi}\rangle_{L^2(\mathbb{R})},
\end{align}
where $\varphi^{\B,\xxi}(\xx)=\frac{1}{\sqrt{2\pi}}\mathcal{M}_{\xxi}\tau_{\B}\mathcal{D}_{\xxi}\varphi(\xx).$ Another important, equivalent form of the definition of the WT \eqref{R2000_Wilczok_DefinitionWaveletTransform} is in terms of the convolution given by
\begin{align}\label{R2009_Molahajloo-Wong_DefnitionST_GeneralWindow_Convolution}
(\mathcal{S}_{\varphi}f)(\B,\xxi)=\frac{1}{\sqrt{2\pi}}\left(\mathcal{M}_{-\xxi}f\star\mathcal{D}_{\xi}\tilde{\varphi}\right)(\B),
\end{align}
where $\tilde{\psi}(t)=\overline{\psi(-t)}.$ It is  a matter of simple verification that $\mathcal{S}_{\varphi}f$ is linear in $f$ and conjugate linear in $\varphi.$ The ST $\mathcal{S}_\varphi$ has been extended to a family of modified ST $\mathcal{S}^s_\varphi$ for $0<s\leq\infty$ in \cite{guo2009modified}(also see \cite{guo2008modified}), which are defined as
\begin{align}
(\mathcal{S}^s_{\varphi}f)(\B,\xxi)&=\frac{1}{\sqrt{2\pi}}|\xxi|^{\frac{1}{s}}\int_{\mathbb{R}}e^{-i\xx\xxi}f(\xx)\overline{\varphi(\xxi(\xx-\B))}d\xx\\
&=\langle f,\varphi_s^{\B,\xxi}\rangle_{L^2(\mathbb{R})}
\end{align}
where $\varphi_s^{\B,\xxi}(\xx)=\frac{1}{\sqrt{2\pi}}\mathcal{M}_{\xxi}\tau_{\B}\mathcal{D}^{s}_{\xxi}\varphi(\xx)$ and ${D}^{s}_{\xxi}\varphi(\xx)=|\xxi|^{\frac{1}{s}}\varphi(\xxi\xx).$ This family include the classical ST \eqref{R1996_Stockwell-Mansinha-Lowe_Definition_STransform} when $s=1$ and a variant of WT when $s=2.$ The modified ST can be represented in the classical one as $(\mathcal{S}^s_{\varphi}f)(\B,\xxi)=|\xxi|^{-\frac{1}{s'}}(\mathcal{S}_{\varphi}f)(\B,\xxi),$ where $\frac{1}{s}+\frac{1}{s'}=1.$ The modified ST satisfies the resolution of identity formula \cite{guo2009modified}
\begin{align}\label{R2009_Molahajloo-Wong_ST_Resolution-of-Identity-formula}
\langle f,g \rangle_{L^2(\mathbb{R})}=\frac{1}{c_\varphi}\int_{\mathbb{R}}(\mathcal{S}^s_{\varphi}f)(\B,\xxi)\overline{(\mathcal{S}^s_{\varphi}g)(\B,\xxi)}\frac{d\B d\xxi}{|\xxi|^{1-\frac{1}{s'}}},
\end{align}
where $\varphi\in L^1(\mathbb{R})\cap L^2(\mathbb{R})$ is such that $\int_{R}\varphi(\xx)d\xx=1$ and 
\begin{align}\label{R2009_Molahajloo-Wong_ST_AdmissibilityCondition}
c_\varphi:=\int_{\mathbb{R}}\frac{|\hat{\varphi}(\xxi-1)|^2}{|\xxi|}d\xxi <\infty,
\end{align}
with $\hat{\varphi}(\xxi)=\frac{1}{\sqrt{2\pi}}\int_{\mathbb{R}}e^{-i\xx\xxi}\varphi(\xx)d\xx,~\xxi\in\mathbb{R}.$ The modified ST $\mathcal{S}^s_\varphi,~0<s\leq \infty,$ for $f\in L^2(\mathbb{R})$ can be reconstructed from the time-frequency representation $\mathcal{S}^s_\varphi f$ using the formula 
\begin{align}\label{R2009_Molahajloo-Wong_ST_Inversion-Formula}
f=\frac{1}{c_\varphi}\int_{\mathbb{R}}\int_{\mathbb{R}}(\mathcal{S}^s_\varphi f)(\B,\xxi)\varphi^{\B,\xxi}_s\frac{d\B d\xxi}{|\xxi|^{1-\frac{2}{s'}}}.
\end{align}
For $s=1,$ the resolution of identity formula \eqref{R2009_Molahajloo-Wong_ST_Resolution-of-Identity-formula} and the inversion formula \eqref{R2009_Molahajloo-Wong_ST_Inversion-Formula} reduces to the one obtained by Du et al. in \cite{du2007continuous}. Moreover, if $\varphi$ is a normalized window function such that $c_\varphi=1,$ then for all $f,g\in L^2(\mathbb{R}),$ $\langle f,g\rangle_{L^2(\mathbb{R})}=\langle \mathcal{S}_\varphi f,\mathcal{S}_\varphi g \rangle_{L^2\left(\mathbb{R}\times\mathbb{R}_{0},\frac{d\B d\xxi}{|\xxi|}\right)},$ which implied that the transform $\mathcal{S}_\varphi f\in L^2\left(\mathbb{R}\times\mathbb{R}_{0},\frac{d\B d\xxi}{|\xxi|}\right),$ for all $f\in L^2(\mathbb{R})$ and also $\|\mathcal{S}_\varphi f\|_{L^2\left(\mathbb{R}\times\mathbb{R}_{0},\frac{d\B d\xxi}{|\xxi|}\right)}=\|f\|_{L^2(\mathbb{R})}.$ Thus, the transform $\mathcal{S}_\varphi$ is an isometry from $L^2(\mathbb{R})$ to the Hilbert space $L^2\left(\mathbb{R}\times\mathbb{R}_{0},\frac{d\B d\xxi}{|\xxi|}\right)$ with its range $\mathcal{S}_\varphi(L^2(\mathbb{R}))$ as a closed subspace. Furthermore, the range $\mathcal{S}_\varphi(L^2(\mathbb{R}))$ is a RKHS with the kernel $K(\B,\xxi;\B',\xxi')=\overline{(\mathcal{S}_\varphi^{\B,\xxi}\varphi^{\B,\xxi})(\B',\xxi')}=\langle \varphi^{\B',\xxi'},\varphi^{\B,\xxi} \rangle,$ for $(\B,\xxi),(\B'\xxi')\in\mathbb{R}\times\mathbb{R}_{0}$ was proved in \cite{hutnikova2013range}.
\\

So far, we have considered the ST in one dimension. We now see the multidimensional ST and explore the results that have been explored in the literature. For all $\bldxi\in\mathbb{R}^n,$ let $A_{\bldxi}=\left(a_{ij}(\bldxi)\right)_{1\leq i,j\leq n}$ be an invertible $n\times n$ real matrix and let $\varphi\in L^1(\mathbb{R}^n)\cap L^2(\mathbb{R}^n)$ be such that $\int_{\mathbb{R}^n}\varphi(\xx)d\xx=1.$ Then the $n-$dimensional ST $\mathcal{S}_\varphi f$ of $f\in L^2(\mathbb{R}^n)$ with respect to the window function $\varphi$ is given by Riba et al.\cite{riba2013continuous} as
\begin{align}\label{R2013_Riba-Wong_MMNP_MultiD-ST_Definition}
(\mathcal{S}_\varphi f)(\bldb,\bldxi)=(2\pi)^{-\frac{n}{2}}|\det A_{\bldxi}|\int_{\mathbb{R}^n}f(\bldx)e^{-i\bldx\cdot\bldxi}\overline{\varphi(A_{\bldxi}(\bldx-\bldb))}d\bldx,~(\bldb,\bldxi)\in \mathbb{R}^n\times\mathbb{R}^n,
\end{align}
with the inversion formula as $f=\mathcal{F}^{-1}A\mathcal{S}_\varphi f,$ where $\mathcal{F}$ is the inverse FT and the operator $A$ is given by $(AF)(\bldxi)=\int_{\mathbb{R}^n}F(\bldb,\bldxi)d\bldb,~\bldxi\in\mathbb{R}^n$ for all measurable function $F$ on $\mathbb{R}^n\times\mathbb{R}^n$ provided the integral exists. It right hand side of equation \eqref{R2013_Riba-Wong_MMNP_MultiD-ST_Definition} can be can be written in the form  $\langle f,\varphi^{\bldb,\bldxi,A_{\bldxi}}\rangle_{L^2(\mathbb{R}^n)}$ with $\varphi^{\bldb,\bldxi,A_{\bldxi}}=(2\pi)^{-\frac{n}{2}}\mathcal{M}_{\bldxi}\tau_{\bldb}\mathcal{D}_{A_{\bldxi}}\varphi,$ where $(\mathcal{D}_{A_{\bldxi}}\varphi)(\bldx)=|\det A_{\bldxi}|\varphi(A_{\bldxi}\bldx),$ $(\tau_{\bldb}\varphi)(\bldx)=\varphi(\bldx-\bldb)$ and $(\mathcal{M}_{\bldxi}\varphi)(\bldx)=e^{i\bldx\cdot\bldxi}\varphi.$ Author in \cite{riba2013continuous}, gave a resolution of identity formula for the $n-$dimensional ST associated with some special class of matrices, which are discussed below:
 
Let $A_{\bldxi}\in GL(n,\mathbb{R}),~\bldxi\in \mathbb{R}^n\setminus\{0\}$ satisfying the conditions
\begin{enumerate}[label=(\roman*)]
\item\label{R2013_Riba-Wong_MMNP_MatricProperty1} $\boldsymbol\eta=\left(A_{\bldxi}^{-1}\right)^t\boldsymbol\zeta$ for some $\boldsymbol\zeta\in \mathbb{R}^n$ then there exists positive functions $f_1,f_2$ on $\mathbb{R}^n$ such that $\left|\det\left(\frac{\partial\boldsymbol\eta}{\partial\bldxi}\right)\right|=\frac{f_1(\boldsymbol\eta)}{f_2(\bldxi)},\bldxi,\boldsymbol\eta\in\mathbb{R}^n,$ where $\frac{\partial\boldsymbol\eta}{\partial\bldxi}$ is the Jacobian matrix of $\boldsymbol\eta$ with respect to $\boldsymbol\bldxi$
\item\label{R2013_Riba-Wong_MMNP_MatricProperty2} there exists a vector $\boldsymbol v\in\mathbb{R}^n$ such that $\left(A_{\bldxi}^{-1}\right)^t\bldxi=\boldsymbol v,$ 
\end{enumerate}
then, for such matrices the following resolution of identity holds 
\begin{align}\label{R2013_Riba-Wong_MMNP_Resolution-Identity_Multi-DST_Definition}
\int_{\mathbb{R}^n}\int_{\mathbb{R}^n}(\mathcal{S}_\varphi f)(\bldb,\bldxi)\overline{(\mathcal{S}_\varphi g)(\bldb,\bldxi)}d\bldb\frac{d\bldxi}{f_2(\bldxi)}=c_{\varphi}\langle f,g \rangle_{L^2(\mathbb{R}^n)},
\end{align}
for all $f,g\in L^2(\mathbb{R}^n),$ where relation $c_\varphi=\int_{\mathbb{R}^n}\left|\hat{\varphi}\left(\left(A_{\bldxi}^{-1}\right)^t(\boldsymbol\zeta-\bldxi)\right)\right|^2\frac{d\bldxi}{f_2(\bldxi)}=\int_{\mathbb{R}^n}\left|\hat{\varphi}(\boldsymbol\eta-\boldsymbol v)\right|^2\frac{d\boldsymbol\eta}{f_1(\boldsymbol\eta)}<\infty.$ Furthermore, if for all $\bldxi=(\xxi_1,\xxi_2,\cdots,\xxi_n)\in\mathbb{R}^n,~\xxi_l\neq 0,l=1,2,\cdots,n,$ we consider the matrices $I_{\bldxi}=(\delta_{ll'})_{n\times n},\delta_{ll}=\xxi_l$ and $\delta_{ll'}=0,$ if $l\neq l',$ the the corresponding transform is called the non-isotropic $n-$dimensional ST given as 
\begin{align}\label{R2013_Riba-Wong_MMNP_MultiD-NonIsotropic-ST_Definition}
(\mathcal{S}_\varphi f)(\bldb,\bldxi)=(2\pi)^{-\frac{n}{2}}\Pi_{l=1}^n|\xxi_l|\int_{\mathbb{R}^n}e^{-i\bldx\cdot\bldxi}f(\bldx)\overline{\varphi\left(I_{\bldxi}(\bldx-\bldb)\right)}d\bldx
\end{align}
and the corresponding resolution of identity is given as equation \eqref{R2013_Riba-Wong_MMNP_Resolution-Identity_Multi-DST_Definition} with $f_2(\bldxi)=\Pi_{l=1}^n|\xxi_l|$ and $c_\varphi=\int_{\mathbb{R}^n}\left|\hat{\varphi}(\bldxi-\boldsymbol 1)\right|^2\frac{d\bldxi}{\Pi_{l=1}^n|\xxi_l|},$ where $\boldsymbol 1=(1,1,\cdots,1)\in\mathbb{R}^n.$ These results, in \cite{riba2013continuous}, extend those for the two-dimensional ST studied by Liu et al. in \cite{liu2007inversion}. Hleili et al. \cite{hleili2021time} proved that the transform $\mathcal{S}_\varphi f$ given in \eqref{R2013_Riba-Wong_MMNP_MultiD-NonIsotropic-ST_Definition} have the property that for window functions $\varphi$ and $\psi,$ $(\mathcal{S}_\varphi f)(\mathcal{S}_\psi g)\in L^p\left(\mathbb{R}^n\times\mathbb{R}^n,\frac{d\bldb d\bldxi}{|\det I_{\bldxi}|}\right)$ for all $p\in [1,\infty]$ and $f,g\in L^2(\mathbb{R}^n).$ Riba et al.\cite{riba2013continuous} also introduced a special class of matrices $A_{\bldxi}$ satisfying
\begin{enumerate}[label=(\Roman*)]
\item\label{R2013_Riba-Wong_MMNP_Special-Class-MatricProperty1} $\frac{1}{|\bldxi|}A_{\bldxi}\in O(n,\mathbb{R})$ for every $\bldxi\in\mathbb{R}^n\setminus\{0\},$
\item\label{R2013_Riba-Wong_MMNP_Special-Class-MatricProperty2} there exists a matrix $P\in O(n,\mathbb{R})$ such that $A_{\bldxi}\boldsymbol\zeta=P A_{\boldsymbol\zeta}\bldxi,$
\item\label{R2013_Riba-Wong_MMNP_Special-Class-MatricProperty3} there  exists a vector $\boldsymbol v\in\mathbb{R}^n$ such that $\left(A_{\bldxi}^{-1}\right)^t\bldxi=\boldsymbol v.$
\end{enumerate}
This class of matrices are such that $A_{\bldxi}\in GL(n,\mathbb{R}),\bldxi\in \mathbb{R}_{0}$ and it satisfies both the property \ref{R2013_Riba-Wong_MMNP_MatricProperty1} and \ref{R2013_Riba-Wong_MMNP_MatricProperty2}. The matrices 
$A_{\bldxi}^{2\times 2}=\begin{pmatrix}
\xxi_1 & \xxi_2 \\
-\xxi_2 & \xxi_1
\end{pmatrix},$
$A_{\bldxi}^{4\times 4}=\begin{pmatrix}
\xxi_1 & \xxi_2 & \xxi_3 & \xxi_4\\
-\xxi_2 & \xxi_1 & \xxi_4 & -\xxi_3\\
-\xxi_3 & -\xxi_4 & \xxi_1 & \xxi_2\\
-\xxi_4 & \xxi_3 & -\xxi_2 & \xxi_1
\end{pmatrix}$
and 
$A_{\bldxi}^{8\times 8}=\begin{pmatrix}
\xxi_1 & \xxi_2 & \xxi_3 & \xxi_4 &\xxi_5 & \xxi_6 & \xxi_7 & \xxi_8\\
-\xxi_2 & \xxi_1 & \xxi_4 & -\xxi_3 & \xxi_6 & -\xxi_5 & -\xxi_8 & \xxi_7\\
-\xxi_3 & -\xxi_4 & \xxi_1 & \xxi_2 & \xxi_7 & \xxi_8 & -\xxi_5 & -\xxi_6\\
-\xxi_4 & \xxi_3 & -\xxi_2 & \xxi_1 & \xxi_8 & -\xxi_7 & \xxi_6 & -\xxi_5\\
-\xxi_5 & -\xxi_6 & -\xxi_7 & -\xxi_8 &\xxi_1 & \xxi_2 & \xxi_3 & \xxi_4\\
-\xxi_6 & \xxi_5 & -\xxi_8 & \xxi_7 & -\xxi_2 & \xxi_1 & -\xxi_4 & \xxi_3\\
-\xxi_7 & \xxi_8 & \xxi_5 & -\xxi_6 & -\xxi_3 & \xxi_4 & \xxi_1 & -\xxi_2\\
-\xxi_8 & -\xxi_7 & \xxi_6 & \xxi_5 & -\xxi_4 & -\xxi_3 & \xxi_2 & \xxi_1
\end{pmatrix}$
of order $2,$ $4$ and $8$ respectively satisfy properties \ref{R2013_Riba-Wong_MMNP_Special-Class-MatricProperty1},\ref{R2013_Riba-Wong_MMNP_Special-Class-MatricProperty2} and \ref{R2013_Riba-Wong_MMNP_Special-Class-MatricProperty3}. Authors in \cite{riba2015continuous} explained that the family $A_{\bldxi}\in GL(n,\mathbb{R}),\bldxi\in \mathbb{R}_{0}$ for which the properties \ref{R2013_Riba-Wong_MMNP_MatricProperty1} and \ref{R2013_Riba-Wong_MMNP_MatricProperty2} hold cannot be satisfied for the dimensions apart from $1,2,4$ and $8.$ They also extended the multidimensional ST to the multidimensional modified ST. For that, they considered the multidimensional dilation operator $\mathcal{D}^s_A\varphi,~1\leq s<\infty,$ defined for $A\in GL(n,\mathbb{R})$ as
\begin{align}\label{RRiba-Wong_InversionFormula_Multi-DModifiedST_DilationOperator}
(\mathcal{D}^s_A\varphi)(\bldx)=|\det A|^{-\frac{1}{s}}\varphi(A^{-1}\bldx),~\bldx\in\mathbb{R}^n
\end{align}
for all measurable function $\varphi$ on $\mathbb{R}^n.$ In particular, $\mathcal{D}^2_A$ is a unitary operator on $L^2(\mathbb{R}^n).$ Also let, $A:\mathbb{R}^n\rightarrow GL(n,\mathbb{R}),$ be given by $\bldxi\mapsto A_{\bldxi}.$ Then the multi-dimensional  modified ST of $f\in L^2(\mathbb{R}^n)$ with respect to $\varphi\in L^2(\mathbb{R}^n)$ was given by Riba et al. \cite{riba2015continuous} as
\begin{align}\label{RRiba-Wong_InversionFormula_Multi-DModifiedST_Definition}
(\mathcal{S}^s_{A,\varphi}f)(\bldb,\bldxi)=(2\pi)^{-\frac{n}{2}}|\det A_{\bldxi}|^{-\frac{1}{s}}\int_{\mathbb{R}^n}f(\bldx)e^{-i\bldx\cdot\bldxi}\overline{\varphi(A_{\bldxi}^{-1}(\bldx-\bldb))}d\bldx,~(\bldb,\bldxi)\in \mathbb{R}^n\times\mathbb{R}^n,
\end{align}
which can be expressed in several equivalent form as 
\begin{align}
(\mathcal{S}^s_{A,\varphi}f)(\bldb,\bldxi)&=|\det A_{\bldxi}|^{\frac{1}{2}-\frac{1}{s}}(\mathcal{S}^2_{A,\varphi}f)(\bldb,\bldxi)\\
&=(2\pi)^{-\frac{n}{2}}|\det A_{\bldxi}|^{\frac{1}{2}-\frac{1}{s}}\langle f,\mathcal{M}_{\bldxi}\tau_{\bldb}\mathcal{D}^s_{A_{\bldxi}}\varphi\rangle_{L^2(\mathbb{R}^n)}.
\end{align}
Equation \eqref{RRiba-Wong_InversionFormula_Multi-DModifiedST_Definition} can also be expressed in terms of the FT $(\mathcal{F}f)(\bldxi)=\hat{f}(\bldxi)=(2\pi)^{-\frac{n}{2}}\int_{\mathbb{R}^n}f(\bldx)e^{-2\pi i\bldx.\bldxi}d\bldt$ as
\begin{align}\label{RRiba-Wong_InversionFormula_Multi-DModifiedST_EquivalentForm_Fourier-Domain}
(\mathcal{S}^s_{A,\varphi}f)(\bldb,\bldxi)=|\det A_{\bldxi}|^{1-\frac{1}{s}}e^{-i\bldb\cdot\bldxi}\left(\mathcal{F}^{-1}_{\boldsymbol\zeta\mapsto\bldb}f_{\bldxi,A_{\bldxi}}\right)(\bldb),~\bldb,\bldxi\in\mathbb{R}^n,
\end{align}
where $f_{\bldxi,A_{\bldxi}}(\boldsymbol \zeta)=\hat{f}(\boldsymbol \zeta)\overline{\hat{\varphi}\left(A_{\bldxi}^t(\boldsymbol \zeta-\bldxi)\right)},~\boldsymbol \zeta\in\mathbb{R}^n.$ Moreover, the inversion formula for the modified ST of $f\in L^2(\mathbb{R}^n)$ and the window function $\varphi\in L^1(\mathbb{R}^2)\cap L^2(\mathbb{R}^n)$ satisfying $\int_{\mathbb{R}^n}\varphi(\bldx)d\bldx=1$ can be obtained from the relation $\hat{f}(\bldxi)=\left|\det A_{\bldxi}\right|^{\frac{1}{s}-1}\int_{\mathbb{R}^n}(\mathcal{S}^s_{A,\varphi}f)(\bldb,\bldxi)d\bldb,~\bldxi\in\mathbb{R}^n.$ Authors in \cite{riba2015continuous} also obtained the resolution of identity formula for $\mathcal{S}^s_{A,\varphi}$ as 
\begin{align}\label{RRiba-Wong_Resolution-of-dentity_Multi-DModifiedST}
\int_{\mathbb{R}^n}\int_{\mathbb{R}^n}(\mathcal{S}^s_{A,\varphi}f)(\bldb,\bldxi)\overline{(\mathcal{S}^s_{A,\varphi}g)(\bldb,\bldxi)}d\bldb\left|\det A_{\bldxi}\right|^{\frac{2}{s}-1}d\bldxi=c_{\varphi}\langle f,g \rangle_{L^2(\mathbb{R}^n)},
\end{align}
for all $f,g,\varphi\in L^2(\mathbb{R}^n)$ with $c_{\varphi}=\int_{\mathbb{R}^n}|\hat{\varphi}(\bldxi)|^2\frac{d\bldxi}{\left|\det \left(F^i_{jk}\xxi^j+\delta^i_k\right)\right|}<\infty,$ where $A:\mathbb{R}^n\rightarrow GL(n,\mathbb{R})$ be given by $\bldxi\mapsto A_{\bldxi}$ is a piecewise differentiable function for which there exist fixed $(1,2)-$tensor $F$ and a fixed $(1,1)-$tensor $G$ such that $\left(A_{\bldxi}^t\right)^{-1}=\left(F^i_{jk}\xxi^j+\delta^i_k\right)_{1\leq i,j\leq n}$ and  $\eta_{\boldsymbol\zeta}(\bldxi):=A_{\bldxi}^t(\boldsymbol\zeta-\bldxi),~\bldxi\in\mathbb{R}^n$ is such that $\eta_{\boldsymbol\zeta}(\mathbb{R}^n)=\mathbb{R}^n$ for all $\boldsymbol\zeta\in\mathbb{R}^n.$ These assumptions on the map $A:\mathbb{R}^n\rightarrow GL(n,\mathbb{R}^n)$ are satisfied by the matrix-valued functions $A_{\bldxi}^{2\times 2},$ $A_{\bldxi}^{4\times 4}$ and $A_{\bldxi}^{8\times 8}$ mentioned earlier. Moreover, the equation \eqref{RRiba-Wong_Resolution-of-dentity_Multi-DModifiedST} is applicable in all the dimensions to the constant and diagonal matrix-valued functions $I_{\bldxi}=(\delta_{ll'})_{n\times n},\delta_{ll}=\xxi_l$ and $\delta_{ll'}=0,$ if $l\neq l'.$ 

The two-dimensional ST given in  \cite{liu2007inversion} is inadequate for capturing anisotropic features like corners and edges of signals due to the directional insensitivity of the window function. So, Shah et al. \cite{shah2021non} proposed a novel transform called the non-isotropic angular ST with respect to a window function $\psi\in L^2(\mathbb{R}^2)$ defined as 
\begin{align}\label{R2021_Shah-Tantary_AA_NonIsotropic_AST_Definition}
(\mathcal{M}_\psi f)(\bldb,\bldxi,\theta)=\frac{1}{2\pi}\langle f,\psi^\theta_{\bldb,\bldxi}\rangle_{L^2(\mathbb{R}^2)}
\end{align}
for any $f\in L^2(\mathbb{R}^2),$ with 
\begin{align}\label{R2021_Shah-Tantary_AA_NonIsotropic_AST_DaughterWindows}
\psi^\theta_{\bldb,\bldxi}(\bldt)=E_{\bldxi}R_{\theta}\mathcal{D}_{A_{\bldxi}}\tau_{\bldb}\psi(\bldt)=|\det A_{\bldxi}|e^{i\bldt\cdot\bldxi}\psi(R_{-\theta}A_{\bldxi}(\bldt-\bldb)),
\end{align}
where $(\mathcal{D}_{A_{\bldxi}}\psi)(\bldt)=|\det A_{\bldxi}|\psi(A_{\bldxi}\bldt),$ $(E_{\bldxi}\psi)(\bldt)=e^{i\bldxi\cdot\bldt}\psi(\bldt),$ $\tau_{\bldb}\psi)(\bldt)=\psi(\bldt-\bldb),$ $(R_\theta \psi)(\bldt)=\psi(R_{-\theta}\bldt),$ $A_{\bldxi}=\begin{pmatrix}
\xxi_1 & 0 \\
 0 & \xxi_2
\end{pmatrix},\xxi_1,\xxi_2\neq 0$  and $R_\theta=\begin{pmatrix}
\cos\theta & \sin\theta \\
 -\sin\theta & \cos\theta
\end{pmatrix},$ $\bldxi,\bldb\in\mathbb{R}^2,\theta\in[0,2\pi).$ Apart from the several important properties like linearity, anti-linearity, translation, scaling, reflection and conjugation, authors expressed the non-isotropic angular ST in the Fourier domain as $(\mathcal{M}_\psi f)(\bldb,\bldxi,\theta)=e^{-i\bldxi\cdot\bldb}\mathcal{F}^{-1}\left[\hat{f}(\boldsymbol\omega)\overline{\hat{\psi}\left((\boldsymbol\omega-\bldxi)A_{\bldxi}^{-1}R_\theta\right)}\right](\bldb)$ and obtained the following resolution of identity formula 
\begin{align}\label{R2021_Shah-Tantary_AA_NonIsotropic_AST_IPR}
\int_{\mathbb{R}^2}\int_{\mathbb{R}^2\setminus\{0\}}\int_0^{2\pi} (\mathcal{M}_\psi f)(\bldb,\bldxi,\theta)\overline{(\mathcal{M}_\psi g)(\bldb,\bldxi,\theta)}d\theta d\bldxi d\bldb=C_\psi\langle f,g \rangle_{L^2(\mathbb{R}^2)},~f,g\in L^2(\mathbb{R}^2)
\end{align}
$C_\psi=\int_{\mathbb{R}^2\setminus\{0\}}\int_0^{2\pi}\left|(\mathcal{F}\psi)\left((\boldsymbol\omega-\bldxi)A_{\bldxi}^{-1}R_\theta\right)\right|^2d\theta d\bldxi<\infty.$ This in particular, for $f=g,$ implies that $\|\mathcal{M}_{\psi}f\|^2_{L^2(\mathbb{R}^2)}=C_\psi\|f\|^2_{L^2(\mathbb{R})},$ which further shows that for the function $\psi\in L^2(\mathbb{R}^2)$ satisfying $C_\psi=1,$ $\mathcal{M}_{\psi}$ is an isometry from $L^2(\mathbb{R}^2)$ into $L^2(\mathbb{R}^2\times\mathbb{R}^2\setminus\{0\}\times [0,2\pi))$ and hence it is one-one. Moreover, $\mathcal{M}_{\psi}(L^2(\mathbb{R}^2))$ is a proper subspace of $L^2(\mathbb{R}^2\times\mathbb{R}^2\setminus\{0\}\times [0,2\pi)).$ Thus, square integrability of $\mathcal{M}_{\psi}f$ is necessary condition but not sufficient.  The range $\mathcal{M}_{\psi}(L^2(\mathbb{R}^2)),$  of  bounded linear operator $\mathcal{M}_{\psi}f:L^2(\mathbb{R}^2)\to L^2(\mathbb{R}^2\times\mathbb{R}^2\setminus\{0\}\times [0,2\pi))$ which is isometric when $C_\psi=1,$ is a closed subspace. Being closed subspace of the Hilbert space $L^2(\mathbb{R}^2\times\mathbb{R}^2\setminus\{0\}\times [0,2\pi)),$ $\mathcal{M}_{\psi}(L^2(\mathbb{R}^2))$ is also a Hilbert space. Furthermore, it is a RKHS, with the kernel given by $K_\psi(\bldb,\bldxi,\theta;\boldsymbol c,\boldsymbol \eta,\theta')=\frac{1}{(2\pi)^2}\langle \psi^{\theta}_{\bldb,\bldxi},\psi^{\theta'}_{\boldsymbol c,\boldsymbol \eta} \rangle_{L^2(\mathbb{R}^2)}.$ The proposed transform \eqref{R2021_Shah-Tantary_AA_NonIsotropic_AST_Definition}, has the following inversion formula 
%\textcolor{red}{(See some thing important lines can be written here, also not the the inversion formula for the window FT and the WT in the earlier section is not written)}
\begin{align}\label{R2021_Shah-Tantary_AA_NonIsotropic_AST_InversionFormula}
f(\bldt)=\frac{1}{2\pi C_\psi}\int_{\mathbb{R}^2}\int_{\mathbb{R}^2\setminus\{0\}}\int_0^{2\pi}(\mathcal{M}_\psi f)(\bldb,\bldxi,\theta)\psi^{\theta}_{\bldb,\bldxi}(\bldt)d\theta d\bldxi d\bldb,~a.e.
\end{align}
where $\psi^{\theta}_{\bldb,\bldxi}(\bldt)$ is given in \eqref{R2021_Shah-Tantary_AA_NonIsotropic_AST_DaughterWindows}.\\

\noindent\textbf{Uncertainty Principles:} Motivated by the idea of Singer \cite{singer1999uncertainty} and Wilczok \cite{wilczok2000new} in the case of WFT and WT, authors in \cite{shah2021non} obtained a class of UPs results for comparing the localization of $f$ and $\mathcal{F}{f}$ with the localization of the non-isotropic angular ST regarded as function of one variable as
$R_{\bldb,\boldsymbol\omega}\geq\frac{C_{\psi}}{2}\|f\|^2_2$ and $R_{\bldxi,\bldt}\geq\frac{C_{\psi}}{2}\|f\|^2_2,$ where 
$$R_{\bldb,\boldsymbol\omega}=\left\{\int_{\mathbb{R}^2}\int_{\mathbb{R}^2\setminus\{0\}}\int_0^{2\pi} |\bldb|^2|(\mathcal{M}_\psi f)(\bldb,\bldxi,\theta)|^2d\theta d\bldxi d\bldb\right\}^{\frac{1}{2}}\left\{\int_{\mathbb{R}^2}|\boldsymbol\omega|^2|(\mathcal{F}f)(\boldsymbol\omega)|^2d\boldsymbol\omega\right\}^{\frac{1}{2}}$$ 
and 
$$R_{\bldxi,\bldt}=\left\{\int_{\mathbb{R}^2}\int_{\mathbb{R}^2\setminus\{0\}}\int_0^{2\pi} |\bldxi|^2|(\mathcal{M}_\psi f)(\bldb,\bldxi,\theta)|^2d\theta d\bldxi d\bldb\right\}^{\frac{1}{2}}\left\{\int_{\mathbb{R}^2}|\bldt|^2|f(\bldt)|^2d\bldt\right\}^{\frac{1}{2}}.$$ 
As a consequence of which we have $R_{\bldb,\boldsymbol\omega}R_{\bldxi,\bldt}\geq\frac{C_{\psi}}{2}\|f\|^2_2,$ which shows that the phase-space localization of the pair $(f,\mathcal{F}f)$ varies inversely with the phase space localization of the transform \eqref{R2021_Shah-Tantary_AA_NonIsotropic_AST_Definition}, i.e., better the phase-space localization of $(f,\mathcal{F}f)$, the worse is the phase-space localization of $(\mathcal{M}_\psi f)(\bldb,\bldxi,\theta)$ and vice versa. In fact, a more general version for $f\in L^p(\mathbb{R}^2)$ is proved by the authors as
\begin{align}\label{R2021_Shah-Tantary_AA_NonIsotropic_AST_GeneralUP_1<=p<=2}
\left\{\int_{\mathbb{R}^2}\int_{\mathbb{R}^2\setminus\{0\}}\int_0^{2\pi} |\bldb|^p|(\mathcal{M}_\psi f)(\bldb,\bldxi,\theta)|^pd\theta d\bldxi d\bldb\right\}^{\frac{1}{p}}\left\{\int_{\mathbb{R}^2}|\boldsymbol\omega|^p|(\mathcal{F}f)(\boldsymbol\omega)|^pd\boldsymbol\omega\right\}^{\frac{1}{p}}\geq\frac{C_{\psi}^{\frac{1}{2}}}{2}\|f\|^2_2,~1\leq p\leq 2
\end{align}
and
\begin{align}\label{R2021_Shah-Tantary_AA_NonIsotropic_AST_GeneralUP_1<=p<=2}
\left\{\int_{\mathbb{R}^2}\int_{\mathbb{R}^2\setminus\{0\}}\int_0^{2\pi} |\bldb|^p|(\mathcal{M}_\psi f)(\bldb,\bldxi,\theta)|^2d\theta d\bldxi d\bldb\right\}^{\frac{1}{p}}\left\{\int_{\mathbb{R}^2}|\boldsymbol\omega|^p|(\mathcal{F}f)(\boldsymbol\omega)|^2d\boldsymbol\omega\right\}^{\frac{1}{p}}\geq\frac{C_{\psi}^{\frac{1}{p}}}{2}\|f\|^{\frac{4}{p}}_2,~p\geq 2
\end{align}
following the idea of Cowling and Price \cite{cowling1984bandwidth}. Furthermore, they also obtained the following logarithmic UP for the transform \eqref{R2021_Shah-Tantary_AA_NonIsotropic_AST_Definition} using the Pitt's inequality \cite{beckner1995pitt}
\begin{align}\label{R2021_Shah-Tantary_AA_NonIsotropic_AST_LUP}
\int_{\mathbb{R}^2}\int_{\mathbb{R}^2\setminus\{0\}}\int_0^{2\pi}\ln |\bldb||(\mathcal{M}_\psi f)(\bldb,\bldxi,\theta)|^2d\theta d\bldxi d\bldb+C_\psi \int_{\mathbb{R}^2}\ln|\boldsymbol\omega||(\mathcal{F}f)(\boldsymbol\omega)|^2d\boldsymbol\omega\geq C_{\psi}\|f\|^2_2\left(\frac{\Gamma'(\frac{1}{2})}{\Gamma(\frac{1}{2})}-\ln\pi\right).
\end{align}

\section{Short time fractional Fourier transform (STFrFT)}
The FrFT only reveals the overall fractional Fourier domain (FrFD) frequency content of the signal. Sometimes, we want to know the FrFD-frequency contents and how they change over time. So, the representation combining the time and FrFD-frequency information should be developed. A localized FrFT was developed in \cite{mendlovic1996signal},\cite{zalevsky1997localized},\cite{ozaktas2000applications} based on which time/space-dependent FrFD filtering approach was obtained. But they cannot be regarded as the time-FrFD-frequency representation since it has no explicit $2-D$ domain. A kind of STFrFT was proposed in \cite{capus2003short}. It modifies the classical STFrFT by utilizing a dechirp window whose chirp rate is obtained by the FrFT in a searching scheme. An adaptive STFrFT was given in \cite{zhang2000chip} for the function $f\in L^2(\mathbb{R})$ with respect to the FrFT domain window function $g(\oomega')\in L^2(\mathbb{R})$ as 
\begin{align}\label{R2000_Zhang-Bi-Chen_Conference_STFrFT_Definition}
(\mathcal{G}^\alpha_gf)(\oomega,\tee)=\int_{\mathbb{R}}(\mathcal{F}_\alpha f)(\oomega')\overline{g(\oomega'-\oomega)}e^{-i\oomega'\tee}d\oomega',
\end{align}  
which rotates the signal in the time-frequency domain by the FrFT before implementing the STFT. They could be classified as a kind of adaptive STFT rather than time-FrFD-frequency representation because they analyze signals in the time-frequency domain but not in the time-FrFD-frequency domain. Such a definition was independently derived by Stankovi\`c et al. \cite{stankovic2003time} for time-frequency analysis of chirp signals. From equation \eqref{R2000_Zhang-Bi-Chen_Conference_STFrFT_Definition}, it follows that the STFrFT, which is given in \cite{stankovic2003time},\cite{zhang2000chip} looses all time-varying information and hence cannot be viewed as a joint time-FrFD-frequency analysis transform. As a kind of time-FrFD-frequency representation, the different STFrFT was in Zhang et al. \cite{zhang1998new} and Tao et al. \cite{tao2009short}, which is defined for functions $f\in L^2(\mathbb{R})$ with respect to the window function $g\in L^2(\mathbb{R})$
\begin{align}\label{R2010_Tao-Li-Wang_IEEETSP_STFrFT_Definition}
\left(^{(0)}\mathcal{G}^\alpha_gf\right)(\oomega,\tee)=\frac{1}{\sqrt{2\pi}}\int_{\mathbb{R}}f(\xx)g(\xx-\tee)K_\alpha(\xx,\oomega)d\xx.
\end{align}
The transform \eqref{R2010_Tao-Li-Wang_IEEETSP_STFrFT_Definition} is linear, which indicates that it satisfies the superposition principle, which is favorable for the analysis of multicomponent signals. It satisfies the time marginal constraints, i.e., integrating the STFrFT  with respect to the time variable gives the LCT multiplied by a constant, i.e., $\int_{\mathbb{R}}\left(^{(0)}\mathcal{G}^\alpha_gf\right)(\oomega,\tee)d\tee=C(\mathcal{F}_\alpha f)(\oomega),$ for some constant $C$ and the additivity of rotation, i.e., $\mathcal{F}_\beta\left(\left(^{(0)}\mathcal{G}^\alpha_gf\right)(\oomega,\tee)\right)(\oomega',\tee)=\left(^{(0)}\mathcal{G}^{\alpha+\beta}_gf\right)(\oomega',\tee).$ Apart from this, author in \cite{tao2009short} proved the covariance properties 
$$\left(^{(0)}\mathcal{G}^\alpha_g[e^{i\oomega_0\xx}f(\xx-\tee_0)]\right)(\oomega,\tee)=e^{i\phi(\oomega)}\left(^{(0)}\mathcal{G}^\alpha_gf\right)(\oomega-(\tee_0\cos\alpha+\oomega_0\sin\alpha),\tee-\tee_0),$$
$\phi(\oomega)=(\oomega_0\cos\alpha-\tee_0\sin\alpha)\oomega+\frac{\tee_0^2-\oomega_0^2}{2}\sin\alpha\cos\alpha+\tee_0\oomega_0\sin^2\alpha.$ With the  assumption that the window function satisfies the conditions $g(0)\neq 0,$ the $1-D$ inversion formula is 
\begin{align}\label{R2009_Tao_STFT-1D-Inversion-Formula}
f(\tee)=\frac{1}{g(0)}\int_{\mathbb{R}}\left(^{(0)}\mathcal{G}^\alpha_gf\right)(\oomega,\tee)\overline{K_\alpha(\tee,\oomega)}d\oomega
\end{align}
and under the assumption that $\|g\|_{L^2(\mathbb{R})}=1,$ the $2-D$ inversion formula is $f(\xx)=\int_{\mathbb{R}}\int_{\mathbb{R}}\left(^{(0)}\mathcal{G}^\alpha_gf\right)(\oomega,\tee)g(\xx-\tee)\overline{K_\alpha(\xx,\oomega)}d\tee d\oomega$ are also obtained in \cite{tao2009short}. Based on the fractional convolution $\star_\alpha,$ Shi et al. \cite{shi2020novel} obtained a novel STFrFT which is defined for function $f\in L^2(\mathbb{R})$ with respect to the window function $g\in L^2(\mathbb{R})$ as 
\begin{align}\label{R2020_Shi-Xiang-Zhang_IEEETSP_NovelSTFT_Definition}
\left(^{(1)}\mathcal{G}^\alpha_gf\right)(\oomega,\tee)&=e^{-i\tee \oomega\csc\alpha}\left[f(\xx)\star_\alpha\left(\overline{g(-\xx)}e^{i\xx \oomega\csc\alpha}\right)\right](\oomega)\notag\\
&=\int_{\mathbb{R}}f(\xx)\overline{g_{\alpha,\tee,\oomega}(\xx)}d\xx,
\end{align}
where $g_{\alpha,\tee,\oomega}(\xx)=g(\xx-\tee)e^{-i\frac{\xx^2-\tee^2}{2}\cot\alpha+i\xx \oomega \csc\alpha}.$ This definition can be equivalently written in fractional Fourier domain as 
\begin{align}\label{R2020_Shi-Xiang-Zhang_IEEETSP_NovelSTFT_Definition-in_FrFoureirDomain}
\left(^{(1)}\mathcal{G}^\alpha_gf\right)(\oomega,\tee)=e^{-i\tee \oomega\csc\alpha}\int_{\mathbb{R}}(\mathcal{F}_\alpha f)(\oomega')\overline{(\mathcal{F}g)((\oomega'-\oomega)\csc\alpha)}\overline{K_\alpha(\oomega',\tee)}d\oomega',
\end{align}
where $(\mathcal{F}g)(\oomega'\csc\alpha)$ denotes the FT of $g$ with its argument scaled by $\csc\alpha.$ Transform \eqref{R2020_Shi-Xiang-Zhang_IEEETSP_NovelSTFT_Definition} is linear in $f$ and anti-linear in $g.$ It satisfy the property called the additivity of rotation i.e., $\mathcal{F}_\beta\left\{A_\alpha e^{\frac{\tee^2+\oomega^2}{2}\cot\alpha}\left(^{(1)}\mathcal{G}^\alpha_gf\right)(\oomega,\tee)\right\}(\oomega',\tee)=\left(^{(1)}\mathcal{G}^{\alpha+\beta}_gf\right)(\oomega',\tee),$ $A_\alpha=\sqrt{\frac{1-i\cot\alpha}{2\pi}}.$ Apart from this, it has the following covariance properties 
\begin{align*}
\left(^{(1)}\mathcal{G}^\alpha_g\left(\tau_{\xx_0}f\right)\right)(\oomega,\tee)=(^{(1)}\mathcal{G}^\alpha_g f)(\oomega,\tee-\tee_0)e^{i\left(\frac{3}{2}\xx_0^2-\xx_0 \xx\right)\cot\alpha-i\tee_0 \oomega\csc\alpha}
\end{align*}
and 
\begin{align*}
\left(^{(1)}\mathcal{G}^\alpha_g\left(e^{i\oomega_0 \xx}f(\xx)\right)\right)(\oomega,\tee)=\left(^{(1)}\mathcal{G}^\alpha_g f\right)(\oomega-\oomega_0\sin\alpha,\tee).
\end{align*}
The covariance property obtained above indeed gives the formula for the WQPFT for the time and frequency shift of a function $f.$ For $f_1,f_2,g\in L^2(\mathbb{R}),$ the WQPFT \eqref{R2020_Shi-Xiang-Zhang_IEEETSP_NovelSTFT_Definition} satisfies the following inner product relation 
\begin{align}\label{R2020_Shi-Xiang-Zhang_IEEETSP_NovelSTFT_IPR}
\csc\alpha\int_{\mathbb{R}}\int_{\mathbb{R}}\left(^{(1)}\mathcal{G}^\alpha_gf_1\right)(\oomega,\tee)\overline{\left(^{(1)}\mathcal{G}^\alpha_g f_2\right)(\oomega,\tee)}d\tee d\oomega=\langle f_1,f_2 \rangle.
\end{align}
In particular, if $f_1=f_2=f$ then $(\csc\alpha)\left\|^{(1)}\mathcal{G}^\alpha_gf\right\|^2_{L^2(\mathbb{R}^2)}=\|f\|^2_{L^2(\mathbb{R})}$ and thus WQPFT is an isometry from $L^2(\mathbb{R})$ to $L^2(\mathbb{R}^2)$ upto some constant, this is called the energy conservation property. With the assumption that the time domain window function $g$ satisfying the condition $g(0)\neq 0,$ author obtained the $1-D$ inversion formula for the transform \eqref{R2020_Shi-Xiang-Zhang_IEEETSP_NovelSTFT_Definition} as $f(\tee)=\frac{\csc\alpha}{2\pi\overline{g(0)}}\int_{\mathbb{R}}\left(^{(1)}\mathcal{G}^\alpha_gf\right)(\oomega,\tee)e^{i\tee\oomega\csc\alpha}d\oomega$ and also for the time domain window function $\gamma(\xx),$ $h(\xx):=\frac{\csc\alpha}{2\pi}\int_{\mathbb{R}}\int_{\mathbb{R}}\left(^{(1)}\mathcal{G}^\alpha_gf\right)(\oomega,\tee)\gamma_{\alpha,\tee,\oomega}(\xx)d\tee d\oomega$ is such that if $\int_{\mathbb{R}}\gamma(\xx)\overline{g(\xx)}d\xx=1$ then $h(\xx)=f(\xx).$ There are several choices  of the window function $\gamma(\xx)$ that satisfies $\int_{\mathbb{R}}\gamma(\xx)\overline{g(\xx)}d\xx=1.$ One such important choice is that $\gamma(\xx)=g(\xx).$ For which $\int_{\mathbb{R}}|g(\xx)|^2d\xx=1$ and the $2-D$ reconstruction formula for the WLCT is given as 
\begin{align}\label{R2020_Shi-Xiang-Zhang_IEEETSP_NovelSTFT_Inversion-Formula}
f(\xx)=\frac{\csc\alpha}{2\pi}\int_{\mathbb{R}}\int_{\mathbb{R}}\left(^{(1)}\mathcal{G}^\alpha_gf\right)(\tee,\oomega)g_{\alpha,\tee,\oomega}(\xx)d\tee d\oomega.
\end{align}
The transform \eqref{R2020_Shi-Xiang-Zhang_IEEETSP_NovelSTFT_Definition} represents a $1-D$ signal $f(\xx)$ using the $2-D$ function $\left(^{(1)}\mathcal{G}^\alpha_gf\right)(\xx,\uu).$ Based on the energy conservation property, it follows that $^{(1)}\mathcal{G}^\alpha_gf\in L^2(\mathbb{R}^2).$  It is not in general true that, a function $F$ in $L^2(\mathbb{R}^2)$ is a WQPFT of some function in $L^2(\mathbb{R}).$ The author in \cite{shi2020novel} characterized the range $^{(1)}\mathcal{G}^\alpha_g(L^2(\mathbb{R}))$ of the transform \eqref{R2020_Shi-Xiang-Zhang_IEEETSP_NovelSTFT_Definition} showing that $F\in~^{(1)}\mathcal{G}^\alpha_g(L^2(\mathbb{R}))$ if and only if it satisfies $F(\tee',\oomega')=\frac{1}{2\pi}\int_{\mathbb{R}}\int_{\mathbb{R}}F(\tee,\oomega)\langle g_{\alpha,\tee,\oomega},g_{\alpha,\tee',\oomega'}\rangle_{L^2(\mathbb{R})}d\tee d\oomega.$ Consequently, it follows that the range  $^{(1)}\mathcal{G}^\alpha_g(L^2(\mathbb{R}))$ is a RKHS with the kernel $\frac{1}{2\pi}\langle g_{\alpha,\tee,\oomega},g_{\alpha,\tee',\oomega'}\rangle_{L^2(\mathbb{R})}.$ In view of equation \eqref{R2020_Shi-Xiang-Zhang_IEEETSP_NovelSTFT_Definition} and its representation in the fractional Fourier domain  \eqref{R2020_Shi-Xiang-Zhang_IEEETSP_NovelSTFT_Definition-in_FrFoureirDomain}, the WQPFT gives the local information of $f(\xx)$ in the time window $\left[E_g+\tee-\Delta_g,E_g+\tee+\Delta_g\right]$ and the local information of $f(\xx)$ in the fractional Fourier domain $\left[E_G\sin\alpha+\oomega-\Delta_G\sin\alpha,E_G\sin\alpha+\oomega+\Delta_G\sin\alpha\right]$, where $G$ is the FT of $g$ having the constant window area $4\Delta_g\Delta_G\sin\alpha$ in the time-fractional-frequency plane with the optimal window function as $g_\sigma(\tee)=\frac{1}{\sqrt{\sqrt{\pi}\sigma}}e^{-\frac{\tee^2}{2\sigma^2}},$ where $\sigma$ is a positive number. Kamalakkannan et al.\cite{kamalakkannan2021short} proposed the short time coupled FrFT for the function $f\in L^2(\mathbb{R}^2)$ with respect to a non-zero $g\in L^2(\mathbb{R}^2)$ and the  parameters $\alpha,\beta\in \mathbb{R}$ satisfying $\alpha+\beta\neq 2n\pi,~n\in\mathbb{Z}$ as
\begin{align}\label{R2021_Kamalakkannan-Roopkumar-Zayed_Coupled_STFT_Definition}
(\mathcal{G}^{\alpha,\beta}_gf)(\bldomega,\bldt)=\int_{\mathbb{R}^2}f(\bldx)g(\bldy-\bldt)K_{\alpha,\beta}(\bldy,\bldomega)d\bldy.
\end{align}
The authors proved the linearity, time delay, frequency shift, and addition property of the transform. On applying the equivalent form $(\mathcal{G}^{\alpha,\beta}_gf)(\bldomega,\bldt)=(\mathcal{F}_{\alpha,\beta}\left(f\tau_{\bldt}\bar{g}\right))(\bldomega),$ $\tau_{\bldt}\bar{g}(\bldy)=\overline{g(\bldy-\bldt)}$ of  \eqref{R2021_Kamalakkannan-Roopkumar-Zayed_Coupled_STFT_Definition}, authors obtained the Parseval's formula $\left\langle \mathcal{G}^{\alpha,\beta}_gf_1,\mathcal{G}^{\alpha,\beta}_gf_2\right\rangle_{L^2(\mathbb{R}^2)}=\frac{1}{2\pi}\|g\|^2_{L^2(\mathbb{R}^2)}\left\langle f_1,f_2\right \rangle_{L^2(\mathbb{R}^2)},$ and the inversion formula 
$$f=\frac{1}{\|g\|^2_{L^2(\mathbb{R}^2)}}\int_{\mathbb{R}^2}\int_{\mathbb{R}^2}\left(\mathcal{G}^{\alpha,\beta}_gf\right)(\bldomega,\bldt) M_{\bldomega,\alpha,\beta}(\tau_{\bldt}g)d\bldomega d\bldt,$$ weakly in $L^2(\mathbb{R}^2),$ where $M_{\bldomega,\alpha,\beta}(h)(\bldy)=h(\bldy)\overline{K_{\alpha,\beta}(\bldy,\bldomega)}.$ They also characterized the range $\mathcal{G}^{\alpha,\beta}_g\left(L^2(\mathbb{R}^2)\right)$ and proved that it is a reproducing kernel Hilbert space with the kernel $\frac{2\pi}{\|g\|^2_{L^2(\mathbb{R}^2)}}\left\langle M_{\bldomega,\alpha,\beta}(\tau_{\bldomega}g), M_{\bldomega',\alpha,\beta}(\tau_{\bldomega'}g)\right\rangle$\\

\noindent\textbf{Uncertainty Principles:} Using relation \eqref{R2022_Verma-Gupta_IJWMIP_nD-FrFT-Heisenberg_UP_2FrFDomain} for $n=1,$ author in \cite{guanlei2009logarithmic} obtained the Heisenberg UP for the STFrFT \eqref{R2000_Zhang-Bi-Chen_Conference_STFrFT_Definition} as
\begin{align}\label{R2009_Guanlei_Xiaotong_Xiaogang_SP_STFrFT_Heisenberg-UP}
\Delta\tee^2_\alpha \Delta\oomega^2_\beta\geq \frac{\left(|\cos(\alpha-\beta)|+1\right)^2}{4},
\end{align}
where $\tee_0=\int_{\mathbb{R}}\int_{\mathbb{R}}\tee\left|(\mathcal{G}^\alpha_gf)(\oomega,\tee)\right|^2d\tee d\oomega$ and $\oomega_0=\int_{\mathbb{R}}\int_{\mathbb{R}}\oomega\left|(\mathcal{G}^\alpha_gf)(\oomega,\tee)\right|^2d\tee d\oomega$ are expected values of the variables $\tee$ and $\oomega$ of the spectrogram $\left|(\mathcal{G}^\alpha_gf)(\oomega,\tee)\right|^2$ respectively. $\Delta\tee^2_\alpha=\int_{\mathbb{R}}\int_{\mathbb{R}}(\tee-\tee_0)^2\left|(\mathcal{G}^\alpha_gf)(\oomega,\tee)\right|^2d\tee d\oomega$ and $\Delta\oomega^2_\alpha=\int_{\mathbb{R}}\int_{\mathbb{R}}(\oomega-\oomega_0)^2\left|(\mathcal{G}^\alpha_gf)(\oomega,\tee)\right|^2d\tee d\oomega$ are the spread for the variables $\tee$ and $\oomega$ respectively. It can be seen from \eqref{R2009_Guanlei_Xiaotong_Xiaogang_SP_STFrFT_Heisenberg-UP} the UP is also connected with the parameters $\alpha$ and $\beta.$ In particular for $\alpha=\beta,$ the UP reads $\Delta\tee^2_\alpha \Delta\oomega^2_\beta\geq 1.$ With the assumptions that $f,g\in L^2(\mathbb{R}^2)$ are real valued and $P_lf,~P_l\hat{f},~P_lg,~P_l\hat{g}\in L^2(\mathbb{R}^2),$ $P_l(\bldx)=\xx_l,~l=1,2,$ authors in \cite{kamalakkannan2021short}, obtained the Heisenberg UP for the short time coupled FrFT \eqref{R2021_Kamalakkannan-Roopkumar-Zayed_Coupled_STFT_Definition} in $\bldomega$ as 
\begin{align*}
\Delta\bldomega^2_{\alpha,\beta}\Delta\bldx^2_f\geq \cos^2\gamma\|g\|^2_{L^2(\mathbb{R}^2)}\left(\Delta\bldx^2_f\right)^2 +\sin^2(\gamma)\left(\|f\|^2_{L^2(\mathbb{R}^2)}+\left(|\bldxi^0_f|^2+|\bldxi^0_g|^2\right)\Delta\bldx^2_f+\frac{\Delta\bldx^2_f}{\Delta\bldx^2_g}\|g\|^2_{L^2(\mathbb{R}^2)}\right)\|g\|^2_{L^2(\mathbb{R}^2)}\|f\|^2_{L^2(\mathbb{R}^2)},
\end{align*}
where $\gamma=\frac{\alpha+\beta}{2},$
$\Delta\bldomega^2_{\alpha,\beta}=\int_{\mathbb{R}^2}\int_{\mathbb{R}^2}|\bldomega-\bldomega^{\alpha,\beta,0}|^2\left|(\mathcal{G}^{\alpha,\beta}_gf)(\bldomega,\bldt)\right|^2d\bldomega d\bldt$ with $\bldomega^{\alpha,\beta,0}=(\oomega_1^{\alpha,\beta,0},\oomega_2^{\alpha,\beta,0}),$
$$\oomega_l^{\alpha,\beta,0}=\frac{1}{\|\mathcal{G}^\alpha_gf\|^2_{L^2(\mathbb{R}^4)}}\int_{\mathbb{R}^2}\int_{\mathbb{R}^2}\oomega_l\left|\left(\mathcal{G}^\alpha_gf\right)(\bldomega,\bldt)\right|^2d\bldomega d\bldt,~l=1,2$$
and $\Delta \bldx^2_f,~\Delta \bldxi^2_f$ are given respectively by \eqref{R2021_Kamalakkannan-Roopkumar-Zayed_Coupled_FrFT_Spread_Time_Domain} and \eqref{R2021_Kamalakkannan-Roopkumar-Zayed_Coupled_FrFT_Spread_Frequency_Domain}. 
\section{Fractional wavelet transform (FrWT)}
The concept of FrWT was proposed by Mendlovic et al.\cite{mendlovic1997fractional}, where fractional FT is firstly used to derive the fractional spectrum of a signal, and the WT is then performed on the obtained fractional spectrum. Since the fractional spectrum derived by the FrFT only represents the FrFD-frequency over the entire duration of the signal, the FRWT defined in \cite{mendlovic1997fractional} actually fails in obtaining the information of the local property of the signal. To overcome this drawback, using the relation \eqref{R2000_Wilczok_DefinitionWaveletTransform_using_InnerProduct}, between the WT and the classical convolution, Shi et al. \cite{shi2012novel} defined a novel FrWT by replacing the classical convolution by the fractional convolution \eqref{RShiFrFTConvolution} as
\begin{align}\label{R2012_Shi-Na-Liu_FrWT_Definition}
(\mathcal{W}^\alpha_{\psi}f)(\A,\B)=&\left(f\star_\alpha\delta_{\A}\tilde{\psi}\right)(\B)\notag\\
=&\int_{\mathbb{R}}f(\tee)\overline{\psi_{\alpha,\A,\B}(\tee)}d\tee,
\end{align}
where 
\begin{align}\label{R2012_Shi-Na-Liu_FrWT_Daughter_Wavelets}
\psi_{\alpha,\A,\B}(\tee)=e^{-\frac{i}{2}(\tee^2-\B^2)\cot\alpha}\psi_{\A,\B}(\tee)
\end{align}
and $\psi_{\A,\B}(\tee)=\frac{1}{\sqrt{a}}\psi\left(\frac{\tee-\B}{\A}\right).$
Equation \eqref{R2012_Shi-Na-Liu_FrWT_Definition} can be expressed in terms of the FrFT of $f$ as 
\begin{align}\label{R2012_Shi-Na-Liu_FrWT_FrFourierDomain}
(\mathcal{W}^\alpha_{\psi}f)(\A,\B)=\int_{\mathbb{R}}\sqrt{2\pi \A}(\mathcal{F}^\alpha f)(\xxi)\overline{\hat{\psi}(\A \xxi\csc \alpha)}\mathcal{K}_{-\alpha}(\xxi,\B)d\xxi
\end{align}
It is to be noted that for $\alpha=1,$ the FrWT \eqref{R2012_Shi-Na-Liu_FrWT_Definition} reduces to the classical WT. It is a matter of simple verification that $\mathcal{W}^\alpha_{\psi}f$ is linear in $f$ and conjugate linear in $\psi.$ The following properties, called the covariance properties of the fraction WT, give the formula for the FrWT of the dilation by $c>0$ of the function $f$ with respect to $\psi$
\begin{align}\label{R2012_Shi-Na-Liu_FrWT_ScalingProperty}
[(\mathcal{W}^\alpha_{\psi}(f(c\tee))](\A,\B)=\frac{1}{\sqrt{c}}(\mathcal{W}^\beta_{\psi}f)\left(c\A,c\B\right),~\A\in\mathbb{R}^+,\B\in\mathbb{R}
\end{align}
where $\varphi=\cot^{-1}(c^2\cot\alpha)$ and $\beta=\varphi/(\frac{\pi}{2}).$ Bahri et al.\cite{bahri2017logarithmic} obtained the fractional convolution and the correlation results associated with the FrWT as 
$$(\mathcal{W}^\alpha_{\psi_1\star_\alpha\psi_2}(f\star_\alpha g))(\A,\B)=\frac{1}{\sqrt{\A}}\int_{\mathbb{R}}(\mathcal{W}^\alpha_{\psi_1}f)(\A,\B-u)(\mathcal{W}^{-\alpha}_{\psi_2}g)(\A,u)du$$
and
$$(\mathcal{W}^\alpha_{\psi_1\circ_\alpha\psi_2}(f\circ_\alpha g))(\A,\B)=\frac{1}{\sqrt{\A}}\int_{\mathbb{R}}\overline{(\mathcal{W}^{-\alpha}_{\psi_1}f)(\A,-\B+u)}(\mathcal{W}^{-\alpha}_{\psi_2}g)(\A,u)du,$$
respectively, where $\circ_\alpha$ is the fractional correlation given as $(f\circ_\alpha g)(x)=\int_{\mathbb{R}}e^{-\frac{i}{2}(x^2+y^2)\cot\alpha}\overline{f(y)}g(x+y)dy.$ Authors in \cite{shi2012novel} obtained the resolution of identity formula for the transform \eqref{R2012_Shi-Na-Liu_FrWT_Definition} as
\begin{align}\label{R2012_Shi-Na-Liu_FrWT_IPR}
\int_{\mathbb{R}}\int_{\mathbb{R}^+}(\mathcal{W}^\alpha_{\psi}f)(\A,\B)\overline{(\mathcal{W}^\alpha_{\psi}g)(\A,\B)}\frac{d\A d\B}{\A^2}=2\pi c_{\psi}\langle f,g\rangle_{L^2(\mathbb{R})},~f,g\in L^2(\mathbb{R})
\end{align}
where $c_\psi$ is a positive constant depending on wavelet  given as $c_\psi=\int_{\mathbb{R}}\frac{|\hat{\psi}(\xxi)|^2}{\xxi}d\xxi<\infty,$ where $\hat{\psi}(\xxi)$ denotes the FT of $\psi(\tee).$ This in particular, for $f=g,$ implies that $\|\mathcal{W}^\alpha_{\psi}f\|^2_{L^2(\mathbb{R}^+\times\mathbb{R},\frac{d\A d\B}{\A^2})}=2\pi c_{\psi}\|f\|^2_{L^2(\mathbb{R})},$ which further shows that the normalized WT $\frac{1}{\sqrt{2\pi c_{\psi}}}\mathcal{W}^\alpha_{\psi}$ is an isometry from $L^2(\mathbb{R})$ into $L^2(\mathbb{R}^+\times\mathbb{R},d\mu_{aff}),$ where $d\mu_{aff}=\frac{d\A d\B}{\A^2}$ denotes the so called affine measure. Furthermore, the corresponding inversion formula is given by 
\begin{align}\label{R2012_Shi-Na-Liu_FrWT_InversionFormula}
f(t)=\frac{1}{2\pi c_\psi}\int_{\mathbb{R}}\int_{\mathbb{R}^+}(\mathcal{W}^\alpha_{\psi}f)(\A,\B)\psi_{\alpha,\A,\B}(\tee)\frac{d\A d\B}{\A^2}
\end{align}

\noindent From equations \eqref{R2012_Shi-Na-Liu_FrWT_Definition} and \eqref{R2012_Shi-Na-Liu_FrWT_FrFourierDomain}, it follows that if the kernel $\psi_{\alpha,\A,\B}(\tee)$ is supported in the time domain, the $(\mathcal{W}^\alpha_{\psi}f)(\A,\B)$ is supported in the time domain. Thus $(\mathcal{W}^\alpha_{\psi}f)(\A,\B)$ contains information about $f(\tee)$ near $\B.$ Similarly, if $\overline{\hat{\psi}}(\A \xxi\csc\alpha)$ is bandpass, i.e., $\hat{\psi}(\xxi)$ satisfies the condition $c_{\psi}<\infty,$
then the multiplication $\hat{f}(\xxi)$ and $\overline{\hat{\psi}}(\A \xxi\csc\alpha)$ can provide fractional-domain local properties of $f(t)$. This implies that the FrWT is capable of providing the time- and fractional-domain information simultaneously, hence giving a time-fractional-frequency representation of the signal. To be specific, both $\psi(t)$ and its FT $\hat{\psi}(\xxi)$ must have sufficiently fast decay so that they can be used as window functions. Suppose that $\psi(t)$ and $\hat{\psi}(\xxi)$ are functions with finite centers $E_\psi$ and $E_{\hat{\psi}}$ and finite radii $\Delta_\psi$ and $\Delta_{\hat{\psi}}$ . Then, the center and radii of the time-domain window function $\psi_{\alpha,\A,\B}(\tee)$ of the FrWT are respectively given by $E[\psi_{\alpha,\A,\B}(\tee)]=\B+\A E_\psi$ and $\Delta[\psi_{\alpha,\A,\B}(\tee)]=\A\Delta_\psi,$ where $E[\cdot]$ and $\Delta[\cdot]$ denote the expectation and deviation operator respectively. Similarly, the center and radii of the frequency-domain window function $\hat{\psi}(\A \xxi\csc\alpha)$ of the WT can be easily derived as $E[\hat{\psi}(\A \xxi\csc\alpha)]=\frac{E_{\hat{\psi}}}{\A}$ and $\Delta[\hat{\psi}(\A \xxi\csc\alpha)]=\frac{\Delta_{\hat{\psi}}}{\A}.$ The $Q-$ factor (or the ratio between the width and the center) of the fractional-domain window function of the FrWT is given by $Q=\frac{\Delta_{\hat{\psi}}}{E_{\hat{\psi}}},$ which is independent of the order $\alpha$ and the scaling parameter $\A.$ This is the constant$-Q$ property of the FrWT. Thus, the FrWT gives the local information of the signal $f(\tee)$ in the time window $[\B+\A E_{\psi}-\A\Delta_\psi,\B+\A E_{\psi}+\A\Delta_\psi]$ and the local FrFD spectrum information in the fractional frequency window $\left[\frac{E_{\hat{\psi}}}{\A}-\frac{\Delta_{\hat{\psi}}}{\A},\frac{E_{\hat{\psi}}}{\A}+\frac{\Delta_{\hat{\psi}}}{\A}\right]\sin\alpha,$ resulting into the time-fractional-frequency window
\begin{align}\label{R2012_Shi-Na-Liu_FrWT_time-fractional-frequency-window}
[\B+\A E_{\psi}-\A\Delta_\psi,\B+\A E_{\psi}+\A\Delta_\psi]\times \left[\frac{E_{\hat{\psi}}}{\A}-\frac{\Delta_{\hat{\psi}}}{\A},\frac{E_{\hat{\psi}}}{\A}+\frac{\Delta_{\hat{\psi}}}{\A}\right]\cdot\sin\alpha
\end{align}
with constant area $4\Delta_\psi\Delta_{\hat{\psi}}\sin\alpha$ in the time-FrFD-frequency plane. The area depends only on the mother wavelet $\psi(\tee)$ and the parameter $\alpha$ but is independent of the parameters $\A$ and $\B.$ However, for a given $\alpha,$ the time-fractional-frequency window varies with the scaling parameter $\A.$ More precisely, the window becomes narrower for detecting high FrFD-frequency  phenomena, i.e., small values of $\A$, and wider for detecting low FrFD-frequency  phenomena, i.e., large $\A.$

In view of \eqref{R2012_Shi-Na-Liu_FrWT_FrFourierDomain}, it follows the FrWT \eqref{R2012_Shi-Na-Liu_FrWT_Definition} analyzes the signal in time-frequency-FrFD domain, its physical meaning requires deeper interpretation. So, Prasad et al.\cite{prasad2014generalized} expressed the FrWT \eqref{R2012_Shi-Na-Liu_FrWT_Definition} in the fractional Fourier domain as 
\begin{align}\label{R2014_Prasad-Manna-Mahato-Singh_WT_in_FrFD}
(\mathcal{W}^\alpha_{\psi}f)(\A,\B)=\sqrt{\A}\int_{\mathbb{R}}e^{-\frac{i}{2}(\B^2+\xxi^2)\cot\alpha+i\B\xxi\csc\alpha+\frac{i}{2}\A^2\xxi^2\cot\alpha}(\mathcal{F}^\alpha f)(\xxi)\overline{\mathcal{F}^\alpha\left[e^{-\frac{i}{2}(\cdot)^2\cot\alpha}\psi\right](\A\xxi)}\mathcal{K}_{-\alpha}(\xxi,\B)d\xxi
\end{align}
as a result, the signal analysis only involves the time-FrFD domain. Moreover, in that case, the wavelet satisfies the following admissibility condition
\begin{align}\label{R2014_Prasad-Manna-Mahato-Singh_Wavelet_Admissibility_Condition}
c_{\psi,\alpha}=\int_{\mathbb{R}^+}\frac{|\mathcal{F}^\alpha\left(e^{-\frac{i}{2}(\cdot)\cot\alpha}\psi\right)(\xxi)|^2}{|\xxi|}d\xxi<\infty.
\end{align}
Authors in \cite{prasad2014generalized}, generalized the inner product relation \eqref{R2012_Shi-Na-Liu_FrWT_IPR} in terms of two wavelets $\psi_1,\psi_2\in L^2(\mathbb{R})$ as 
\begin{align}\label{R2014_Prasad-Manna-Mahato-Singh_2Wavelet_IPR}
\int_{\mathbb{R}}\int_{\mathbb{R}^+}(\mathcal{W}^\alpha_{\psi_1}f)(\A,\B)\overline{(\mathcal{W}^\alpha_{\psi_2}g)(\A,\B)}\frac{d\A d\B}{\A^2}=2\pi \sin\alpha c_{\psi_1,\psi_2,\alpha}\langle f,g\rangle_{L^2(\mathbb{R})},~f,g\in L^2(\mathbb{R}),
\end{align}
where $\psi_1$ and $\psi_2$ are such that $C_{\psi_1,\psi_2,\alpha}=\int_{\mathbb{R}^+}\overline{\mathcal{F}^\alpha\left[e^{-\frac{i}{2}(\cdot)^2\cot{\alpha}}\psi_1\right](\A)}\mathcal{F}^\alpha\left[e^{-\frac{i}{2}(\cdot)^2\cot{\alpha}}\psi_2\right](\A)\frac{d\A}{\A}<\infty.$ From \eqref{R2014_Prasad-Manna-Mahato-Singh_2Wavelet_IPR} it follows that, if $\psi_1=\psi_2=\psi$ then we have $\int_{\mathbb{R}}\int_{\mathbb{R}^+}(\mathcal{W}^\alpha_{\psi}f)(\A,\B)\overline{(\mathcal{W}^\alpha_{\psi}g)(\A,\B)}\frac{d\A d\B}{\A^2}=2\pi \sin\alpha c_{\psi,\alpha}\langle f,g\rangle_{L^2(\mathbb{R})}.$ Furthermore, if $f=g$ then $\int_{\mathbb{R}}\int_{\mathbb{R}^+}\left|(\mathcal{W}^\alpha_{\psi}f)(\A,\B)\right|^2\frac{d\A d\B}{\A^2}=2\pi\sin\alpha c_{\psi,\alpha}\|f\|_{L^2(\mathbb{R})}^2.$ For fixed $\A$ treating $(\mathcal{W}^\alpha_{\psi}f)(\A,\B)$ as a function of the variable $\B,$ authors also obtained the following
$$\int_{\mathbb{R}}\left[(\mathcal{W}^\alpha_{\psi}f)(\A,\B)\overline{(\mathcal{W}^\alpha_{\psi}g)(\A,\B)}\right]d\B=2\pi \A\sin\alpha\langle U_\alpha,V_\alpha\rangle_{L^2(\mathbb{R})},$$ 
where $U_\alpha(\xxi)=(\mathcal{F}^\alpha f)(\xxi)e^{\frac{i}{2}\A^2\xxi^2\cot\alpha}\overline{\mathcal{F}^\alpha\left[e^{-\frac{i}{2}(\cdot)^2\cot\alpha}\psi\right](\A\xxi)}$ and 
$V_\alpha(\xxi)=(\mathcal{F}^\alpha g)(\xxi)e^{-\frac{i}{2}\A^2\xxi^2\cot\alpha}\overline{\mathcal{F}^\alpha\left[e^{-\frac{i}{2}(\cdot)^2\cot\alpha}\psi\right](\A\xxi)}.$ Motivated by the convolution operator \eqref{RConvolutionZayedFrFT} given by Zayed, Dai et al. \cite{dai2017new} introduced a new definition of the fractional WT
\begin{align}\label{R2017_Dai-Zheng-Wang_FrWT_Definition}
\left(\mathcal{W}_\psi^{\alpha}f\right)(\A,\B)&=e^{-\frac{i}{2}\B^2\cot\alpha}\left\{\left[f(\tee)e^{\frac{i}{2}\tee^2\cot\alpha}\right]\star\overline{\left[\frac{1}{\sqrt{\A}}\psi\left(-\frac{\tee}{\A}\right)e^{\frac{i}{2}\left(\frac{\tee}{\A}\right)^2\cot\alpha}\right]}\right\}(\B)\notag\\
&=\int_{\mathbb{R}}f(\tee)
\overline{\psi_{\alpha,\A,\B}(\tee)}d\tee
\end{align}
where $\psi_{\alpha,\A,\B}(\tee)=e^{-\frac{i}{2}\left(\tee^2-\B^2-\left(\frac{\tee-\B}{\A}\right)\right)\cot\alpha}\psi_{\A,\B}(\tee).$ Equation \eqref{R2017_Dai-Zheng-Wang_FrWT_Definition}, can be expressed in the fractional Fourier domain as
\begin{align}
(\mathcal{W}^\alpha_{\psi}f)(\A,\B)=\sqrt{\frac{2\pi a}{1+j\cot\alpha}}\int_{\mathbb{R}}e^{-\frac{i}{2}\A^2\xxi^2\cot\alpha}(\mathcal{F}^\alpha f)(\xxi)\overline{(\mathcal{F}^\alpha\psi)(\A \xxi)}\mathcal{K}_{-\alpha}(\xxi,\B)d\xxi.
\end{align}
As a result of which, the authors in \cite{dai2017new} obtained the following inner product relation
\begin{align}\label{R2017_Dai-Zheng-Wang_FrWT_IPR}
\int_{\mathbb{R}^+}\int_{\mathbb{R}}(\mathcal{W}^\alpha_{\psi}f)(\A,\B)\overline{(\mathcal{W}^\alpha_{\psi}f)(\A,\B)}\frac{d\A}{\A}d\B=2\pi\sin\alpha c_\psi\langle f,g\rangle 
\end{align}
followed by the inversion formula
\begin{align}\label{R2017_Dai-Zheng-Wang_FrWT_InversionFormula}
f(\tee)=\frac{1}{2\pi\sin\alpha c_\psi}\int_{\mathbb{R}}\int_{\mathbb{R}}\frac{1}{\A^2}(\mathcal{W}^\alpha_{\psi}f)(\A,\B)\psi_{\A,\B,\alpha}(\tee)d\A d\B,
\end{align}
where $c_\psi=\int_{\mathbb{R}^+}\frac{|\mathcal{F}^\alpha\psi(s)|^2}{|s|}ds$ is a positive constant. From equation \eqref{R2017_Dai-Zheng-Wang_FrWT_IPR}, the FrWT satisfies the Parseval identity $\int_{\mathbb{R}}|f(\tee)|^2d\tee=\frac{1}{2\pi\sin\alpha C_\psi}\int_{\mathbb{R}^+}\int_{\mathbb{R}}|(\mathcal{W}^\alpha_{\psi}f)(\A,\B)|^2\frac{d\A}{\A^2}d\B.$ As a consequence of which it can be shown that the range $\mathcal{W}^\alpha_{\psi}(L^2(\mathbb{R}))$ of the FrWT given by \eqref{R2017_Dai-Zheng-Wang_FrWT_Definition} is a RKHS with the kernel $K_{\phi_\alpha}(\A_0,\B_0;\A,\B)=\frac{1}{2\pi \sin\alpha C_\psi}\int_{\mathbb{R}}\psi_{\A,\B,\alpha}(\tee)\overline{\psi_{\A_0,\B_0,\alpha}(\tee)}d\tee.$ The FrWT \eqref{R2017_Dai-Zheng-Wang_FrWT_Definition} gives the local information of the signal $f(\tee)$ in the time window $[\B+\A E_{\psi}-\A\Delta_\psi,\B+\A E_{\psi}+\A\Delta_\psi]$ and the local FrFD spectral information in the FrFD as $[kE_{\mathcal{F}^\beta\psi}-k\Delta_{\mathcal{F}^\beta\psi},kE_{\mathcal{F}^\beta\psi}+k\Delta_{\mathcal{F}^\beta\psi}],$ where $k=\sqrt{\frac{1-i\cot\alpha}{\A^2-i\cot\alpha}},\beta=\arctan(\A^2\tan\alpha),$ with the constant window area $4\A k\Delta_\psi\Delta_{\mathcal{F}^\beta\psi}$ in the time-FrFD-frequency plane. The area depends only on the mother wavelet $\psi$ and the angle $\alpha$ and is independent of the parameters $\A$ and $\B.$ Guo et al. \cite{guo2021novel} defined a general FrWT, using some particular type of convolutions given in \eqref{RShiFrFTConvolutionII}, as
\begin{align}\label{R2021_Guo-Li-Yang_FrWT_Definition}
(\mathcal{W}^{\alpha,\beta}_{\psi}f)(\A,\B)=&f(\tee)\star_{\alpha,\alpha,\beta}\frac{1}{\sqrt{\A}}\overline{\psi\left(-\frac{\tee}{\A}\right)}\notag\\
=&e^{-\frac{i}{2}\B^2\cot\alpha}\left[f(\tee)e^{\frac{i}{2}\tee^2\cot\alpha}\star\frac{1}{\sqrt{\A}}\psi\left(-\frac{\tee}{\A}\right)e^{\frac{i}{2}\tee^2\cot\beta}\right](\B).
\end{align}
This definition of the FrWT has the equivalent form 
\begin{align}\label{R2021_Guo-Li-Yang_FrWT_Definition_Time_Frequency_Analysis_Form}
(\mathcal{W}^{\alpha,\beta}_{\psi}f)(\A,\B)=\int_{\mathbb{R}}f(\tee)\overline{\psi_{\A,\B,\alpha,\beta}(\tee)}d\tee,
\end{align}
where 
\begin{align}\label{R2021_Guo-Li-Yang_FrWT_Daughter_wavelet}
\psi_{\A,\B,\alpha,\beta}(\tee)=\frac{1}{\sqrt{\A}}\psi\left(\frac{\tee-\B}{\A}\right)e^{-\frac{i}{2}\left[(\tee^2-\B^2)\cot\alpha-\left(\frac{\tee-\B}{\A}\right)^2\cot\beta\right]}.
\end{align}
The window function $\psi_{\alpha,\beta,\A,\B}(\tee)$ of proposed FrWT has one more free parameter $\beta,$ and thus the proposed transform has more degree of freedom of parameter selection. The transform have the scaling property $(\mathcal{W}^{\alpha,\beta}_{\psi}f(c\tee))(\A,\B)=\frac{1}{\sqrt{c}}(\mathcal{W}^{\alpha',\beta}_{\psi}f)(c\A,c\B),$ where $c>0$ and $\alpha'$ satisfies $\cot\alpha'=\frac{\cot\alpha}{c^2}.$ Moreover, it reduces to the classical WT when $\alpha=\beta=\frac{\pi}{2}$ and to the FrWT \eqref{R2012_Shi-Na-Liu_FrWT_Definition} and \eqref{R2017_Dai-Zheng-Wang_FrWT_Definition} when $\beta=\frac{\pi}{2}$ and $\alpha=\beta$ respectively. Also, it can be expressed in the fractional Fourier domain with $\lambda=\frac{\csc\beta}{\csc\alpha}$ as
\begin{align}\label{R2021_Guo-Li-Yang_FrWT_Definition-FrFD}
(\mathcal{W}^{\alpha,\beta}_{\psi}f)(\A,\B)=\frac{\sqrt{\A}}{A_\beta}\int_{\mathbb{R}}e^{\frac{i}{2}\left(\frac{\A\xxi}{\lambda}\right)^2\cot\beta}\left(\mathcal{F}^\alpha f\right)(\xxi)\overline{\left(\mathcal{F}^\beta \psi\right)\left(\frac{\A \xxi}{\lambda}\right)}K_{-\alpha}(\xxi,\B)d\xxi.
\end{align}
In term of fractional domain, the window and signal are restricted to same fractional domain (i.e., $\mathcal{F}^\alpha\psi$ and $\mathcal{F}^\alpha f$) for the transform \eqref{R2017_Dai-Zheng-Wang_FrWT_Definition} defined by \cite{dai2017new}, whereas for the transform \eqref{R2021_Guo-Li-Yang_FrWT_Definition}, it appears from \eqref{R2021_Guo-Li-Yang_FrWT_Definition-FrFD}, that the window and signal can be transformed freely to different fractional domains (i.e., $\mathcal{F}^\beta\psi$ and $\mathcal{F}^\alpha f$). The following resolution of the identity formula holds for the transform \eqref{R2017_Dai-Zheng-Wang_FrWT_Definition}
\begin{align}\label{R2021_Guo-Li-Yang_FrWT_IPR}
\int_{\mathbb{R}}\int_{\mathbb{R}}(\mathcal{W}^{\alpha,\beta}_{\psi}f)(\A,\B)\overline{(\mathcal{W}^{\alpha,\beta}_{\psi}g)(\A,\B)}\frac{1}{\A^2}d\A d\B=2\pi\sin\beta c_{\psi,\beta}\int_{\mathbb{R}}f(\tee)\overline{g(\tee)}d\tee,
\end{align}
followed by the inversion formula
\begin{align}\label{R2021_Guo-Li-Yang_FrWT_InversionFormula}
f(\tee)=\frac{1}{2\pi\sin\beta c_{\psi,\beta}}\int_{\mathbb{R}}\int_{\mathbb{R}}(\mathcal{W}^{\alpha,\beta}_{\psi}f)(\A,\B)\psi_{\A,\B,\alpha,\beta}(\tee)\frac{1}{\A^2}d\A d\B.
\end{align}
From \eqref{R2021_Guo-Li-Yang_FrWT_IPR}, it follows that $\int_{\mathbb{R}}\int_{\mathbb{R}}|(\mathcal{W}^{\alpha,\beta}_{\psi}f)(\A,\B)|^2\frac{d\A d\B}{\A^2}=2\pi\sin\beta c_{\psi,\beta}\int_{\mathbb{R}}|f(\tee)|^2d\tee,$ as a result of which it follows that the range $\mathcal{W}^{\alpha,\beta}_{\psi}(L^2(\mathbb{R}))$ of the FrWT given by \eqref{R2021_Guo-Li-Yang_FrWT_Definition} is a RKHS with the kernel given by $K_{\psi,\alpha,\beta}(\A_0,\B_0,\A,\B)=\frac{1}{2\pi\sin\beta c_{\psi,\beta}}\int_{\mathbb{R}}\psi_{\A,\B,\alpha,\beta}(\tee)\overline{\psi_{\A_0,\B_0,\alpha,\beta}(\tee)}d\tee.$ In view of equations \eqref{R2021_Guo-Li-Yang_FrWT_Definition_Time_Frequency_Analysis_Form} and \eqref{R2021_Guo-Li-Yang_FrWT_Definition-FrFD}, the FrWT \eqref{R2021_Guo-Li-Yang_FrWT_Definition} gives the local information of the signal $f(\tee)$ in the time window $[\B+\A E_{\psi}-\A\Delta_\psi,\B+\A E_{\psi}+\A\Delta_\psi]$ and the local FrFD spectral information in the FrFD as $\left[\frac{\lambda}{\A}E[\mathcal{F}^\beta\psi]-\frac{\lambda}{\A}\Delta\mathcal{F}^\beta\psi,\frac{\lambda}{\A}E[\mathcal{F}^\beta\psi]+\frac{\lambda}{\A}\Delta\mathcal{F}^\beta\psi\right],$ with the constant window area $4\lambda\Delta\psi\Delta\mathcal{F}^\beta\psi$ in the time-FrFD-frequency plane, where $\lambda=\frac{\csc\beta}{\csc\alpha}.$ The area depends only on the mother wavelet $\psi$ and the parameter $\alpha,\beta$ and is independent of the parameters $\A$ and $\B.$ Motivated by the work of Pandey et al. (\cite{pandey2015continuous}, \cite{pandey2016continuous}), Verma et al. \cite{verma2022note} extended the FrWT \eqref{R2012_Shi-Na-Liu_FrWT_Definition} to $n-$dimensional Euclidean space with the dilation parameter $\blda=(\A_1,\A_2,\cdots,\A_n)\in\mathbb{R}^n$ such that $\A_i\neq 0,i=1,2,\cdots,n,$ given as
\begin{align}\label{R2022_Verma-Gupta_nD_FrWT_Definition}
(\mathcal{W}^\alpha_{\psi}f)(\blda,\bldb)=\int_{\mathbb{R}}f(\bldt)\overline{\psi_{\blda,\bldb,\alpha}(\bldt)}d\bldt,
\end{align}
where $f,\psi\in L^2(\mathbb{R}^n),$ $\psi_{\blda,\bldb,\alpha}(\bldt)=\frac{1}{\sqrt{|\blda|_{p}}}\psi\left(\frac{\bldt-\bldb}{\blda}\right)e^{-\frac{i}{2}(|\bldt|^2-|\bldb|^2)\cot\alpha}$ and $\psi$ satisfies the admissibility condition
$\ C_{\psi ,\alpha}:=\int_{\mathbb{R}_{0}^n}\left|\left(\mathcal{F}_\alpha\left\{e^{-\frac{i}{2}|\cdot|^2\cot\alpha}\psi(\cdot)\right\}\right)(\bldu)\right|^2\frac{1}{|\bldu|_{p}} d\bldu<\infty.$ They considered two functions $\phi,\psi\in L^2(\mathbb{R})$ satisfying this admissibility condition with the assumption that
\begin{equation}
\int_{\mathbb{R}_{0}^n}\left|\left(\mathcal{F}_\alpha\left\{e^{-\frac{i}{2}|\cdot|^2\cot\alpha}\phi(\cdot)\right\}\right)(\bldxi)\right|\left|\left(\mathcal{F}_\alpha\left\{e^{-\frac{i}{2}|\cdot|^2\cot\alpha}\psi(\cdot)\right\}\right)(\bldxi)\right|\frac{1}{|\bldxi|_{p}}d\bldxi<\infty
\end{equation}
and obtained the following inner product relation 
\begin{align}\label{R2022_Verma-Gupta_nDFrWT_IPR2Wavlet}
\int_{\mathbb{R}_{0}^n}\int_{\mathbb{R}^n}(W^\alpha_\phi f)(\blda,\bldb)\overline{(W^\alpha_\psi g)(\blda,\bldb)}\frac{d\blda d\bldb}{|\blda|_{p}^2}=\frac{C_{\phi,\psi,\alpha}}{|C_{\alpha}|^2}\langle f,g\rangle_{L^2(\mathbb{R}^n)},
\end{align}
followed by the inversion formula 
$$f(\bldt)=\frac{|C_{\alpha}|^2}{C_{\phi,\psi,\alpha}}\int_{\mathbb{R}_{0}^n}\int_{\mathbb{R}^n}\phi_{\blda,\bldb,\alpha}(\bldt)(W^\alpha_\psi f)(\blda,\bldb)\frac{d\blda d\bldb}{|\blda|_{p}^2},$$
where $C_{\phi,\psi,\alpha}=\int_{\mathbb{R}_{0}^n}\overline{\left(\mathcal{F}_\alpha\left\{e^{-\frac{i}{2}|\cdot|^2\cot\alpha}\phi(\cdot)\right\}\right)(\bldxi)}\left(\mathcal{F}_\alpha\left\{e^{-\frac{i}{2}|\cdot|^2\cot\alpha}\psi(\cdot)\right\}\right)(\bldxi)\frac{1}{|\bldxi|_{p}}d\bldxi<\infty.$ From \eqref{R2022_Verma-Gupta_nDFrWT_IPR2Wavlet}, it follows that, if $f=g$ and $\phi=\psi$ then the range $W^\alpha_\psi(L^2(\mathbb{R}^n))$ is not only the subspace of $L^2(\mathbb{R}_{0}^n\times\mathbb{R}^n,\frac{d\blda d\bldb}{|\blda|^2_{p}})$ but also the RKHS with the kernel given by $K_{\psi,\alpha}(\blda_{0},\bldb_{0};\blda,\bldb)=\frac{|C_{\alpha}|^2}{C_{\psi,\alpha}}\int_{\mathbb{R}^n}\psi_{\blda,\bldb,\alpha}(\bldt)\overline{{\psi_{\blda_{0},\bldb_{0},\alpha}}(\bldt)}d\bldt.$ Kaur et al.\cite{kaur2022multidimensional}, introduced a more general definition of  FrWT based on the convolution associated with the multi-dimensional FrFT \eqref{R2020_Kamalakkannan-Roopkumar_ITSF_nD-FrFT_Convolution} introduced in \cite{kamalakkannan2020multidimensional}, which is defined for the $f\in L^2(\mathbb{R}^n)$ and the function $\psi\in L^2(\mathbb{R}^n)$ satisfying the admissibility condition $c_{\boldsymbol\alpha,\lambda}=(2\pi)^n\int_{\mathbb{R}^n_0}|(\mathfrak{F}_{\boldsymbol\alpha_\lambda}\psi)(\bldxi)|^2\frac{d\bldxi}{|\bldxi|_m}<\infty,$ where $\mathfrak{F}_{\boldsymbol\alpha_\lambda}$ is given by \eqref{R2020-Kamalakkanan-Roopkumar_ITSF_nD-FrFT_nalpha-Definition}, as 
\begin{align}\label{R2022_Kaur-Gupta-Verma_nD_FrWT_Definition_Convolution}
(W_{\psi}^{\boldsymbol\alpha,\lambda}f)(\blda,\bldb)=c(\boldsymbol\alpha_\lambda)e_{\boldsymbol\alpha,-\lambda^2}(\tilde{f} \ast \breve{g})(\bldb), \quad \blda\in \mathbb{R}^{n}_0,\quad \bldb\in\mathbb{R}^{n},
\end{align}
where $e_{\boldsymbol\alpha,\lambda^2}(\bldx)=e^{i\lambda^2 \sum_{k=1}^{n}a(\alpha_{k}){x_{k}}^2},$ $a(\alpha_k)=\frac{\cot(\alpha_k)}{2},$ $\tilde{\psi}(\bldx)=\psi(\bldx)e_{\boldsymbol\alpha,\lambda^2}(\bldx),$ $\tilde{f}(\bldx)=f(\bldx)e_{\boldsymbol\alpha,\lambda^2}(\bldx),$ ${\breve{g}}(\bldt)=\overline{\frac{1}{\sqrt{|\blda|_m}}\tilde{\psi}\left(\frac{-\bldt}{\blda}\right)},$ $c(\boldsymbol\alpha_\lambda)=c(\alpha_1)c(\alpha_2) \cdots c(\alpha_{n}),$ $ \boldsymbol\alpha_\lambda = \left(\cot^{-1}(\lambda^2\cot\alpha_1),\cot^{-1}(\lambda^2\cot\alpha_2),\cdots,\cot^{-1}(\lambda^2\cot\alpha_{n})\right)$ and $\ast$ is the convolution given by
$$(\mu\ast \nu)(\mathbf{x})=\int_{\mathbb{R}^{n}}\mu(\bldy)\nu(\bldx-\bldy)d\bldy,$$ 
which can equivalently be written as $(W_{\psi}^{\boldsymbol\alpha,\lambda}f)(\blda,\bldb)=\left\langle f(\bldt),\mathbf{\Psi}_{\boldsymbol\alpha,\lambda,\blda,\bldb}(\bldt)\right\rangle,$ where 
$$\mathbf{\Psi}_{\boldsymbol\alpha,\lambda,\blda,\bldb}(\bldt)=\overline{c(\boldsymbol\alpha_\lambda)}e_{\boldsymbol\alpha,\lambda^2}(\bldb)e_{\boldsymbol\alpha,-\lambda^2}(\bldt){\frac{1}{\sqrt{|\blda|_m}}\tilde{\psi}}\left(\frac{\bldt-\bldb}{\blda}\right).$$ 
They  proved basic properties like linearity, anti-linearity, dilation, parity, conjugation, and translation of the FrWT \eqref{R2022_Kaur-Gupta-Verma_nD_FrWT_Definition_Convolution}, followed by the inner product relation
\begin{align}
\int_{\mathbb{R}^n_0}\int_{\mathbb{R}^n}(W_{\psi}^{\boldsymbol\alpha,\lambda}f)(\blda,\bldb)\overline{(W_{\psi}^{\boldsymbol\alpha,\lambda}g)(\blda,\bldb)}\frac{d\blda d\bldb}{|\blda|_m}=c_{\boldsymbol\alpha,\lambda}\langle f,g\rangle_{L^2(\mathbb{R}^n)},
\end{align}
inversion formula $f(\bldt)=\frac{1}{c_{\boldsymbol\alpha,\lambda}}\int_{\mathbb{R}^n_0}\int_{\mathbb{R}^n}(W_{\psi}^{\boldsymbol\alpha,\lambda}f)(\blda,\bldb)\mathbf{\Psi}_{\boldsymbol\alpha,\lambda,\blda,\bldb}(\bldt)\frac{d\blda d\bldb}{|\blda|_m}.$ In addition, authors characterized the range $W_{\psi}^{\boldsymbol\alpha,\lambda}(L^2(\mathbb{R}^n))$ and proved that it is a RKHS with the kernel $K_{\boldsymbol\alpha,\lambda,\psi}(\blda_0,\bldb_0;\blda,\bldb)=\int_{\mathbb{R}^n}\mathbf{\Psi}_{\boldsymbol\alpha,\lambda,\blda,\bldb}(\bldt)\overline{\mathbf{\Psi}_{\boldsymbol\alpha,\lambda,\blda_0,\bldb_0}(\bldt)}d\bldt.$  Definition \eqref{R2022_Kaur-Gupta-Verma_nD_FrWT_Definition_Convolution} not only combine the advantages of the multidimensional FrFT but also inherit the characteristics of the classical WT in $n-$dimension studied in \cite{pandey2015continuous}. It reduces to \eqref{R2017_Dai-Zheng-Wang_FrWT_Definition} with a constant multiple $\sqrt{1-i\cot\alpha},$ when $n=1,$ $\lambda=1$ and $\alpha_1=\alpha$ and to the $n-$dimensional WT introduced in \cite{pandey2015continuous}, when $\alpha_i=\frac{\pi}{2},i=1,2,\cdots,n$ and $\lambda=1.$
\\

Bahri et al.\cite{bahri2017logarithmic} studied the FrWT \eqref{R2012_Shi-Na-Liu_FrWT_Definition} on the generalized Sobolev spaces $H^s_{r,\alpha},~s,\alpha\in\mathbb{R}$ associated with the FrFT which is defined as the space of all measurable functions $f$ such that 
$$\|f\|^2_{H^s_{r,\alpha}}=\int_{\mathbb{R}}(1+|\xxi\csc\alpha|)^{-s}|(\xxi\csc\alpha)^r|^2|(\mathcal{F}^\alpha f)(\xxi)|^2d\xxi<\infty.$$ 
They obtained the boundedness of the operator $W_\psi^\alpha:H^s_{r,\alpha}\rightarrow H^s_{r,\alpha}$ given by $f\rightarrow (W_\psi^\alpha f)(\A,\cdot)$ with the assumption that $\psi$ is integrable and proved that $\int_{\mathbb{R}^+}\|(\mathcal{W}^\alpha_{\psi}f)(\A,\cdot)\|^2_{H^s_{r,\alpha}}\frac{d\A}{\A^2}=2\pi\sin\alpha c_{\psi,\alpha}\|f\|^2_{H^s_{r,\alpha}}.$ The result analogous to \eqref{RRieder_WT_Argument-Dependence-Property} for the FrWT with the norm $\|\cdot\|_{H^s_{r,\alpha}}$ as a counterpart of $\|\cdot\|_{H^\alpha(\mathbb{R})}$ is also obtained. They  derived the convergence for the convolution for $\mathcal{W}^\alpha_{\psi},$ i.e., $(\psi_\A\star_\alpha f)(\cdot)\rightarrow f(\cdot)$ in $H^s_{r,\alpha}$ as $\A\rightarrow 0,$ where $\psi\in L^1(\mathbb{R})$ is a real valued function satisfying $\int_{\mathbb{R}}\psi(\tee)d\tee=1$ and $\psi_\A(x)=\frac{1}{\A}\psi(\frac{x}{\A}),\A\in \mathbb{R}^+.$ Prasad et al.\cite{prasad2016continuous} studied the similar result for the generalized Sobolev space $B^{\eta,\alpha}_{p,k}(\mathbb{R}),$ which is defined as 
%$L^p\left(\left(\mathbb{R}_+,\frac{d\A}{\A^2}\right),B^{\xi,\alpha}_{p.k}(\mathbb{R})\right),~1\leq p<\infty$ abbreviated by $W^p_{k,\alpha}$ consists of all measurable  functions $f$ on $\mathbb{R}\times\mathbb{R}_+$ such that 
%$$\|f\|_{W^p_{k,\alpha}}=\left(\int_{\mathbb{R}_+}\|f(\cdot,\A)\|^p_{p,k,\alpha}\frac{d\A}{\A^{\frac{p}{2}+1}}\right)<\infty$$
%$$c_{\psi,p,\alpha}=(2\pi\sin\alpha)^{\frac{p}{2}}\int_{\mathbb{R}^+}\left|\mathcal{F}^\alpha[e^{\frac{i}{2}(\cdot)^2\cot\alpha}\psi](\xxi)\right|^p\frac{d\xxi}{\xxi}<\infty.$$
\begin{definition}\label{RGeneralized_Ultra-distribution_LCDomain}
Let $\eta\in\mathfrak{M}_c,k\in K_\eta,$ then the generalized Sobolev space $B^{\eta,\alpha}_{p,k}(\mathbb{R})$ is defined as the space of all the generalized ultra-distribution $f$ in $(S^\alpha_\eta)'$ such that 
$\|f\|_{p,k,\alpha}=\left(\int_{\mathbb{R}}|k(\xxi)(\mathcal{F}^\alpha f)(\xxi)|^pd\xxi\right)^{\frac{1}{p}}<\infty$ and $\|f\|_{\infty,k,\alpha}=ess\sup|k(\xxi)(\mathcal{F}^\alpha f)(\xxi)|,$ where $K_\eta$ the set of positive functions $k$ in $\mathbb{R}$ satisfying $k(\xxi_1+\xxi_2)\leq e^{\lambda\eta(-\xxi)}k(\xxi_2),$ for some $\lambda>0.$ The space $B^{\eta,\alpha}_{p,k}(\mathbb{R})$ reduces to the  H\"ormander space $B_{p,k}(\mathbb{R}),$ when $\eta(\xxi)=\log(1+|\xxi|)$ and $\alpha=\frac{\pi}{2}.$ 
\end{definition}
Here, $\mathfrak{M}_c$ denotes the set of all real valued continuous functions $\eta$ on $\mathbb{R}$ satisfying $\eta(0)=0,$ for all $\xxi_1,\xxi_2,\xxi,$ $0\leq \eta(\xxi_1+\xxi_2)\leq \eta(\xxi_1)+\eta(\xxi_2),$ $\int_{\mathbb{R}}\frac{\eta(\xxi)}{(1+|\xxi|)^2}d\xxi<\infty$ and $\eta(\xxi)>b+a\log (1+|\xxi|),~a>0,~b\in\mathbb{R}$ and $\eta(\xxi)=\Omega(|\xxi|),$ where $\Omega$ is a concave function on $[0,\infty).$ 
%Let $\mathfrak{M}$ denotes the class of all real valued continuous functions $\eta$ on $\mathbb{R}$ such that $\eta(0)=0,$ and for all $\xxi_1,\xxi_2,\xxi,$ $0\leq \eta(\xxi_1+\xxi_2)\leq \eta(\xxi_1)+\eta(\xxi_2),$ $\int_{\mathbb{R}}\frac{\eta(\xxi)}{(1+|\xxi|)^2}d\xxi<\infty$ and $\eta(\xxi)>b+a\log (1+|\xxi|),~a>0,~b\in\mathbb{R}.$
%Also let  $\mathfrak{M}_c$ denotes the set of all $\eta$ in $\mathfrak{M}$ satisfying $\eta(\xxi)=\Omega(|\xxi|),$ where $\Omega$ is a concave function on $[0,\infty).$ 
The space $S^\alpha_\xxi$ is a collection of $\phi\in L^1(\mathbb{R})$ such that  $\phi$ and $\mathcal{F}^\alpha\phi$ are smooth functions and for every $l\in\mathbb{N}_0$ and $\lambda\geq 0,$ $p_{l,\lambda,\alpha}(\phi)=\sup_{t\in\mathbb{R}}e^{\lambda\eta(t)}|\Delta_t^\alpha\phi(t)|<\infty$ and $r_{l,\lambda,\alpha}(\phi)=\sup_{\xxi\in\mathbb{R}}e^{\lambda\eta(\xxi)}|\Delta_\xxi^\alpha\phi(\xxi)|<\infty,$ where $\Delta_t=\frac{d}{dt}-it\cot\alpha.$ The topology on $S^\alpha_\eta$ is defined by the semi norms $p_{l,\lambda,\alpha}(\phi)$ and $r_{l,\lambda,\alpha}(\phi).$ The elements of the dual $(S^\alpha_\eta)'$ are called the generalized ultra-distributions. Note that for $\eta(\xxi)=\log(1+|\xxi|),\alpha=\frac{\pi}{2}$ and $\eta(\xxi)=\xxi^\sigma(0<\sigma<1),\alpha=\frac{\pi}{2}$ the space $S^\alpha_\eta$ reduces to the Schwartz space and the Gevrey space of test functions, respectively.

Authors in \cite{verma2022note}, studied the properties of FrWT \eqref{R2022_Verma-Gupta_nD_FrWT_Definition}, with the assumption that $\psi$ is compactly supported, on Morrey space $L^{1,\nu}_{M}(\mathbb{R}^n),$ define as 
\begin{definition}
Let $1\leq p<\infty$ and $\nu\geq 0.$ The Morrey space, denoted by $L^{p,\nu}_{M}(\mathbb{R}^n),$ is defined as
$$ L^{p,\nu}_{M}(\mathbb{R}^n)=\left\{f\in L^p_{loc}(\mathbb{R}^n): \sup_{\bldx\in\mathbb{R}^n,r>0}\frac{1}{r^{\nu}}\int_{B(\bldx,r)}|f(\bldt)|^pd\bldt<\infty\right\},$$ which is a Banach space with respect to the norm $\|f\|_{L^{p,\nu}_{M}(\mathbb{R}^n)}=\sup_{\bldx\in\mathbb{R}^n,r>0}\left(\frac{1}{r^{\nu}}\int_{B(\bldx,r)}|f(\bldt)|^p\right)^{\frac{1}{p}}.$
\end{definition}
\noindent They proved  the boundedness of the operator $W_\psi^\alpha:L^{1,\nu}_{M}(\mathbb{R}^n)\rightarrow L^{1,\nu}_{M}(\mathbb{R}^n)$ given by $f\rightarrow (W_\psi^\alpha f)(\blda,\cdot)$ and 
%proved that 
%$\int_{\mathbb{R}^+}\|(\mathcal{W}^\alpha_{\psi}f)(\blda,\cdot)\|^2_{L^{1,\nu}_{M}(\mathbb{R}^n)}\frac{d\blda}{\blda^2}=2\pi\sin\alpha c_{\psi,\alpha}\|f\|^2_{L^{1,\nu}_{M}(\mathbb{R}^n)}.$ 
the result analogous to \eqref{RRieder_WT_Argument-Dependence-Property} for the FrWT with the norm $|\cdot|_{L^{1,\nu}_{M}(\mathbb{R}^n)}$ as a counterpart of $|\cdot|_{H^\alpha(\mathbb{R})}.$\\

\noindent\textbf{Uncertainty Principles:} Using the logarithmic UP for the FrFT given by equation \eqref{R2008_Guanlei-Xiaotong-Xiaogang_SP_FrFT_Logarithmic-UP}, obtained in \cite{guanlei2009logarithmic}, authors in \cite{bahri2017logarithmic} obtained the logarithmic UP for the FrWT \eqref{R2012_Shi-Na-Liu_FrWT_Definition} as 
\begin{align}
\int_{\mathbb{R}^+}\int_{\mathbb{R}}\ln|\xxi|\left|\left(\mathcal{F}^\beta_\B\{(\mathcal{W}^\alpha_\psi f)(\A,\B)\}\right)(\xxi)\right|^2\frac{d\xxi d\A}{\A^2}+2\pi\sin\alpha &c_{\psi,\alpha}\int_{\mathbb{R}}\ln|\xxi|\left|(\mathcal{F}^\alpha f)(\xxi)\right|^2d\xxi\notag\\
&\hspace{-2.2cm}\geq 2\pi\sin\alpha c_{\psi,\alpha}\left(\ln|\sin(\beta-\alpha)|+\frac{\Gamma'(\frac{1}{4})}{\Gamma(\frac{1}{4})}\right)\int_{\mathbb{R}}|f(\xx)|^2d\xx.
\end{align}
Using the Heisenberg UP  \eqref{R2022_Verma-Gupta_IJWMIP_nD-FrFT-Heisenberg_UP_2FrFDomain} in two fractional Fourier domain, Verma et al.\cite{verma2022note} obtained the the Heisenberg type UPs for the FrWT \eqref{R2022_Verma-Gupta_nD_FrWT_Definition} as
\begin{align}
 \left(\int_{\mathbb{R}^n}|\bldt|^2|(\mathcal{F}^\beta_{\bldb}\{(W^\alpha_\psi f)(\blda,\bldb)\})(\bldt)|^2d\bldt\frac{d\blda}{|\blda|^2_{p}}\right)\left(\int_{\mathbb{R}^n}|\bldxi|^2|(\mathcal{F}^\alpha f)(\bldxi)|^2d\bldxi\right)\geq\frac{n^2C_{\psi,\alpha}}{4|C_{\alpha}|^2}|\sin(\alpha-\beta)|^2\|f\|^4_{L^2(\mathbb{R}^n)}. 
 \end{align}
and also its local UP, using local UP \eqref{R2022_Verma-Gupta_IJWMIP_nD-FrFT-Local_UP_2FrFDomain} for the FrFT in two fractional Fourier domain, as
\begin{align*}
\displaystyle\int_E|(\mathcal{F}_{\alpha}f)(\bldxi)|^2d\bldxi
\leq
\begin{cases}
\frac{A_{\theta}{|C_{\alpha}|^2}}{|\sin(\alpha-\beta)|^{2\theta}C_{\psi,\alpha}}(\lambda(E))^\frac{2\theta}{n}\displaystyle\int_{\mathbb{R}_{0}^n}\int_{\mathbb{R}^n}|\bldx|^{2\theta}|(\mathcal{F}^\beta_{\bldb}\{(W^\alpha_\psi f)(\blda,\bldb)\})(\bldx)|^2d\bldx\frac{d\blda}{|\blda|^2_{p}},\ \mbox{if}\ 0<\theta<\frac{n}{2},\\
\frac{A_{\theta}|C_{\alpha}|^2}{|\sin(\alpha-\beta)|^{n}C_{\psi,\alpha}}\lambda(E)\displaystyle\int_{\mathbb{R}_{0}^n}\|(W^\alpha_\psi f)(\cdot,\blda)\|^{2-\frac{n}{\theta}}_{L^2(\mathbb{R}^n)}\left\||\bldx|^{\theta}(\mathcal{F}^\beta_{\bldb}\{(W^\alpha_\psi f)(\blda,\bldb)\})\right\|^{\frac{n}{\theta}}_{L^2(\mathbb{R}^n)}\frac{d\blda}{|\blda|^2_{p}},\ \mbox{if}\ \theta>\frac{n}{2}.
\end{cases}
\end{align*} 
Kaur et al.\cite{kaur2022multidimensional} obtained the same for the FrWT \eqref{R2022_Kaur-Gupta-Verma_nD_FrWT_Definition_Convolution}, using the corresponding UPs \eqref{R2022-Kaur-Gupta-Verma_nD-FrFT_nalpha-Heisenberg-UP} and \eqref{R2022-Kaur-Gupta-Verma_nD-FrFT_nalpha-Logarithmic-UP} for the multidimensional FrFT \cite{kamalakkannan2020multidimensional} as
\begin{align*}
\left(\int_{\mathbb{R}_0^n}\int_{\mathbb{R}^n}|\bldb|^2|(W_{\psi}^{\boldsymbol\alpha,\lambda}{f})(\blda,\bldb)|^2d\bldb\frac{d\blda}{|\blda|_m^2}\right)\left(\int_{\mathbb{R}^n} |\bldxi|^2|(\mathcal{F}^{\boldsymbol\alpha,\lambda}f)(\bldxi)|^2d\bldxi\right) \geq P_{\boldsymbol\alpha,\lambda}C_{\boldsymbol\alpha,\lambda}\frac{n^2}{4}\|f\|^2_{L^2(\mathbb{R}^n)},
\end{align*}
where $P_{\boldsymbol\alpha,\lambda}=\frac{|c(\boldsymbol\alpha_\lambda)|^2|\sin\boldsymbol\alpha|_m}{M^2(\lambda^2)^{n+2}},$ $M=\mbox{max}\left\{\frac{1}{|\sin\alpha_l|},l=1,2,\cdots,n\right\}$ and 
\begin{align*}
\int_{\mathbb{R}_0^n}\int_{\mathbb{R}^n}\ln|\bldb||(W_{\psi}^{\boldsymbol\alpha,\lambda}{f})(\blda,\bldb)|^2d\bldb\frac{d\blda}{|\blda|_m^2}+C_{\boldsymbol\alpha,\lambda}P'_{\boldsymbol\alpha,\lambda}\int_{\mathbb{R}^n}\ln|\bldxi||(\mathcal{F}^{\boldsymbol\alpha,\lambda}f)(\bldxi)|^2d\bldxi\geq (D-P'_{\boldsymbol\alpha,\lambda}\ln(\lambda^2M))C_{\boldsymbol\alpha,\lambda}\|f\|^2_{L^2(\mathbb{R}^n)},
\end{align*}
where $P'_{\boldsymbol\alpha,\lambda}=\frac{\lambda^{2n}}{|c(\boldsymbol\alpha)|^2|\sin\boldsymbol\alpha|_m}.$ They also obtained the local UP for the FrWT using the same for the multidimensional FrFT \eqref{R2022-Kaur-Gupta-Verma_nD-FrFT_nalpha-Local-UP}.

\subsection{Multiresolution Analysis (MRA)}
Multiresolution analysis (MRA) has become an important tool for the construction of a discrete wavelet system, i.e., an orthonormal basis for $L^2(\mathbb{R}).$ It was first formulated by Mallet \cite{mallat1989multiresolution}. Some of the important wavelets that have been obtained from the MRA are Haar wavelet, Shannon wavelet, Meyer wavelet, Daubechies wavelet, etc. \cite{daubechies1992ten},\cite{debnath2002wavelet}. Not all wavelets are derived via MRA; there are some wavelets, like the Mexican hat wavelet and Morlet wavelet, that do not come from an MRA (see\cite{pereyra2012harmonic}). 
\begin{definition}
The MRA consists of a sequence of closed subspaces $\{V_m\}_{m\in\mathbb{Z}}$ in $L^2(\mathbb{R})$ satisfying the properties
\begin{enumerate}
\item $V_m\subset V_{m+1},~\forall~m\in\mathbb{Z}$
\item $f(\tee)\in V_m$ if and only if $f(2\tee)\in V_{m+1},~\forall~m\in\mathbb{Z}$
\item $\cap_{m\in\mathbb{Z}}V_m=\{0\},$ $\overline{\cup_{m\in\mathbb{Z}}V_m}=L^2(\mathbb{R})$
\item there exists a function $\phi\in V_0$ such that $\{\phi_{0,n}(\tee):n\in\mathbb{Z}\}$ form an orthonormal basis of $V_0,$ where 
$$\phi_{m,n}(\tee)=2^{\frac{m}{2}}\phi(2^m\tee-n),~m,n\in\mathbb{Z}.$$
\end{enumerate}
\end{definition}\label{RMRA_Classical_Definition}

Shi et al. \cite{shi2015multiresolution} developed the MRA associated with the FrWT \eqref{R2012_Shi-Na-Liu_FrWT_Definition} by discretizing the parameters $\A$ and $\B$ by dyadic grids, i.e., $\A=\A_0^m$ and $\B=n\A_0^m\B_0,~m,n\in\mathbb{Z}$ to the family \eqref{R2012_Shi-Na-Liu_FrWT_Daughter_Wavelets}, results in the following, which contains one extra parameter $\alpha.$ As a result of which, a new class of orthonormal basis of $L^2(\mathbb{R})$ can be obtained
\begin{enumerate}[label=(\Roman*)]
\item \label{R2012_Shi-Na-Liu_FrWT_MRA_Con1} $V^\alpha_m\subset V^\alpha_{m+1},~\forall~m\in\mathbb{Z}$
\item \label{R2012_Shi-Na-Liu_FrWT_MRA_Con2}$f(\tee)\in V^\alpha_m$ if and only if $f(2\tee)e^{\frac{i}{2}((2\tee)^2-\tee^2)\cot\alpha}\in V^\alpha_{m+1},~\forall~m\in\mathbb{Z}$
\item \label{R2012_Shi-Na-Liu_FrWT_MRA_Con3}$\cap_{m\in\mathbb{Z}}V^\alpha_m=\{0\},$ $\overline{\cup_{m\in\mathbb{Z}}V^\alpha_m}=L^2(\mathbb{R})$
\item \label{R2012_Shi-Na-Liu_FrWT_MRA_Con4}there exists a function $\phi\in V^\alpha_0$ such that $\{\phi_{\alpha,0,n}(\tee):n\in\mathbb{Z}\}$ form an orthonormal basis of $V^\alpha_0,$ where 
$$\phi_{\alpha,0,n}(\tee)=2^{\frac{m}{2}}\phi(2^m\tee-n)e^{-i\frac{\tee^2-n^2}{2}\cot\alpha},~n\in\mathbb{Z}.$$
\end{enumerate}
Indeed condition \ref{R2012_Shi-Na-Liu_FrWT_MRA_Con4} can be relaxed by assuming that the set of functions  $\{\phi_{\alpha,0,n}(\tee):n\in\mathbb{Z}\}$ is a Riesz basis of $V^\alpha_0=\{\sum_{n\in\mathbb{Z}}c[n]\phi_{0,n,\alpha}(\tee):c[n]\in l^2(\mathbb{Z})\}.$ The fractional MRA so defined consists of one extra parameter $\alpha$ compared to the classical and reduces to the later for $\alpha=\frac{\pi}{2}.$ With the scaling function $\phi$ given by 
$$
(\mathcal{F}\phi)(\xxi\csc\alpha)=
\begin{cases}
1,&|\xxi|<\xxi_\alpha\\
\cos\left[\frac{\pi}{2}v\left(\frac{\xxi}{2\xxi_\alpha}-1\right)\right],&\xxi_\alpha\leq|\xxi|<2\xxi_\alpha\\
0,&\mbox{otherwise},
\end{cases}
$$
where $\xxi_\alpha=\frac{2\pi}{3|\csc\alpha|}$ and $v(\tee)$ is a smooth function satisfying $v(\tee)=
\begin{cases}
0,&\tee\leq 0\\
1,&\tee\geq 1,
\end{cases}
$
satisfying the additional property $v(\tee)+v(1-\tee)=1,$ authors obtained the wavelet $\psi$ by taking the inverse FT in
$$\left(\mathcal{F}\psi\right)(\xxi\csc\alpha)=\left[\left(\mathcal{F}\psi\right)(\xxi\csc\alpha+2\pi)+\left(\mathcal{F}\psi\right)(\xxi\csc\alpha-2\pi)\right]e^{-\frac{i}{2}\xxi\csc\alpha}\left(\mathcal{F}\psi\right)\left(\frac{\xxi\csc\alpha}{2}\right).$$
%$$(\mathcal{F}^\alpha\psi)(\xxi\csc\alpha)=\sqrt{3}e^{\frac{i}{2}\xxi\csc\alpha}\sin^2\left(\frac{\xxi\csc\alpha}{2}\right)\left(\mathcal{F}^\alpha\phi\right)\left(\frac{\xxi\csc\alpha}{2}\right)\left\{\frac{1+2\sin^2\left(\frac{\xxi\csc\alpha}{4}\right)}{\left[1+2\cos^2\left(\frac{\xxi\csc\alpha}{2}\right)\right]\left[1+2\cos^2\left(\frac{\xxi\csc\alpha}{2}\right)\right]}\right\}^\frac{1}{2}.$$
Also for the scaling function $\phi(\tee)=\chi_{[0,1)}(\tee),$ the corresponding $\psi$ is given by $\psi(\tee)=\chi_{[0,\frac{1}{2})}(\tee)-\chi_{[\frac{1}{2},1)}(\tee).$ Guo et al.\cite{guo2021novel} developed the fractional MRA associated with the family \eqref{R2021_Guo-Li-Yang_FrWT_Daughter_wavelet} by discretizing the dilation parameter $\A$ and the translation parameter $\B,$ with $\A=\A_0^m$ and $\B=n\A_0^m\B_0,$ where $m,n\in\mathbb{Z},$ $\A_0>0,\B_0\neq 0.$ Of particular interest, is the dyadic grid, which occur for $\A_0=2,\B_0=1.$ Authors considered a sequence of closed sub-spaces $\{V^{\alpha,\beta}_m\}_{m\in\mathbb{Z}}$ of $L^2(\mathbb{R})$ satisfies the following properties.  
\begin{enumerate}[label=(\roman*)]
\item \label{R2021_Guo-Li-Yang_FrWT_MRA_Con1} $V^{\alpha,\beta}_m\subset V^{\alpha,\beta}_{m+1},~\forall~m\in\mathbb{Z}$
\item \label{R2021_Guo-Li-Yang_FrWT_MRA_Con2}$f(\tee)\in V^{\alpha,\beta}_m$ if and only if $f(2\tee)e^{\frac{i}{2}((2\tee)^2-\tee^2)\cot\alpha}\in V^{\alpha,\beta}_{m+1},~\forall~m\in\mathbb{Z}$
\item \label{R2021_Guo-Li-Yang_FrWT_MRA_Con3}$\cap_{m\in\mathbb{Z}}V^{\alpha,\beta}_m=\{0\},$ $\overline{\cup_{m\in\mathbb{Z}}V^{\alpha,\beta}_m}=L^2(\mathbb{R})$
\item \label{R2021_Guo-Li-Yang_FrWT_MRA_Con4}there exists a function $\phi\in V^{\alpha,\beta}_0$ such that $\{\phi_{\alpha,\beta,0,n}(\tee):n\in\mathbb{Z}\}$ form an orthonormal basis of $V^{\alpha,\beta}_0,$ where 
$$\phi_{\alpha,\beta,m,n}(\tee)=2^{\frac{m}{2}}\phi(2^m\tee-n)e^{-\frac{i}{2}\left[\left(\tee^2-\left(\frac{n}{2^m}\right)^2\right)\cot\alpha-(2^m\tee-n)^2\cot\beta\right]},~m,n\in\mathbb{Z}.$$
\end{enumerate}
If the orthogonality in \ref{R2021_Guo-Li-Yang_FrWT_MRA_Con4} is relaxed by the assumption that $\{\phi_{\alpha,\beta,0,n}\}_{n\in\mathbb{Z}}$ is a Riesz basis of $V_0^{\alpha,\beta},$ the $\phi$ generates a generalized MRA associated with the FrWT. Authors obtained a result showing a way to construct an orthonormal basis from Riesz basis $\{\phi_{\alpha,\beta,0,n}\}_{n\in\mathbb{Z}}$ of $V_0^{\alpha,\beta}.$ The MRA so developed reduces to the MRA proposed by Dai et al.\cite{dai2017new} when $\alpha=\beta$ and to definition \ref{RMRA_Classical_Definition} if $\alpha=\beta=\frac{\pi}{2}.$  The subspace $W^{\alpha,\beta}_m$ is defined as the orthogonal complement of $V^{\alpha,\beta}_m$ in $V^{\alpha,\beta}_{m+1},$ i.e., 
$$W^{\alpha,\beta}_m\perp V^{\alpha,\beta}_m,~V^{\alpha,\beta}_{m+1}=V^{\alpha,\beta}_m\oplus V^{\alpha,\beta}_m,~\forall~m\in\mathbb{Z}.$$
Then according to the definition, $\{W^{\alpha,\beta}_m\}_{m\in\mathbb{Z}}$ have the following properties
\begin{enumerate}
\item $W^{\alpha,\beta}_k\perp V^{\alpha,\beta}_l,~\forall~k\neq l;$
\item $\oplus_{m\in\mathbb{Z}}W^{\alpha,\beta}_m=L^2(\mathbb{R})$
\item $g(\tee)\in W^{\alpha,\beta}_m$ if and only if $g(2\tee)e^{\frac{i}{2}((2\tee)^2-\tee^2)\cot\alpha}\in W^{\alpha,\beta}_{m+1},~\forall~m\in\mathbb{Z}.$
\end{enumerate}
From the second property, it follows that an orthonormal basis of $L^2(\mathbb{R})$ can be constructed by finding out an orthonormal basis of $W^{\alpha,\beta}_m.$ Third property says that the problem of finding this basis is centered to finding an orthonormal basis of $W^{\alpha,\beta}_0.$ Thus the crucial point is to construct a function $\psi\in L^2(\mathbb{R})$ such that $\{\psi_{\alpha,\beta,0,n}(\tee)\}_{n\in\mathbb{Z}}$ form an orthonormal basis of $W^{\alpha,\beta}_0,$ where $\psi_{\alpha,\beta,0,n}(\tee)=\psi(\tee-n)e^{-\frac{i}{2}((\tee^2-n^2)\cot\alpha-(\tee-n)^2\cot\beta)}.$ The two free parameters $\alpha$ and $\beta$ involved in the MRA helps in the construction of several new classes of orthonormal basis of $L^2(\mathbb{R}).$ The presence of these two free parameters makes the MRA developed in \cite{guo2021novel} more flexible as compared to that in \cite{shi2015multiresolution} and \cite{dai2017new}. Thus, it results in developing a new class of orthonormal basis of $L^2(\mathbb{R}).$ For example, for some scaling function $\phi$ taken in the fractional Fourier domain the  
%(This can be added later)
%if the scaling function $\phi$ is taken in the fractional Fourier domain as 
%\textcolor{purple}{
%$$(\mathcal{F}^\beta\phi)(\frac{\xxi}{\lambda})=
%\begin{cases}
%(2\pi\sin\beta)^{-\frac{1}{2}},&|\xxi|\leq \xxi_\alpha\\
%(0,(2\pi\sin\beta)^{-\frac{1}{2}}),&\xxi_\alpha<|-1/2\xxi|<2\xxi_\alpha\\
%0,&|\xxi|>2\xxi_\alpha,
%\end{cases}
%$$
%}
%where $\xxi_\alpha=\frac{2\pi|\sin\alpha|}{3},$ then the
%(This can be added later)
corresponding fractional wavelet is given in the fractional Fourier domain as 
$\left(\mathcal{F}^\beta\psi\right)\left(\frac{\xxi}{\lambda}\right)=e^{\frac{i}{2}\frac{3}{4}\left(\frac{\xxi}{\lambda}\right)^2\cot\alpha}(2\pi\sin\alpha)^{\frac{1}{2}}e^{\frac{i}{2}\frac{3}{4}\left(\frac{\xxi+2\pi\sin\alpha}{\lambda}\right)^2\cot\beta}(\mathcal{F}^\beta\phi)(\frac{\xxi+2\pi\sin\alpha}{\lambda})(\mathcal{F}^\beta\phi)(\frac{\xxi}{2\lambda}),$ where $\lambda=\frac{\csc\beta}{\csc\alpha}.$ Authors in \cite{ahmad2021fractional} proved that the conditions in the definition of the fractional MRA given in \cite{dai2017new}, which is a particular case when $\alpha=\beta,$ are not independent. More precisely, the intersection triviality condition $\cap_{m\in\mathbb{Z}}V^{\alpha,\alpha}=\{0\}$ can be obtained from the conditions \ref{R2021_Guo-Li-Yang_FrWT_MRA_Con1}, \ref{R2021_Guo-Li-Yang_FrWT_MRA_Con2} and \ref{R2021_Guo-Li-Yang_FrWT_MRA_Con4} for $\alpha=\beta.$ In addition to these assumptions if $(\mathcal{F}^\alpha\phi)(\xi)$ is continuous at $\xxi=0,$ they proved that the union density condition $\overline{\cup_{m\in\mathbb{Z}}V^{\alpha,\alpha}}=L^2(\mathbb{R})$ is equivalent to the condition $(\mathcal{F}^\alpha\phi)(0)\neq 0.$ These results are in fact the counter part of those given in \cite{hernandez1996first}.
\section{Fractional Stockwell Transform (FrST)}
In 2012, Xu et al.\cite{xu2012fractional} first proposed a FrST based on the idea of of the FrFT and the ST with Gaussian window as 
\begin{align}\label{R2012_Xu-Guo_FrST_Definition}
(\mathcal{S}^\alpha f)(\tau,\uu)=\int_{-\infty}^\infty f(\tee)g(\tau-\tee,\uu)R_\alpha(\tee,\uu)d\tee,
\end{align}
where the Gaussian window $g(\tau-\tee,\uu)$ is a scalable function of $\tee$ and the fractional Fourier frequency (FrFfr) $\uu$
$g(\tee,\uu)=\frac{|\uu\csc\alpha|^p}{\sqrt{2\pi}q}e^{-\frac{\tee^2(\uu\csc\alpha)^{2p}}{2q^2}}$ and $R_\alpha(\tee,\uu)=\sqrt{1-i\cot\alpha}e^{i\pi(\uu^2\cot\alpha-2\uu\tee\csc\alpha+\tee^2\cot\alpha)},~\alpha\nin\pi\mathbb{Z}.$ It is to be noted that, when $\alpha=\frac{\pi}{2}$ and $p=q=1,$ then the FrST \eqref{R2012_Xu-Guo_FrST_Definition} is the classical ST \eqref{R1996_Stockwell-Mansinha-Lowe_Definition_STransform} given in \cite{stockwell1996localization}. Observe that the window $g(\tee,\uu)$ width is directed by $\uu,$ i.e., the window width is wider in the time domain at lower FrFfr $\uu,$ whereas it is narrower at higher $\uu.$ Authors proved the linearity of the transform \eqref{R2012_Xu-Guo_FrST_Definition} followed by the scaling formula 
$$
(\mathcal{S}^\alpha f(c\tee))(\tau,\uu)=\sqrt{\frac{c^2(1-i\cot\alpha)}{c^2-i\cot\alpha}}e^{i\pi\uu^2\cot\alpha\left(1-\frac{\cos^2\beta}{\cos^2\alpha}\right)}(\mathcal{S}^\beta f)\left(c\tau,\uu\frac{\sin\beta}{c\sin\alpha}\right),
$$
where $\beta=\arctan(c^2\tan\alpha)$ followed by the inversion formula 
$$f(\tee)=\int_{-\infty}^\infty [\int_{-\infty}^\infty (\mathcal{S}^\beta f)(\tau,\uu)d\tau]\overline{R_\alpha(\tee,\uu)}d\uu.$$
For more properties and the applications of the fractional Stockwell transform \eqref{R2012_Xu-Guo_FrST_Definition}, we refer the reader to a survey paper  \cite{ranjan2020fractional}. Recently, Srivastava et al. \cite{srivastava2020family} introduced the FrST of a function $f\in L^2(\mathbb{R})$ with respect to the window function $\psi\in L^2(\mathbb{R})$ as 
\begin{align}\label{R2020_Srivastave-Shah-Tantary_FrST_Convolution_Definition}
(\mathcal{S}^\alpha_{\psi}f)(\A,\B)=\frac{1}{\sqrt{2\pi}}\left(\mathcal{M}_{-\xxi}f\star_\alpha\mathcal{D}_{\xi}\tilde{\psi}\right)(\B),
\end{align}
where $\tilde{\psi}(t)=\overline{\psi(-t)},$ and the convolution $\star_\alpha$ is the fractional convolution given by \eqref{RShiFrFTConvolution}, after multiplying it by the factor $\sqrt{2\pi}$. Equation \eqref{R2020_Srivastave-Shah-Tantary_FrST_Convolution_Definition} can be written in the form 
\begin{align}\label{R2020_Srivastave-Shah-Tantary_FrST_Definition_Inner-Product-Form}
(\mathcal{S}^\alpha_{\psi}f)(\A,\B)=\frac{1}{\sqrt{2\pi}}\int_{\mathbb{R}}f(\tee)\overline{\psi^\alpha_{\A,\B}(\tee)}d\tee,
\end{align}
where $\psi^\alpha_{\A,\B}(\tee)=|\A|\psi(\A(\tee-\B))e^{i\A\tee-\frac{i(\tee^2-\B^2)\cot\alpha}{2}},~\A\in\mathbb{R}^*,\B\in\mathbb{R}.$ It is  a matter of simple verification that $\mathcal{S}^\alpha_{\psi}f$ is linear in $f$ and conjugate linear in $\psi.$ The following properties called the covariance properties of the FrST give respectively the formula for the FrST of the translation by $\xx_0$ and dilation by $c\in \mathbb{R}^{\ast}$ of function $f$ with respect to the $\psi$
$$[(\mathcal{S}^\alpha_{\psi}(\tau_{\xx_0}f)](\A,\B)=e^{i\xx_0((\xx_0-\B)\cot\alpha-\A)}(\mathcal{S}^\alpha_{\psi}(\mathcal{M}_{\xx_0\cot\alpha}f))(\A,\B-\xx_0)$$
and
$$[(\mathcal{S}^\alpha_{\psi}(f(c\tee))](\A,\B)=(\mathcal{S}^\beta_{\psi}f)\left(\frac{\A}{c},c\B\right),~\beta=arc\cot\left(\frac{\cot\alpha}{c}\right).$$
If the function $\psi$ is such that it satisfied the admissibility condition, i.e., $c_\psi=\int_{\mathbb{R}^*}\left|\mathcal{F}[e^{iz}\psi(z)]\left(\frac{\xxi\csc\alpha}{\A}\right)\right|^2d\A$ is non-zero finite, then $\mathcal{S}^\alpha_\psi$ satisfies the resolution of identity formula 
\begin{align}\label{R2020_Srivastave-Shah-Tantary_FrST_IPR}
\int_{\mathbb{R}^*\times\mathbb{R}}(\mathcal{S}^\alpha_{\psi}f)(\A,\B)\overline{(\mathcal{S}^\alpha_{\psi}g)(\A,\B)}d\A d\B=c_{\psi}\langle f,g\rangle_{L^2(\mathbb{R})},
\end{align}
followed by the inversion formula 
$$f(\tee)=\frac{1}{\sqrt{2\pi}c_\psi}\int_{\mathbb{R}\times\mathbb{R}^*}(\mathcal{S}^\alpha_{\psi}f)(\A,\B)\psi^\alpha_{\A,\B}(\tee)d\A d\B.$$
For $f=g,$ \eqref{R2020_Srivastave-Shah-Tantary_FrST_IPR}, yields the energy preserving relation $\|(\mathcal{S}^\alpha_{\psi}g)(\A,\B)\|^2_{L^2(\mathbb{R}^*\times\mathbb{R})}=c_\psi\|f\|^2_{L^2(\mathbb{R})}.$ If further $c_\psi=1,$ the FrST  \eqref{R2020_Srivastave-Shah-Tantary_FrST_Convolution_Definition} becomes an isometry from $L^2(\mathbb{R})$ to $L^2(\mathbb{R}^*\times\mathbb{R}).$ Moreover, the range $\mathcal{S}^\alpha_{\psi}(L^2(\mathbb{R}))$ is a RKHS  with the kernel $K^\alpha_\psi(\A,\B;\A_0,\B_0)=\frac{1}{2\pi c_\psi}\langle \psi^\alpha_{\A,\B},\psi^\alpha_{\A_0,\B_0}\rangle.$ In view of equation \eqref{R2020_Srivastave-Shah-Tantary_FrST_Definition_Inner-Product-Form} and its representation in the FrFD, the FrST  \eqref{R2020_Srivastave-Shah-Tantary_FrST_Convolution_Definition} gives the local information of $f(\tee)$ in the time window $\left[\B+\frac{E_\psi}{\A}-\frac{\Delta_\psi}{\A},\B+\frac{E_\psi}{\A}+\frac{\Delta_\psi}{\A}\right]$ and the local FrFD spectral information in the FrFD as $[\A\sin\alpha E_\Psi-\A\sin\alpha\Delta_\Psi,\A\sin\alpha E_\Psi+\A\sin\alpha\Delta_\Psi]$, where $\Psi(\xxi)=\mathcal{F}^\alpha[e^{iz}\psi(z)](\xxi)$ with the constant window area $4\sin\alpha\Delta_\psi\Delta_\Psi$ in the time-fractional-frequency plane. The area depends only on $\psi$ and the parameter $\alpha$ and is independent of the parameter $\A$ and $\B.$ Also, it is worth noting that the fractional parameter $\alpha$ can be invoked to optimize the concentration of the Stockwell spectrum in the time-fractional-frequency plane. Wei et al. \cite{wei2021fractional} used the definition of fractional convolution $\star_\alpha$ \eqref{RShiFrFTConvolution}, to give a new definition of FrST for the function $f\in L^2(\mathbb{R})$ with respect to the window function $\psi$ as
\begin{align}\label{R2021_Wei-Zhang_FrST_Definition}
(\mathcal{S}_\psi^\alpha f)(\A,\B)=&\frac{1}{\sqrt{2\pi}}[e^{-i\A\tee}f(\tee)\star_\alpha \A\tilde{
\psi}(\A\tee)](\B)\notag\\
=&\frac{1}{\sqrt{2\pi}}\int_{\mathbb{R}}f(\tee)\overline{\psi_{\alpha,\A,\B}(\tee)}d\tee,
\end{align}
where $\psi_{\alpha,\A,\B}(\tee)=\A\psi(\A(\tee-\B))e^{i\A\tee\csc\alpha-\frac{i}{2}(\tee^2-\B^2)\cot\alpha}.$ They proved that the transform is linear in $f$ and anti-linear in $\psi$ and obtained the formula for the FrST of the translation and dilation of $f$ as 
$$(\mathcal{S}_\psi^\alpha f(\tee-k))(\A,\B)=e^{-i\A k\csc\alpha+i(k^2-k\B)\cot\alpha}(\mathcal{S}_\psi^\alpha (f(\tee)e^{ik\tee\cot\alpha}))(\A,\B-k)$$
and 
$$(\mathcal{S}_\psi^\alpha f(\lambda\tee))(\A,\B)=(\mathcal{S}_\psi^\beta f)\left(\frac{\A}{\lambda},\lambda\B\right),~\lambda\in\mathbb{R}^+,~\beta=\mbox{arc}\cot\left(\frac{\alpha}{\lambda^2}\right)$$
If the window function $\psi$ satisfies the admissible conditions $c_\psi=\int\left|\mathcal{F}(\psi(\tee)e^{j\tee\csc\alpha})(u\csc\alpha/\A)\right|d\A<\infty$ then $\mathcal{S}_\psi^\alpha$ satisfies the inner product relations
$$\int_0^\infty\int_{-\infty}^\infty (\mathcal{S}_\psi^\alpha f)(\A,\B)\overline{(\mathcal{S}_\psi^\alpha g)(\A,\B)}d\A d\B=\frac{1}{2\pi}c_\psi\langle f,g\rangle$$
followed by the inversion formula 
$$f(\tee)=\frac{2\pi}{c_\psi}\int_0^\infty\int_{-\infty}^\infty(\mathcal{S}_\psi^\alpha f)(\A,\B)\psi_{\alpha,\A,\B}(\tee)d\A d\B$$
In view of \eqref{R2021_Wei-Zhang_FrST_Definition} and its representation in the fractional Fourier domain, it follows that the FrST \eqref{R2021_Wei-Zhang_FrST_Definition} is capable of providing the time and spectral information simultaneous and representing signal in the time-fractional-frequency domain. It gives local information of the signal $f(\tee)$ with the time window $\left[\B+\frac{E_\psi}{\A}-\frac{\Delta_\psi}{\A},\B+\frac{E_\psi}{\A}+\frac{\Delta_\psi}{\A}\right]$ and the local information about the fractional Fourier domain spectrum of $f(\tee)$ as $\left[\A E_G-\A\Delta_G,\A  E_G+\A\Delta_G\right]\sin\alpha,$ with the constant window area $4\sin\alpha\Delta_\psi\Delta_G$ in the time-fractional-frequency plane where $G(\uu)=\mathcal{F}(\psi(\tee)e^{i\tee\csc\alpha})(\uu)$. This area depends on the function $\psi$ and the rotation angle $\alpha$ and is independent of the parameters $\A$ and $\B.$ Moreover, the window narrows for the small value and widens for larger values of $\A.$ 

Kamalakkannan et al. \cite{kamalakkannan2022two} proposed a two-dimensional FrST, associated with the coupled FrFT \eqref{R2022_Kamalakkannan-Roopkumar-Zayed_Coupled-FrFT-Definition}, of a function $f\in L^2(\mathbb{R}^2)$ with respect to a non-zero window function $\psi\in L^1(\mathbb{R}^2)\cap L^2(\mathbb{R}^2)$ as 
\begin{align}\label{R2022_Kamalakkannan-Roopkumar_Coupled-FrSFT-Definition}
(\mathcal{S}_{\psi,\alpha,\beta}f)(\boldsymbol z,\bldv)=(E_\gamma D_{\boldsymbol\rho}\breve{\overline{\psi}}\bigotimes\mathcal{M}_{\boldsymbol\rho}f)(\boldsymbol z),
\end{align}
$\boldsymbol z\in\mathbb{R}^2,$ $\bldv\in (\mathbb{R}^*)^2=\{(\xx_1,\xx_2)^t\in\mathbb{R}^2:\xx_1\neq 0,\xx_2\neq 0\}$ and $\boldsymbol\rho=(\rho_1,\rho_2)^t=\boldsymbol\zeta\bldv$ (see \eqref{R2022_Kamalakkannan-Roopkumar-Zayed_Coupled-FrFT-Definition}), where $\bigotimes$ denotes the convolution associated with the coupled FrFT given in \eqref{R2022_Kamalakkannan-Roopkumar-Zayed_Coupled-FrFT-Convolution}, $E_\gamma$ is same as that in \eqref{R2022_Kamalakkannan-Roopkumar-Zayed_Coupled-FrFT-FT_Relation} and $(D_{\boldsymbol\rho}f)(\bldx)=|\det(A_{\boldsymbol\rho})|f(A_{\boldsymbol\rho}\bldx),$ $A_{\boldsymbol\rho}=
\begin{bmatrix}
\rho_1 & 0\\
0&\rho_2
\end{bmatrix}
.$ With the assumption that the window function $\psi$ satisfies the admissibility condition  
$$0\neq c_\psi=\frac{1}{2\pi}\int_{(\mathbb{R}^*)^2}\left|\hat{\psi}(\bldxi-(1,1)^t)\right|^2\frac{d\bldxi}{|\det(A_{\bldxi})|}<\infty$$
authors obtained the inner product relation of the FrST
\begin{align}\label{R2022_Kamalakkannan-Roopkumar_Coupled-FrSFT-IPR}
c_\psi\langle f,g\rangle=\int_{(\mathbb{R}^*)^2}(\mathcal{S}_{\psi,\alpha,\beta}f)(\boldsymbol z,\bldv)\overline{(\mathcal{S}_{\psi,\alpha,\beta}g)(\boldsymbol z,\bldv)}d\boldsymbol z\frac{d\bldv}{|\det(A_{\boldsymbol \rho})|},~\boldsymbol\rho=\boldsymbol\zeta\bldv,
\end{align}
followed by the inversion formula
\begin{align}\label{R2022_Kamalakkannan-Roopkumar_Coupled-FrSFT-Inversion_Formula}
f(\bldx)=\frac{d}{2\pi\bar{d}c_\psi}\int_{(\mathbb{R}^*)^2}\mathcal{M}_{-\boldsymbol\rho}\left((\mathcal{S}_{\psi,\alpha,\beta}f)(\cdot,\bldv)\bigotimes E_\gamma D_{\boldsymbol\rho}\psi\right)(\bldx)\frac{d\bldv}{|\det(A_{\boldsymbol\rho})|}.
\end{align}
As a consequence of equation \eqref{R2022_Kamalakkannan-Roopkumar_Coupled-FrSFT-IPR}, it follows that $\mathcal{S}_{\psi_\alpha,\beta}f\in L^2\left(\mathbb{R}^2\times L^2(\mathbb{R}^2),\frac{d\boldsymbol zd\bldv}{|\det(A_{\boldsymbol\rho})|}\right).$ Authors in \cite{kamalakkannan2022two}, characterized the range of the transform and proved that $\mathcal{S}_{\psi,\alpha,\beta}(L^2(\mathbb{R}^2))$ is a RKHS with the kernel $K^{\alpha,\beta}_\psi(\boldsymbol z,\bldv,\boldsymbol z',\bldv')=\frac{1}{2\pi}\int_{\mathbb{R}^2}\psi^{\alpha,\beta}_{\boldsymbol z,\bldv}(\bldx)\overline{\psi^{\alpha,\beta}_{\boldsymbol z',\bldv'}(\bldx)}d\bldx,$ where $\psi^{\alpha,\beta}_{\boldsymbol z,\bldv}(\bldx)=2\pi\bar{d}E^{-1}_\gamma(\boldsymbol z)(E_\gamma\mathcal{M}_{\boldsymbol\rho}\tau_{\boldsymbol z}D_{\boldsymbol\rho}\psi)(\bldx).$ For $f\in L^p(\mathbb{R}^2),~p=1,2$ and $g\in L^1(\mathbb{R}^2),$ the formula for the FrST of convolution $\bigotimes$ is obtained as
\begin{align}
\mathcal{S}_{\psi,\alpha,\beta}(f\bigotimes g)(\boldsymbol z,\bldv)=\left((\mathcal{S}_{\psi,\alpha,\beta}f)(\cdot,\bldv)\bigotimes\mathcal{M}_{\boldsymbol\rho}g\right)(\boldsymbol z).
\end{align}

\noindent\textbf{Uncertainty Principles:} Srivastava et al. \cite{srivastava2020family} obtained the Heisenberg's uncertainty inequality associated with the FrST \eqref{R2020_Srivastave-Shah-Tantary_FrST_Convolution_Definition} as
\begin{align}\label{R2020_Srivastave-Shah-Tantary_FrST_Heisenberg's_UP}
\left\{\int_{\mathbb{R}}|\xxi|^2|(\mathcal{F}^\alpha f)(\xxi)|^2d\xxi\right\}^\frac{1}{2}\left\{\int_{\mathbb{R}\times\mathbb{R}^*}|\omega|^2|(\mathcal{F}^\beta((\mathcal{S}^\alpha_{\psi}f)(\A,\B)))(\omega)|^2d\A d\omega\right\}^\frac{1}{2}\geq \frac{\sqrt{c_\psi}|\sin(\alpha-\beta)|}{2}\|f\|^2_{L^2(\mathbb{R})},
\end{align}
where $\alpha$ and $\beta$ are chosen such that $\beta=\alpha-\gamma$ and $\sin\alpha,\sin\beta,\sin\gamma\neq 0.$ The fractional parameters govern these relations and are fruitful in understanding mutual relations among different transform domains. Using the Heisenberg UP for the coupled FrFT, the author in \cite{kamalakkannan2022two} obtained the Heisenberg UP for the FrST \eqref{R2022_Kamalakkannan-Roopkumar_Coupled-FrSFT-Definition}, which is given by 
$$\Delta\boldsymbol z^2_{\alpha,\beta}\Delta\bldu^2_{\alpha,\beta}\geq c_\psi\sin^2\gamma\|f\|^4_2,$$
where $\Delta\boldsymbol z^2_{\alpha,\beta}=\frac{1}{4\pi^2}\int_{(\mathbb{R}^*)^2}\int_{\mathbb{R}^2}|\boldsymbol z|^2|(\mathcal{S}_{\psi,\alpha,\beta}f)(\boldsymbol z,\bldv)|^2d\boldsymbol z\frac{d\bldv}{|\det(A_{\boldsymbol \rho})|}$ and $\Delta\bldu^2_{\alpha,\beta}=\frac{1}{2\pi}\int_{\mathbb{R}^2}|\bldu|^2|(\mathcal{F}_{\alpha,\beta}f)(\bldu)|^2d\bldu.$

\section{Windowed Linear Canonical Transform (WLCT)}

Compared to the FT and the FrFT, LCT is more flexible with more extra degrees of freedom but with a similar computational cost as the conventional FT. Because of the mentioned advantages, the LCT has found many applications in filter designing, signal synthesis, optics, radar analysis, pattern recognition, etc. For example, filtering in the linear canonical domain, given in \cite{barshan1997optimal}, can achieve better performance than in the fractional Fourier domain owing to more degree of freedom. Especially when multi-component chirp signals interfere with the desired signal, only one filter is used in the LCT domain usually, but several filters are required in the FRFT domain \cite{erden1999repeated}. Despite these benefits, it has limitations on its own. The LCT cannot reveal the local LCT-frequency component because of its global kernel. As a result, Kou et al. \cite{kou2012windowed} introduced the WLCT by substituting the FT kernel with the LCT kernel in the WFT definition, which is defined for the function $f$ with respect to $g$ as

\begin{align}\label{R2012_Kou-Xu_SP_WLCT_Definition}
(\mathcal{V}^M_gf)(\xx,\uu)=\int_{\mathbb{R}}f(\tee)\overline{\tau_{\xx}g(\tee)}K_M(\tee,\uu)d\tee,
\end{align}
where $M=(\Ma,\Mb;\Mc,\Md)$ with $\det(M)=1,$ and $\Mb>0.$ Observe that, if $f\in L^2(\mathbb{R}),~g\in L^2(\mathbb{R})$ and the window function $f\in L^1(\mathbb{R}),~g\in L^\infty(\mathbb{R})$ then the WLCT \eqref{R2012_Kou-Xu_SP_WLCT_Definition} is well defined. Similarly, if $f\in L^q(\mathbb{R})$ and $f\in L^p(\mathbb{R}), 1\leq p,q<\infty$ and $\frac{1}{p}+\frac{1}{q}=1,$ then by H\"older's inequality it follows that the WLCT \eqref{R2012_Kou-Xu_SP_WLCT_Definition} is well defined. Bahri et al.\cite{bahri2016some} proved that the WLCT satisfied the inequality $\|(\mathcal{V}^M_gf)(\cdot,\uu)\|_{L^p(\mathbb{R})}\leq \frac{1}{\sqrt{2\pi\Mb}}\|g\|_{L^p(\mathbb{R})}\|f\|_{L^1(\mathbb{R})}$ for $g\in L^p(\mathbb{R})$ and $f\in L^1(\mathbb{R}).$ They also established the relation between the WLCT and the WFT and obtained  the Lieb's inequality for the WLCT 
\begin{align*}
\int_{\mathbb{R}}\int_{\mathbb{R}}\left|(\mathcal{V}^M_gf)(\xx,\uu)\right|^pd\xx d\uu\leq\frac{2}{p}\left(\frac{1}{\sqrt{2\pi}}|\Mb|^{\frac{1}{p}-\frac{1}{2}}\right)^p\|f\|^p_{L^2(\mathbb{R})}\|g\|^p_{L^2(\mathbb{R})}
\end{align*}
using the same \eqref{RDonoho-Stark_UP_Strong_Version_STFT} for the WFT \cite{grochenig2001foundations}. The proof is simple as compared to the one given in \cite{kou2012paley}. Transform \eqref{R2012_Kou-Xu_SP_WLCT_Definition} is linear in $f,$ anti-linear in $g.$ It satisfies the time marginal constraints, i.e., integrating the WLCT  with respect to the time variable gives the LCT multiplied by a constant, i.e., $\int_{\mathbb{R}}(\mathcal{V}^M_gf)(\xx,\uu)d\xx=(\int_{\mathbb{R}}g(\tee)d\tee)\mathcal{L}^M(f)(\uu)$ and the additivity of rotation, i.e., $\mathcal{L}^{M_2}((\mathcal{V}^{M_1}_gf)(\xx,\uu))(\xx,\uu')=(\mathcal{V}^{M_2M_1}_gf)(\xx,\uu'),$ where $M_l=(\Ma_l,\Mb_l;\Mc_l,\Md_l),~l=1,2$ and $\det(M_l)=1,~B_l>0$ \cite{kou2012windowed}. Apart from this, it has the following covariance properties 
\begin{align*}
\left(\mathcal{V}^M_g\left(\tau_{\xx_0}\mathcal{M}^{M}_\mu f\right)\right)(\xx,\uu)=\frac{1}{\sqrt{i2\pi\Mb}}e^{-\frac{i}{\Mb}\left(\frac{\Md}{2}\uu^2-\frac{\Md}{2}\mu^2+\frac{\Md}{2}\xx_0^2-2t\xx_0\mu-\xx_0\uu\right)}\mathcal{V}^M_gf(\xx-\xx_0,\frac{1}{\Mb}[\uu-(\mu+\Ma\xx_0)])
\end{align*}
and 
\begin{align*}
\left(\mathcal{V}^M_g\left(\mathcal{M}^{M}_\mu\tau_{\xx_0} f\right)\right)(\xx,\uu)=\frac{1}{\sqrt{i2\pi\Mb}}e^{-\frac{i}{\Mb}\left(\frac{\Md}{2}\uu^2-\frac{\Md}{2}\mu^2+\frac{\Md}{2}\xx_0^2-2t\xx_0\mu-\xx_0\uu\right)}\mathcal{V}^M_gf(\xx-\xx_0,\frac{1}{\Mb}[\uu-(\mu+\Ma\xx_0)]),
\end{align*}
where $\mathcal{M}^{M}_\mu f(\tee)=e^{-i\left(\frac{\Ma}{2\Mb}\tee^2-\frac{1}{\Mb}\tee\mu+\frac{\Md}{2\Mb}\mu^2\right)}f(\tee)$ is generalized modulation by $\mu.$ The covariance property obtained above indeed gives the formula for the WLCT for the time-frequency shift $\tau_{\xx_0}\mathcal{M}^{M}_\mu$ or $\mathcal{M}^{M}_\mu\tau_{\xx_0}$ of the function $f.$ Specifically, the formula for the parity, conjugation,  $\left(\mathcal{V}^M_g\left(\tau_{\xx_0} f\right)\right)(\xx,\uu)$ and $\left(\mathcal{V}^M_g\left(e^{i\uu_0\xx}f\right)\right)(\xx,\uu)$ can be found in \cite{bahri2016some}. For $f_1,f_2,g_1,g_2\in L^2(\mathbb{R}),$ the WLCT \eqref{R2012_Kou-Xu_SP_WLCT_Definition} satisfies the following inner product relation \cite{kou2012windowed},\cite{bahri2016some}
\begin{align}\label{R2012_Kou-Xu_SP_WLCT_IPR}
\left\langle(\mathcal{V}^M_{g_1}f_1)(\xx,\uu),(\mathcal{V}^M_{g_2}f_2)(\xx,\uu)\right\rangle_{L^2(\mathbb{R}^2)}=\langle f_1,f_2 \rangle\overline{\langle g_1,g_2 \rangle}.
\end{align}
In particular, if $f_1=f_2=f$ and $\|g_1\|_{L^2(\mathbb{R})}=\|g_2\|_{L^2(\mathbb{R})}=1$ then $\|\mathcal{V}^M_gf\|_{L^2(\mathbb{R}^2)}=\|f\|_{L^2(\mathbb{R})}$ and thus WLCT is an isometry from $L^2(\mathbb{R})$ to $L^2(\mathbb{R}^2)$. The inversion formula for the WLCT is given as 
\begin{align}\label{R2012_Kou-Xu_SP_WLCT_Inversion-Formula}
f(\tee)=\frac{1}{\langle g_1,g_2\rangle}\int_{\mathbb{R}}\int_{\mathbb{R}}(\mathcal{V}^M_{g_2}f)(\xx,\uu)\overline{K_M(\tee,\uu)}\tau_{\xx}g_1(\tee)d\uu d\xx
\end{align}
where $g_1,g_2\in L^2(\mathbb{R}),$ $\langle g_1,g_2\rangle\neq 0$  and $f\in L^2(\mathbb{R}).$ In fact, a stronger version of the inversion formula \eqref{R2012_Kou-Xu_SP_WLCT_Inversion-Formula} is proved by the authors in \cite{kou2012windowed}. Its formulation resembles an $L^2-$function by an approximation procedure. For that a nested sequence of compact sets $K_n\subset \mathbb{R}^2$ satisfying $\cup_{n\in\mathbb{N}}K_n=\mathbb{R}^2$ and $K_n\subset \mbox{int} K_{n+1}$ is considered and that for each $n\in\mathbb{N},$ $f_n(\tee)=\frac{1}{\langle g_1,g_2\rangle}\int\int_{K_n}(\mathcal{V}^M_{g_2}f)(\xx,\uu)\overline{K_M(\tee,\uu)}\tau_{\xx}g_1(\tee)d\uu d\xx$ then $\lim_{n\rightarrow\infty}\|f-f_n\|_{L^2(\mathbb{R})}=0.$ Bahri et al.\cite{bahri2016some} characterized the range $\mathcal{V}^M_g(L^2(\mathbb{R})$ of the transform \eqref{R2012_Kou-Xu_SP_WLCT_Definition} showing that $h\in\mathcal{V}^M_g(L^2(\mathbb{R}))$ if and only if it satisfies $h(\xx',\uu')=\int_{\mathbb{R}}\int_{\mathbb{R}}h(\xx,\uu)\langle g^M_{\xx,\uu},g^M_{\xx',\uu'}\rangle_{L^2(\mathbb{R})}d\xx d\uu.$ Consequently, it follows that the range  $\mathcal{V}^M_g(L^2(\mathbb{R})$ is a RKHS with the kernel $\langle g^M_{\xx,\uu},g^M_{\xx',\uu'}\rangle_{L^2(\mathbb{R})},$ where $g^M_{\xx,\uu}(\tee)=\frac{1}{\sqrt{i2\pi\Mb}}e^{-\frac{i}{2}\left(\frac{\Ma}{\Mb}\tee^2-\frac{2}{\Mb}\tee\uu+\frac{\Md}{\Mb}\uu^2\right)}\phi(\tee-\xx).$ Bahri in \cite{bahri2022windowed} obtained a relation between the WLCT \eqref{R2012_Kou-Xu_SP_WLCT_Definition} and the WFT, using which he gave an alternative proof of the properties like orthogonality property and inversion formula. Han et al.\cite{han2022inversion} studied the inversion formula for recovering a signal from its WLCT different from the inversion formula \eqref{R2012_Kou-Xu_SP_WLCT_Inversion-Formula} by invoking the double integral and showing that every signal can be recovered from its WLCT, with a univariate integral, with a suitable choice of window function $g$ as
\begin{align}\label{R_1DInversion_Formula_WLCT}
f(\xx)=\frac{1}{\overline{g(0)}}\int_{\mathbb{R}}(\mathcal{V}^M_{g}f)(\xx,\uu)\overline{K_M(\xx,\uu)}d\uu,
\end{align}
which generalizes a similar result for the WFT \cite{sun2012inversion}. More precisely, they proved that for $f\in L^p(\mathbb{R}),$ the integral in \eqref{R_1DInversion_Formula_WLCT} is convergent in $L^p(\mathbb{R})$ for all $1<p<\infty$ by showing that 
$\lim_{u_1,u_2\rightarrow \infty}\|R_{u_1,u_2}f-\overline{g(\xx_0)}f\|_{L^p(\mathbb{R})}=0$
and they also proved $\lim_{\uu\rightarrow \infty}(R_{u,u}f)(\xx)=\overline{g(\xx_0)}f(\xx),~a.e.$, by applying  the Carleson-Hunt theorem, showing that the integral \eqref{R_1DInversion_Formula_WLCT} converges a.e. on $\mathbb{R},$ where where $R_{u_1,u_2}f(\xx)=\int_{-u_1}^{u_2}(\mathcal{V}^M_gf)(\xx-\xx_0,\uu)\overline{K_M(\xx,\uu)}d\uu,$ $u_1,u_2>0,$ $\xx_0$ is a fixed real number and $g$ is a continuous function such that $g,\hat{g}\in L^1(\mathbb{R}).$ With the same assumption on $g,$ authors in \cite{han2022inversion} obtained the recovery formula of function in $L^1(\mathbb{R})$ from its WLCT by using the Ces\`aro summability method by showing that 
$\lim_{h\rightarrow\infty}\|\sigma_hf-\overline{g(\xx_0)}f\|_{L^p(\mathbb{R})}=0,~1\leq p<\infty$
and $\lim_{h\rightarrow\infty}(\sigma_hf)(\xx)=\overline{g(\xx_0)}f(\xx),~a.e.,$
where  $\sigma_hf(\xx)=\frac{1}{h}\int_0^h(R_{u,u}f)(\xx)d\uu.$ Using the convolution \eqref{RConvolutionGeneralizedLCT} given by Shi et al.\cite{shi2014generalizedconvolution}, Zhang et al.\cite{zhang2022novel} introduced a novel WLCT, which is defined for $f\in L^2(\mathbb{R})$ and the parameter $M_l=(\Ma_l,\Mb_l;\Mc_l,\Md_l),~l=1,2,3$ as
\begin{align}\label{R2022_Zhang-Guo-Mao_DSP_WLCT_Definition}
(\mathcal{V}^{M_1,M_2,M_3}_gf)(\tee,\uu)=\int_{\mathbb{R}}f(\tau)\overline{g_{M_1,M_2,M_3}(\tee,\uu,\tau)}d\tau,
\end{align}
where $g_{M_1,M_2,M_3}(\tee,\uu,\tau)=g(\tau-\tee)e^{i\frac{\uu\tau}{\Mb_1}-i\frac{\Ma_1}{2\Mb_1}\tau^2-i\frac{\Ma_2}{2\Mb_2}(\tee-\tau)^2+i\frac{\Ma_3}{2\Mb_3}\tee^2}.$ Equation \eqref{R2022_Zhang-Guo-Mao_DSP_WLCT_Definition} can be written in the generalized canonical convolution \eqref{RConvolutionGeneralizedLCT} form as 
\begin{align*}
\left(\mathcal{V}^{M_1,M_2,M_3}_gf\right)(\tee,\uu)=\left[\left(f(v)e^{-i\frac{\uu}{\Mb_1}v}\right)\star_{M_1,M_2,M_3}\overline{g(-v)}\right](\tee).
\end{align*}
With a particular choice of the matrices i.e.,  $M_1=M_3$ and $M_2=(0,1,-1,0),$ the transform reduces to the WLCT mentioned in \cite{shi2020novel},\cite{wei2021theory} and for the choices $M_1=M_3=(\cos\alpha,\sin\alpha,-\sin\alpha,\cos\alpha),$ $M_2=(0,1,-1,0),$ the transform \eqref{R2022_Zhang-Guo-Mao_DSP_WLCT_Definition} reduces to the novel WFrFT \eqref{R2020_Shi-Xiang-Zhang_IEEETSP_NovelSTFT_Definition} which analyze the signal in one linear canonical and one Fourier domain, whereas it can be seen using the convolution theorem \eqref{RConvolutionTheoremGeneralizedLCT} that the transform analyzes the signal in two distinct linear canonical domains. Authors in \cite{shi2020novel} only defined the WLCT, for which the mathematical properties were recently explored by Wei et al.\cite{wei2021theory}. The linear transformation \eqref{R2022_Zhang-Guo-Mao_DSP_WLCT_Definition} have the covariance property
$$(\mathcal{V}^{M_1,M_2,M_3}_g(\tau_{\tee_0}f))(\tee,\uu)=e^{-i\frac{\uu\tee_0}{\Mb_1}+i\frac{\Ma_1}{2\Mb_1}\tee_0^2+i\frac{\Ma_3}{2\Mb_3}\tee_0^2-i\frac{\Ma_3}{\Mb_3}\tee\tee_0}(\mathcal{V}^{M_1,M_2,M_3}_gf)(\tee-\tee_0,\uu-\Ma_1\tee_0)$$
for classical time shift operator and for the linear canonical time shift $\tau^{M_1}_{\tee_0}f$ the novel WLCT is given by
$$(\mathcal{V}^{M_1,M_2,M_3}_g(\tau^{M_1}_{\tee_0}f))(\tee,\uu)=e^{-i\frac{\uu\tee_0}{\Mb_1}+i\frac{\Ma_3}{2\Mb_3}\tee_0^2-i\frac{\Ma_3}{\Mb_3}\tee\tee_0}(\mathcal{V}^{M_1,M_2,M_3}_gf)(\tee-\tee_0,\uu),$$
where the operator is $\tau^{M_1}_{\tee_0}$ is given by $\tau^{M_1}_{\tee_0}f(t)=f(\tee-\tee_0)e^{-i\frac{\Ma_1}{\Mb_1}\left(\tee-\frac{\tee_0}{2}\right)\tee_0}.$ Also the transform for the classical modulation of $f$ is given as $(\mathcal{V}^{M_1,M_2,M_3}_g(\mathcal{M}_{\uu_0}f))(\tee,\uu)=(\mathcal{V}^{M_1,M_2,M_3}_gf)(\tee,\uu-\Mb_1\uu_0).$ Authors in \cite{zhang2022novel} obtained the inner product relation 
\begin{align}\label{R2022_Zhang-Guo-Mao_DSP_WLCT_IPR}
\int_{\mathbb{R}}\int_{\mathbb{R}}(\mathcal{V}^{M_1,M_2,M_3}_gf_1)(\tee,\uu)\overline{(\mathcal{V}^{M_1,M_2,M_3}_gf_2)(\tee,\uu)}d\tee d\uu=2\pi|\Mb_1|\|g\|^2_{L^2(\mathbb{R})}\langle f_1,f_2\rangle_{L^2(\mathbb{R})} 
\end{align}
followed by the inversion formula given by $f(\xx)=\frac{1}{2\pi|\Mb_1|}\int_{\mathbb{R}}\int_{\mathbb{R}}(\mathcal{V}^{M_1,M_2,M_3}_gf)(\tee,\uu)g_{M_1,M_2,M_3}(\tee,\uu,\xx)d\tee d\uu.$ With the assumption that $f_1=f_2=f$ and $\|g\|_{L^2(\mathbb{R})}=\frac{1}{2\pi},$ it follows from \eqref{R2022_Zhang-Guo-Mao_DSP_WLCT_IPR}, that the novel WLCT \eqref{R2022_Zhang-Guo-Mao_DSP_WLCT_Definition} is an isometry. Moreover, the range $\mathcal{V}^{M_1,M_2,M_3}_g(L^2(\mathbb{R}))$ is a reproducing kernel with the kernel given by $K(\tee,\uu;\tee_0,\uu_0)=\int_{\mathbb{R}}g_{M_1,M_2,M_3}(\tee,\uu,\tau)\overline{g_{M_1,M_2,M_3}(\tee_0,\uu_0,\tau)}d\tau.$ In view of equation \eqref{R2022_Zhang-Guo-Mao_DSP_WLCT_Definition} and its representation in the linear canonical domain obtained using \eqref{RConvolutionTheoremGeneralizedLCT}, the WLCT gives the local information of $f(\tau)$ in the time window $\left[E_g+\tee-\Delta_g,E_g+\tee+\Delta_g\right]$ and the local information of $f(\tau)$ in the linear canonical domain as 
$$\left[\Mb_3E_{\tilde{\tilde{G}}}+\Mb_3\left(\frac{\uu}{\Mb_1}-\frac{\Ma_1}{\Mb_1}\tee+\frac{\Ma_3}{\Mb_3}\tee\right)-\Mb_3\Delta_{\tilde{\tilde{G}}},\Mb_3E_{\tilde{\tilde{G}}}+\Mb_3\left(\frac{\uu}{\Mb_1}-\frac{\Ma_1}{\Mb_1}\tee+\frac{\Ma_3}{\Mb_3}\tee\right)+\Mb_3\Delta_{\tilde{\tilde{G}}}\right],$$ 
where $\tilde{\tilde{G}}(\xi)$ is the FT of $\tilde{\tilde{g}}(\yy)=g(\yy)e^{-i\frac{\Ma_1}{2\Mb_1}\yy^2-i\frac{\Ma_2}{2\Mb_2}\yy^2+i\frac{\Ma_3}{2\Mb_3}\yy^2},$ having the constant window area $4\Mb_3\Delta_g\Delta_{\tilde{\tilde{G}}}$ in the time-linear-canonical-frequency plane with the optimal window function as 
$$g(\tee)=\frac{1}{\sqrt{2\pi}\sigma}e^{-\frac{\tee^2}{2\sigma^2}}e^{i\frac{\Ma_1}{2\Mb_1}\tee^2+i\frac{\Ma_2}{2\Mb_2}\tee^2-i\frac{\Ma_3}{2\Mb_3}\tee^2},$$ 
where $\sigma$ denoted the window width in time domain.\\

\noindent\textbf{Uncertainty Principles:} Bahri et al. \cite{bahri2016some} obtained the Heisenberg's UP for the WLCT \eqref{R2012_Kou-Xu_SP_WLCT_Definition} using the Heisenberg's UP  \cite{guanlei2009new} for the LCT and adopting the method used by Wilczok \cite{wilczok2000new}
$$\left(\int_{\mathbb{R}}\int_{\mathbb{R}}\uu^2|(\mathcal{V}^M_gf)(\xx,\uu)|^2d\xx d\uu\right)^{\frac{1}{2}}\left(\int_{\mathbb{R}}\xx^2|f(\xx)|^2d\xx\right)^{\frac{1}{2}}\geq \frac{\Mb}{2}\|f\|^2_{L^2(\mathbb{R})}\|g\|_{L^2(\mathbb{R})}$$
and also gave the logarithmic UP for the transform
\begin{align}\label{R2022_Bahri_IJWMIP_WLCT_Logarithmic_UP}
\int_{\mathbb{R}}\int_{\mathbb{R}}\ln|\uu||(\mathcal{V}^M_gf)(\xx,\uu)|^2d\xx d\uu+\|g\|^2_{L^2(\mathbb{R})}\int_{\mathbb{R}}\ln|\xx||f(\xx)|^2d\xx\geq (D+\ln|\Mb|)\|f\|^2_{L^2(\mathbb{R})}\|g\|^2_{L^2(\mathbb{R})}.
\end{align}
Bahri in \cite{bahri2022windowed} proved that if $\mathcal{V}^M_gf$ is $\epsilon-$concentrated on a measurable set $V\subset\mathbb{R}\times\mathbb{R},$ then for $r>0,$ $\mu(V)\geq (1-\epsilon^2)^{\frac{r}{r-2}}(2\pi \Mb).$ Moreover, the author obtained the local UP for the WLCT \eqref{R2012_Kou-Xu_SP_WLCT_Definition}, using its relation with the FT and the local UP for the FT, as   
\begin{align}\label{R2022_BahriJoIA__WLCT_Local-UP1}
\int_{\mathbb{R}}\int_{\Mb X}|(\mathcal{V}^M_gf)(\xx,\uu)|^2d\xx d\uu\leq 
\begin{cases}
C_\alpha \|g\|^2_{L^2(\mathbb{R})}[\mu(X)]^{2\alpha}\int_{\mathbb{R}}|\xx|^{2\alpha}|f(\xx)|^2d\xx,~0<\alpha<\frac{1}{2}\\
C_\alpha \mu(X)\|g\|^2_{L^2(\mathbb{R})}\left(\int_{\mathbb{R}}|f(\xx)|^2d\xx\right)^{\frac{2\alpha-1}{2\alpha}}\left(\int_{\mathbb{R}}|\xx|^{2\alpha}|f(\xx)|^2d\xx\right)^{\frac{1}{2\alpha}},~\alpha>\frac{1}{2}
\end{cases}
\end{align}
where $C_\alpha$ is some constant. Hleili \cite{hleili2022windowed} obtained the entropy UP for the WLCT, 
$$E\left(|\mathcal{V}^M_gf|^2\right)\geq \|f\|^2_{L^2(\mathbb{R})}\|g\|^2_{L^2(\mathbb{R})}\left(\ln\left(2\pi|\Mb|\right)-\ln\left(\|f\|^2_{L^2(\mathbb{R})}\|g\|^2_{L^2(\mathbb{R})}\right)\right).$$ Using the logarithmic Sobolev-type UP for the FT give by Kubo et al.\cite{kubo2019beckner}, author in \cite{bahri2022windowed} also obtained the logarithmic Sobolev-type inequality for the WLCT
\begin{align}\label{R2022_BahriJoIA__WLCT_Logarithmic-Sobolev-type-UP}
\|g\|^2_{L^2(\mathbb{R})}\int_{\mathbb{R}}|f(\xx)|^2\ln\left(\frac{1+|\xx|^2}{2}\right)d\xx+\int_{\mathbb{R}}\int_{\mathbb{R}}\ln|\uu||(\mathcal{V}^M_gf)(\xx,\uu)|^2d\xx d\uu\geq \left(\frac{\Gamma'(\frac{1}{2})}{\Gamma(\frac{1}{2})}+\ln|\Mb|\right)\|g\|^2_{L^2(\mathbb{R})}\int_{\mathbb{R}}|f(\xx)|^2d\xx
\end{align}
for the window function $g\in L^2(\mathbb{R})$ and $f\in H^1(\mathbb{R})\cap L_1^2(\mathbb{R}),$ where the Sobolev space $H^1(\mathbb{R})=\{f\in L^2(\mathbb{R}):Df\in L^2(\mathbb{R})\},$ $D$ stands for differential operator and the weighted Lebesgue space $L_s^p(\mathbb{R})=\{f\in L_{loc}^p(\mathbb{R}):(1+|x|^2)^{\frac{s}{2}}f\in L^p(\mathbb{R})\},~s>0.$ Using the property that $\ln$ function is an increasing function, equation \eqref{R2022_Bahri_IJWMIP_WLCT_Logarithmic_UP} can be obtained from \eqref{R2022_BahriJoIA__WLCT_Logarithmic-Sobolev-type-UP}. Using the method adopted by Wilczok \cite{wilczok2000new}, Hleili \cite{hleili2022windowed} proved that $\mathcal{V}^M_gf$ cannot have a support of finite measure. The author also proved that the portion of $\mathcal{V}^M_gf$ lying outside some set $S\subset \mathbb{R}^2$ cannot be arbitrarily small. For sufficiently small $S,$ this fact is justified by the following theorem 
\begin{theorem}
(Concentration of $\mathcal{V}^M_gf$ on a small set) Let $g$ be a window function and $S\subset\mathbb{R}^2$ with $\frac{\|g\|^2_{L^2(\mathbb{R})}}{2\pi|\Mb|}\lambda^{(2)}(S)<1.$ Then, for $f\in L^2(\mathbb{R})$ arbitrary
\begin{align}
\|\mathcal{V}^M_gf-\chi_S\cdot \mathcal{V}^M_gf\|_{L^2(\mathbb{R}^2)}\geq \left(\sqrt{1-\frac{\|g\|^2_{L^2(\mathbb{R})}}{2\pi|\Mb|}\lambda^{(2)}(S)}\right)\|f\|_{L^2(\mathbb{R})}\|g\|_{L^2(\mathbb{R})}.
\end{align}
\end{theorem}
\noindent Apart from this, the author also studied the UP for the orthonormal sequence. They proved that if the WLCT of each member of an orthonormal sequence is $\epsilon-$concentrated in a set of finite measure then the sequence is necessarily finite. More precisely, if $s>0,0<\epsilon<\frac{1}{\|g\|_{L^2(\mathbb{R})}},$ $G_s=\{(\xx,\uu)\in\mathbb{R}^2:\xx^2+\uu^2\leq s^2\},$ $\{\phi_{n}\}_{1\leq n\leq N}\subset L^2(\mathbb{R})$ is an orthonormal sequence and $\mathcal{V}^M_g\phi_n$ is $\epsilon-$concentrated in $G_s$ then $N$ is bounded above by $\frac{s^2\|\psi\|^2_{L^2(\mathbb{R})}}{2|\Mb|\left(1-\epsilon\|g\|_{L^2(\mathbb{R})}\right)}.$ As a consequence  the Shapiro's dispersion theorem for the WLCT is obtained, which states that if $\{\phi_{n}\}_{1\leq n\leq N}$ is an orthonormal sequence in $L^2(\mathbb{R}),$ then for every $p>0$ and $N\in \mathbb{N},$ 
$$\sum_{n=1}^N\left\{\rho_{p}\left(\mathcal{V}^M_g\phi_n\right)\right\}^p\geq N^{\frac{p}{2}+1}\left(\frac{3|\Mb|}{2^{\frac{4}{p}}\|g\|^2_{L^2(\mathbb{R})}}\right)^{\frac{p}{2}},$$
where $\rho_{p}\left(\mathcal{V}^M_gf\right)=\left(\int_{\mathbb{R}\times\mathbb{R}}|(\xx,\uu)|^p|(\mathcal{V}^M_gf)(\xx,\uu)|^2d\xx d\uu\right)^{\frac{1}{p}},~|(\xx,\uu)|=\sqrt{\xx^2+\uu^2},$ is the generalized $p^{th}$ time-frequency dispersion of $\mathcal{V}^M_gf.$

\section{Linear canonical wavelet transform (LCWT)}
Compared to the FrFT with one extra degree of freedom, LCT is more flexible for its extra degree of three degrees of freedom and has been used frequently in time-frequency analysis and non-stationary signal processing. Inspired by the FrWT, Wei et al. \cite{wei2014generalized} introduced a new transform by combining the idea of LCT and the WT, called the LCWT. This definition is obtained by replacing the classical convolution in \eqref{R2000_Wilczok_DefinitionWaveletTransform_using_InnerProduct} by the linear canonical convolution \eqref{R2011_Wei-Ran-Li-LCT-Convolution} and is defined for a function $f\in L^2(\mathbb{R})$ and the function $\psi\in L^2(\mathbb{R})$ satisfying the admissibility condition $c_\psi:=\int_{\mathbb{R}}\frac{|\hat{\psi}(\xxi)|^2}{\xxi}d\xxi=\int_{\mathbb{R}}\frac{|\hat{\psi}(-\xxi)|^2}{\xxi}d\xxi<\infty,$ where $\hat{\psi}(\xxi)=\frac{1}{\sqrt{2\pi}}\int_{\mathbb{R}}f(\tee)e^{-i\tee\xxi}d\tee,$ as 
%\textcolor{red}{(check the condition and the definition of FT used, we have used a nice one i.e., $\frac{1}{\sqrt{2\pi}}$ from \cite{guo2018linear})} 
\begin{align}\label{R2014_Wei-Li_Optik_LCWT_Definition_COnvolution}
(\mathcal{W}^M_{\psi}f)(\A,\B)=&\left(f\star'_M\delta_{\A}\tilde{\psi}\right)(\B)\notag\\
=&\int_{\mathbb{R}}f(\tee)\overline{\psi^M_{\A,\B}(\tee)}d\tee,
\end{align}
where $\psi^M_{\A,\B}(\tee)=e^{-i\frac{\Ma(\tee^2-\B^2)}{2\Mb}}\psi_{\A,\B}(\tee),~\A\in\mathbb{R}^+,\B\in\mathbb{R}.$ In particular, for $M=(\cos\alpha,\sin\alpha,-\sin\alpha,\cos\alpha),$ it reduces to the FrWT \eqref{R2012_Shi-Na-Liu_FrWT_Definition}. The LCWT \eqref{R2014_Wei-Li_Optik_LCWT_Definition_COnvolution} can be expressed as $(\mathcal{W}^M_{\psi}f)(\A,\B)=e^{-i\frac{\Ma\B^2}{2\Mb}}\left(\mathcal{W}_{\psi}[f(\tee)e^{-\frac{\Ma\B^2}{2\Mb}}]\right)(\A,\B),$ in term of the classical WT. The transform is linear and satisfy the following convolution formula $(\mathcal{W}^M_{\psi}(f\star_Mg))(\A,\B)=(\mathcal{W}^M_{\psi}f)(\A,\B)(\mathcal{W}_{\psi}g)(\A,\B),$ for $f\in L^2(\mathbb{R})$ and $g\in L^1(\mathbb{R})\cap L^2(\mathbb{R}).$ Based on the convolution theorem \eqref{R2011_Wei-Ran-Li-LCT-Convolution-Theorem} and the inverse of LCT, the transform \eqref{R2014_Wei-Li_Optik_LCWT_Definition_COnvolution} can be expressed as 
\begin{align}\label{R2014_Wei-Li_Optik_LCWT_LCD}
(\mathcal{W}^M_{\psi}f)(\A,\B)=\sqrt{2\pi\A}\int_{\mathbb{R}}\left(\mathcal{L}^Mf\right)(\xxi)\overline{\hat{\psi}\left(\frac{\A\xxi}{\Mb}\right)}K_{M^{-1}}(\xxi,\B)d\xxi.
\end{align}
As a consequence of which the authors in \cite{wei2014generalized},\cite{guo2018linear} obtained the inner product theorem 
\begin{align}\label{R2014_Wei-Li_Optik_LCWT_IPR}
\int_{\mathbb{R}^+}\int_{\mathbb{R}}(\mathcal{W}^M_{\psi}f)(\A,\B)\overline{(\mathcal{W}^M_{\psi}g)(\A,\B)}\frac{d\A d\B}{\A^2}=2\pi|\Mb| c_\psi\langle f,g\rangle,
\end{align}
followed by the inversion formula 
$$f(\tee)=\frac{1}{2\pi|\Mb|c_\psi}\int_{\mathbb{R}^+}\int_{\mathbb{R}}(\mathcal{W}^M_{\psi}f)(\A,\B)\psi^M_{\A,\B}(\tee)\frac{d\A d\B}{\A^2}.$$
given in \cite{guo2018linear}. In view of \eqref{R2014_Wei-Li_Optik_LCWT_Definition_COnvolution} and \eqref{R2014_Wei-Li_Optik_LCWT_LCD}, it follows that the linear canonical wavelet transform is capable of providing the time and spectral information simultaneous and representing signal in the time-fractional-frequency domain. It gives local information of the signal $f(\tee)$ with the time window $\left[\B+\A E_\psi-\A\Delta_\psi,\B+\A E_\psi+\A \Delta_\psi\right]$ and the local information about the linear canonical spectrum of $f(\tee)$ as $\left[\frac{E_{\hat{\psi}}}{\A}-\frac{\Delta_{\hat{\psi}}}{\A},\frac{E_{\hat{\psi}}}{\A}+\frac{\Delta_{\hat{\psi}}}{\A}\right]\Mb,$ with the constant window area $4\Mb\Delta_\psi\Delta_{\hat{\psi}}$ in the time-LCT-frequency plane. This area depends on the function $\psi$ and $\Mb$ and is independent of the parameters $\A$ and $\B.$ Thus, the LCWT \eqref{R2014_Wei-Li_Optik_LCWT_Definition_COnvolution} is capable of providing time and linear canonical domain information simultaneously and representing the signal in the time-LCT-frequency plane with the transform parameter $M.$ Prasad et al. \cite{prasad2019continuous} explored further the LCWT \eqref{R2014_Wei-Li_Optik_LCWT_Definition_COnvolution}, with the alternative admissible condition of $\psi$ as $C^M_\psi=\int_{\mathbb{R}^+}\frac{\left|\left(\mathcal{L}^M[e^{-\frac{\i\A}{2\B}(\cdot)}\psi]\right)(\xxi)\right|^2}{|\xxi|}d\xxi<\infty.$ They studied the properties like dilation property, parity and also obtained the formula $(\mathcal{W}^M_\psi f)(\lambda\A,\lambda\B)=\lambda^{n+\frac{1}{2}}(\mathcal{W}^{M'}_\psi f)(\A,\B),~M'=(\lambda^2\Ma,\Mb,\Mc,\Md)$ when $\psi$ is a homogeneous function of degree $n.$ They also generalized the inner product relation \eqref{R2014_Wei-Li_Optik_LCWT_IPR} in the context of two linear canonical wavelet as 
\begin{align*}
\int_{\mathbb{R}^+}\int_{\mathbb{R}}(\mathcal{W}^M_{\psi_1}f)(\A,\B)\overline{(\mathcal{W}^M_{\psi_2}g)(\A,\B)}\frac{d\A d\B}{\A^2}=2\pi|\Mb| c^M_{\psi_1,\psi_2}\langle f,g\rangle,
\end{align*}
where $c^M_{\psi_1,\psi_2}=\int_{\mathbb{R}^+}\overline{\left(\mathcal{L}^M[e^{-\frac{\i\Ma}{2\Mb}(\cdot)}\psi_1]\right)(\A)}\left(\mathcal{L}^M[e^{-\frac{\i\Ma}{2\Mb}(\cdot)}\psi_2]\right)(\A)\frac{d\A}{\A}<\infty.$ With the assumption that $c^M_{\psi_1,\psi_2}=1,$ Ansari et al. \cite{ansari2021abelian} obtained the Calde\`ron's reproducing formula given as
\begin{align*}
f(\tee)=\frac{1}{2\pi|B|}\int_{\mathbb{R}^+}\left(\left(f(\cdot)\star_M(\psi_1)_{\A}(\cdot)\right)\star_M\overline{(\psi_2)}(-(\cdot))\right)(\tee)\frac{d\A}{\A^2},
\end{align*}
where $\psi_{\A}(\tee)=\frac{1}{\sqrt{\A}}\overline{\psi}\left(\frac{-\tee}{\A}\right).$ Authors also established the Abelian theorem, which relates the behavior of LCWT $(\mathcal{W}^M_\psi f)(\A,\B)$ as $\A\rightarrow \infty (\A\rightarrow 0)$ to the behavior of function $(\mathcal{L}^M f)(\xxi)$ as $\xxi\rightarrow \infty (\xxi\rightarrow 0).$ These results are also referred to as initial value and final value Abelian theorems. Guo et al. \cite{guo2018linear} considered the test function space $\mathcal{S}_M(\mathbb{R})$ defined as 
\begin{align}\label{RGuo-Generalized-Sobolev-Space(3)_LCDomain}
\mathcal{S}_M(\mathbb{R})=\{f:f\in C^\infty,\sup_{\tee\in\mathbb{R}}|\tee^\alpha\Delta_\tee^\beta f(\tee)|<\infty,~\alpha,\beta\in \mathbb{N}\},
\end{align}
where $\Delta_\tee=\frac{d}{d\tee}-i\tee\frac{\Ma}{\Mb}.$ It is equipped with the topology generated by the collection of seminorms $\gamma^M_{\alpha,\beta}$ defined by $\gamma^M_{\alpha,\beta}(f)=\sup_{\tee\in\mathbb{R}}|\tee^\alpha\Delta_\tee^\beta f(\tee)|.$ It is a Fr\`echet space and reduces to the Schwartz space $\mathcal{S}(\mathbb{R}),$ when $M=(0,1,-1,0)$ and proved that the LCT is a homeomorphism from  $\mathcal{S}_M(\mathbb{R})$ onto itself. As a subset of the dual $\mathcal{S}'_M(\mathbb{R})$ of $\mathcal{S}_M(\mathbb{R})$ authors in \cite{guo2018linear} also considered a generalized Sobolev space $H^{p,s}_M$ defined by 
$$H^{p,s}_M=\{f:f\in\mathcal{S}'_M(\mathbb{R}),\xxi^{\frac{s}{p}}(\mathcal{L}^Mf)(\xxi)\in L^p(\mathbb{R})\},~1\leq p<\infty,~s\in \mathbb{R},$$ 
with norm given as $\|f\|_{H^{p,s}_M}=\|\xxi^{\frac{s}{p}}(\mathcal{L}^Mf)(\xxi)\|_{L^p}$ and proved that the operator $\mathcal{W}^M_\psi$ \eqref{R2014_Wei-Li_Optik_LCWT_Definition_COnvolution} is a bounded linear operator from $H^{p,s}_M$ into $H^{p,s}_{q,M}(\mathbb{R}\times\mathbb{R}^+)$  for $1\leq p\leq 2,~\frac{1}{p}+\frac{1}{q}=1,$ with $\psi$ satisfying $c^{s,q}_{\psi,M}:=\int_{\mathbb{R}^+}\frac{|\mathcal{L}^M[e^{-\frac{-i\Ma\xxi^2}{2\Mb}}\psi(\xxi)](\A \xxi)|^q}{\A|\xxi|^s}\frac{d\A}{\A}$ is a constant independent of $\xxi,$ where $H^{p,s}_{q,M}(\mathbb{R}\times\mathbb{R}^+)$  is the set of all measurable functions $\varphi(\A,\B)$ satisfying $\|\varphi(\A,\B)\|_{H^{p,s}_{q,M}}=\left[\int_{\mathbb{R}^+}\left(|\varphi(\A,\B)|^p d\B\right)^{\frac{q}{p}}\A^{-s-1}d\A\right]^{\frac{1}{q}}<\infty.$. With the notation $\mathfrak{M}_c$ and $K_\eta$ given in definition \ref{RGeneralized_Ultra-distribution_LCDomain}, Prasad et al. \cite{prasad2019approximation} considered the space $S^M_\eta,$ the collection of $\phi\in L^1(\mathbb{R})$ such that  $\phi$ and $\hat{\phi}$ are smooth functions and for every $l\in\mathbb{N}_0$ and $\lambda\geq 0,$ $p_{l,\lambda,M}(\phi)=\sup_{t\in\mathbb{R}}e^{\lambda\eta(t)}|\Delta_t^l\phi(t)|<\infty$ and $r_{l,\lambda,M}(\phi)=\sup_{\xxi\in\mathbb{R}}e^{\lambda\eta(\xxi)}|\Delta_\xxi^l\hat{\phi}(\xxi)|<\infty,$ where $\Delta_\tee=\frac{d}{d\tee}-i\frac{\Ma}{\Mb}\tee.$ The topology on $S^M_\eta$ is defined by the semi norms $p_{l,\lambda,M}(\phi)$ and $r_{l,\lambda,M}(\phi).$ The elements of the dual $(S^M_\eta)'$ are called the generalized ultra-distributions. Note that for $\eta(\xxi)=\log(1+|\xxi|),M=(0,1,-1,0)$ and $\eta(\xxi)=\xxi^\sigma(0<\sigma<1),M=(0,1,-1,0)$ the space $S^M_\eta$ reduces to the Schwartz space and the Gevrey space of test functions, respectively. Author also considered generalized Sobolev space $B^{\eta,M}_{p,k}(\mathbb{R}),$ defined as the space of all the generalized ultra-distribution $f$ in $(S^M_\eta)'$ such that  $\|f\|_{p,k,M}=\left(\int_{\mathbb{R}}|k(\xxi)(\mathcal{L}^M f)(\xxi)|^pd\xxi\right)^{\frac{1}{p}}<\infty$ and $\|f\|_{\infty,k,M}=ess\sup|k(\xxi)(\mathcal{L}^M f)(\xxi)|.$ The space $B^{\eta,M}_{p,k}(\mathbb{R})$ reduces to the  H\"ormander space $B_{p,k}(\mathbb{R}),$ when $\eta(\xxi)=\log(1+|\xxi|)$ and $M=(0,1,-1,0).$ They determine the $B^{\eta,M}_{p,k}(\mathbb{R})-$distance of two LCWT with different argument functions and the wavelet to study the dependence of the transform on its wavelets and its argument and determining the LCWT in terms of linear canonical convolution they proves that if $f\in B^{\eta,M}_{p,k}(\mathbb{R})$ and $\psi\in L^1(\mathbb{R})$ with $\int_{\mathbb{R}}\psi(\tee)d\tee=1,$ then $\psi_{\A}\star_M f\rightarrow f$ in $B^{\eta,M}_{p,k}(\mathbb{R})$  as $\A\rightarrow 0.$ Apart from this authors in \cite{prasad2019approximation} defined the convolution of LCWTs $(\mathcal{W}^M_{\psi_1}f)(\A,\B)$ and $(\mathcal{W}^M_{\psi_1}f)(\A',\B')$ for $\psi_1\in L^2(\mathbb{R}),$  $\psi_2\in L^1(\mathbb{R})$ denoted by $\left[\left(\mathcal{W}^M_{\psi_1}\circ \mathcal{W}^M_{\psi_2}\right)f\right](\A,\A',\B)=\left(\mathcal{W}^M f\right)(\A,\A',\B)$ is defined by
\begin{align}\label{R2019_Prasad-Ansari_JoPDOA_LCWT_Convolution-LCWT}
\left(\mathcal{W}^M f\right)(\A,\A',\B)=\left(\mathcal{W}^M_{\psi_1}\left(\mathcal{W}^M_{\psi_2}\right)(\A',\B')\right)(\A,\A',\B),
\end{align}
which can equivalently written as
$$\left(\mathcal{W}^M f\right)(\A,\A',\B)=\sqrt{\frac{2\pi\Mb \A \A'}{i}}\int_{\mathbb{R}}e^{-i\frac{\Md}{2\Mb}\xxi^2+\frac{i}{\Mb}\B\xxi-i\frac{\Ma}{2\Mb}\B^2+i\frac{\Md}{2\Mb}(\A^2+\A'^2)\xxi^2}\overline{\mathcal{L}^M[e^{-i\frac{\Ma}{2\Mb}(\cdot)^2}\psi_1](\A\xxi)}\overline{\mathcal{L}^M[e^{-i\frac{\Ma}{2\Mb}(\cdot)^2}\psi_2](\A'\xxi)}(\mathcal{L}^M_{\B'}f)(\xxi)d\xxi,$$
where $\mathcal{L}^M_{\B'}$ denoted the LCWT with respect to the variable $\B'.$ Authors obtained several important properties like translation, parity, dilation, etc for the composition \eqref{R2019_Prasad-Ansari_JoPDOA_LCWT_Convolution-LCWT} of the LCWTs along with the inner product relation 
\begin{align*}
\int_{\mathbb{R}}\int_{\mathbb{R}^+}\int_{\mathbb{R}^+}\left(\mathcal{W}^M f\right)(\A,\A',\B)\overline{\left(\mathcal{W}^M g\right)(\A,\A',\B)}\frac{d\A d\A' d\B}{\A \A'}=4\pi^2\Mb c^M_{\psi_1,\psi_2}\langle f,g\rangle,
\end{align*}
with the assumption that $\Ma=\Mb$ and $c^M_{\psi_1,\psi_2}:=\int_{\mathbb{R}^+}\int_{\mathbb{R}^+}|\mathcal{L}^M[e^{-i\frac{\Ma}{2\Mb}(\cdot)^2}\psi_1](\A\xxi)|^2|\mathcal{L}^M[e^{-i\frac{\Ma}{2\Mb}(\cdot)^2}\psi_2](\A'\xxi)|^2\frac{d\A d\A'}{\A \A'}<\infty.$ Furthermore, if $f\in L^2(\mathbb{R})$ is a continuous function at $\tee\in\mathbb{R},$ then $f$ can be reconstructed by the formula 
$$f(\tee)=\frac{1}{4\pi^2\Mb^2c^M_{\psi_1,\psi_2}}\int_{\mathbb{R}}\int_{\mathbb{R}^+}\int_{\mathbb{R}^+}\left(\mathcal{W}^M f\right)(\A,\A',\B)\overline{\left(\mathcal{W}^M (\tau_{\Mb}\delta)\right)(\A,\A',\B)}\frac{d\A d\A' d\B}{\A \A'},$$
where $\tau_{\Mb}\delta$ denote the translation of Dirac delta function. Prasad et al. \cite{prasad2020composition} considered that the generalized Sobolev space $H^{p,s}_{q,M}(\mathbb{R}\times\mathbb{R}^+\times\mathbb{R}^+)$  is the set of all measurable functions $\varphi(\A,\A',\B)$  satisfying $\|\varphi(\A,\A',\B)\|_{H^{p,s}_{q,M}}=\left[\int_{\mathbb{R}^+}\left(|\varphi(\A,\A',\B)|^p d\B\right)^{\frac{q}{p}}(\A \A')^{-s-1}d\A d\A'\right]^{\frac{1}{q}}<\infty,$ where $1\leq p,~q<\infty,~s\in\mathbb{R}$ and proved that the composition operator $\mathcal{W}^M$ \eqref{R2019_Prasad-Ansari_JoPDOA_LCWT_Convolution-LCWT} is a bounded linear operator from $H^{p,s}_M$ into 
$H^{p,s}_{q,M}(\mathbb{R}\times\mathbb{R}^+\times\mathbb{R}^+)$ for $1\leq p\leq 2,~q=\frac{p}{p-1},$ with the linear canonical wavelet $\psi_1$ and $\psi_2$ satisfying $c^{s,q}_{\psi_1,\psi_2,M}:=\int_{\mathbb{R}^+}\int_{\mathbb{R}^+}\frac{\left|\left(\mathcal{L}^M\left[e^{-\frac{i\Ma\tee^2}{2\Mb}}\psi_1(\tee)\right]\right)(\A \xxi)\left(\mathcal{L}^M\left[e^{-\frac{i\Ma\tee^2}{2\Mb}}\psi_2(\tee)\right]\right)(\A' \xxi)\right|^q}{\left(\A\A'|\xxi|\right)^s}\frac{d\A d\A'}{\A \A'}$ is independent of $\xxi$ where $H^{p,s}_M$ is defined in \eqref{RGuo-Generalized-Sobolev-Space(3)_LCDomain}. Gupta et al. \cite{gupta2024new} introduced a novel LCWT for a function $f\in L^2(\mathbb{R})$ with respect to the linear canonical wavelet $\psi,$ i.e., $\psi\in L^2(\mathbb{R})$ satisfying the admissible condition $c^M_\psi:=\int_{\mathbb{R}^+}\left|(\mathcal{L}^M\psi)\left(\frac{\xxi}{\A}\right)\right|^2\frac{d\A}{\A}$ is a positive constant independent of $\xxi$ satisfying $|\xxi|=1,$ as 
\begin{align}\label{RGupta_LCWT-Definition}
(\mathcal{W}^M_\psi f)(\A,\B)=&e^{-\frac{i\Ma}{2\Mb}\tee}\left\{f(\tee)e^{\frac{i\Ma}{2\Mb}\tee}\star\overline{\sqrt{\A}\psi(-\A\tee)e^{-]\frac{i\Ma}{2\Mb}\tee}}\right\}(\B)\notag\\
=&\int_{\mathbb{R}}f(\tee)\overline{\psi^M_{\A,\B}(\tee)}d\tee,
\end{align}
where $\psi^M_{\A,\B}(\tee)=e^{-\frac{i\Ma}{2\Mb}\{(\tee^2-\B^2)-(\A(\tee-\B))^2\}}\sqrt{\A}\psi(\A(\tee-\B))$
which can equivalently be written in the linear canonical domain as  
\begin{align}\label{RGupta_LCWT-Definition_LCDomain}
(\mathcal{W}^M_\psi f)(\A,\B)=\frac{\sqrt{-2\pi i\Mb}}{\sqrt{\A}}\int_{\mathbb{R}}e^{-\frac{i\Md}{2\Mb}\left(\frac{\xxi}{\A}\right)^2}(\mathcal{L}^Mf)(\xxi)\overline{(\mathcal{L}^M\psi)(\frac{\xxi}{\A})}K_{M^{-1}}(\B,\xxi)d\xxi.
\end{align}
Apart from the basic properties of the  \eqref{RGupta_LCWT-Definition} like linearity, anti-linearity and the for the translation and dilation formula of the transform author obtain the inner product relation $\langle \mathcal{W}^M_{\psi_1}f,\mathcal{W}^M_{\psi_2}g\rangle_{L^2(\mathbb{R}^+\times\mathbb{R})}=2\pi c^M_{\psi_1,\psi_2}\langle f,g\rangle_{L^2(\mathbb{R})},$ followed by the inversion formula $f(\tee)=\frac{1}{2\pi|\Mb|c^M_{\psi_1,\psi_2}}\int_{\mathbb{R}^+\times\mathbb{R}}(\mathcal{W}^M_{\psi_1}f)(\A,\B)(\psi_2)^M_{\A,\B}(\tee)d\A d\B,$ a.e. $\tee\in\mathbb{R},$ where $c^M_{\psi_1,\psi_2}:=\int_{\mathbb{R}^+}\overline{(\mathcal{L}^M\psi_1)\left(\frac{\xxi}{\A}\right)}(\mathcal{L}^M\psi_2)\left(\frac{\xxi}{\A}\right)\frac{d\A}{\A}$ is a non-zero complex constant independent of $\xxi$ satisfying $|\xxi|=1.$ 
In view of \eqref{RGupta_LCWT-Definition} and \eqref{RGupta_LCWT-Definition_LCDomain}, it follows that the LCWT \eqref{RGupta_LCWT-Definition} localizes the signal $f$ in the time window 
$\left[\frac{1}{\A}E_\psi+\B-\frac{1}{\A}\Delta_
\psi,\frac{1}{\A}E_\psi+\B+\frac{1}{\A}\Delta_
\psi\right]$ and the linear canonical domain window $\left[\A E_{\mathcal{L}^M\psi}-\A\Delta_{\mathcal{L}^M\psi},\A E_{\mathcal{L}^M\psi}+\A\Delta_{\mathcal{L}^M\psi}\right],$ respectively, with the joint resolution in the time and linear canonical domain  window as
$$\left[\frac{1}{\A}E_\psi+\B-\frac{1}{\A}\Delta_
\psi,\frac{1}{\A}E_\psi+\B+\frac{1}{\A}\Delta_
\psi\right]\times\left[\A E_{\mathcal{L}^M\psi}-\A\Delta_{\mathcal{L}^M\psi},\A E_{\mathcal{L}^M\psi}+\A\Delta_{\mathcal{L}^M\psi}\right],$$
with constant area $4\Delta_{\psi}\Delta_{\mathcal{L}^M\psi}$ in the time-LCT-frequency plane. Thus, for a given parameter $M,$ the window area depends on the linear canonical admissible wavelets and is independent of the parameters $\A$ and $\B.$ But it is noted that the window gets narrower for a large value of $\A$ and wider for a small value of $\A.$ Thus, the window given by the transform is flexible, and hence, it can simultaneously provide the time-linear canonical domain information.\\

\noindent\textbf{Uncertainty Principles:} Bahri et al. \cite{bahri2022novel} studied the UPs associated with the LCWT \eqref{R2014_Wei-Li_Optik_LCWT_Definition_COnvolution}. Adopting the method in \cite{wilczok2000new}, they obtained the generalized Heisenberg and logarithmic UP for the transform using the same for the LCT given, respectively as 
\begin{align*}
\left(\int_{\mathbb{R}^+}\int_{\mathbb{R}}\B^p\left|(\mathcal{W}^M_\psi f)(\A,\B)\right|^pd\B\frac{d\A}{\A^2}\right)^{\frac{1}{p}}\left(\int_{\mathbb{R}}\xxi^p\left|(\mathcal{L}^Mf)(\xxi)\right|^pd\xxi\right)^{\frac{1}{p}}\geq \Mb^{1+\frac{1}{p}}\sqrt{\frac{\pi c^M_\psi}{2}}\|f\|^2_{L^2(\mathbb{R})}
\end{align*}
and 
\begin{align*}
\int_{\mathbb{R}^+}\int_{\mathbb{R}}\ln|\B|\left|(\mathcal{W}^M_\psi f)(\A,\B)\right|^2d\B\frac{d\A}{\A^2}+2\pi\Mb c^M_\psi\int_{\mathbb{R}}\ln|\xxi|\left|(\mathcal{L}^Mf)(\xxi)\right|^2d\xxi\geq (D+\ln|\Mb|)2\pi\Mb c^M_\psi\|f\|^2_{L^2(\mathbb{R})}.
\end{align*}
followed by the local UP for which states that if $E\subset \mathbb{R}^+\times\mathbb{R}$ and $B_{r_0}=\{(\A,\B)\in\mathbb{R}^+\times\mathbb{R}:|(\A,\B)|<r_0\},$ $r_0\in (0,1]$ such that $0<\frac{m(E)}{2\pi\Mb c^M_\psi},\frac{m(B_{r_0})}{2\pi\Mb c^M_\psi}<1$ and $\|\psi\|_{L^2(\mathbb{R})}=1$ then for all $f\in L^2(\mathbb{R})$ and $l>0$ then 
$$\|f\|_{L^2(\mathbb{R})}\leq \frac{1}{\sqrt{2\pi\Mb c^M_\psi-m(E)}}\left(\int\int_{E^c}\left|(\mathcal{W}^M_\psi f)(\A,\B)\right|^2d\B\frac{d\A}{\A^2}\right)^\frac{1}{2}$$ 
and 
$$\|f\|_{L^2(\mathbb{R})}\leq \frac{1}{r_0^l\sqrt{2\pi\Mb c^M_\psi-m(B_{r_0})}}\left(\int_{\mathbb{R}^+}\int_{\mathbb{R}}|(\A,\B)|^{2l}\left|(\mathcal{W}^M_\psi f)(\A,\B)\right|^2d\B\frac{d\A}{\A^2}\right)^\frac{1}{2}$$
where $m(E)$ represents the Lebesgue measure of $E$ and $E^c=\mathbb{R}^+\times\mathbb{R}\setminus E.$ Authors in \cite{gupta2024new} obtained the Donoho-Stark's UP for the LCWT \eqref{RGupta_LCWT-Definition}, which gives a lower bound of the measure of essential support of the $\mathcal{W}^M_\psi f,$, i.e., if  $\mathcal{W}^M_\psi f$ is $\epsilon-$ concentrated on a measurable set $\Omega\subset\mathbb{R}^+\times\mathbb{R}$, then 
$|\Omega|\|\psi\|^2_{L^2(\mathbb{R})}\geq 2\pi|\Mb|c^M_\psi(1-\epsilon^2),~0\leq \epsilon<1.$ Using the method adopted by Wilczok \cite{wilczok2000new}, they proved that $\mathcal{W}^M_\psi f$ cannot have a support of finite measure. The author also proved that the portion of $\mathcal{W}^M_\psi f$ \eqref{RGupta_LCWT-Definition} lying outside some set $\Omega\subset \mathbb{R}^2$ cannot be arbitrarily small. For sufficiently small $\Omega,$ this fact is justified by the following theorem 
\begin{theorem}(Concentration of $\mathcal{W}^M_\psi f$ on a small set)
Let $f\in L^2(\mathbb{R}),$ $\psi$ is an ALCW and $\Omega\subset\mathbb{R}^+\times\mathbb{R}$ such that  $|\Omega|<\frac{2\pi |\Mb|C_{\psi,M}}{\|\psi\|^2_{L^2(\mathbb{R})}}.$ Then
$$\|\chi_{\Omega^c}W_{\psi}^Mf\|_{L^2(\mathbb{R}^+\times\mathbb{R})}\geq\sqrt{2\pi|\Mb|C^M_\psi-|\Omega|\|\psi\|^2_{L^2(\mathbb{R})}}\|f\|_{L^2(\mathbb{R})}.$$
\end{theorem}
\noindent Apart from this, the author also studied the UP for the orthonormal sequence. They proved that if the LCWT \eqref{RGupta_LCWT-Definition} of each member of an orthonormal sequence is $\epsilon-$concentrated in a set of finite measure, then the sequence is necessarily finite. More precisely, if $s>0,0<\epsilon<1,$ $G_s=\{(\A,\B)\in\mathbb{R}^+\times\mathbb{R}:\A^2+\B^2\leq s^2\},$ $\{\phi_{n}\}_{n\in\wedge}\subset L^2(\mathbb{R}),~\wedge\subset\mathbb{N},$ is an orthonormal sequence and $\mathcal{W}^M_\psi\left(\frac{\phi_n}{\sqrt{2\pi|\Mb|c^M_\psi}}\right)$ is $\epsilon-$concentrated in $G_s$ for all $n\in\wedge,$ then the cardinality of $\wedge$ is bounded above by $\frac{s^2\|\psi\|^2_{L^2(\mathbb{R})}}{4|\Mb|c^M_\psi\left(1-\epsilon\right)}.$ As a consequence of which the Shapiro's dispersion theorem for the LCWT \eqref{RGupta_LCWT-Definition} is obtained, which states that if $\{\phi_{n}\}_{n\in\mathbb{N}}$ is an orthonormal sequence in $L^2(\mathbb{R}),$ then for every $p>0$ and non-empty finite $\wedge\subset\mathbb{N},$ 
$$\sum_{n\in\wedge}\left\{\rho_{p}\left(W_{\psi}^M\left(\frac{\phi_{n}}{\sqrt{2\pi|B|C_{\psi,M}}}\right)\right)\right\}^p\geq \frac{(Card(\wedge))^{\frac{p}{2}+1}}{2^{p+1}}\left(\frac{3|\Mb|C_{\psi,M}}{2^{\frac{4}{p}+2}\|\psi\|^2_{L^2(\mathbb{R})}}\right)^{\frac{p}{2}}.$$
where $\rho_{p}\left(\mathcal{W}^M_\psi f\right)=\left(\int_{\mathbb{R}^+\times\mathbb{R}}|(\A,\B)|^p|(\mathcal{W}^M_\psi f)(\A,\B)|^2d\A d\B\right)^{\frac{1}{p}},~|(\A,\B)|=\sqrt{\A^2+\B^2},$ is the generalized $p^{th}$ time-frequency dispersion of $\mathcal{W}^M_\psi f.$

\section{Linear Canonical Stockwell Transform (LCST)}
Bahri et al. \cite{bahri2019generalized} replaced the Fourier kernel in \eqref{R2009_Molahajloo-Wong_DefnitionST_GeneralWindow} by the kernel $K^M$ of the LCT to define the LCST for the  function $f\in L^2(\mathbb{R})$ with respect to the window function $\psi\in L^2(\mathbb{R})$ as 
\begin{align}\label{R2019_Bahri-Toaha-Lande_LCST_Definition}
(\mathcal{S}^M_\psi f)(\A,\B)=|\A|\int_{\mathbb{R}}K^M(\xx,\A)f(\xx)\overline{\psi(\A(\xx-\B))}d\xx,~\xx\in\mathbb{R},~\A\in\mathbb{R}_{0}.
\end{align}
and obtained the formula for the translation and modulation of the function, including the inner product relation and the reconstruction formula. Later, Bahri et al. \cite{bahri2019two} established the Heisenberg uncertainty and logarithmic UPs for the transform \eqref{R2019_Bahri-Toaha-Lande_LCST_Definition}. Replacing the classical convolution in \eqref{R2009_Molahajloo-Wong_DefnitionST_GeneralWindow_Convolution} by the linear canonical convolution \eqref{RCovolutionLCT2}, Shah et al.\cite{shah2020linear} defined the LCST of  a function $f\in L^2(\mathbb{R})$ with respect to the window function $\psi\in L^2(\mathbb{R})$ as 
\begin{align}\label{R2020_Shah-Tantary_LCST_Definition_Convolution}
(\mathcal{S}^M_{\psi}f)(\A,\B)=\frac{1}{\sqrt{2\pi}}\left(\mathcal{M}_{-\A}f\star_M\mathcal{D}_{\xi}\tilde{\psi}\right)(\B),
\end{align}
where $\tilde{\psi}(t)=\overline{\psi(-t)}$ and $\star_M$ is the convolution given by \eqref{RCovolutionLCT2}. It is  a matter of simple verification that $\mathcal{S}_{\psi}f$ is linear in $f$ and conjugate linear in $\psi.$ Transform \eqref{R2020_Shah-Tantary_LCST_Definition_Convolution} has an equivalent form 
\begin{align}\label{R2020_Shah-Tantary_LCST_Definition_Integral-Form}
(\mathcal{S}^M_{\psi}f)(\A,\B)=\frac{1}{\sqrt{2\pi}}\int_{\mathbb{R}}f(\tee)\overline{\psi^M_{\A,\B}(\tee)}d\tee,
\end{align}
where $\psi^M_{\A,\B}(\tee)=\A \psi\left(\A(\tee-\B)\right)e^{i\A\tee+\frac{i\Ma\tee (\B-\tee)}{\Mb}},$ $\A\in\mathbb{R}^+,~\B\in\mathbb{R}.$ It reduces to the FrST for $M=(\cos\alpha,\sin\alpha,-\sin\alpha,\cos\alpha),$ $\alpha\neq n\pi,~n\in\mathbb{Z}$, Fresnel-Stockwell transform for $M=(1,\Mb,0,1),~\Mb\neq 0$ and the classical ST for $M=(0,1,-1,0).$ Authors obtained the formula for the LCST of the translation and dilation of $f$ as 
$$(\mathcal{S}_\psi^M f(\tee-k))(\A,\B)=e^{-i\left(\frac{\Ma k}{\Mb}(k+\B)+\A k\right)}(\mathcal{S}_\psi^M (\mathcal{M}_{-\frac{\Ma k}{\Mb}}f))(\A,\B+k)$$
and 
$$(\mathcal{S}_\psi^M f(\lambda\tee))(\A,\B)=(\mathcal{S}_\psi^{\tilde{M}} f)\left(\frac{\A}{\lambda},\lambda\B\right),~\lambda\in\mathbb{R}^+,~\tilde{M}=(\Ma,\lambda^2\Mb,\Mc,\Md)$$
If further, the function $\psi$ satisfies the admissible conditions $c_\psi=\int_0^\infty\left|\left(\mathcal{L}^M\left[e^{i\left(\tee-\frac{\Ma\tee^2}{2\Mb}\left(1+\frac{1}{\A^2}\right)\right)}\psi(\tee)\right]\right)(\frac{\xxi}{\A})\right|d\A<\infty$ then $\mathcal{S}_\psi^M$ satisfies the inner product relations
$$\int_0^\infty\int_{-\infty}^\infty (\mathcal{S}_\psi^M f)(\A,\B)\overline{(\mathcal{S}_\psi^M g)(\A,\B)}d\A d\B=c_\psi\langle f,g\rangle$$
followed by the inversion formula 
$f(\tee)=\frac{1}{c_\psi}\int_{-\infty}^\infty\int_0^\infty(\mathcal{S}_\psi^M f)(\A,\B)\psi^M_{\A,\B}(\tee)d\A d\B,~a.e.$ As a consequence of the inner product relation it follows that, if $c_\psi=1$ then the operator $\mathcal{S}^M_\psi$ is an isometry from $L^2(\mathbb{R})$ to $L^2(\mathbb{R}^+\times\mathbb{R})$ and the range $\mathcal{S}^M_\psi(L^2(\mathbb{R}))$ is a RKHS with the kernel $K^M_\psi(\A,\B,\A',\B')=\langle\psi^M_{\A,\B},\psi^M_{\A',\B'}\rangle.$ In view of \eqref{R2020_Shah-Tantary_LCST_Definition_Integral-Form} and its representation in the linear canonical domain, it follows that the LCST is capable of providing the time and spectral information simultaneous and representing signal in the time-LCT-frequency domain. It gives local information of the signal $f(\tee)$ with the time window $\left[\B+\frac{E_\psi}{\A }-\frac{\Delta_\psi}{\A },\B+\frac{E_\psi}{\A }+\frac{\Delta_\psi}{\A}\right]$ and the local information about the linear canonical spectrum of $f(\tee)$ as $\left[\A E_H-\A \Delta_H,\A E_H+\A \Delta_H\right],$ where $H(\xxi)=\left(\mathcal{L}^M\left[e^{i\left(\tee-\frac{\Ma\tee^2}{2\Mb}\left(1+\frac{1}{\A^2}\right)\right)}\psi(\tee)\right]\right)(\xxi)$ with the constant window area $4\Mb\Delta_\psi\Delta_H$ in the time-LCT-frequency plane. This area depends on the the function $\psi$ and $\Mb$ and is independent of the parameters $\A$ and $\B.$ Thus the LCST is capable of providing time and linear canonical domain information simultaneously and representing signal in the time-LCT-frequency plane with the transform parameter $M.$ Wei et al. \cite{wei2022linear} extended the Gaussian function $g(\tee,\vv)$ to the linear canonical domain as $g(\tee,\vv)=\frac{|\frac{\vv}{\Mb}|^p}{\sqrt{2\pi}q}e^{-\frac{\tee^2\left(\frac{\vv}{\Mb}\right)^{2p}}{2q^2}}$ and used this function  as the window function and the linear canonical convolution \eqref{R2011_Wei-Ran-Li-LCT-Convolution} to define the LCST of a function $f\in L^2(\mathbb{R})$ as 
\begin{align}\label{R2022_IEEETSP_LCST_Definition}
(\mathcal{S}^Mf)(\tau,\vv)=&[e^{-i2\pi\frac{\vv}{\Mb}\tee}f(t)]\star'_Mg(-\tee,\vv)\notag\\
=&\int_{\mathbb{R}}f(\tee)g_M(\tee-\tau,\vv)e^{-i2\pi\frac{\vv}{\Mb}\tee}d\tee
\end{align} 
where the window function $g_M(\tee-\tau,\vv)$ is given by $g_M(\tee-\tau,\vv)=\frac{1}{\sqrt{2\pi}q}|\frac{\vv}{\Mb}|^pe^{-\frac{(\tee-\tau)^2}{2q^2}\left(\frac{\vv}{\Mb}\right)^{2p}}e^{i\frac{\Ma}{2\Mb}(\tee^2-\tau^2)}.$ This transform is linear and have the translation and scaling formula as $(\mathcal{S}^Mf(\tee-\tee_0))(\tau,\vv)=e^{i\frac{\Ma}{\Mb}(-\tau\tee_0+\tee_0^2)-i2\pi\frac{\vv}{\Mb}\tee_0}(\mathcal{S}^M[f(\tee)e^{i\frac{\Ma}{\Mb}\tee_0\tee}])(\tau-\tee_0,\vv)=$ and $(\mathcal{S}^{M'}f(k\tee))(\tau,\vv)=(\mathcal{S}^M[f(\tee)e^{-i2\pi\frac{\vv}{\Mb}\left(\frac{1}{k}-1\right)\tee}])(k\tau,\vv),$ $M'=(\frac{\Ma}{k^2},\Mb,\Mc,\Md).$ Moreover, it satisfies the inner product relation $\int_0^\infty\int_{-\infty}^\infty(\mathcal{S}^Mf)(\tau,\vv)\overline{(\mathcal{S}^Mg)(\tau,\vv)}=C_g\langle f,g\rangle$ with the $1-D$ inversion formula 
$$f(\tee)=\frac{1}{C_g}\int_0^\infty\int_{-\infty}^\infty\frac{1}{g(0,\vv)}(\mathcal{S}^Mf)(\tau,\vv)e^{i2\pi\frac{\vv}{\Mb}\tee}d\frac{\vv}{\Mb},~g(0,\vv)\neq 0$$
and the 2-D inversion formula as 
$$f(\tee)=\frac{1}{C_g}\int_0^\infty\int_{-\infty}^\infty(\mathcal{S}^Mf)(\tee-\tau,\vv)e^{-i\frac{\Ma}{2\Mb}(\tee^2-\tau^2)}e^{i2\pi\frac{\vv}{\Mb}\tee}d\tau\frac{d\vv}{\vv},$$
where $C_g=\int_0^\infty\left|G\left(\frac{\uu-2\pi\vv}{\Mb}\right)\right|^2\frac{d\vv}{\vv}$ is independent of $\uu$ in $\mathbb{R}_{0}$ with $G\left(\frac{\uu-2\pi\vv}{\Mb}\right)=e^{-2\pi^2q^2\frac{(\uu-2\pi\vv)^2}{\vv^2}}.$ Authors in \cite{wei2022linear} obtained the convolution theorem and the cross-correlation theorem for the transform \eqref{R2022_IEEETSP_LCST_Definition} along with that they also showed that the range $\mathcal{S}^{M}(L^2(\mathbb{R}))$ is a RKHS with the kernel $K_{g_M}(\tau_0,\vv_0;\tau,\vv)=\int_{-\infty}^\infty\overline{g_M(\tee-\tau,\vv)}g_M(\tee-\tau_0,\vv_0)e^{i2\pi\frac{(\vv-\vv_0)}{\Mb}\tee}d\tee$ and proved that the transform \eqref{R2022_IEEETSP_LCST_Definition} is capable of representing the signal in the time-linear canonical domain plane by jointly time and linear canonical domain information with the time-linear canonical-frequency window $[E_g+\tau-\Delta_g,E_g+\tau+\Delta_g]\times[\Mb E_G+2\pi\vv-\Mb \Delta_G,\Mb E_G+2\pi\vv+\Mb \Delta_G],$ with constant window area $2\Mb\Delta_e\Delta_G$ in the time-LCT domain-frequency plane. Gupta et al. \cite{gupta2024newST} gave a more general definition of a LCST. They considered the matrix parameters $M_{l}=(\Ma_l,\Mb_l;\Mc_l,\Md_l),~l=1,2$ with $\Ma_l\Md_l-\Mb_l\Mc_l=1~\mbox{and}~\Mb_{l}\neq 0,~l=1,2$ and defined the LCST for $f\in L^2(\mathbb{R})$ with respect to the window function $\psi\in L^2(\mathbb{R})$ as
\begin{align}\label{RGupta-Verma_LCST(2M)_Definition}
\left(\mathcal{S}^{M_1,M_2}_{\psi}f\right)(\A,\B)=&e^{-\frac{i\Ma_1}{2\Mb_1}\B^2}\left\{\tilde{f}(\tee)e^{\frac{i\Ma_1}{2\Mb_1}\tee^2}\star\overline{\A\psi(-\A\tee)e^{\frac{i\Ma_2}{2\Mb_2}(\A \tee)^2}}\right\}(\B),~\A\in\mathbb{R^+},\B\in\mathbb{R}\\
=&\int_{\mathbb{R}}f(\tee)\overline{\psi^{M_1,M_2}_{\A,\B}(\tee)}d\tee,
\end{align}
where $\tilde{f}(\tee)=e^{-i\frac{\A \tee}{B_{1}}}f(\tee),$ $\psi^{M_1,M_2}_{\A,\B}(\tee)=e^{-\frac{i\Ma_1}{2\Mb_1}(\tee^2-\Mb^2)+\frac{i\Ma_2}{2\Mb_2}\{\A(\tee-\B)\}^2+\frac{i\A \tee}{\Mb_1}}\A\psi(\A(\tee-\B))$ and $\star$ denotes the classical convolution. Depending upon the different choices of the matrices $M_1$ and $M_2$ we have different families of novel integral transforms, i.e., 
\begin{itemize}
\item For $M_1=(\cos\alpha,\sin\alpha;-\sin\alpha,\cos\alpha),~M_2=(\cos\beta,\sin\beta;-\sin\beta,\cos\beta),~\alpha,\beta\neq n\pi~(n\in\mathbb{Z}),$ \eqref{RGupta-Verma_LCST(2M)_Definition} give rise to a  novel FrST $\left(\mathcal{S}^{\alpha,\beta}_\psi f\right)(\A,\B)=\A\int_{\mathbb{R}}f(\tee)\overline{\psi\left(\A(\tee-\B)\right)}e^{\frac{i}{2}(\tee^2-\B^2)\cot\alpha-\frac{i}{2}\left\{\A(\tee-\B)\right\}^2\cot\beta-i\A \tee\csc\alpha}d\tee.$
\item For $M_1=(1,\Mb_1;0,1),~M_2=(1,\Mb_2;0,1),~\Mb_1,\Mb_2\neq 0,$ \eqref{RGupta-Verma_LCST(2M)_Definition} give rise to novel Fresnel Stockwell transform $\left(\mathcal{S}^{\Mb_1,\Mb_2}_\psi f\right)(\A,\B)=\A\int_{\mathbb{R}}f(\tee)\overline{\psi\left(\A(\tee-\B)\right)}e^{\frac{i\Ma_{1}}{2\Mb_{1}}(\tee^2-\B^2)-\frac{i\Ma_{2}}{2\Mb_{2}}\left\{\A(\tee-\B)\right\}^2-\frac{i\A \tee}{\Mb_{1}}}d\tee.$
\item For $M_1=M_2=(0,1;-1,0),$ \eqref{RGupta-Verma_LCST(2M)_Definition} reduces to the classical ST 
$\left(\mathcal{S}_\psi f\right)(\A,\B)=\A\int_{\mathbb{R}}f(\tee)\overline{\psi\left(\A(\tee-\B)\right)}e^{-i\A \tee}d\tee.$
\end{itemize}
Authors obtained the LCST of the translation and dilation of $f$ as 
$\left(\mathcal{S}^{M_1,M_2}_{\psi} f(\lambda\tee)\right)(\A,\B)=\left(\mathcal{S}^{\tilde{M_1},M_2}_{\psi}f\right)\left(\frac{\A}{\lambda},\B\lambda\right)$
and $\left(\mathcal{S}^{M_1,M_2}_{\psi} f(\tee-\yy)\right)(\A,\B)=e^{-\frac{i\A \yy}{\Mb_1}-\frac{i\Ma_1}{\Mb_1}\yy(\B-\yy)}\left(\mathcal{S}^{\tilde{M_1},M_2}_{\psi} \left\{e^{\frac{i\Ma_1}{\Mb_1}\yy \tee}f(\tee)\right\}\right)(\A,\B-\yy)$ where $\tilde{M_1}=\left(\frac{\Ma_1}{\lambda^2},\Mb_1;\Mc_1,\Md_1\lambda^2\right).$ If further, the window function $\psi$ satisfies the admissible conditions $c^{M_1,M_2}_\psi=\int_0^\infty\left|\left(\mathcal{L}^{M_2}\left[e^{\frac{i\tee}{\Mb_1}}\psi(\tee)\right]\right)(\frac{\Mb_2\xxi}{\Mb_1\A})\right|\frac{d\A}{\A}$ is a non zero positive constant then $\mathcal{S}_\psi^{M_1,M_2}$ satisfies the inner product relations
$$\int_0^\infty\int_{-\infty}^\infty (\mathcal{S}_\psi^{M_1,M_2} f)(\A,\B)\overline{(\mathcal{S}_\psi^M g)(\A,\B)}\frac{d\A d\B}{\A}=2\pi |\Mb_1|c^{M_1,M_2}_\psi\langle f,g\rangle_{L^2(\mathbb{R})}$$
followed by the inversion formula 
$f(\tee)=\frac{1}{2\pi |\Mb_1|c^{M_1,M_2}_\psi}\int_{-\infty}^\infty\int_0^\infty(\mathcal{S}_\psi^{M_1,M_2} f)(\A,\B)\psi^{M_1,M_2}_{\A,\B}(\tee)\frac{d\A d\B}{\A},~a.e.$ As a consequence of the inner product relation it follows that, if $c^{M_1,M_2}_\psi=\frac{1}{2\pi|\Mb_1|}$ then the operator $\mathcal{S}^{M_1,M_2}_\psi$ is an isometry from $L^2(\mathbb{R})$ to $L^2(\mathbb{R}^+\times\mathbb{R},\frac{d\A d\B}{\A})$ and the range $\mathcal{S}^{M_1,M_2}_\psi(L^2(\mathbb{R}))$ is a RKHS with the kernel $K^{M_1,M_2}_\psi(\A,\B,\A',\B')=\langle\psi^{M_1,M_2}_{\A,\B},\psi^{M_1,M_2}_{\A',\B'}\rangle_{L^2(\mathbb{R})}.$ Author studied the time-LCT-frequency resolution and proved that the joint resolution of the LCST in the time and linear canonical domain is 
\begin{align*}
\left[\frac{1}{\A}E_\psi+\B-\frac{1}{\A}\Delta_\psi,\frac{1}{A}E_\psi+\B+\frac{1}{\A}\Delta_\psi\right]\times\left[\frac{\Mb_1\A}{\Mb_2}E_{\Psi}-\frac{\Mb_1\A}{\Mb_2}\Delta_{\Psi},\frac{\Mb_1\A}{\Mb_2}E_{\Psi}+\frac{\Mb_1\A}{\Mb_2}\Delta_{\Psi}\right],
\end{align*}
with constant window area $4\frac{\Mb_1}{\Mb_2}\Delta_{\psi}\Delta_{\Psi}$ in the time-LCT-frequency plane, where $\Psi(\xxi)=\left(\mathcal{L}^{M_2}\{e^{\frac{i\tee}{\Mb_1}}\psi(\tee)\}\right)(\xxi).$ Author also developed the MRA $\{V^{M_1,M_2}_m\}_{m\in\mathbb{Z}}$ which satisfies the properties 
\begin{enumerate}[label=(\roman*)]
\item $V^{M_1,M_2}_m\subset V^{M_1,M_2}_{m+1},~\forall~m\in\mathbb{Z}$
\item $f(\tee)\in V^{M_1,M_2}_m$ if and only if $f(2\tee)e^{\frac{i\Ma_1}{2\Mb_1}((2\tee)^2-\tee^2)}\in V^{M_1,M_2}_{m+1},~\forall~m\in\mathbb{Z}$
\item $\cap_{m\in\mathbb{Z}}V^{M_1,M_2}_m=\{0\},$ $\overline{\cup_{m\in\mathbb{Z}}V^{M_1,M_2}_m}=L^2(\mathbb{R})$
\item there exists a function $\phi\in L^2(\mathbb{R})$, called the scaling function $ V^{M_1,M_2}_0$ such that $\{\phi_{M_1,M_2,0,n}(\tee):n\in\mathbb{Z}\}$ form an orthonormal basis of $V^{M_1,M_2}_0,$ where 
$$\phi_{M_1,M_2,m,n}(\tee)=2^{\frac{m}{2}}\phi(2^m\tee-n)e^{-\frac{i}{2}\left[\left(\tee^2-\left(\frac{n}{2^m}\right)^2\right)\frac{\Ma_1}{\Mb_1}-(2^m\tee-n)^2\frac{\Ma_2}{\Mb_2}-\frac{2(2^m\tee+n)}{\Mb_1}\right]},~m,n\in\mathbb{Z}.$$
\end{enumerate}
and obtained a function $\psi\in L^2(\mathbb{R})$ corresponding to the function $\phi(\tee)$ such that $\{\psi_{M_1,M_2,m,n}(\tee)\}_{m,n\in\mathbb{Z}}$ forms a novel orthonormal basis of $L^2(\mathbb{R}).$
%%%%%%%%%%%%%%%%%%%%%%%%%%%%%%%%%%%%%%%%%%%%%%%%%%%%%%%%%
%%%%%%%%%%%%%%%%%%%%%%%%%%%%%%%%%%%%%%%%%%%%%%%%%%%%%%%%%
%%%%%%%%%%%%%%%%%%%%%%%%%%%%%%%%%%%%%%%%%%%%%%%%%%%%%%%%%
%%%%%%%%%%%%%%%%%%%%%%%%%%%%%%%%%%%%%%%%%%%%%%%%%%%%%%%%%
%%%%%%%%%%%%%%%%%%%%%%%%%%%%%%%%%%%%%%%%%%%%%%%%%%%%%%%%%
%%%%%%%%%%%%%%%%%%%%%%%%%%%%%%%%%%%%%%%%%%%%%%%%%%%%%%%%%
%%%%%%%%%%%%%%%%%%%%%%%%%%%%%%%%%%%%%%%%%%%%%%%%%%%%%%%%%
\section{Quadratic Phase Windowed Fourier Transform (QPWFT)}

The importance of arbitrary real parameters involved in the definition of QPFT \eqref{RDefinitionQPFT} lies in the fact that an appropriate choice of the parameters can be used to inculcate a sense of rotation and shift into both the time and frequency axes, resulting in an efficient analysis of the chirp-like signals which are ubiquitous in nature. Thus, the chirp-like signals arising in radar and other communication systems reveal their characteristics better in the quadratic-phase Fourier domain than in the usual Fourier domain, fractional Fourier, and linear canonical domain. Despite several merits of the QPFT, one of its severe limitations is that it is inept for a localized analysis of the quadratic-phase spectrum of non-transient signals. To circumvent this limitation, Shah et al. \cite{shah2021short} introduced the notion of short-time QPFT, which can effectively localize the quadratic phase spectrum of a non-transient signal in the time-frequency plane for the function $f\in L^2(\mathbb{R})$ with respect to the window function $\phi$ as
\begin{align}\label{R2021_Shah-Lone-Tantary_Optik_WQPFT_Definition}
(\mathcal{V}^\wedge_\phi f)(\xx,\uu)=\frac{1}{\sqrt{2\pi}}\int_{\mathbb{R}}f(\tee)\overline{\tau_{\xx}g(\tee)}e^{i(\LmA\tee^2+\LmB\tee\uu+\LmC\uu^2+\LmD\tee+\LmE\uu)}d\tee,
\end{align}
where $\wedge=(\LmA,\LmB,\LmC,\LmD,\LmE)$ with $\LmB\neq 0.$ For the parameter $\wedge=\left(\frac{\LmA}{2\LmB},-\frac{1}{\LmB},\frac{\LmC}{2\LmB},0,0\right),$ $\wedge=\left(\frac{\cot\alpha}{2},-\csc\alpha,\frac{\cot\alpha}{2},0,0\right),$ $\alpha\neq n\pi,~n\in\mathbb{Z}$ and $\wedge=\left(0,1,0,0,0\right)$ and multiplying the right of the equation \eqref{R2021_Shah-Lone-Tantary_Optik_WQPFT_Definition} by $\frac{1}{\sqrt{i\LmB}},$ $\sqrt{1-i\cot\alpha}$ and $1$ respectively, it can be seen that the WQPFT \eqref{R2021_Shah-Lone-Tantary_Optik_WQPFT_Definition} reduces to the WLCT, WFrFT and WFT. Observe that, if $f\in L^2(\mathbb{R}),~g\in L^2(\mathbb{R})$ or $f\in L^1(\mathbb{R}),~g\in L^\infty(\mathbb{R})$ then the WQPFT \eqref{R2021_Shah-Lone-Tantary_Optik_WQPFT_Definition} is well defined. Similarly, if $f\in L^q(\mathbb{R})$ and $g\in L^p(\mathbb{R}), 1\leq p,q<\infty$ and $\frac{1}{p}+\frac{1}{q}=1,$ then by H\"older's inequality it follows that the transform \eqref{R2021_Shah-Lone-Tantary_Optik_WQPFT_Definition} is well defined. Based on the relation between the WQPFT and the WFT, they obtained  the Lieb's inequality for the WQPFT 
\begin{align*}
\int_{\mathbb{R}}\int_{\mathbb{R}}\left|(\mathcal{V}^\wedge_\phi f)(\xx,\uu)\right|^pd\xx d\uu\leq\frac{2}{p|\LmB|}\|f\|^p_{L^2(\mathbb{R})}\|\phi\|^p_{L^2(\mathbb{R})}.
\end{align*}
Transform \eqref{R2021_Shah-Lone-Tantary_Optik_WQPFT_Definition} is linear in $f,$ anti-linear in $g$ and has the following covariance properties 
\begin{align*}
\left(\mathcal{V}^\wedge_\phi\left(\mathcal{M}^{M}_\mu \tau_{\xx_0}f\right)\right)(\xx,\uu)=e^{i\left(\LmB\xx_0(\uu-\mu)+\LmC(\uu^2-\mu^2)+\LmE(\uu-\mu)\right)}\mathcal{V}^\wedge_\phi f(\xx-\xx_0,\LmB(\uu-\mu))
\end{align*}
and 
\begin{align*}
\left(\mathcal{V}^\wedge_\phi\left(\tau_{\xx_0}\mathcal{M}^{M}_\mu f\right)\right)(\xx,\uu)=e^{i\left(\LmA\xx_0^2+\LmB \xx_0(\mu+\uu)+\LmC\left(\uu^2-\mu^2\right)+\LmD\xx_0+\LmE(\uu-\mu^2)+2\LmA\xx_0-\LmA\LmB\right)}\mathcal{V}_\phi f(\xx-\xx_0,2\Ma\xx_0-\LmB(\uu-\mu)),
\end{align*}
where $\mathcal{M}^\wedge_\mu f(\tee)=e^{i(\LmA\tee^2+\LmB\tee\mu+\LmC\uu^2+\LmD\tee+\LmE\mu)}f(\tee)$ is generalized modulation by $\mu.$ The covariance property obtained above indeed gives the formula for the WQPFT for the time-frequency shift $\tau_{\xx_0}\mathcal{M}^\wedge_\mu$ or $\mathcal{M}^\wedge_\mu\tau_{\xx_0}$ of the function $f.$ Specifically, the formula for the parity, conjugation,  $\left(\mathcal{V}^\wedge_\phi\left(\tau_{\xx_0} f\right)\right)(\xx,\uu)$ and $\left(\mathcal{V}^\wedge_\phi\left(e^{i\uu_0\xx}f\right)\right)(\xx,\uu)$ can be found in \cite{shah2021short}. For $f_1,f_2,\phi_1,\phi_2\in L^2(\mathbb{R}),$ the WQPFT \eqref{R2021_Shah-Lone-Tantary_Optik_WQPFT_Definition} satisfies the following inner product relation
\begin{align}\label{R2021_Shah-Lone-Tantary_Optik_WQPFT_IPR}
\left\langle(\mathcal{V}^\wedge_{\phi_1}f_1)(\xx,\uu),(\mathcal{V}^\wedge_{\phi_2}f_2)(\xx,\uu)\right\rangle_{L^2(\mathbb{R}^2)}=\frac{1}{|\LmB|}\langle f_1,f_2 \rangle\overline{\langle \phi_1,\phi_2 \rangle}.
\end{align}
In particular, if $f=g,$ $\phi_1=\phi_2=\phi$ and $\|\phi_1\|_{L^2(\mathbb{R})}=\|\phi_2\|_{L^2(\mathbb{R})}=1$ then $\|\mathcal{V}^\wedge_\phi f\|_{L^2(\mathbb{R}^2)}^2=\frac{1}{|\Mb|}\|f\|_{L^2(\mathbb{R})}^2$ and thus WQPFT is an isometry from $L^2(\mathbb{R})$ to $L^2(\mathbb{R}^2)$. The inversion formula for the WQPFT is given as 
\begin{align}\label{R2021_Shah-Lone-Tantary_Optik_WQPFT_Inversion-Formula}
f(\tee)=\frac{|\Mb|}{\sqrt{2\pi}~\overline{\langle \phi_1,\phi_2\rangle}}\int_{\mathbb{R}}\int_{\mathbb{R}}(\mathcal{V}^\wedge_{\phi_2}f)(\xx,\uu)e^{-i(\LmA\tee^2+\LmB\tee\uu+\LmC\uu^2+\LmD\tee+\LmE\uu)}\tau_{\xx}\phi_1(\tee)d\uu d\xx
\end{align}
where $\phi_1,\phi_2\in L^2(\mathbb{R}),$ $\langle \phi_1,\phi_2\rangle\neq 0$  and $f\in L^2(\mathbb{R}).$ Authors in \cite{shah2021short} also characterized the range of the transform \eqref{R2021_Shah-Lone-Tantary_Optik_WQPFT_Definition} and obtained the Heisenberg's uncertainty  principle for the transform, for the functions $f\in L^2(\mathbb{R}),$ using corresponding UP \eqref{R2021_Shah-Nisar-Lone-Tantary_MMAS_QPFT_Heisenberg-UP} for the QPFT as
\begin{align}\label{R2021_Shah-Lone-Tantary_Optik_WQPFT_Heisenberg-UP}
\left(\int_{\mathbb{R}}\int_{\mathbb{R}}\uu^2\left|(\mathcal{V}^\wedge_{\phi}f)(\xx,\uu)\right|^2d\xx d\uu\right)^{\frac{1}{2}}\left(\int_{\mathbb{R}}\tee^2|f(\tee)|^2d\tee\right)^{\frac{1}{2}}\geq \frac{1}{2|\LmB|}\|f\|^2_{L^2(\mathbb{R})}\|\phi\|_{L^2(\mathbb{R})}.
\end{align}
For the choice of the parameter $\wedge=\left(\frac{\LmA}{2\LmB},-\frac{1}{\LmB},\frac{\LmC}{2\LmB},0,0\right),$ $\wedge=\left(\frac{\cot\alpha}{2},-\csc\alpha,\frac{\cot\alpha}{2},0,0\right),~\alpha\neq n\pi,~n\in\mathbb{Z}$ and $\wedge=\left(0,1,0,0,0\right)$ the inequality \eqref{R2021_Shah-Lone-Tantary_Optik_WQPFT_Heisenberg-UP} reduces to the corresponding Heisenberg uncertainty inequality for the WLCT, WFrFT and WFT respectively.\\
 
\noindent\textbf{Uncertainty Principles:} They also obtained the logarithmic UP for the transform \eqref{R2021_Shah-Lone-Tantary_Optik_WQPFT_Definition}, for the function $f\in \mathcal{S}(\mathbb{R}),$ using the same \eqref{R2021_Shah-Nisar-Lone-Tantary_MMAS_QPFT_Logarithmic-UP} for the QPFT  as
\begin{align}
|\LmB|^2\|\phi\|^2_{L^2(\mathbb{R})}\int_{\mathbb{R}}\ln|\tee||f(\tee)|^2d\tee+|\LmB|\int_{\mathbb{R}}\int_{\mathbb{R}}\ln|\uu|\left|(\mathcal{V}^\wedge_{\phi}f)(\xx,\uu)\right|^2d\xx d\uu\geq \left[\frac{\Gamma'\left(\frac{1}{4}\right)}{\Gamma\left(\frac{1}{4}\right)}-\ln\pi-\ln|\LmB|\right]\|f\|^2_{L^2(\mathbb{R})}\|\phi\|^2_{L^2(\mathbb{R})}.
\end{align}
%%%%%%%%%%%%%%%%%%%%%%%%%%%%%%%%%%%%%%%%%%%%%%%%%%%%%%%%%
%%%%%%%%%%%%%%%%%%%%%%%%%%%%%%%%%%%%%%%%%%%%%%%%%%%%%%%%%
%%%%%%%%%%%%%%%%%%%%%%%%%%%%%%%%%%%%%%%%%%%%%%%%%%%%%%%%%
%%%%%%%%%%%%%%%%%%%%%%%%%%%%%%%%%%%%%%%%%%%%%%%%%%%%%%%%%
\section{Quadratic Phase Wavelet Transform (QPWT)}
The QPFT  is a generalization of classical WT, continuous FrWT, and the LCWT. The transform is a more beneficial tool than the QPFT and the QPWFT. Prasad et al. \cite{prasad2020quadratic} first introduced the definition, which was based on the quadratic phase convolution \eqref{R2018_Prasad-Sharma_MMAS_QPWT_QP-Convolution} and is defined for the function $f\in L^2(\mathbb{R})$ with respect to the function $\psi\in L^2(\mathbb{R})$ and the parameter $\wedge=(\LmA,\LmB,\LmC,\LmD,\LmE)$ as 
\begin{align}\label{R2018_Prasad-Sharma_MMAS_Quadratic-Phase-WT_Definition}
(\mathcal{W}^\wedge_\psi f)(\A,\B)=&\left(f\star_{\Omega^{\LmA}_{\LmD}}\psi_\A\right)(\B)\notag\\
=&\int_{\mathbb{R}}f(\xx)\overline{\psi^{\Omega^{\LmA}_{\LmD}}_{\A,\B}(\xx)}d\xx,~\A\in \mathbb{R}^+,\B\in\mathbb{R},
\end{align}
where $\psi_{\A}(\tee)=\frac{1}{\sqrt{\A}}\overline{\psi\left(-\frac{\tee}{\A}\right)},$ $\psi^{\Omega^{\LmA}_{\LmD}}_{\A,\B}(\xx)=\frac{1}{\sqrt{\A}}\psi\left(\frac{\tee-\B}{\A}\right)e^{-i\Ma(\tee^2-\B^2)-i\Md(\tee-\B)},~\B\in\mathbb{R},~\A\in\mathbb{R}^+.$ The transform \eqref{R2018_Prasad-Sharma_MMAS_Quadratic-Phase-WT_Definition} is linear in $f,$ anti-linear in $\psi.$ With the definition of the translation operator $\tau'_{\LmA r}f(\tee)=e^{-2\LmA ir\tee}f(\tee-r),$ author in \cite{prasad2020quadratic} obtained the translational formula  $(\mathcal{W}^\wedge_\psi \tau_{\LmA r} f)(\A,\B)=e^{-2\LmA i\B r}(\mathcal{W}^\wedge_\psi f)(\A,\B-r)$ and the scaling formula  $(\mathcal{W}^\wedge_\psi  [\frac{1}{\sqrt{m}}f\left(\frac{\tee}{m}\right)])(\A,\B)=e^{-2\LmA i\B r}(\mathcal{W}^\wedge_\psi f)(\frac{\A}{m},\frac{\B}{m}),~m>0.$ With the assumption that the functions $\psi_1,\psi_2\in L^2(\mathbb{R})$ satisfying the condition
$$c_{\psi_1,\psi_2}^{\Omega^{\LmA}_{\LmD}}:=\int_{\mathbb{R}^+}\overline{\left(\mathcal{Q}^{\LmA,\LmB,\LmC}_{\LmD,\LmE}[e^{-i\LmA(\cdot)^2-i\LmD(\cdot)}\psi_1]\right)(\A)}\left(\mathcal{Q}^{\LmA,\LmB,\LmC}_{\LmD,\LmE}[e^{-i\LmA(\cdot)^2-i\LmD(\cdot)}\psi_2]\right)(\A)\frac{d\A}{\A}$$ 
is non-zero complex number, they proved the inner product relation  for the functions $f_1,f_2\in L^2(\mathbb{R})$
\begin{align}\label{R2021_Shah-Lone_MMAS_Quadratic-Phase-WT_IPR}
\int_{\mathbb{R}}\int_{\mathbb{R}^+}(\mathcal{W}^\wedge_{\psi_1} f)(\A,\B)\overline{(\mathcal{W}^\wedge_{\psi_2} g)(\A,\B)}\frac{d\B d\A}{\A^2}=\frac{2\pi}{\Mb}c_{\psi_1,\psi_2}^{\Omega^{\LmA}_{\LmD}}\langle f,g\rangle,
\end{align}
followed by the reconstruction formula $f(\xx)=\frac{\LmB}{2\pi c^{\Omega^{\LmA}_{\LmD}}_{\psi,\psi}}\int_{\mathbb{R}}\int_{\mathbb{R}^+}(\mathcal{W}^\wedge_\psi f)(\A,\B)\psi^{\Omega^{\LmA}_{\LmD}}_{\A,\B}(\xx)\frac{d\B d\A}{\A^2},~a.e.$ Motivated by the Shannon sampling theorem \eqref{RLeonSamplingTheoremQPFT}, Gupta et al. \cite{gupta2023quadratic} studied the MRA associated with QPWT given by  \eqref{R2018_Prasad-Sharma_MMAS_Quadratic-Phase-WT_Definition} and obtain a novel class of orthonormal bases. Sharma et al.\cite{sharma2023abelian} established the Abelian theorem, which relates the behavior of the QPWT $(\mathcal{W}^\wedge_\psi f)(\A,\B)$ as $\A\rightarrow \infty (\A\rightarrow 0)$ to the behavior of function $(\mathcal{Q}^{\LmA,\LmB,\LmC}_{\LmD,\LmE} f)(\xxi)$ as $\xxi\rightarrow \infty (\xxi\rightarrow 0).$ These results are also referred to as initial value and final value Abelian theorems. Based on the convolution \eqref{R2023_SIVP_Shah-Tantary_QPFT_Convolution}, of the QPFT introduced in \cite{shah2022sampling}, Shah et al.\cite{shah2022quadratic} introduced a new definition of the QPWT for a function $f\in L^2(\mathbb{R}),$ with respect to the function $\psi\in L^2(\mathbb{R})$ and the parameter set $\wedge=(\LmA,\LmB,\LmC,\LmD,\LmE),$ which we denote by $^{(1)}\mathcal{W}^\wedge_\psi$ and is defined  as
\begin{align}\label{R2021_Shah-Lone_MMAS_QuadraticPhaseWT_Definition}
\left(^{(1)}\mathcal{W}^\wedge_\psi f\right)(\A,\B)=&\left(f(\xx)\star'_\wedge\right(\mathcal{D}_{\A}\hat{\psi}\left)(-\xx)\right)(\B)\notag\\
=&\int_{\mathbb{R}}f(\xx)\overline{\psi^\wedge_{\A,\B}(\xx)}d\xx,~\A\in \mathbb{R}^+,\B\in\mathbb{R}
\end{align}  
where $\mathcal{D}_\A$ is the usual dilation operator and the family $\psi^\wedge_{\A,\B}(\xx)=\frac{1}{\sqrt{2\pi\A}}\psi\left(\frac{\xx-\B}{\A}\right)e^{2i\LmA \xx(\xx-\B)}.$ For the choice of the parameter $\wedge=\left(\frac{\LmA}{2\LmB},-\frac{1}{\LmB},\frac{\LmC}{2\LmB},0,0\right),$ $\wedge=\left(\frac{\cot\alpha}{2},-\csc\alpha,\frac{\cot\alpha}{2},0,0\right),~\alpha\neq n\pi,~n\in\mathbb{Z}$ and $\wedge=\left(0,1,0,0,0\right)$ the QPWT \eqref{R2021_Shah-Lone_MMAS_QuadraticPhaseWT_Definition} reduces to the novel LCWT, novel FrWT and the classical WT respectively. The transform \eqref{R2021_Shah-Lone_MMAS_QuadraticPhaseWT_Definition} is linear in $f$ and anti linear in $\psi$ and has the following covariance property 
$$\left(^{(1)}\mathcal{W}^\wedge_\psi f(\xx-k)\right)(\A,\B)=e^{2i\LmA k(\B-k)}\left(^{(1)}\mathcal{W}^\wedge_\psi \left[e^{-2i\LmA\xx k}f(\xx)\right]\right)(\A,\B)$$
and
$$\left(^{(1)}\mathcal{W}^\wedge_\psi f(\mu\xx)\right)(\A,\B)= \frac{1}{\mu}\left(^{(1)}\mathcal{W}^\wedge_\psi f(\mu\xx)\right)(\mu\A,\mu\B),~\wedge'=\left(\frac{\LmA}{\mu^2},\LmB,\LmC,\LmD,\LmE\right),$$
which gives the formula for the QPWT of the translation and scaling of a function $f$, respectively. Apart from this, using the machinery of the QPFT, authors in \cite{prasad2020quadratic} obtained orthogonality relation
\begin{align}\label{R2021_Shah-Lone_MMAS_Quadratic-Phase-WT_IPR}
\int_{\mathbb{R}}\int_{\mathbb{R}^+}\left(^{(1)}\mathcal{W}^\wedge_\psi f\right)(\A,\B)\overline{\left(^{(1)}\mathcal{W}^\wedge_\psi g\right)(\A,\B)}\frac{d\B d\A}{\A^2}=c_\psi\langle f,g\rangle,
\end{align}
followed by the reconstruction formula $f(\xx)=\frac{1}{c_\psi}\int_{\mathbb{R}}\int_{\mathbb{R}^+}\left(^{(1)}\mathcal{W}^\wedge_\psi f\right)(\A,\B)\psi^\wedge_{\A,\B}(\xx)\frac{d\B d\A}{\A^2},~a.e,$ with the assumption that the function $\psi\in L^2(\mathbb{R})$ satisfies the admissibility condition, i.e.,  $0<c_\psi:=\int_{\mathbb{R}^+}\frac{\left|\left(\mathcal{Q}_\wedge \Psi\right)(\A \oomega)\right|^2}{\A}d\A<\infty,$ where $\Psi(\xx,\A)=e^{-i\left(\LmA^2\A^2+\LmD \xx \A^2-\LmD \xx\right)}\tilde{\psi}(\xx),$ $\tilde{\psi}(\xx)=\bar{\psi}(-\xx).$ For $f=g$ and $c_\psi=1,$ it follows from \eqref{R2021_Shah-Lone_MMAS_Quadratic-Phase-WT_IPR} that $\int_{\mathbb{R}^+}\int_{\mathbb{R}}\left|\left(^{(1)}\mathcal{W}^\wedge_\psi f\right)(\A,\B)\right|^2\frac{d\A d\B}{\A^2}=\|f\|^2_{L^2(\mathbb{R})},$ which means that the QPWT is an isometry from $L^2(\mathbb{R})$ into $L^2(\mathbb{R}^+\times\mathbb{R}).$ They also obtained the reproducing kernel for the range $\mathcal{W}^\wedge_\psi\left(L^2(\mathbb{R})\right)$ to be $\mathbb{K}_\psi(\A,\B;\A',\B')=\left\langle \psi^\wedge_{\A,\B},\psi^\wedge_{\A',\B'} \right\rangle.$\\

\noindent\textbf{Uncertainty Principles:} Using the properties of the QPFT and its relation with the QPWT, authors in \cite{castrouncertainty},\cite{prasad2020quadratic} obtained the Heisenberg's UP for the transform \eqref{R2021_Shah-Lone_MMAS_QuadraticPhaseWT_Definition} using equation \eqref{R2021_Shah-Nisar-Lone-Tantary_MMAS_QPFT_Heisenberg-UP} as
\begin{align*}
\left(\int_{\mathbb{R}}\int_{\mathbb{R}^+}\B^2\left|\left(^{(1)}\mathcal{W}^\wedge_\psi f\right)(\A,\B)\right|^2\frac{d\A d\B}{\A^2}\right)^{\frac{1}{2}}\left(\int_{\mathbb{R}}\oomega^2\left|(\mathcal{Q}_\wedge f)(\oomega)\right|^2d\oomega\right)^{\frac{1}{2}}\leq \frac{\sqrt{c_\psi}}{2|\LmB|}\|f\|^2_{L^2(\mathbb{R})}.
\end{align*}
Moreover, using the logarithmic UP  \eqref{R2021_Shah-Nisar-Lone-Tantary_MMAS_QPFT_Logarithmic-UP} authors in \cite{castrouncertainty} obtained the logarithmic UP for the QPWT \eqref{R2018_Prasad-Sharma_MMAS_Quadratic-Phase-WT_Definition} as
\begin{align*}
\int_{\mathbb{R}}\int_{\mathbb{R}^+}\log|\B|\left|\left(\mathcal{W}^\wedge_\psi f\right)(\A,\B)\right|^2\frac{d\A d\B}{\A^2}+|\LmB|\int_{\mathbb{R}}\oomega^2\left|(\mathcal{Q}_\wedge f)(\oomega)\right|^2d\oomega\leq \left(\frac{\Gamma'\left(\frac{1}{4}\right)}{\Gamma\left(\frac{1}{4}\right)}-\ln\pi-\ln|\LmB|\right)\|f\|^2_{L^2(\mathbb{R})}.
\end{align*} 
In addition to these, using the local and entropy-based UP  for the QPFT \cite{shah2021uncertainty}, authors also obtained the corresponding UP for the QPWT. 
%%%%%%%%%%%%%%%%%%%%%%%%%%%%%%%%%%%%%%%%%%%%%%%%%%%%%%%%%
%%%%%%%%%%%%%%%%%%%%%%%%%%%%%%%%%%%%%%%%%%%%%%%%%%%%%%%%%
%%%%%%%%%%%%%%%%%%%%%%%%%%%%%%%%%%%%%%%%%%%%%%%%%%%%%%%%%
%%%%%%%%%%%%%%%%%%%%%%%%%%%%%%%%%%%%%%%%%%%%%%%%%%%%%%%%%
%%%%%%%%%%%%%%%%%%%%%%%%%%%%%%%%%%%%%%%%%%%%%%%%%%%%%%%%%
%%%%%%%%%%%%%%%%%%%%%%%%%%%%%%%%%%%%%%%%%%%%%%%%%%%%%%%%%
%%%%%%%%%%%%%%%%%%%%%%%%%%%%%%%%%%%%%%%%%%%%%%%%%%%%%%%%%
\section{Quaternion Windowed Fourier transform (QWFT)}
The QFT generalizes the real and complex FT to the quaternion case. However, the QFT is ineffective in studying the local behavior of the quaternion signal. To overcome this, B\"ulow and Sommer constructed one such extension by studying a particular case of QWFT with quaternionic Gabor filters \cite{bulow1999hypercomplex},\cite{bulow2001non}. Their generalization was obtained using the inverse (two-sided) quaternion Fourier kernel.
Bahri et al.\cite{bahri2010windowed} generalized the classical WFT to quaternion-valued functions of two variables using the right-sided QFT. Before we present here the definition of the QWFT, let us recall that a quaternion window function is a function $\phi\in L^2_{\mathbb{H}}(\mathbb{R}^2)$ such that $|\bldx|^{\frac{1}{2}}\phi(\bldx)\in L^2_{\mathbb{H}}(\mathbb{R}^2).$ Moreover, $\phi_{\bldomega,\bldb}(\bldx)=\frac{1}{(2\pi)^2}e^{j\omega_2\xx_2}e^{i\omega_1\xx_1}\phi(\bldx-\bldb)$ will be called as a quaternionic window daughter function. Thus the QWFT of a function $f\in L^2_{\mathbb{H}}(\mathbb{R}^2)$ with respect to the quaternion window function $\phi$ is defined as 
\begin{align}\label{R2010_Bahri-Hitzer-Ashino-Vaillancourt_AMC_RsQWFT_Definition}
(\mathcal{G}_\phi f)(\bldomega,\bldb)=&(f,\phi_{\bldomega,\bldb})\notag\\
=&\frac{1}{(2\pi)^2}\int_{\mathbb{R}^2}f(\bldx)\overline{\phi(\bldx-\bldb)}e^{-i\omega_1\xx_1}e^{-j\omega_2\xx_2}d\bldx.
\end{align}
where the inner product $(\cdot,\cdot)$ is defined in \eqref{RQuaternionInnerProduct}. The QWFT  \eqref{R2010_Bahri-Hitzer-Ashino-Vaillancourt_AMC_RsQWFT_Definition} is different from \cite{bahri2010two} in the definition of the exponential kernel. It is worth mentioning that the order of the exponentials in \eqref{R2010_Bahri-Hitzer-Ashino-Vaillancourt_AMC_RsQWFT_Definition} is fixed because of the non-commutativity of the product of quaternions. An example of a quaternion window function is the two-dimensional Haar function 
$$\phi(\bldx)=
\begin{cases}
1, &\mbox{for}~0\leq\xx_1<\frac{1}{2}~\mbox{and}~ 0\leq\xx_2<\frac{1}{2}\\
-1, &\mbox{for}~\frac{1}{2}\leq\xx_1<1~\mbox{and}~ \frac{1}{2}\leq\xx_1<1\\
0, &\mbox{otherwise.}
\end{cases}
$$
The transform \eqref{R2010_Bahri-Hitzer-Ashino-Vaillancourt_AMC_RsQWFT_Definition} is linear in $f$ i.e., $(\mathcal{G}_\phi (\lambda f+\mu g))(\bldomega,\bldb)=\lambda(\mathcal{G}_\phi f)(\bldomega,\bldb)+\mu (\mathcal{G}_\phi f)(\bldomega,\bldb),~\lambda,\mu\in \mathbb{H},$ $(\mathcal{G}_{P\phi} Pf)(\bldomega,\bldb)=(\mathcal{G}_\phi f)(-\bldomega,-\bldb)$ and for specific choice of functions $f=f_0+if_2,~\phi=\phi_0+i\phi_1$ the transform have the formula $(\mathcal{G}_\phi \tau_{\xx_0}f)(\bldomega,\bldb)=e^{-i\omega_1\xx_0}\left((\mathcal{G}_\phi f)(\bldomega,\bldb-\bldx_0)\right)e^{-j\omega_2\xx_0},~\bldx_0=(\xx_0,\xx_1).$ Apart from this, using the machinery of the right-sided QFT \cite{hitzer2007quaternion},\cite{bahri2008uncertainty}, authors in \cite{bahri2010windowed} obtained the reconstruction formula $f(\bldx)=\frac{(2\pi)^2}{C_\phi}\int_{\mathbb{R}^2}\int_{\mathbb{R}^2}(\mathcal{G}_\phi f)(\bldomega,\bldb)\phi_{\bldomega,\bldb}(\bldx)d\bldb d\bldomega,$ with the assumption that $C_\phi=\|\phi\|^2_{L^2_{\mathbb{H}}(\mathbb{R}^2)}$ and $0<C_\phi<\infty$ and proved the orthogonality relation
\begin{align*}
\int_{\mathbb{R}^2}\int_{\mathbb{R}^2}(\mathcal{G}_\phi f)(\bldomega,\bldb)\overline{(\mathcal{G}_\phi g)(\bldomega,\bldb)}d\bldomega d\bldb=\frac{C_\phi}{(2\pi)^2}(f,g).
\end{align*}  
In particular, for $f=g,$ and $C_\phi=(2\pi)^2$ gives $\int_{\mathbb{R}^2}\int_{\mathbb{R}^2}\left|(\mathcal{G}_\phi f)(\bldomega,\bldb)\right|^2d\bldb d\bldomega=\|f\|^2_{L^2_{\mathbb{H}}(\mathbb{R}^2)},$ this shows that the QWFT is an isometry from $L^2_{\mathbb{H}}(\mathbb{R}^2)$ into $L^2_{\mathbb{H}}(\mathbb{R}^2\times\mathbb{R}^2),$ i.e., the energy of the quaternion-valued signal computed in the spatial domain is equal to the total energy computed in the QWFT domain. They also obtain the reproducing kernel for the range to be $\mathbb{K}_\phi(\bldomega,\bldb;\bldomega',\bldb')=\frac{(2\pi)^2}{C_\phi}(\phi_{\bldomega,\bldb},\phi_{\bldomega',\bldb'}).$ Bahri et al. in \cite{bahri2020uncertainty} explored further properties and some important inequalities by considered the transform \eqref{R2010_Bahri-Hitzer-Ashino-Vaillancourt_AMC_RsQWFT_Definition} with an extra factor $(2\pi)^2,$ we call it for our convenience $\mathcal{G}'_\phi f=(2\pi)^2\mathcal{G}_\phi f.$ Using the relation between the $\mathcal{G}'_\phi f$ and the right-sided Fourier transform authors in \cite{bahri2020uncertainty} gave the alternate proofs of these properties. They obtained the orthogonality relation and the inversion formula depending on two window functions. They also obtained the Lieb's inequality as $\int_{\mathbb{R}^2}\int_{\mathbb{R}^2}\left|\left(\mathcal{G}'_\phi f\right)(\bldomega,\bldb)\right|^p_Qd\bldomega d\bldb~
\begin{cases}
\leq (\|f\|_{L^2_{\mathbb{H}}(\mathbb{R}^2)}\|\phi\|_{L^2_{\mathbb{H}}(\mathbb{R}^2)})^p &\mbox{for}~2\leq p<\infty \\
\geq(\|f\|_{L^2_{\mathbb{H}}(\mathbb{R}^2)}\|\phi\|_{L^2_{\mathbb{H}}(\mathbb{R}^2)})^p &\mbox{for}~1\leq p<2
\end{cases}
.$  Authors in \cite{bahri2010windowed} also pointed out that the extension of the WFT to the quaternion case by means of the two-sided QFT is complicated because of the non-commutativity. As a result, a major problem for the representation of the kernel of the two-sided QWFT  arises. Fu et al. \cite{fu2012balian} defined the two-sided QWFT of $f\in L^2_{\mathbb{H}}(\mathbb{R})$ with respect to a non-zero windowed function $g\in L^2(\mathbb{R}^2,\mathbb{R})$ as 
\begin{align}\label{R2012_Fu-Kahler-Cerejeiras_AACA_QWFT_Definition}
(\mathtt{Q}_gf)(\bldomega,\bldb)=\int_{\mathbb{R}^2}e^{-2\pi i\xx_1\omega_1}f(\bldx)g(\bldx-\bldb)e^{-2\pi j\xx_2\omega_2}d\bldx.
\end{align}
It can be seen from the above definition that the classical way of expressing it as a linear inner product between function and kernel is not possible. So authors defined the right multiplicative operator $\mathcal{C}[\lambda,\cdot]$ and $\mathcal{C}[\cdot,\lambda],$ $\lambda$ being a quaternion-valued function, as $\mathcal{C}[\lambda,f]=f(\bldx)\lambda(\bldx)$ and $\mathcal{C}[f,\lambda]=\lambda(\bldx)f(\bldx),$ respectively and expressed the kernel of the QWFT as an atom operator $g_{\bldomega,\bldb}=e^{2\pi j\omega_2\xx_2}g(\bldx-\bldb)\mathcal{C}[e^{2\pi i\omega_1\xx_1},\cdot],$ $\bldomega,\bldb,\bldx\in\mathbb{R}^2.$ The kernel operator is such that for an arbitrary function $f:\mathbb{R}^2\rightarrow \mathbb{H},$ 
$$g_{\bldomega,\bldb}\overline{f(\bldx)}=e^{2\pi j\omega_2 \xx_2}g(\bldx-\bldb)\overline{f(\bldx)}e^{2\pi i\omega_1\xx_1}~\mbox{and}~f(\bldx)\overline{g_{\bldomega,\bldb}}=e^{-2\pi i\omega_1\xx_1}g(\bldx-\bldb)f(\bldx)e^{-2\pi j\omega_2\xx_2}$$
and so $g$ being real valued, $f(\bldx)\overline{g_{\bldomega,\bldb}}=\overline{g_{\bldomega,\bldb}\overline{f(\bldx)}}.$ Thus the integral representation formula for the transform \eqref{R2012_Fu-Kahler-Cerejeiras_AACA_QWFT_Definition} as 
$$(\mathtt{Q}_gf)(\bldomega,\bldb)=(f,g_{\bldomega,\bldb}):=\int_{\mathbb{R}^2}f(\bldx)\overline{g_{\bldomega,\bldb}}d\bldx~\mbox{and}~(\mathtt{Q}_gf)(\bldomega,\bldb)=(f,g_{\bldomega,\bldb})=\overline{\int_{\mathbb{R}^2}g_{\bldomega,\bldb}\overline{f(\bldx)}d\bldx}.$$
Authors in \cite{fu2012balian} obtained several properties of the QWFT \eqref{R2012_Fu-Kahler-Cerejeiras_AACA_QWFT_Definition} like $(\mathtt{Q}_{Pg}Pf)(\bldomega,\bldb)=(\mathtt{Q}_gf)(-\bldomega,-\bldb),$ $(\mathtt{Q}_g(\tau_{\bldx_0}f))(\bldomega,\bldb)=e^{-2\pi i\xx_1\omega_1}(\mathtt{Q}_gf)(\bldomega,\bldb-\bldx_0)e^{-2\pi j\xx_2\omega_2}$ $(\mathtt{Q}_g[e^{2\pi i\xx_1\omega_{01}}f(\bldx)e^{2\pi j\xx_2\omega_{02}}])(\bldomega,\bldb)=(\mathtt{Q}_gf)(\bldomega-\bldomega_0,\bldb),\bldomega_0=(\omega_{01},\omega_{02}).$ Apart from this author also obtained the reconstruction formula 
\begin{align}\label{R2012_Fu-Kahler-Cerejeiras_AACA_QWFT_Inversion-Formula}
f(\bldx)=\frac{1}{C_g}\int_{\mathbb{R}^2}\int_{\mathbb{R}^2}(\mathtt{Q}_gf)(\bldomega,\bldb)\overline{g_{-\bldomega,\bldb}}d\bldomega d\bldb,
\end{align}
and the orthogonality relation $\langle \mathtt{Q}_gf_1,\mathtt{Q}_gf_2)\rangle=C_g\langle f_1,f_2\rangle,$ where $\langle \cdot,\cdot\rangle$ is defined in \eqref{RQuaternionScalerProduct} and $C_g=\|g\|^2_{L^2(\mathbb{R}^2,\mathbb{R})}.$ For more properties and examples of QWFT \eqref{R2012_Fu-Kahler-Cerejeiras_AACA_QWFT_Definition} one can refer \cite{el2020two}. With the window function, $g\in L^2_{\mathbb{H}}(\mathbb{R}^2),$ recently Kamel et al. \cite{kamel2019uncertainty} considered a two sided QWFT for the function $f\in L^2_{\mathbb{H}}(\mathbb{R}^2)$ as
\begin{align}\label{R2019_Kamel-Tefjeni_ITSF_2sQWFT_Definition}
\left(\mathit{G}_gf\right)(\bldomega,\bldb)=\frac{1}{2\pi}\int_{\mathbb{R}^2}e^{-i\omega_1\xx_1}f(\bldx)\overline{g(\bldx-\bldb)}e^{-j\omega_2\xx_2}d\bldx
\end{align}
and expressed the transform in term of the two-sided QFT $\mathcal{F}_s$ \eqref{R2019Chen-KouDefinitionTsQFT} as $\left(\mathit{G}_gf\right)(\bldomega,\bldb)=\left(\mathcal{F}_s\left(f\overline{\tau_{\bldb}g}\right)\right)(\bldomega)$ using which they proved the orthogonality relation for QWFT \eqref{R2019_Kamel-Tefjeni_ITSF_2sQWFT_Definition} depending on two windowed functions $g_1$ and $g_2$ as $\langle\mathit{G}_{g_1}f_1,\mathit{G}_{g_2}f_2\rangle=\langle f_1(\bar{g_1},\bar{g_2}),f_2\rangle,$ where $\langle\cdot,\cdot\rangle$ and $(\cdot,\cdot)$ are respectively given by \eqref{RQuaternionScalerProduct} and \eqref{RQuaternionInnerProduct}. With the assumption that the windowed function $g$ is real valued authors obtained the reconstruction formula for the transform in \eqref{R2019_Kamel-Tefjeni_ITSF_2sQWFT_Definition} as $f(\bldx)=\frac{1}{2\pi\|g\|^2_{L^2(\mathbb{R}^2,\mathbb{R})}}\int_{\mathbb{R}^2}\int_{\mathbb{R}^2}e^{i\xx_1\omega_1}\left(\mathit{G}_gf\right)(\bldomega,\bldb)g(\bldx-\bldb)e^{j\xx_2\omega_2}d\bldomega d\bldb$ and proved that the range $\mathit{G}_gf(L^2_{\mathbb{H}}(\mathbb{R}^2))$ is a RKHS with the kernel $K(\bldomega,\bldb,\bldomega',\bldb')=\frac{2\pi}{\|g\|^2_{L^2(\mathbb{R}^2,\mathbb{R})}}(\mathit{G}_g(\tau_{\bldb}f))(\bldb', \bldomega'-\bldomega).$ Li et al. \cite{li2021some} expressed the two sided QWFT \eqref{R2012_Fu-Kahler-Cerejeiras_AACA_QWFT_Definition}, also given in \cite{el2020two}, in terms of the convolution \eqref{R2014AkilaConvolution1DQFT}  for the function defined on $\mathbb{R}^2$ and proved the linearity and anti-linearity with respect to the real scalars, translation, scaling and parity properties along with the Parseval's formula and with the quaternion window they also obtained the inversion formula and characterized the range of the transform. Brahim et al. \cite{brahim2020uncertainty} further explored the properties of the transform \eqref{R2019_Kamel-Tefjeni_ITSF_2sQWFT_Definition} and proved the inequality 
\begin{align}\label{R2020_Brahim-Tefjeni_JoPDOA_2sQWFT_Lieb-Inequalilty}\|\mathit{G}_gf\|_{L^p_{\mathbb{H}}(\mathbb{R}^2\times\mathbb{R}^2)}\leq (2\pi)^{\frac{2}{p}-1}C_{p,q}\|f\|_{L^2_{\mathbb{H}}(\mathbb{R}^2)}\|g\|_{L^2_{\mathbb{H}}(\mathbb{R}^2)},~C_{p,q}=\left(\frac{4}{p}\right)^{\frac{1}{p}}\left(\frac{1}{q}\right)^{\frac{1}{q}},~\frac{1}{p}+\frac{1}{q}=1,~p\geq 2
\end{align}
called the Lieb's inequality for the QWFT.
\\
In the Cayley-Dickson form representation of the quaternion $q=(a_1+ib_1)+j(a_2+ib_2)$,  two imaginary numbers $i$ and $j$ are involved. In view of the above definitions of the QWFT, it is observed that for each $i$ and $j$, the Fourier kernels $e^{-i\xx_1\omega_1}$ and $e^{-j\xx_2\omega_2}$ are used to define it. Thus, it is necessary to assume that $f$ must be a function of two variables (or even dimension in general). Akila et al.\cite{akila2016multidimensional} gave two definitions of the QWFT for the function defined on $\mathbb{R}^d,$ which are given by
\begin{align}\label{R2016_Akila-Roopkumar_QWFT_Definition1-Convolution}
\left(\mathcal{G}^{(1)}_gf\right)(\bldomega,\bldb)=\left(f_{\bldomega}\ostar\check{\tilde{g}}\right)(\bldb),~(\bldomega,\bldb)\in \mathbb{R}^{2d},
\end{align} 
and 
\begin{align}\label{R2016_Akila-Roopkumar_QWFT_Definition2-Splited}
\left(\mathcal{G}^{(2)}_gf\right)(\bldomega,\bldb)=\left(G_{g_1}f_1+G_{\bar{g}_2}\bar{f}_2+j\left(G_{g_1}f_2+G_{\bar{g}_2}\bar{f}_1\right)\right)(\bldomega,\bldb)
\end{align}
where $\ostar$ denotes the quaternion convolution given by \eqref{R2014AkilaConvolution1DQFT},  $f=f_1+jf_2\in L^2_{\mathbb{H}}(\mathbb{R}^d)$ with $f_{\bldomega}(\bldx)=f(\bldx)e^{-i\bldomega\cdot\bldx},$ $\check{f}(\bldx)=f(-\bldx)$ and $g=g_1+jg_2\in L^2_{\mathbb{H}}(\mathbb{R}^d)$ is window function with $\tilde{g}(\bldx)=\hat{g}_1(\bldx)-j\check{g}_2(\bldx).$ Although both definitions can be defined for the case when $d=2,$ none of these coincides with those discussed earlier. Moreover, they are both consistent with the WFT on $L^2(\mathbb{R}^d,\mathbb{C}).$ The properties like linearity, scaling and pairity given respectively by $\mathcal{G}^{(l)}_g(\alpha f_1+\beta f_2)=\alpha\mathcal{G}^{(l)}_gf_1+\beta\mathcal{G}^{(l)}_gf_2,$ $\mathcal{G}^{(l)}_g(f(\lambda\bldx))=\frac{1}{|\lambda|^d}\left(\mathcal{G}^{(l)}_{g(\lambda\bldx)}f(\bldx)\right)(\frac{\bldomega}{\lambda},\lambda\bldb)$ and $\left(\mathcal{G}^{(l)}_{\check{g}}\check{f}\right)(\bldomega,\bldb)=\left(\mathcal{G}^{(l)}_gf\right)(-\bldomega,-\bldb),$ for $\alpha,\beta\in\mathbb{H},~\lambda\in\mathbb{R}_{0},$ are common to both the transforms \eqref{R2016_Akila-Roopkumar_QWFT_Definition1-Convolution} and \eqref{R2016_Akila-Roopkumar_QWFT_Definition2-Splited}, but the translation formula are different and are given respectively by $\left(\mathcal{G}^{(1)}_g(f(\bldx-\textbf{a}))\right)(\bldomega,\bldb)=\left(\mathcal{G}^{(1)}_{g_1}f\right)(\bldomega,\bldb-\textbf{a})e^{-i\textbf{a}\cdot\bldomega}+\left(\mathcal{G}^{(1)}_{jg_2}f\right)(\bldomega,\bldb+\textbf{a})e^{i\textbf{a}\cdot\bldomega}$ and $\left(\mathcal{G}^{(2)}_g(f(\bldx-\textbf{a}))\right)(\bldomega,\bldb)=\left(\mathcal{G}^{(2)}_gf\right)(\bldomega,\bldb-\textbf{a})e^{-i\textbf{a}\cdot\bldomega}.$ The following Parseval's formula holds for both the transforms 
$$\int_{\mathbb{R}^d}\int_{\mathbb{R}^d}\left(\mathcal{G}^{(l)}_gf\right)(\bldomega,\bldb)\overline{\left(\mathcal{G}^{(l)}_gh\right)(\bldomega,\bldb)}d\bldomega d\bldb=(2\pi)^d\|g\|^2_{L^2_{\mathbb{H}}(\mathbb{R}^d)}(f,h),~l=1,2$$
thus the transforms maps $L^2_{\mathbb{H}}(\mathbb{R}^d)$ into $L^2_{\mathbb{H}}(\mathbb{R}^{2d})$ and is one-one. Furthermore, if the window function $g$ satisfies the condition $(2\pi)^d\int_{\mathbb{R}^d}\left(\hat{g}_1\check{g}_1+\hat{g}_2\check{g}_2\right)=C\neq 0,$ then the following reconstruction formula holds for the transform \eqref{R2016_Akila-Roopkumar_QWFT_Definition1-Convolution}
$$f(\bldx)=\frac{1}{C}\int_{\mathbb{R}^d}\int_{\mathbb{R}^d}\left(\mathcal{G}^{(1)}_gf\right)(\bldomega,\bldb)e^{i\bldomega\cdot\bldx}
\bullet\check{g}(\bldx-\bldb)d\bldomega d\bldb,$$
where $(f_1+jf_2)\bullet(g_1+jg_2)=f_1g_1-\check{\bar{f}}_2g_2+j\left(\check{\bar{f}}_1g_2+f_2g_1\right)$ and for the transform \eqref{R2016_Akila-Roopkumar_QWFT_Definition2-Splited}, the reconstruction formula is given by $f(\bldx)=\frac{1}{(2\pi)^d\|g\|^2_{L^2_{\mathbb{H}}(\mathbb{R}^d)}}\int_{\mathbb{R}^d}\int_{\mathbb{R}^d}\left(\mathcal{G}^{(2)}_gf\right)(\bldomega,\bldb)e^{i\bldomega\cdot\bldx}
g(\bldx-\bldb)d\bldomega d\bldb.$ Authors in \cite{akila2016multidimensional}, also characterized the range of the transform \eqref{R2016_Akila-Roopkumar_QWFT_Definition2-Splited} and proved that $H\in\mathcal{G}^{(2)}_g(L^2_{\mathbb{H}}(\mathbb{R}^d))$ if and only if it satisfies the relation $H(\bldomega',\bldb')=\frac{1}{(2\pi)^d\|g\|^2_{L^2_{\mathbb{H}}(\mathbb{R}^d)}}\int_{\mathbb{R}^d}\int_{\mathbb{R}^d}H(\bldomega,\bldb)\left(g_{\bldomega,\bldb},g_{\bldomega',\bldb'}\right)d\bldomega d\bldb,$ where $g_{\bldomega,\bldb}(\bldx)=e^{i\bldomega\cdot\bldx}g(\bldx-\bldb).$\\

\noindent\textbf{Uncertainty Principles:} Using the Heisenberg UP for the right sided QFT\cite{bahri2008uncertainty} and following the methods adopted by Wilczok \cite{wilczok2000new}, authors in \cite{bahri2010windowed} obtained the Heisenberg UP for the QWFT \eqref{R2010_Bahri-Hitzer-Ashino-Vaillancourt_AMC_RsQWFT_Definition} as
\begin{align*}
\left(\int_{\mathbb{R}^2}\int_{\mathbb{R}^2}\omega_l^2\left|(\mathcal{G}_\phi f)(\bldomega,\bldb)\right|^2d\bldomega\bldb\right)^{\frac{1}{2}}\left(\int_{\mathbb{R}^2}\xx_l|f(\bldx)|^2d\bldx\right)^{\frac{1}{2}}\geq \frac{1}{4\pi}\|f\|^2_{L^2_{\mathbb{H}}(\mathbb{R}^2)}\|\phi\|_{L^2_{\mathbb{H}}(\mathbb{R}^2)},~l=1,2.
\end{align*}
Bahri et al. in \cite{bahri2020uncertainty} considered the function $\phi,f\in \mathcal{S}(\mathbb{R}^2,\mathbb{H})$ and using the Pitt's inequality and logarithmic UP for the right-sided QFT, obtained, respectively, the Pitt's inequality 
$$\int_{\mathbb{R}^2}\int_{\mathbb{R}^2}|\bldomega|^{-\alpha}\left|(\mathcal{G}'_\phi f)(\bldomega,\bldb)\right|^2_Q d\bldomega d\bldb\leq C_\alpha\|\phi\|^2_{L^2(\mathbb{R}^2)}\int_{\mathbb{R}^2}|\bldx|^{\alpha}|f(\bldx)|^2d\bldx,~0\leq\alpha<2$$ 
and the logarithmic UP for the QWFT \eqref{R2010_Bahri-Hitzer-Ashino-Vaillancourt_AMC_RsQWFT_Definition} as 
\begin{align*}
\int_{\mathbb{R}^2}\int_{\mathbb{R}^2}\ln|\bldomega|\left|(\mathcal{G}'_\phi f)(\bldomega,\bldb)\right|^2d\bldomega d\bldb+\|\phi\|^2_{L^2_{\mathbb{H}}(\mathbb{R}^2)}\int_{\mathbb{R}^2}\ln|\bldx||f(\bldx)|^2d\bldx\geq \left(\frac{\Gamma'(\frac{1}{2})}{\Gamma(\frac{1}{2})}-\ln\pi\right)\|f\|^2_{L^2_{\mathbb{H}}(\mathbb{R}^2)}\|\phi\|^2_{L^2_{\mathbb{H}}(\mathbb{R}^2)},
\end{align*}
where $C_\alpha$ and the module $|\cdot|_Q$ are same as in \eqref{R2015ChenPitt'sInequalityQuaternion}.
Fu et al. \cite{fu2012balian} (see also \cite{el2020two}) obtained the Heisenberg UP for the two sided QWFT \eqref{R2012_Fu-Kahler-Cerejeiras_AACA_QWFT_Definition} as 
\begin{align*}
\left(\int_{\mathbb{R}^2}\int_{\mathbb{R}^2}\omega^2_l\left|(\mathtt{Q}_gf)(\bldomega,\bldb)\right|^2d\bldomega d\bldb\right)\left(\int_{\mathbb{R}^2}\xx^2_l|f(\bldx)|^2d\bldx\right)\geq \frac{1}{16\pi^2}\|f\|^4_{L^2_{\mathbb{H}}(\mathbb{R}^2)}\|g\|^2_{L^2(\mathbb{R}^2,\mathbb{R})}
\end{align*}
and also considering the discrete version of the kernel of the inverse QWFT \eqref{R2012_Fu-Kahler-Cerejeiras_AACA_QWFT_Inversion-Formula}, authors established a new Gabor system $\overline{g_{-\boldsymbol m,\boldsymbol n}}:=\mathcal{C}[\cdot,e^{2\pi i m_1\xx_1}]g(\bldx-\boldsymbol n)e^{2\pi j m_2\xx_2},~\bldx\in\mathbb{R}^2,\boldsymbol m=(m_1,m_2),\boldsymbol n\in\mathbb{Z}^2,$ which induces an orthonormal basis with respect to the scalar inner product. Furthermore, they also proved the corresponding Balian-Low theorem, which can be regarded as the strong version of the Heisenberg UP for the Gabor system. More precisely, if the system $\{\overline{g_{-\boldsymbol m,\boldsymbol n}}:\boldsymbol m,\boldsymbol n\in\mathbb{Z}^2\}$ form a Gabor orthonormal basis, then 
$$\left(\int_{\mathbb{R}^2}\xx^2_l|g(\bldx)|^2d\bldx\right)\left(\int_{\mathbb{R}^2}\omega_l|(\mathcal{F}_Qg)(\bldomega)|^2d\bldomega\right)=\infty,~l=1,2,$$ where $\mathcal{F}_Q$ is given by \eqref{RChenDefinitionQFT}. Moreover, the author also pointed out that the Balian-Low theorem also holds for the right-sided QWFT \eqref{R2010_Bahri-Hitzer-Ashino-Vaillancourt_AMC_RsQWFT_Definition} discussed in \cite{bahri2010windowed}. Kamel et al. \cite{kamel2019uncertainty} studied the UPs for the QWFT \eqref{R2019_Kamel-Tefjeni_ITSF_2sQWFT_Definition}. Using the logarithmic UP \eqref{R2015ChenLogarithmicUPQFT} for the two-sided QFT $\mathcal{F}_s$, authors obtained the logarithmic UP for the QWFT as
\begin{align}\label{R2019_Kamel-Tefjeni_ITSP_2sQWFT_LogarithminUP}
\int_{\mathbb{R}^2}\int_{\mathbb{R}^2}\ln|\bldomega||\left(\mathit{G}_gf\right)(\bldomega,\bldb)|^2_Q d\bldomega d\bldb+\|g\|^2_{L^2_{\mathbb{H}}(\mathbb{R}^2)}\int_{\mathbb{R}^2}\ln|\bldx||f(\bldt)|^2d\bldx\geq \left(\ln(2\pi^2)-2\frac{\Gamma'(\frac{1}{2})}{\Gamma(\frac{1}{2})}\right)\|f\|^2_{L^2_{\mathbb{H}}(\mathbb{R}^2)}\|g\|^2_{L^2_{\mathbb{H}}(\mathbb{R}^2)}.
\end{align}
They also proved that for a subset $\Omega\subset \mathbb{R}^2\times\mathbb{R}^2$ satisfying the condition $|\Omega| <(2\pi)^2,$ 
$$\|f\|^2_{L^2_{\mathbb{H}}(\mathbb{R}^2)}\|g\|^2_{L^2_{\mathbb{H}}(\mathbb{R}^2)}\sqrt{1-\frac{|\Omega|}{(2\pi)^2}}\leq \left(\int_{\mathbb{R}^2}\int_{\mathbb{R}^2}\chi_{\Omega^c}(\bldomega,\bldb)\left|\left(\mathit{G}_gf\right)(\bldomega,\bldb)\right|^2 d\bldomega d\bldb\right)^{\frac{1}{2}},$$
which says that the portion of $\mathit{G}_gf$ lying outside some set $\Omega$ cannot be arbitrary small, followed by the local Price's inequality which says that for $p>1,$ $0<\epsilon<2$ and a measurable subset $\Omega_0$ of $\mathbb{R}^2\times\mathbb{R}^2$ there exists a constant $M_{\epsilon,p}$ depending on $p$ and $\epsilon$ such that 
$$
\|\chi_{\Omega_0}\mathit{G}_gf\|^{p(p+1)}_{L^p_{\mathbb{H}}(\mathbb{R}^2\times\mathbb{R}^2)}\leq M_{\epsilon,p}|\Omega_0|\||(\bldomega,\bldb)|^\epsilon\mathit{G}_gf\|^{p(p+1)}_{L^p_{\mathbb{H}}(\mathbb{R}^2\times\mathbb{R}^2)}\left(\|f\|_{L^2_{\mathbb{H}}(\mathbb{R}^2)}\|g\|_{L^2_{\mathbb{H}}(\mathbb{R}^2)}\right)^{p\left(p-\frac{2-\epsilon}{2+\epsilon}\right)}.
$$
The logarithmic UP for the QWFT \eqref{R2012_Fu-Kahler-Cerejeiras_AACA_QWFT_Definition}, with the module $|\cdot|_Q$ replaced with the quaternion module $|\cdot|$ in \eqref{R2019_Kamel-Tefjeni_ITSP_2sQWFT_LogarithminUP} can be found in \cite{el2020two}. Using Lieb's inequality \eqref{R2020_Brahim-Tefjeni_JoPDOA_2sQWFT_Lieb-Inequalilty}, authors in \cite{brahim2020uncertainty}, proved that if $f\neq 0,$ $\Omega \subset \mathbb{R}^2\times\mathbb{R}^2,~\epsilon \geq 0$ and $ $ is $\epsilon-$ concentrated on $\Omega$ then for every $p>2,$ $\left(C_{p,q}\right)^{\frac{2p}{2-p}}(1-\epsilon^2)^{\frac{p}{p-2}}\leq \frac{|\Omega|}{(2\pi)^2},$ which is Lieb's UP for the QWFT \eqref{R2019_Kamel-Tefjeni_ITSF_2sQWFT_Definition}. They also proved the Beckner's UP in terms of entropy as $E\left(\left|\mathit{G}_gf\right|^2\right)\geq \|f\|^2_{L^2_{\mathbb{H}}(\mathbb{R}^2)}\|g\|^2_{L^2_{\mathbb{H}}(\mathbb{R}^2)}\left(2\ln(2)-\ln\left(\frac{\|f\|^2_{L^2_{\mathbb{H}}(\mathbb{R}^2)}\|g\|^2_{L^2_{\mathbb{H}}(\mathbb{R}^2)}}{(2\pi)^2}\right)\right)$
using which they obtained the generalized Heisenberg's UP for the QWFT \eqref{R2019_Kamel-Tefjeni_ITSF_2sQWFT_Definition} as
\begin{align}
\left(\int_{\mathbb{R}^2}\int_{\mathbb{R}^2}|\bldomega|^p\left|\left(\mathit{G}_gf\right)(\bldomega,\bldb)\right|^2d\bldomega d\bldb\right)^{\frac{q}{p+q}}\left(\int_{\mathbb{R}^2}\int_{\mathbb{R}^2}|\bldb|^q\left|\left(\mathit{G}_gf\right)(\bldomega,\bldb)\right|^2d\bldomega d\bldb\right)^{\frac{p}{p+q}}\geq D_{p,q}\|f\|^2_{L^2_{\mathbb{H}}(\mathbb{R}^2)}\|g\|^2_{L^2_{\mathbb{H}}(\mathbb{R}^2)},
\end{align}
where $D_{p,q}=\frac{2}{p^{\frac{q}{p+q}}q^{\frac{p}{P+q}}}e^{\frac{pq\left(2\ln(2)+\ln\left(\frac{pq}{\Gamma(\frac{2}{p})\Gamma(\frac{2}{q})}\right)\right)}{2(p+q)}-1}$ and $p,q$ are positive real numbers.. Akila et al. \cite{akila2016multidimensional} obtained the Heisenberg UPs for the QWFT \eqref{R2016_Akila-Roopkumar_QWFT_Definition1-Convolution} and \eqref{R2016_Akila-Roopkumar_QWFT_Definition2-Splited} as
$$\left(\int_{\mathbb{R}^d}\int_{\mathbb{R}^d}\omega_k^2\left|\left(\mathcal{G}^{(l)}_gf\right)(\bldomega,\bldb)\right|^2d\bldomega d\bldb\right)^{\frac{1}{2}}\left(\int_{\mathbb{R}^d}\xx_k^2|f(\bldx)|^2\right)^{\frac{1}{2}}\geq \frac{(2\pi)^{\frac{d}{2}}}{2}\|f\|^2_{L^2_{\mathbb{H}}(\mathbb{R}^d)}\|g\|_{L^2_{\mathbb{H}}(\mathbb{R}^d)},~l=1,2$$

\section{Quaternion Wavelet transform (QWT)}
Let $\mathcal{G}$ denotes the similitude group $SIM(2)$ on $\mathbb{R}^2$ defined as $\mathcal{G}=\mathbb{R}^+\times SO(2)\times\mathbb{R}^2=\{(\A,r_{\theta},\bldb):\A\in\mathbb{R}^+,r_{\theta}\in SO(2),\bldb\in\mathbb{R}^2\},$ $SO(2)$ being special orthogonal group of $\mathbb{R}^2.$ The group law associated with $\mathcal{G}$ is given by 
$(\A,\bldb,r_{\theta})(\A',\bldb',r_{\theta'})=(\A\A',\bldb+\A r_{\theta}\bldb',r_{\theta+\theta'}).$ The rotation $r_{\theta}\in SO(2)$ acts on $\bldx\in\mathbb{R}^2$ as $r_{\theta}(\bldx)=(\xx_1\cos\theta-\xx_2\sin\theta,\xx_1\sin\theta+\xx_2\cos\theta),~0\leq \theta\leq 2\pi$ and the left Haar measure on $\mathcal{G}$ is given by $d\lambda(\A,\theta,\bldb)=d\mu(\A,\theta)d\bldb,$ where $d\mu(\A,\theta)=\frac{d\A d\theta}{\A^3}$ and $d\theta$ is the Haar measure on $SO(2).$ 

A function $\psi\in L^2_{\mathbb{H}}(\mathbb{R}^2)$ is said to be an admissible quaternion wavelet if it satisfies the admissibility condition $c_\psi:=\int_{SO(2)}\int_{\mathbb{R}^+}\left|\hat{\psi}\left(\A r_{-\theta}(\bldomega)\right)\right|^2\frac{d\A d\theta}{\A}$ is a real positive constant independent of $\bldomega$ satisfying $|\bldomega|=1,$ and hence independent of $\bldomega\in \mathbb{R}^2.$ Bahri et al. \cite{bahri2011two}, defined the unitary linear operator $U_{\A,\theta,\bldb}:L^2_{\mathbb{H}}(\mathbb{R}^2)\rightarrow L^2_{\mathbb{H}}(\mathcal{G})$ as $U_{\A,\theta,\bldb}(\psi)=\psi_{\A,\theta,\bldb}(\bldx)=\frac{1}{\A}\psi\left(r_{-\theta}\left(\frac{\bldx-\bldb}{\A}\right)\right),~\A\in\mathbb{R}^+,\bldb\in\mathbb{R}^2$ and $r_{-\theta}\in SO(2)$ and called the family $\psi_{\A,\theta,\bldb}$ as the quaternion daughter wavelets. They defined the QWT of $f\in L^2_{\mathbb{H}}(\mathbb{R}^2)$ with respect to the quaternion admissible wavelet $\psi$ as 
\begin{align}\label{R2011_Bahri-Ashino-Vaillancourt_AMC_QWT_Definition}
(\mathcal{W}_\psi f)(\A,\theta,\bldb)=(f,\psi_{\A,\theta,\bldb})=\int_{\mathbb{R}^2}f(\bldx)\frac{1}{\A}\overline{\psi\left(r_{-\theta}\left(\frac{\bldx-\bldb}{\A}\right)\right)}d\bldx.
\end{align} 
Here, $\A$ is a dilation parameter, $\bldb$ a translation parameter and $\theta$ an $SO(2)$ rotation parameter. It is to be noted that the order of the terms on \eqref{R2011_Bahri-Ashino-Vaillancourt_AMC_QWT_Definition} is fixed because of the non-commutativity of the quaternion multiplication. Moreover, with the assumption that $(\mathfrak{F}_q\psi_l)(\bldomega),~l=0,1,2,3$ are real-valued, the transform \eqref{R2011_Bahri-Ashino-Vaillancourt_AMC_QWT_Definition}, has the quaternion Fourier representation of the form 
\begin{align}\label{R2011_Bahri-Ashino-Vaillancourt_AMC_QWT_Definition-QFDomain}
(\mathcal{W}_\psi f)(\A,\theta,\bldb)=\frac{\A}{(2\pi)^2}\int_{\mathbb{R}^2}(\mathfrak{F}_qf)(\bldomega)e^{i\B_2\omega_1}\left[\overline{(\mathfrak{F}_q\psi_{0l})(\A r_{-\theta}(\bldomega))}e^{i\B_1\omega_1}+\overline{(\mathfrak{F}_q\psi_{1l})(\A r_{-\theta}(\bldomega))}e^{-i\B_1\omega_1}\right]d\bldomega,
\end{align}
where $(\mathfrak{F}_q\psi_{0l})(\A r_{-\theta}(\bldomega))=(\mathfrak{F}_q\psi_0)(\A r_{-\theta}(\bldomega))+i(\mathfrak{F}_q\psi_1)(\A r_{-\theta}(\bldomega)),$ $(\mathfrak{F}_q\psi_{1l})(\A r_{-\theta}(\bldomega))=j(\mathfrak{F}_q\psi_2)(\A r_{-\theta}(\bldomega))+k(\mathfrak{F}_q\psi_3)(\A r_{-\theta}(\bldomega))$ and $(\mathfrak{F}_qf)(\bldomega)=\int_{\mathbb{R}^2}f(\bldx)e^{-i\xx_1\omega_1}e^{-j\xx_2\omega_2}d\bldx.$ The transform \eqref{R2011_Bahri-Ashino-Vaillancourt_AMC_QWT_Definition} is linear in $f$ and anti-linear in $\psi$ with respect to the quaternion scalars and satisfies the following covariance properties $(\mathcal{W}_\psi f(\cdot-\bldx_0))(\A,\theta,\bldb)=(\mathcal{W}_\psi f)(\A,\theta,\bldb-\bldx_0),$ $(\mathcal{W}_\psi f(c\cdot))(\A,\theta,\bldb)=\frac{1}{c}(\mathcal{W}_\psi f)(\A c,\theta,\bldb c),~c>0$ and $(\mathcal{W}_\psi f(r_{\theta_0}\cdot))(\A,\theta,\bldb)=(\mathcal{W}_\psi f)(\A,\theta',r_{\theta_0}\bldb),~r_{\theta'}=r_{\theta_0}r_{\theta}$ which are respectively called the translation, dilation and rotational covariance properties. Similarly, results can be obtained with changes in the admissible quaternion wavelet $\psi.$ Apart from these basic properties authors in \cite{bahri2011two} obtained the inner product relation 
\begin{align}\label{R2011_Bahri-Ashino-Vaillancourt_AMC_QWT_IPR}
\int_{SO(2)}\int_{\mathbb{R}^+}\left(\int_{\mathbb{R}^2}(\mathcal{W}_\psi f)(\A,\theta,\bldb)\overline{(\mathcal{W}_\psi g)(\A,\theta,\bldb)}d\bldb\right)\frac{d\A d\theta}{\A^3}=c_\psi(f,g)
\end{align}
for $f,g\in L^1_{\mathbb{H}}(\mathbb{R}^2)\cap L^2_{\mathbb{H}}(\mathbb{R}^2)$ with the assumption that $\mathfrak{F}_q\psi\in L^2(\mathbb{R}^2,\mathbb{R}).$ Bahri \cite{bahri2011quaternion} gave an alternative proof for \eqref{R2011_Bahri-Ashino-Vaillancourt_AMC_QWT_IPR}. In particular if $f=g,$ then $\|\mathcal{W}_\psi f\|^2_{L^2_{\mathbb{H}}(\mathcal{G})}=c_\psi\|f\|^2_{L^2_{\mathbb{H}}(\mathbb{R}^2)},$ as a result of which it follows that, except for the factor $c_\psi,$ the transform $\mathcal{W}_\psi$ is an isometry from $L^2_{\mathbb{H}}(\mathbb{R}^2)$ to $L^2_{\mathbb{H}}(\mathcal{G}).$ With similar assumptions on $\psi$ as in \eqref{R2011_Bahri-Ashino-Vaillancourt_AMC_QWT_IPR}, the reconstruction formula for the transform \eqref{R2011_Bahri-Ashino-Vaillancourt_AMC_QWT_Definition} is given by
$f(\bldx)=\frac{1}{c_\psi}\int_{\mathcal{G}}(\mathcal{W}_\psi f)(\A,\theta,\bldb)\psi_{\A,\theta,\bldb}(\bldx)d\lambda(\A,\theta,\bldb),$
where the integral converges in the weak sense. Moreover, the range of the transform is a RKHS with the kernel $K_{\psi}(\A,\theta,\bldb;\A',\theta',\bldb')=\frac{1}{c_\psi}(\psi_{\A,\theta,\bldb},\psi_{\A',\theta',\bldb'}).$ Bahri et al. \cite{bahri2012two},\cite{bahri2014continuous} studied the properties of the QWT \eqref{R2011_Bahri-Ashino-Vaillancourt_AMC_QWT_Definition} using the machinery of the QFT given by $(\mathbf{F}f)(\bldomega)=\int_{\mathbb{R}^2}f(\bldx)e^{\frac{I+j+k}{\sqrt{3}}\bldomega\cdot\bldx}d\bldomega$ (see \cite{bahri2014continuous} for more details). With the assumption that the functions $\{\psi,\phi\}$ in $L^2_{\mathbb{H}}(\mathbb{R}^2)$ is an admissible quaternion wavelet pair, i.e., $C_{\{\psi,\phi\}}:=\int_{SO(2)}\int_{\mathbb{R}^+}\overline{(\mathbf{F}\psi)(\A r_{\theta}(\bldomega))}(\mathbf{F}\phi)(\A r_{\theta}(\bldomega))\frac{d\A d\theta}{\A}$ is a non-zero quaternion constant independent of $|\bldomega|=1,$ and $\phi$ and $\psi$ satisfies the conditions $e^{\frac{i+j+k}{\sqrt{3}}}(\mathbf{F}\phi)(\A r_{\theta}(\bldomega))=(\mathbf{F}\phi)(\A r_{\theta}(\bldomega))=e^{\frac{i+j+k}{\sqrt{3}}}$ and $e^{\frac{i+j+k}{\sqrt{3}}}(\mathbf{F}\psi)(\A r_{\theta}(\bldomega))=(\mathbf{F}\psi)(\A r_{\theta}(\bldomega))=e^{\frac{i+j+k}{\sqrt{3}}}$ respectively,  authors obtained the following inner product relation 
$$\int_{SO(2)}\int_{\mathbb{R}^+}\left(\int_{\mathbb{R}^2}(\mathcal{W}_\psi f)(\A,\theta,\bldb)\overline{(\mathcal{W}_\phi g)(\A,\theta,\bldb)}d\bldb\right)\frac{d\A d\theta}{\A^3}=(fC_{\{\psi,\phi\}},g)$$
followed by the inversion formula $f(\bldx)=\int_{\mathcal{G}}(\mathcal{W}_\psi f)(\A,\theta,\bldb)\phi_{\A,\theta,\bldb}(\bldx)C_{\{\psi,\phi\}}^{-1}\frac{d\A d\theta}{\A^3}d\bldb,$ where the integral converges in the weak sense. Tefjeni et al. \cite{tefjeni2020uncertainty} considered an admissible quaternion wavelet pair as function $\phi,\psi\in L^2_{\mathbb{H}}(\mathbb{R}^{2d})$ such that 
$$C'_{\psi,\phi}:=\int_0^\infty\overline{(\mathfrak{F}_q\psi)(\A\bldomega)}(\mathfrak{F}_q\phi)(\A\bldomega)\frac{d\A}{\A}$$
is a non-zero quaternion constant independent of $\bldomega$ satisfying $|\bldomega|=1.$ In particular, a function $\psi\in L^2_{\mathbb{H}}(\mathbb{R}^{2d})$ is said to be admissible quaternion wavelet if $C'_{\psi}=C'_{\psi,\psi}$ is a positive constant independent of $\bldomega$ satisfying $|\bldomega|=1.$ With the assumption that the admissible quaternion wavelet $\psi$ (and also $\phi$) satisfies the conditions
\begin{align}\label{R2020_Tefjini-Brahim_ITSF_RsQWT_Condition1}
\psi=\psi_0+i\psi_1~\mbox{and}~ e^{j\B_2\omega_2}\overline{(\mathfrak{F}_q\psi)(\A\bldomega)}=\overline{(\mathfrak{F}_q\psi)(\A\bldomega)}e^{j\B_2\omega_2}
\end{align}
or 
\begin{align}\label{R2020_Tefjini-Brahim_ITSF_RsQWT_Condition2}
\psi=j\psi_2+k\psi_3~\mbox{and}~ e^{j\B_2\omega_2}\overline{(\mathfrak{F}_q\psi)(\A\bldomega)}=\overline{(\mathfrak{F}_q\psi)(\A\bldomega)}e^{-j\B_2\omega_2}
\end{align}
authors in \cite{tefjeni2020uncertainty} proved that the QWT for the function $f\in L^2_{\mathbb{H}}(\mathbb{R}^2),$ given by   
\begin{align}\label{R2020_Tefjini-Brahim_ITSF_RsQWT_Definition}
(T_\psi f)(\A,\bldb):=\int_{\mathbb{R}^{2d}}f(\bldx)\frac{1}{\A^d}\overline{\psi\left(\frac{\bldx-\bldb}{\A}\right)}d\mu_{2d}(\bldx),
\end{align}
where $d\mu_d(\bldx)=\frac{d\bldx}{(2\pi)^{\frac{d}{2}}},$ satisfies the Parseval's theorem 
%$$\langle T_\psi f,T_\phi g \rangle_{2,G}=\langle f C'_{\psi,\phi},g \rangle_{2,2d}$$
$$\int_0^\infty\int_{\mathbb{R}^{2d}} (T_\psi f)(a,\bldb)\overline{(T_\phi g)(\A,\bldb)}\frac{d\A}{\A^{2d+1}}d\mu_{2d}(\bldb) =\int_{\mathbb{R}^{2d}}f(\bldx) C'_{\psi,\phi}\overline{g(\bldx)} d\mu_{2d}(\bldx)$$
followed by the reconstruction formula $f(\bldx)=\{C'_{\psi}\}^{-1}\int_0^\infty\int_{\mathbb{R}^{2d}}(T_\psi f)(\A,\bldb)\frac{1}{\A^d}\psi\left(\frac{\bldx-\bldb}{\A}\right)\frac{d\A}{\A^{2d+1}}d\mu_{2d}(\bldb).$ 
%$$\|(T_\psi f)(T_\phi g)\|_{p,G}\leq \left(\sqrt{C'_\phi C'_\psi}\right)^{\frac{1}{p}}\left(\|\phi\|_{2,2d}\|\psi\|_{2,2d}\right)^{1-\frac{1}{p}}\|f\|_{2,2d}\|g\|_{2,2d}$$
If $\phi,\psi$ are two admissible quaternion wavelet satisfying \eqref{R2020_Tefjini-Brahim_ITSF_RsQWT_Condition1} or \eqref{R2020_Tefjini-Brahim_ITSF_RsQWT_Condition2} and $p\geq 1,$ then for every $f,g\in L^2_{\mathbb{H}}(\mathbb{R}^{2d})$ the function $(\A,\bldb)\rightarrow (T_\psi f)(a,\bldb)(T_\phi g)(a,\bldb)$ belongs to $L^p_{\mathbb{H}}(\mathbb{R}^+\times\mathbb{R}^{2d})$ and satisfies the Lieb's inequality 
\begin{align*}
&\left(\int_0^\infty\int_{\mathbb{R}^{2d}}|(T_\psi f)(\A,\bldb)(T_\phi g)(\A,\bldb)|^p\frac{d\A}{\A^{2d+1}}d\mu_{2d}(\bldb)\right)^{\frac{1}{p}}\notag\\
&\leq \left(\sqrt{C'_\phi C'_\psi}\right)^{\frac{1}{p}}\left\{\left(\int_{\mathbb{R}^{2d}}|\psi(\bldx)|^2d\mu_{2d}(\bldx)\right)\left(\int_{\mathbb{R}^{2d}}|\phi(\bldx)|^2d\mu_{2d}(\bldx)\right)\right\}^{\frac{1}{2}\left(1-\frac{1}{p}\right)}\left(\int_{\mathbb{R}^{2d}}|f(\bldx)|^2d\mu_{2d}(\bldx)\right)^{\frac{1}{2}}\left(\int_{\mathbb{R}^{2d}}|g(\bldx)|^2d\mu_{2d}(\bldx)\right)^{\frac{1}{2}}.
\end{align*}
Moreover, if $\psi$ be an admissible quaternion wavelet satisfying the assumptions \eqref{R2020_Tefjini-Brahim_ITSF_RsQWT_Condition1} or \eqref{R2020_Tefjini-Brahim_ITSF_RsQWT_Condition2}, then for every $f\in L^2_{\mathbb{H}}(\mathbb{R}^{2d})$ the transform $T_\psi f$ given by \eqref{R2020_Tefjini-Brahim_ITSF_RsQWT_Definition} belong to $L^p_{\mathbb{H}}(\mathbb{R}^+\times\mathbb{R}^{2d})$ for $2\leq p\leq\infty$ and for $p\geq 2,$ it satisfies
%$$\|T_\psi f\|_{p,G}\leq (C'_\psi)^{\frac{1}{p}}\|\psi\|^{1-\frac{2}{p}}_{2,2d}\|f\|_{2,2d}$$
$$\left(\int_0^\infty\int_{\mathbb{R}^{2d}}|(T_\psi f)(\A,\bldb)|^p\frac{d\A}{\A^{2d+1}}d\mu_{2d}(\bldb)\right)^{\frac{1}{p}}
\leq (C'_\psi)^{\frac{1}{p}}\left(\int_{\mathbb{R}^{2d}}|\psi(\bldx)|^2d\mu_{2d}(\bldx)\right)^{1-\frac{2}{p}}\left(\int_{\mathbb{R}^{2d}}|f(\bldx)|^2d\mu_{2d}(\bldx)\right)^{\frac{1}{2}},$$
as a consequence of which it follows that,  if $T_\psi f,$ for $f\neq 0,$ is $\epsilon-$concentrated on a measurable set $\Omega\subset\mathbb{R}^+\times\mathbb{R}^{2d}=G$ of non zero finite measure $\mu_G(\Omega),$ where $\mu_G$ is the measure defined on $G$ by $dw(a,\bldb)=\frac{da}{a^{2d+1}}d\mu_{2d}(\bldb),$ then for every $p>2,$ $\mu_G(\Omega)\int_{\mathbb{R}^{2d}}|\psi(\bldx)|^2d\mu_{2d}(\bldx)\geq C'_\psi(1-\epsilon^2)^{\frac{p}{p-2}}.$
%$$\mu_G(\Omega)\|\psi\|^2_{2,2d}\geq C'_\psi(1-\epsilon^2)^{\frac{p}{p-2}},~p>2$$
This is the Lieb's UP associated to the QWT \eqref{R2020_Tefjini-Brahim_ITSF_RsQWT_Definition}. Apart from this author in \cite{tefjeni2020uncertainty} also obtained the associated Beckner's UP in terms of entropy as 
$$E\left(|T_\psi f|^2\right)\geq C'_\psi\left(\int_{\mathbb{R}^{2d}}|f(\bldx)|^2d\mu_{2d}(\bldx)\right)\ln\left(\frac{1}{\left(\int_{\mathbb{R}^{2d}}|f(\bldx)|^2d\mu_{2d}(\bldx)\right)\left(\int_{\mathbb{R}^{2d}}|f(\bldx)|^2d\mu_{2d}(\bldx)\right)}\right),$$
%$$E\left(|T_\psi f|^2\right)\geq C'_\psi\|f\|^2_{2,2d}\ln\left(\frac{1}{\|f\|_{2,2d}^2\|\psi\|_{2,2d}^2}\right),$$
where $E\left(|T_\psi f|^2\right)=-\int_0^\infty\int_{\mathbb{R}^{2d}}\ln\left(|(T_\psi f)(\A,\bldb)|^2\right)|(T_\psi f)(\A,\bldb)|^2\frac{d\A}{\A^{2d+1}}d\mu_{2d}(\bldb).$ Using the Heisenberg's and logarithmic UPs for the QFT $\mathfrak{F}_q,$ Bahri et al. \cite{bahri2017variation} obtained the Heisenberg and logarithmic UPs for the QWT \eqref{R2020_Tefjini-Brahim_ITSF_RsQWT_Definition} for the case when $d=1.$ So far, we have discussed the QWT for the function with the variable in 2 or even dimension. The QWT for the function $f\in L^2_{\mathbb{H}}(\mathbb{R})$ can be found in \cite{he2004continuous}. He et al.\cite{he2004continuous} decomposed the space $L^2_{\mathbb{H}}(\mathbb{R})=H^{(+,+)}\oplus H^{(-,-)},$ where $H^{(+,+)}=\{f=f_1+jf_2:f_1,f_2\in H^2_+(\mathbb{R})\}$ and $H^{(-,-)}=\{f=f_1+jf_2:f_1,f_2\in H^2_-(\mathbb{R})\}$ with the Hardy space $H^2_+(\mathbb{R})=\{f\in L^2(\mathbb{R}):\hat{f}(\xxi)=0, if~\xxi\nin\mathbb{R}^+\}$ and the conjugate Hardy space $H^2_-(\mathbb{R})=\{f\in L^2(\mathbb{R}):\hat{f}(\xxi)=0, if~\xxi\nin\mathbb{R}^-\}.$ They obtained the characterization of the quaternion admissible wavelet for the function $\phi\in H^{(+,+)}\left(\mbox{or}~H^{(-,-)}\right)$ as $c_\phi=\int_{\mathbb{R}}\frac{|\hat{\phi}(\xxi)|^2}{|\xxi|}d\xxi<\infty$ and studied the properties of the QWT $(\mathcal{W}_\phi f)(\A,\B)=\int_{\mathbb{R}}f(\xx)\left(\frac{1}{\sqrt{\A}}\phi\left(\frac{\xx-\B}{\A}\right)\right)^cd\xx$ for the function $f\in H^{(+,+)}\left(\mbox{or}~H^{(-,-)}\right)$ and the admissible quaternion wavelet $\phi\in H^{(+,+)}\left(\mbox{or}~H^{(-,-)}\right).$  With any one of the assumptions, for the admissible quaternion wavelets $\phi=\phi_1+j\phi_2,\psi=\psi_1+j\psi_2\in H^{(+,+)}$ and $f=f_1+jf_2,g=g_1+jg_2\in H^{(+,+)},$ among the following
\begin{enumerate}
\item $\bar{\hat{\phi}}_l(\xxi)\bar{\hat{\psi}}_l(\xxi),~l=1,2$ are real functions.
\item $f=f_1,~g=g_1$ or $g=jg_2$
\item $f=jf_2,~g=jg_2$ or $g=g_1,$
\end{enumerate}
authors in \cite{he2004continuous} obtained the inner product relation $\int_{\mathbb{R}^+}\int_{\mathbb{R}}(\mathcal{W}_\phi f)(\A,\B)\left((\mathcal{W}_\phi f)(\A,\B)\right)^c\frac{d\A d\B}{\A^2}=c_{\phi,\psi}\int_{\mathbb{R}}f(\xx)(g(\xx))^cd\xx,$ where  $c_{\phi,\psi}=\int_{\mathbb{R}}\left(\bar{\hat{\phi}}_1(\xxi)\hat{\psi}(\xxi)+\bar{\hat{\phi}}_2(\xxi)\hat{\psi}_2(\xxi)\right)\frac{d\xxi}{|\xxi|},$ which was not in general true for the an arbitrary quaternion wavelet $\phi$ and functions $f$ in $H^{(+,+)}.$ They also obtained the following reproducing formula for the admissible quaternion wavelet $\phi$ and $f\in H^{(+,+)}$ 
$$f(\xx)=\frac{1}{c_\phi}\int_{\mathbb{R}}\int_{\mathbb{R}}^+(\mathcal{W}_\phi f)(\A,\B)\frac{1}{\sqrt{\A}}\phi\left(\frac{\xx-\B}{\A^2}\right)\frac{d\xx d\A}{\A^2}$$ 
with the convergence of the integral in the weak sense. Akila et al. \cite{akila2014natural} introduced the QWT with respect to $\phi=\phi_1+j\phi_2\in L^1_{\mathbb{H}}(\mathbb{R})\cap L^2_{\mathbb{H}}(\mathbb{R})$ satisfying $c_\phi=\int_{\mathbb{R}}\frac{|\hat{\phi_1}(\xxi)|^2+|\hat{\phi_2}(\xxi)|^2}{\xxi}d\xxi<\infty$ as
\begin{align}\label{R2014_Akila-Roopkumar_1DQWT_Definition}
(\mathcal{W}_\phi f)(\A,\B)=(f\ostar \check{\tilde{\phi}}_{\A})(\B),~(\A,\B)\in \mathbb{R}^+\times\mathbb{R}
\end{align}
where $\tilde{\phi}=\bar{\phi}_1-j\check{\phi}_2,$ $\phi_{\A}(\xx)=\frac{1}{\sqrt{\A}}\phi\left(\frac{1}{\A}\right),$ $f\in L^2_{\mathbb{H}}(\mathbb{R})$ and $\ostar$ denotes the quaternion convolution given by \eqref{R2014AkilaConvolution1DQFT}. The definition \eqref{R2014_Akila-Roopkumar_1DQWT_Definition} is consistent with the definition of the classical WT on $L^2_{\mathbb{C}}(\mathbb{R}).$ Using the convolution theorem associated with the convolution $\ostar$ \eqref{R2014AkilaConvolution1DQFT}, authors obtained the Parseval's formula for the QWT \eqref{R2014_Akila-Roopkumar_1DQWT_Definition} for the function $f,g\in L^2_{\mathbb{H}}(\mathbb{R})$ as 
\begin{align}\label{R2014_Akila-Roopkumar_1DQWT_IPR}
\int_{\mathbb{R}}\int_{\mathbb{R}^+}(\mathcal{W}_\phi f)(\A,\B)\left((\mathcal{W}_\phi g)(\A,\B)\right)^c\frac{d\A d\B}{\A^2}=c_\phi\int_{\mathbb{R}}f(\xx)(g(\xx))^cd\xx,
\end{align}
followed by the inversion formula $f(\xx)=\frac{1}{c_\phi}\int_{\mathbb{R}^+}\left(\left(f\ostar\check{\hat{\phi}}_{\A}\right)\ostar\phi_{\A}\right)(\xx)\frac{d\A}{\A^2}$ in $L^2_{\mathbb{H}}(\mathbb{R}).$ If the WT of $f\in L^2_{\mathbb{H}}(\mathbb{R})$ is defined as $( \check{\tilde{\phi}}_{\A}\ostar f)(\B),$ then the definition results in the property $\mathcal{W}_\phi(f\ostar g)=\left(\mathcal{W}_\phi f\right)\ostar g$ for $f,g\in L^2_{\mathbb{H}}(\mathbb{R})$ which was not the case when we take the definition as \eqref{R2014_Akila-Roopkumar_1DQWT_Definition}. But for definition $(\check{\tilde{\phi}}_{\A}\ostar f)(\B),$ the Parseval's formula and the inversion formula cannot be achieved as in the case of the transform \eqref{R2014_Akila-Roopkumar_1DQWT_Definition}.

\section{Quaternion Stockwell Transform (QST)}
Akila et al. \cite{akila2016quaternionic} introduced the QST for a function $f=f_1+jf_2\in L^2_{\mathbb{H}}(\mathbb{R})$ with respect to $\phi=\phi_1+j\phi_2\in L^1_{\mathbb{H}}(\mathbb{R})\cap L^2_{\mathbb{H}}(\mathbb{R})$ as
\begin{align}\label{R2016_Akila-Roopkumar_1DQST_Definition}
(\mathcal{S}_\phi f)(\A,\B)=\left([f]_{\A}\ostar \mathcal{M}_{\A}\check{\tilde{\phi}}\right)(\B),~(\A,\B)\in \mathbb{R}^+\times\mathbb{R}
\end{align}
where $[f]_{\A}(\xx)=f(\xx)e^{-i\A\xx},$ $(\mathcal{M}_{\A}f)(\xx)=|\A|f(\A \xx),$ $\tilde{\phi}=\bar{\phi}_1-j\check{\phi}_2,$ $f\in L^2_{\mathbb{H}}(\mathbb{R})$ and $\ostar$ denotes the quaternion convolution given by \eqref{R2014AkilaConvolution1DQFT}. The definition \eqref{R2016_Akila-Roopkumar_1DQST_Definition} can be equivalently written as 
\begin{align}\label{R2016_Akila-Roopkumar_1DQST_Definition_Splitting-Form}
(\mathcal{S}_\phi f)(\A,\B)=\left(\mathcal{S}_{\phi_1} f_1+\mathcal{S}_{\check{\bar{\phi}}_2} \check{\bar{f}}_2+j\left(\mathcal{S}_{\phi_1} f_2-\mathcal{S}_{\check{\bar{\phi}}_2} \check{\bar{f}}_1\right)\right)(\A,\B)
\end{align}
In view of equation \eqref{R2016_Akila-Roopkumar_1DQST_Definition_Splitting-Form}, it follows that the definition \eqref{R2016_Akila-Roopkumar_1DQST_Definition} is consistent with the definition of the classical WT on $L^2_{\mathbb{C}}(\mathbb{R}).$ The transform is linear in $f$ with the scalars in $\mathbb{H}$ and has the translation property
$$(\mathcal{S}_\phi (\tau_{\B_0}f))(\A,\B)=\left(\mathcal{S}_{\phi_1}(f_1+jf_2)\right)(\A,\B-\B_0)-j\left(\mathcal{S}_{\check{\hat{\phi}}_2}(\check{\hat{f}}_2+j\check{\hat{f}}_2)\right)(\A,\B+\B_0)$$
and the scaling property $(\mathcal{S}_\phi (f(\lambda\xx)))(\A,\B)=(\mathcal{S}_\phi f)(\frac{\A}{\lambda},\lambda\B),$ for $\lambda\in \mathbb{R}_{0}.$ Using the convolution theorem associated with the convolution $\ostar$ \eqref{R2014AkilaConvolution1DQFT} and the condition that $c_\phi:=\int_{\mathbb{R}^\star}\frac{|\hat{\phi}(\xxi-1)|^2}{|\xxi|}d\xxi$ in non-zero finite, authors obtained the Parseval's formula for the QST \eqref{R2016_Akila-Roopkumar_1DQST_Definition} for the function $f,g\in L^2_{\mathbb{H}}(\mathbb{R})$ as 
\begin{align}\label{R2014_Akila-Roopkumar_1DQWT_IPR}
\int_{\mathbb{R}}\int_{\mathbb{R}^\star}(\mathcal{S}_\phi f)(\A,\B)\left((\mathcal{S}_\phi g)(\A,\B)\right)^c\frac{d\A d\B}{|\A|}=2\pi c_\phi\int_{\mathbb{R}}f(\xx)(g(\xx))^cd\xx,
\end{align}
using which they obtained the inversion formula $f(\xx)=\frac{1}{2\pi c_\phi}\int_{\mathbb{R}^\star}\int_{\mathbb{R}}(\mathcal{S}_\phi f)(\A,\B)\Phi_{\A,\B}(\xx)\frac{d\B d\A}{|\A|}$ in $L^2_{\mathbb{H}}(\mathbb{R}),$ where $\hat{\Phi}_{\A,\B}(\eta)=e^{i\B(\A-\eta)\hat{\phi}\left(\frac{\eta-\A}{\A}\right)}$ and characterized the range and proved that $\mathcal{S}_\phi\left(L^2_{\mathbb{H}}(\mathbb{R})\right)$ is a RKHS with the kernel $K_\phi(\A,\B;\A',\B')=\int_{\mathbb{R}}\Phi_{\A,\B}(\xx)\left(\Phi_{\A',\B'}(\xx)\right)^cd\xx$

Kamel et al.\cite{kamel2019continuous} extended the definition \eqref{R2013_Riba-Wong_MMNP_MultiD-ST_Definition} of the multidimensional ST to the quaternion valued function $f\in L^2_{\mathbb{H}}(\mathbb{R}^2)$ with respect to the window $\phi\in L^1_{\mathbb{H}}(\mathbb{R}^2)\cap L^2_{\mathbb{H}}(\mathbb{R}^2).$ They called the transform the two-dimensional QST, which is defined as 
\begin{align}\label{R2020_Brahim-Tefjini_arXiv_TsQST_Definition}
\left(\mathcal{S}_\phi f\right)(\bldxi,\bldb)=|\det A_{\bldxi}|^{\frac{1}{2}}\int_{\mathbb{R}^2}e^{-i\xx_1\xxi_1}f(\bldx)e^{-j\xx_2\xxi_2}\overline{\phi\left(A_{\bldxi}(\bldx-\bldb)\right)}d\mu_2(\bldx),
\end{align}
where $d\mu_d(\bldx)=\frac{d\bldx}{(2\pi)^\frac{d}{2}},$ $A_{\bldxi}$ is a $2\times 2$ matrix $diag(\xxi_1,\xxi_2),~\xxi_1,\xxi_2\in \mathbb{R}.$ Authors in \cite{kamel2019continuous} obtained some basic properties like linearity in $f$ with respect to the complex scalars, anti-linearity in $\phi$ with respect to the quaternion scalars and also the scaling formula $\left(\mathcal{S}_\phi f\right)(\bldxi,\bldb)=\frac{1}{|\lambda|}\left(\mathcal{S}_\phi f\right)(\frac{\bldxi}{\lambda},\lambda\bldb),~\lambda\in \mathbb{R}_{0}.$ Furthermore, motivated by \cite{riba2013continuous}, as in the case of ST \eqref{R2013_Riba-Wong_MMNP_MultiD-ST_Definition}, author obtained the inversion formula for the transform \eqref{R2020_Brahim-Tefjini_arXiv_TsQST_Definition}, with the assumption that the window function $\phi$ satisfies $\int_{\mathbb{R}}\phi(\bldx)d\mu_2(\bldx)=1,$ as $f=\mathcal{F}_s^{-1}A\mathcal{S}_\phi f,$ where $\mathcal{F}_s$ is given by \eqref{R2019Chen-KouDefinitionTsQFT}. With the assumption that the quaternion window $\psi$ (and also $\phi$) satisfies the conditions
\begin{align}\label{R2020_Brahim-Tefjini_arXiv_TsQST_Condition1}
(\mathcal{F}_s\psi)(\bldomega)e^{-j\B_1\omega_1}=e^{-j\B_2\omega_2}(\mathcal{F}_s\psi)(\bldomega)
\end{align}
or 
\begin{align}\label{R2020_Brahim-Tefjini_arXiv_TsQST_Condition2}
(\mathcal{F}_s(j\psi))(\bldomega)e^{-j\B_1\omega_1}=je^{-j\B_2\omega_2}(\mathcal{F}_s\psi)(\bldomega)
\end{align}
authors also proved the Parseval's formula for the QST \eqref{R2020_Brahim-Tefjini_arXiv_TsQST_Definition} as $\langle \mathcal{S}_\phi f,\mathcal{S}_\psi g \rangle_{2,\mathbb{R}^4}=\langle f C'_{\phi,\psi},g \rangle_{2,\mathbb{R}^2},$ where $\langle F,G\rangle_{2,\mathbb{R}^d}=Sc\left[\int_{\mathbb{R}^d}F(\bldt)(G(\bldt))^cd\mu_d(\bldt)\right],$ $C'_{\phi,\psi}:=\int_{\mathbb{R}^2}\left(\mathcal{F}_s\bar{\phi}\right)(1-\bldxi)\overline{\left(\mathcal{F}_s\bar{\psi}\right)(1-\bldxi)}\frac{d\mu_2(\bldxi)}{|\det A_{\bldxi}|}$ is a non-zero quaternion constant. With the notation $C'_\phi=C'_{\phi,\phi}$ and $C'_\psi=C'_{\psi,\psi}$ associated respectively with 
%$$\int_{\mathbb{R}^2}\int_{\mathbb{R}^2} (\mathcal{S}_\psi f)(\A,\bldb)\overline{(\mathcal{S}_\phi g)(\A,\bldb)}\frac{d\A}{\A^{2d+1}}d\mu_{2d}(\bldb) =\int_{\mathbb{R}^2}f(\bldx) C'_{\psi,\phi}\overline{g(\bldx)} d\mu_2(\bldx)$$
two quaternion window $\phi$ and $\psi$ satisfying \eqref{R2020_Brahim-Tefjini_arXiv_TsQST_Condition1} or \eqref{R2020_Brahim-Tefjini_arXiv_TsQST_Condition2}, it follows that for every $f,g\in L^2_{\mathbb{H}}(\mathbb{R}^2)$ the function $(\bldxi,\bldb)\rightarrow (\mathcal{S}_\psi f)(\bldxi,\bldb)(\mathcal{S}_\phi g)(\bldxi,\bldb)$ belongs to $L^p_{\mathbb{H}}(\mathbb{R}^2\times\mathbb{R}^2),~p\geq 1$ and satisfies the Lieb's inequality 
\begin{align*}
&\bigg(\int_{\mathbb{R}^2}\int_{\mathbb{R}^2}\left |(\mathcal{S}_\phi f)(\bldxi,\bldb)(\mathcal{S}_\psi g)(\bldxi,\bldb)\right|^p d\mu_4(\bldxi,\bldb)\bigg)^{\frac{1}{p}}\\
&\hspace{1cm}\leq \left(C'_\phi C'_{\psi}\right)^{\frac{1}{p}}\left(\int_{\mathbb{R}^2}|f(\bldx)|^2d\mu_2(\bldx)\right)\left(\int_{\mathbb{R}^2}|g(\bldx)|^2d\mu_2(\bldx)\right)\left(\int_{\mathbb{R}^2}|\phi(\bldx)|^2d\mu_2(\bldx)\right)^{1-\frac{1}{p}}\left(\int_{\mathbb{R}^2}|\psi(\bldx)|^2d\mu_2(\bldx)\right)^{1-\frac{1}{p}}.
\end{align*}
Moreover, if $\phi$ be an quaternion window satisfying the assumptions \eqref{R2020_Brahim-Tefjini_arXiv_TsQST_Condition1} or \eqref{R2020_Brahim-Tefjini_arXiv_TsQST_Condition2}, then for every $f\in L^2_{\mathbb{H}}(\mathbb{R}^2)$ the transform $\mathcal{S}_\phi f$ given by \eqref{R2020_Brahim-Tefjini_arXiv_TsQST_Definition} belong to $L^p_{\mathbb{H}}(\mathbb{R}^2\times\mathbb{R}^2)$ for $p\geq 2$ and satisfies
$$\bigg(\int_{\mathbb{R}^2}\int_{\mathbb{R}^2}\left |(\mathcal{S}_\phi f)(\bldxi,\bldb)\right|^p d\mu_4(\bldxi,\bldb)\bigg)^{\frac{1}{p}}\leq (C'_\phi)^{\frac{1}{p}}\left(\int_{\mathbb{R}^2}|f(\bldx)|^2d\mu_2(\bldx)\right)^{\frac{1}{2}}\left(\int_{\mathbb{R}^2}|\phi(\bldx)|^2d\mu_2(\bldx)\right)^{\frac{1}{2}\left(1-\frac{1}{p}\right)}$$
%$$\left(\int_0^\infty\int_{\mathbb{R}^{2d}}|(T_\psi f)(\A,\bldb)|^p\frac{d\A}{\A^{2d+1}}d\mu_{2d}(\bldb)\right)^{\frac{1}{p}}
%\leq (C'_\psi)^{\frac{1}{p}}\left(\int_{\mathbb{R}^{2d}}|\psi(\bldx)|^2d\mu_{2d}(\bldx)\right)^{1-\frac{2}{p}}\left(\int_{\mathbb{R}^{2d}}|f(\bldx)|^2d\mu_{2d}(\bldx)\right)^{\frac{1}{2}},$$
as a consequence of which it follows that,  if $\mathcal{S}_\psi f,$ for $f\neq 0,$ is $\epsilon-$concentrated on a measurable set $\Omega\subset\mathbb{R}^2\times\mathbb{R}^2$ of non-zero finite measure $\mu_{\mathbb{R}^4}(\Omega),$ where $\mu_4$ is a measure defined on $\mathbb{R}^2\times\mathbb{R}^2$ given by $d\mu_4(\bldxi,\bldb)=d\mu_2(\bldxi)d\mu_2(\bldb),$ then for every $p>2,$ $\mu_{\mathbb{R}^4}(\Omega)\int_{\mathbb{R}^2}|\phi(\bldx)|^2d\mu_2(\bldx)\geq C'_\psi(1-\epsilon^2)^{\frac{p}{p-2}}.$
%$$\mu_G(\Omega)\|\psi\|^2_{2,2d}\geq C'_\psi(1-\epsilon^2)^{\frac{p}{p-2}},~p>2$$
This is the Lieb's UP associated with the QST \eqref{R2020_Brahim-Tefjini_arXiv_TsQST_Definition}.\\

\noindent\textbf{Uncertainty Principles:} Using the UP for the QFT \cite{akila2016multidimensional}, authors in \cite{akila2016quaternionic} obtained the Heisenberg's UP for the QST \eqref{R2016_Akila-Roopkumar_1DQST_Definition}
$$\left(\int_{\mathbb{R}^\star}\int_{\mathbb{R}}|\B (\mathcal{S}_\phi f)(\A,\B)|^2d\B d\A\right)\left(\int_{\mathbb{R}}\xxi^2|\hat{f}(\xxi)|^2d\xxi\right)\geq c_\phi\|f\|^2_{L^2_{\mathbb{H}}(\mathbb{R})},$$
for $f\in L^2_{\mathbb{H}}(\mathbb{R}).$ Author in \cite{kamel2019continuous} also obtained the associated Beckner's UP in terms of entropy as 
%$$E\left(|T_\psi f|^2\right)\geq C'_\psi\left(\int_{\mathbb{R}^{2d}}|f(\bldx)|^2d\mu_{2d}(\bldx)\right)\ln\left(\frac{1}{\left(\int_{\mathbb{R}^{2d}}|f(\bldx)|^2d\mu_{2d}(\bldx)\right)\left(\int_{\mathbb{R}^{2d}}|f(\bldx)|^2d\mu_{2d}(\bldx)\right)}\right),$$
$$E\left(|\mathcal{S}_\phi f|^2\right)\geq C'_\phi\|f\|^2_{2,\mathbb{R}^2}\ln\left(\frac{1}{\left(\int_{\mathbb{R}^2}|f(\bldx)|^2d\mu_2(\bldx)\right)\left(\int_{\mathbb{R}^2}|\phi(\bldx)|^2d\mu_2(\bldx)\right)}\right),$$
where $E\left(|\mathcal{S}_\phi f|^2\right)=-\int_{\mathbb{R}^2}\int_{\mathbb{R}^2}\ln\left(|(\mathcal{S}_\phi f)(\bldxi,\bldb)|^2\right)|(\mathcal{S}_\phi f)(\bldxi,\bldb)|^2 d\mu_2(\bldxi)d\mu_2(\bldb)$ followed by the generalized Heisenberg UP for the QST 
\begin{align*}
\left(\int_{\mathbb{R}^2}\int_{\mathbb{R}^2}|\bldomega|^p\left|\left(\mathcal{S}_\phi f\right)(\bldomega,\bldb)\right|^2d\mu_4(\bldxi,\bldb)\right)^{\frac{q}{p+q}}&\left(\int_{\mathbb{R}^2}\int_{\mathbb{R}^2}|\bldb|^q\left|\left(\mathcal{S}_\phi f\right)(\bldomega,\bldb)\right|^2d\mu_4(\bldxi,\bldb)\right)^{\frac{p}{p+q}}\\
&\geq E_{p,q}\left(\int_{\mathbb{R}^2}|f(\bldx)|^2d\mu_2(\bldx)\right)\left(\int_{\mathbb{R}^2}|g(\bldx)|^2d\mu_2(\bldx)\right),
\end{align*}
where $E_{p,q}=\frac{2}{ep}\left(\frac{p}{q}\right)^{\frac{p}{p+q}}\left(\frac{pq}{\Gamma(\frac{2}{p})\Gamma(\frac{2}{q})}\right)^{\frac{pq}{2(p+q)}}$ and $p,q$ are positive real numbers.

\section{Quaternion Windowed Fractional Fourier Transform (QWFrFT)}
Using the idea given in \cite{akila2016multidimensional} for defining the QWFT (see \eqref{R2016_Akila-Roopkumar_QWFT_Definition2-Splited}), Rajkumar \cite{rajakumar2021quaternionic} defined the QWFrFT for the function $f=f_1+jf_2\in L^2_{\mathbb{H}}(\mathbb{R})$ with respect to a non-zero window function $g=g_1+jg_2\in L^2_{\mathbb{H}}(\mathbb{R})$ is as 
\begin{align}\label{R2021_Akila-Roopkumar_IJoACM_QWFrFT_Definition-Splited}
\left(\mathcal{G}_{g,\alpha}f\right)(\oomega,\B)=\left(G_{g_1,\alpha}f_1+G_{\bar{g_2},\alpha}\bar{f_2}+j\left(G_{g_1,\alpha}f_2+G_{\bar{g_2},\alpha}\bar{f_1}\right)\right)(\oomega,\B).
\end{align}
From equation \eqref{R2021_Akila-Roopkumar_IJoACM_QWFrFT_Definition-Splited}, it follows that the transform is consistent with the classical WFrFT \eqref{R2010_Tao-Li-Wang_IEEETSP_STFrFT_Definition} whenever $f_2=g_2=0$ and is $\mathbb{H}-$linear on $L^2_{\mathbb{H}}(\mathbb{R}).$ Using the translation and dilation formula for the WFrFT of complex-valued functions, author obtain the corresponding formula as 
$$\left(\mathcal{G}_{g,\alpha}(\tau_{t_0}f)\right)(\oomega,\B)=\left(\mathcal{G}_{g,\alpha}f\right)(\oomega-t_0\cos\alpha,\B-t_0)e^{i\frac{t_0^2}{2}\sin\alpha\cos\alpha-it_0\oomega\sin\alpha}$$
and
$$\left(\mathcal{G}_{g,\alpha}(\mathcal{M}_{\oomega_0}f)\right)(\oomega,\B)=\left(\mathcal{G}_{g,\alpha}f\right)(\oomega-\oomega_0\sin\alpha,\B)e^{i\frac{\oomega_0^2}{2}\sin\alpha\cosh\alpha+i\oomega_0\oomega\cos\alpha},$$
where $(\mathcal{M}_{\oomega_0}f)(\tee)=e^{-i\oomega_0\tee}f(\tee).$ Based on the relation $\left(\mathcal{G}_{g,\alpha}f\right)(\oomega,\B)=(\mathcal{F}_\alpha\{f(\tau)g(\tau-\tee)\})(\xxi),$ author obtained the Parseval's formula $\left(\mathcal{G}_{g,\alpha}f,\mathcal{G}_{g,\alpha}h\right)=\|g\|^2_{L^2_{\mathbb{H}}(\mathbb{R})}(f,h),$ for $f,h\in L^2_{\mathbb{H}}(\mathbb{R}),$ as a consequence of which it follows that $\|\mathcal{G}_{g,\alpha}f\|_{L^2_{\mathbb{H}}(\mathbb{R}\times\mathbb{R})}=\|f\|_{L^2_{\mathbb{H}}(\mathbb{R})}\|g\|_{L^2_{\mathbb{H}}(\mathbb{R})},$ here the inner product $(\cdot,\cdot)$ is given by \eqref{RQuaternionInnerProduct}. Looking forward for the inversion formula for the transform \eqref{R2021_Akila-Roopkumar_IJoACM_QWFrFT_Definition-Splited} author obtain the $2-D$ inversion formula for $f\in L^2_{\mathbb{H}}(\mathbb{R})$ and $g\neq 0$ as 
$$f(\xx)=\frac{1}{\|g\|^2_{L^2_{\mathbb{H}}(\mathbb{R})}}\int_{\mathbb{R}}\int_{\mathbb{R}}\left(\mathcal{G}_{g,\alpha}f\right)(\oomega,\B)K_{-\alpha}(\xx,\oomega)[g(\xx-\B)]^cd\oomega d\B,~a.e.$$
Motivated by the idea of the $1-D$ inversion formula for the classical WFT \cite{tao2009short}, they also obtained the $1-D$  inversion formula for the QWFT for the function $f\in L^2_{\mathbb{H}}(\mathbb{R})$ with assumption that $g=g_1+jg_2\in L^2_{\mathbb{H}}(\mathbb{R})$ is continuous in a neighborhood of $0$ and $g(0)\neq 0$ as
$$f(\B)=\left[\int_{\mathbb{R}}\left(\mathcal{G}_{g,\alpha}f\right)(\oomega,\B)K_{-\alpha}(\B,\oomega)d\oomega\right]\frac{1}{g_1(0)+j\bar{g_2}(0)},~a.e.$$
It is to be noted that although the inversion formula for the classical case \eqref{R2009_Tao_STFT-1D-Inversion-Formula} in \cite{tao2009short} is similar, but it lacks rigor, because as $g\in L^2_{\mathbb{C}}(\mathbb{R})$ its value at the single point cannot be determined. Thus to overcome this issue Rajkumar \cite{rajakumar2021quaternionic} while giving the inversion in $1-D$ case for the transform \eqref{R2021_Akila-Roopkumar_IJoACM_QWFrFT_Definition-Splited} considered the function $g\in L^2_{\mathbb{H}}(\mathbb{R})$ in such a way that $g$ is continuous in a neighborhood of $0$  and $g(0)\neq 0.$\\

\noindent\textbf{Uncertainty Principles:} Based on  the UP for the QFrFT $\mathcal{F}_s$ given by 
$$\|\xxi\left(\mathcal{F}_sf\right)(\xxi)\|_{L^2_{\mathbb{H}}(\mathbb{R})}\|\oomega\left(\mathcal{F}_sf\right)(\oomega)\|_{L^2_{\mathbb{H}}(\mathbb{R})}\geq \frac{1}{2\sqrt{2}}|\sin(\alpha-\beta)|\|f\|^2_{L^2_{\mathbb{H}}(\mathbb{R}^2)},$$ 
author also obtained the Heisenberg's UP for the QWFT \eqref{R2021_Akila-Roopkumar_IJoACM_QWFrFT_Definition-Splited} as 
$$\|\oomega\left(\mathcal{G}_{g,\alpha}f\right)(\oomega,\B)\|_{L^2_{\mathbb{H}}(\mathbb{R}\times\mathbb{R})}\|\tee f(\tee)\|_{L^2_{\mathbb{H}}(\mathbb{R})}\geq \frac{1}{2\sqrt{2}}|\sin\alpha|\|f\|^2_{L^2_{\mathbb{H}}(\mathbb{R})}\|g\|_{L^2_{\mathbb{H}}(\mathbb{R})}.$$

\section{Quaternion Fractional Wavelet Transform (QFrWT)}
A function $\psi\in L^1_{\mathbb{H}}(\mathbb{R})\cap L^2_{\mathbb{H}}(\mathbb{R})$ is said to be a quaternion admissible wavelet if it satisfies the condition that $0<C_{\psi,\alpha}:=\int_{\mathbb{R}}\frac{|(\mathcal{F}_\alpha\psi)(\xxi)|^2}{|\xxi|}d\xxi<\infty.$ Roopkumar \cite{roopkumar2018quaternionic} introduced a QWT  for the function $f\in L^2_{\mathbb{H}}(\mathbb{R})$ with respect to the quaternion admissible wavelet $\psi$ as
\begin{align}\label{R2018_Roopkumar_JoA_QFrWT_Definition}
(\mathcal{W}_{\psi,\alpha}f)(\A,\B)=(f\otimes_\alpha [U_\alpha\psi_\A])(\B),~\A,\B\in\mathbb{R}
\end{align}
where $\otimes$ is the convolution given by \eqref{R2016RoopkumarConvolution1DQFrFT}, $\psi_\A(\tee)=\frac{1}{\sqrt{|\A|}}\psi\left(-\frac{\tee}{\A}\right)$ and $U_\alpha$ is defined by $U_\alpha (h)=\mathcal{F}_{-2\alpha}\overline{h_1}-jh_2,~h=h_1+jh_2.$
Using the machinery of the QFrFT \eqref{R2016Roopkumar_1DFrFT_Definition}, author proved that the $\mathbb{H}-$linear transform \eqref{R2018_Roopkumar_JoA_QFrWT_Definition} satisfies the Parseval's identity $(\mathcal{W}_{\psi,\alpha}f,\mathcal{W}_{\psi,\alpha}g)=2\pi C_{\psi,\alpha}(f,g)$ for $f,g\in L^2_{\mathbb{H}}(\mathbb{R})$ followed by the inversion formula $f(\B)=\frac{1}{2\pi C_{\psi,\alpha}}\int_{\mathbb{R}}\left((\mathcal{W}_{\psi,\alpha}f)(\A,\cdot)\otimes_\alpha\psi_\A\right)(\B)\frac{d\A}{\A^2},~a.e.$
%\cite{gupta2023quaternionic},
%$$C_{\psi,\theta}=\int_{\mathbb{R}}\frac{|(\mathfrak{F}_\theta \psi)(\xxi)|^2}{|\xxi|}d\xxi$$
With respect to $\psi\in L^1_{\mathbb{H}}(\mathbb{R})\cap L^2_{\mathbb{H}}(\mathbb{R})$ satisfying the admissibility condition i.e., $0<C_{\psi,\theta}:=\int_{\mathbb{R}}\frac{|(\mathfrak{F}_\theta \psi)(\xxi)|^2}{|\xxi|}d\xxi<\infty,$ where $\mathfrak{F}_\theta$ denotes the QFrFT of type II given by equation \eqref{R2023_BookChapter_Gupta-Verma_Cattani_1DQFrFT_Definition}, Gupta et al.\cite{gupta2023quaternionic} introduced the QFrWT of type II of $f\in L^2_{\mathbb{H}}(\mathbb{R})$ as
\begin{align}\label{R2023_Gupta-Verma-Cattani-BookChapter_QFrWT_Definition}
(W^\theta_\psi f)(\A,\B)=\left(f\otimes(U(\psi_{\A,\theta}))\right)(\B),~\A,\B\in \mathbb{R},
\end{align}
where $\psi_{\A,\theta}(\tee)=\frac{1}{|\A|^{\frac{1}{2\theta}}}\psi\left(\frac{\tee}{(\sgn \A)|\A|^{\frac{1}{\theta}}}\right),$ $\otimes$ denotes the convolution given by equation \eqref{R2014AkilaConvolution1DQFT} and $U(f)=\check{\overline{f_1}}-jf_2,$ for $f=f_1+jf_2.$ Moreover, if $\psi=\psi_1+j\psi_2,$ then the definition \eqref{R2023_Gupta-Verma-Cattani-BookChapter_QFrWT_Definition} can be equivalently written as
$$W^\theta_\psi f=W^\theta_{\psi_1} f_1+W^\theta_{\check{\overline{\psi_2}}} \check{\overline{f_2}}+j\left(W^\theta_{\psi_1} f_2-W^\theta_{\check{\overline{\psi_2}}} \check{\overline{f_1}}\right).$$
Thus, as a consequence of which the translation and scaling formula of the transform is given respectively by
$$(W^\theta_\psi (\tau_\yy f))(\A,\B)=\left(W^\theta_\psi (f_1+jf_2)\right)(\B-\yy,\A)-j\left(W^\theta_{\check{\overline{\psi_2}}}\left(\check{\overline{f_1}}+j\check{\overline{f_1}}\right)\right)(\B+\yy,\A)$$
and $(W^\theta_\psi f_\lambda)(\A,\B)=(W^\theta_\psi f)(\lambda^\theta\A,\lambda\B),~\lambda>0,~f_\lambda(\tee)=\sqrt{\lambda}f(\lambda\tee).$ Apart from this, authors in \cite{gupta2023quaternionic} obtained the inner product relation, using the  convolution theorem and the Parseval's formula of the FrFT \eqref{R2023_BookChapter_Gupta-Verma_Cattani_1DQFrFT_Definition}, which is given by 
$$\int_{\mathbb{R}}\int_{\mathbb{R}}(W^\theta_\psi f)(\A,\B)\left((W^\theta_\psi g)(\A,\B)\right)^c\frac{d\A d\B}{|\A|^{\frac{1}{\theta}+1}}=C_{\psi,\theta}\int_{\mathbb{R}}f(\tee)(g(\tee))^cd\tee$$
followed by the inversion formula $f(\tee)=\frac{1}{C_{\psi,\theta}}\int_{\mathbb{R}}\left((\mathcal{W}_{\psi,\alpha}f)(\A,\cdot)\otimes_\alpha\psi_{\A,\theta}\right)(\B)\frac{d\A}{\A^{\frac{1}{\theta}+1}}$ in $L^2_{\mathbb{H}}(\mathbb{R}).$ With the assumption that the function $\psi$ is complex valued, authors also proved that the range $W^\theta_\psi\left(L^2_{\mathbb{H}}(\mathbb{R})\right)$ is a RKHS with the kernel $K_{\psi,\theta}(\A_0,\B_0;\A,\B)=\frac{1}{C_{\psi,\theta}}\int_{\mathbb{R}}\psi_{\A,\B,\theta}(\tee)\overline{\psi_{\A_0,\B_0,\theta}(\tee)}d\tee,$ where $\psi_{\A,\B,\theta}(\tee)=\frac{1}{|\A|^{\frac{1}{2\theta}}}\psi\left(\frac{\tee-\B}{(\sgn\A)|\A|^{\frac{1}{\theta}}}\right).$ Using the  UP \eqref{R2023_Gupta-Verma-Cattani_BookChapter_QFrFT_UP}, author obtained the Heisenberg's type UP for the  QWT \eqref{R2023_Gupta-Verma-Cattani-BookChapter_QFrWT_Definition} as
$$\left(\int_{\mathbb{R}}\int_{\mathbb{R}}|b|^2|(W^{\theta}_{\psi}f)(b,a)|^2\frac{dbda}{|a|^{\frac{1}{\theta}+1}}\right)\left(\int_{\mathbb{R}}|\xi|^{\frac{3}{\theta}-1}|(\mathfrak{F}_{\theta}f)(\xi)|^2d\xi\right)\geq\frac{\pi\theta}{4}C_{\psi,\theta}\|f\|^4_{L^2_{\mathbb{H}}(\mathbb{R})}$$
for $f=f_1+jf_2,~\psi=\psi_1+j\psi_2$ are such that $f_l,\psi_l,~l=1,2$ are in the Schwartz class.

%\section{Quaternion Fractional Stockwell Transform}
%\textcolor{red}{(No, paper reported till date in this specific topic. However, a more general paper can be found that will be discussed later. For special choice of the parameters we get results for this specific section.)}

\section{Quaternion Windowed Linear Canonical transform (QWLCT)}
Due to the non-commutativity of the multiplication of quaternions, there are three different types of QWLCT depending on the position of the kernel of the LCT. Gao et al. \cite{gao2020quaternion} introduced the two-sided QWLCT of the function $f\in L^2_{\mathbb{H}}(\mathbb{R}^2)$ with respect to the window function $\phi\in L^2_{\mathbb{H}}(\mathbb{R}^2)\setminus\{0\}$ is defined by 
\begin{align}\label{R2020_Gao-Li_AACA_2sQWLCT_Defnition}
\left(\mathcal{G}^{ij,\phi}_{M_1,M_2}f\right)(\bldomega,\bldu)=\int_{\mathbb{R}^2}K^i_{M_1}(\xx_1,\oomega_1)f(\bldx)\overline{\phi(\bldx-\bldu)}K^j_{M_2}(\xx_2,\oomega_2)d\bldx
\end{align} 
where the kernels $K^i_{M_1}(\xx_1,\oomega_1)$ and $K^j_{M_1}(\xx_2,\oomega_2)$ are given respectively by \eqref{R2016Kou-Ou-MoraisDefinitionTsQLCT-KernelI} and \eqref{R2016Kou-Ou-MoraisDefinitionTsQLCT-KernelII} with
%\begin{align}
%K^i_{M_1}(\xx_1,\oomega_1)=
%\begin{cases}
%\frac{1}{\sqrt{i2\pi \Mb_1}}e^{i(\frac{\Ma_1}{2\Mb_1}\xx_1^2-\frac{1}{\Mb_1}\xx_1\oomega_1+\frac{\Md_1}{2\Mb_1}\oomega_1^2-\frac{\pi}{4})},&~\mbox{for}~\Mb_1\neq 0 \\
%\sqrt{\Md_1}e^{i\left(\frac{\Mc_1\Md_1}{2}\right)\oomega_1^2}\delta(\xx_1-\Md_1\oomega_1),&~\mbox{for}~\Mb_1=0,
%\end{cases}
%\end{align}
%\begin{align}
%K^j_{M_2}(\xx_2,\oomega_2)=
%\begin{cases}
%\frac{1}{\sqrt{i2\pi \Mb_2}}e^{i(\frac{\Ma_2}{2\Mb_2}\xx_2^2-\frac{1}{\Mb_2}\xx_2\oomega_2+\frac{\Md_2}{2\Mb_2}\oomega_2^2-\frac{\pi}{4})},&~\mbox{for}~\Mb_2\neq 0 \\
%\sqrt{\Md_2}e^{i\left(\frac{\Mc_2\Md_2}{2}\right)\oomega_2^2}\delta(\xx_2-\Md_2\oomega_2),&~\mbox{for}~\Mb_2=0,
%\end{cases}
%\end{align}
$M_l=(\Ma_l,\Mb_l,\Mc_l,\Md_l)\in\mathbb{R}^{2\times 2}$ are matrix parameters such that $det(M_l)=1,$ for $l=1,2$ and $f\in L^1_{\mathbb{H}}(\mathbb{R}^2).$ The transform \eqref{R2020_Gao-Li_AACA_2sQWLCT_Defnition} is $\mathbb{R}-$linear and is such that $\left|\left(\mathcal{G}^{ij,\phi}_{M_1,M_2}f\right)(\bldomega,\bldu)\right|^2\leq \frac{1}{2\pi\sqrt{|\Mb_1\Mb_2|}}\|f\|_{L^2_{\mathbb{H}}(\mathbb{R}^2)}\|\phi\|_{L^2_{\mathbb{H}}(\mathbb{R}^2)}.$ Apart from the formula of the transform \eqref{R2020_Gao-Li_AACA_2sQWLCT_Defnition}, for the  translation and  modulation of a function, given by 
\begin{align*}
\left(\mathcal{G}^{ij,\phi}_{M_1,M_2}\tau_{\boldsymbol r}f\right)(\bldomega,\bldu)=e^{ir_1\oomega_1\Mc_1}e^{-i\frac{\Ma_1 r_1^2}{2}\Mc_1}\left(\mathcal{G}^{ij,\phi}_{M_1,M_2}f\right)(\boldsymbol m,\boldsymbol n)e^{jr_2\oomega_2\Mc_2}e^{-j\frac{\Ma_2 r_2^2}{2}\Mc_2},
\end{align*}
where $\tau_{\boldsymbol r}f(\bldx)=f(\bldx-\boldsymbol r),~\boldsymbol r=(r_1,r_2),~\boldsymbol m=(m_1,m_2),~\boldsymbol n=(n_1,n_2)\in \mathbb{R}^1, m_l=\oomega_l-\Ma_l r_l,n_l=\uu_l-r_l,~l=1,2$ and 
\begin{align*}
\left(\mathcal{G}^{ij,\phi}_{M_1,M_2}\mathbb{M}_{\boldsymbol s} f\right)(\bldomega,\bldu)=e^{i\oomega_1 s_1\Md_1}e^{-i\frac{\Mb_1\Md_1 s_1^2}{2}}\left(\mathcal{G}^{ij,\phi}_{M_1,M_2}f\right)(\boldsymbol v,\bldu)e^{j\oomega_2 s_2\Md_2}e^{-j\frac{\Mb_2\Md_2 s_2^2}{2}},
\end{align*}
where $\mathbb{M}_{\boldsymbol s} f(\bldx)=e^{i\xx_1 s_1}f(\bldx)e^{j\xx_2 s_2},$ $\boldsymbol s=(s_1,s_2),$ $\boldsymbol v=(v_1,v_2)\in \mathbb{R}^2,$ $v_l=\oomega_l-s_l\Mb_l,~l=1,2$ authors in \cite{gao2020quaternion}, obtianed the inversion formula for the QWLCT \eqref{R2020_Gao-Li_AACA_2sQWLCT_Defnition} as
$$f(\bldx)=\frac{1}{\|\phi\|^2_{L^2_{\mathbb{H}}(\mathbb{R})}}\int_{\mathbb{R}^2}\int_{\mathbb{R}^2}K^{-i}_{M_1}(\xx_1,\uu_1)\left(\mathcal{G}^{ij,\phi}_{M_1,M_2}\mathbb{M}_{\boldsymbol s} f\right)(\bldomega,\bldu)K^{-j}_{M_2}(\xx_2,\oomega_2)\phi(\bldx-\bldu)d\bldomega d\bldu.$$
For the window functions $\phi,\psi$ and $f,g\in L^2_{\mathbb{H}}(\mathbb{R}^2),$ then the transform \eqref{R2020_Gao-Li_AACA_2sQWLCT_Defnition}, satisfies the orthogonality relation $\langle \mathcal{G}^{ij,\phi}_{M_1,M_2}f, \mathcal{G}^{ij,\psi}_{M_1,M_2}g\rangle=Sc[(f,g)(\phi,\psi)].$ In particular, if $f=g$ and $\phi=\psi$ then $\int_{\mathbb{R}^2}\int_{\mathbb{R}^2}\left|\left( \mathcal{G}^{ij,\phi}_{M_1,M_2}f, \mathcal{G}\right)(\bldomega,\bldu)\right|^2d\bldomega d\bldu=\|f\|^2_{L^2_{\mathbb{H}}(\mathbb{R}^2)}\|\phi\|^2_{L^2_{\mathbb{H}}(\mathbb{R}^2)},$ i.e., $\|\mathcal{G}^{ij,\phi}_{M_1,M_2}f\|^2_{L^2_{\mathbb{H}}(\mathbb{R}^2\times\mathbb{R}^2)}=\|f\|^2_{L^2_{\mathbb{H}}(\mathbb{R}^2)}\|\phi\|^2_{L^2_{\mathbb{H}}(\mathbb{R}^2)}.$ Thus for  a given window function $\phi,$ the transform $\mathcal{G}^{ij,\phi}_{M_1,M_2}$ maps $L^2_{\mathbb{H}}(\mathbb{R}^2)$ into $L^2_{\mathbb{H}}(\mathbb{R}^2\times\mathbb{R}^2).$ If $\mathcal{G}^{ij,\phi}_{M_1,M_2}f=0,$ then $\|f\|_{L^2_{\mathbb{H}}(\mathbb{R}^2)},$ so $f=0.$ This shows that the transform is one-one. Prasad et al. \cite{prasad2023uncertainty} characterized the range $\mathcal{G}^{ij,\phi}_{M_1,M_2}(L^2_{\mathbb{H}}(\mathbb{R}^2))$ of the transform \eqref{R2020_Gao-Li_AACA_2sQWLCT_Defnition}. Using the Hausdorff-Young inequality \eqref{R2022Kundu-PrasadSHYITsQLCT}, authors in \cite{prasad2023uncertainty}, \cite{gao2022uncertainty} proved that for a  non-zero function $f\in L^2_{\mathbb{H}}(\mathbb{R}^2),$ the transform \eqref{R2020_Gao-Li_AACA_2sQWLCT_Defnition} satisfies the Lieb's inequality
\begin{align}\label{R2023_Prasad-Kundu_Optik_2sQWLCT_Lieb-Inequality}
\|\mathcal{G}^{ij,\phi}_{M_1,M_2}f\|_{L^p_{\mathbb{H}}(\mathbb{R}^2\times\mathbb{R}^2)}\leq \mathcal{K}^{M_1,M_2}_{p,q}\|f\|_{L^2_{\mathbb{H}}(\mathbb{R}^2)}\|\phi\|_{L^2_{\mathbb{H}}(\mathbb{R}^2)},~p\geq 2
\end{align}
where $\mathcal{K}^{M_1,M_2}_{p,q}=|\Mb_1 \Mb_2|^{\frac{1}{2}-\frac{1}{p}}A_l^{\frac{4}{q}}A_{m'}^{\frac{2}{q}}(2\pi)^{\frac{1}{q}-\frac{1}{p}},~ A_p=\left(\frac{p^{\frac{1}{p}}}{q^{\frac{1}{q}}}\right)^{\frac{1}{2}}, ~l=\frac{p}{q},~m=\frac{2}{q}$ and $\frac{1}{m}+\frac{1}{m'}=1,$
%$D_{s,s'}=\left(\frac{4}{s}\right)^{\frac{1}{s}}\left(\frac{4}{s'}\right)^{\frac{1}{s'}}$ and $\frac{1}{s}+\frac{1}{s'}=1.$
%$$\frac{|\boldsymbol B|^{\frac{1}{s}-\frac{1}{2}}}{2\pi}D_{s,s'}$$
using which authors showed that, for a non-zero function $f\in L^2_{\mathbb{H}}(\mathbb{R}^2)$ and $\epsilon \geq 0,$ if $\mathcal{G}^{ij,\phi}_{M_1,M_2}f$ is $\epsilon-$concentrated on a measurable subset $\Omega$ of $\mathbb{R}^2\times\mathbb{R}^2,$ then $(1-\epsilon^2)|\Omega|^{\frac{pq}{p-q}}\leq \left(\mathcal{K}^{M_1,M_2}_{p,q}\right)^2$ for every $p>2.$ Using the Pitt's inequality for the QLCT, authors in \cite{prasad2023uncertainty}, obtained the inequality for the QWLCT $\int_{\mathbb{R}^2}\int_{\mathbb{R}^2}\left|\left( \mathcal{G}^{ij,\phi}_{M_1,M_2}f\right)(\bldomega,\bldu)\right|^2d\bldomega d\bldu\leq \frac{M_\lambda}{|\Mb_1 \Mb_2|^\lambda}\left\||\bldx|^{\frac{\lambda}{2}}f(\bldx)\right\|^2_{L^2_{\mathbb{H}}(\mathbb{R}^2)}\|\phi\|^2_{L^2_{\mathbb{H}}(\mathbb{R}^2)},$
where $M_\lambda=\left(\frac{\Gamma\left(\frac{1-\lambda}{4}\right)}{\Gamma\left(\frac{1+\lambda}{4}\right)}\right)^2,$ $f\in \mathcal{S}(\mathbb{R}^2,\mathbb{H}),~0\leq\lambda<2.$\\

\noindent\textbf{Uncertainty Principles:} Apart from this author also gave the local Price's inequality for the transform as
$$\left\|\chi_{\Omega}\mathcal{G}^{ij,\phi}_{M_1,M_2}f\right\|_{L^p_{\mathbb{H}}(\mathbb{R}^2\times\mathbb{R}^2)}\leq \frac{|\Omega|^{\frac{1}{p(p+1)}}r^{\frac{2-\epsilon}{p+1}}}{2^2(2-\epsilon)^{\frac{1}{2p+2}}}\left\||(\bldomega,\bldu)|^\epsilon\mathcal{G}^{ij,\phi}_{M_1,M_2}f\right\|^{\frac{1}{p+1}}_{L^2_{\mathbb{H}}(\mathbb{R}^2\times\mathbb{R}^2)}+|\Omega|^{\frac{1}{p(p+1)}}\left\||(\bldomega,\bldu)|^\epsilon\mathcal{G}^{ij,\phi}_{M_1,M_2}f\right\|^{\frac{2}{p+1}}_{L^2_{\mathbb{H}}(\mathbb{R}^2\times\mathbb{R}^2)}r^{-\frac{2\epsilon}{p+1}},$$ 
where $r>0,$ $0<\epsilon<2$ and $p>1.$
%\textbf{Heisenberg's UP}
%
%$$\left(\int_{\mathbb{R}^2}\int_{\mathbb{R}^2}\oomega^2_l\left|\left( \mathcal{G}^{ij,\phi}_{M_1,M_2}f\right)(\bldomega,\bldu)\right|^2d\bldomega d\bldu\right)^{\frac{1}{2}}\left(\int_{\mathbb{R}^2}\xx^2_l|f(\bldx)|^2d\bldx\right)^{\frac{1}{2}}\geq \frac{\Mb_l}{2}\|f\|^2_{L^2_{\mathbb{H}}(\mathbb{R}^2)}\|\phi\|^2_{L^2_{\mathbb{H}}(\mathbb{R}^2)},~l=1,2$$
Gao et al. \cite{gao2022uncertainty} gave the logarithmic UP for the QWLCT \eqref{R2020_Gao-Li_AACA_2sQWLCT_Defnition}. Using the Jensen's inequality and the UP, they also obtained its Heisenberg-Weyl's UP.
%%%%%%%%%%%%%%%%%%%%%%%%%%%%%%%%%%%%%%%%%%%%%%%%%%%%%%%%%%%%%%%%%%%%%%%%%%%%%%%%%%%%%%%%%%%%%%%%%%%
%%%%%%%%%%%%%%%%%%%%%%%%%%%%%%%%%%%%%%%%%%%%%%%%%%%%%%%%%%%%%%%%%%%%%%%%%%%%%%%%%%%%%%%%%%%%%%%%%%%
%%%%%%%%%%%%%%%%%%%%%%%%%%%%%%%%%%%%%%%%%%%%%%%%%%%%%%%%%%%%%%%%%%%%%%%%%%%%%%%%%%%%%%%%%%%%%%%%%%%
%%%%%%%%%%%%%%%%%%%%%%%%%%%%%%%%%%%%%%%%%%%%%%%%%%%%%%%%%%%%%%%%%%%%%%%%%%%%%%%%%%%%%%%%%%%%%%%%%%%
%%%starting with that of the QLCT 
%%%$$\frac{\|\phi\|^2_{L^2_{\mathbb{H}}(\mathbb{R}^2)}}{4\pi^2}\int_{\mathbb{R}^2}\ln|\bldx|\left|f(\bldx)\right|^2d\bldx+\int_{\mathbb{R}^2}\int_{\mathbb{R}^2}\ln|\oomega|\left|\left(\mathcal{G}^{ij,\phi}_{M_1,M_2}f\right)(\bldomega,\bldu)\right|^2d\bldomega d\bldu\geq \textcolor{purple}{\frac{D+\ln|\boldsymbol B|}{4\pi^2}}\|f\|^2_{L^2_{\mathbb{H}}(\mathbb{R}^2)}\|\phi\|^2_{L^2_{\mathbb{H}}(\mathbb{R}^2)},$$
%%%$D=\phi(\frac{1}{2})-\ln(\pi),~\phi(\tee)=\frac{\Gamma'(\tee)}{\Gamma(\tee)}.$ \textcolor{red}{(what is this $\boldsymbol B$, need to check, it appears something wrong here)} 
%Using the Jensen's inequality, the logarithmic UP for the QWLCT \eqref{R2020_Gao-Li_AACA_2sQWLCT_Defnition} implies the Heisenberg-Weyl's UP for the transform. 
%%%%%%%%%%%%%%%%%%%%%%%%%%%%%%%%%%%%%%%%%%%%%%%%%%%%%%%%%%%%%%%%%%%%%%%%%%%%%%%%%%%%%%%%%%%%%%%%%%%
%%%%%%%%%%%%%%%%%%%%%%%%%%%%%%%%%%%%%%%%%%%%%%%%%%%%%%%%%%%%%%%%%%%%%%%%%%%%%%%%%%%%%%%%%%%%%%%%%%%
%%%%%%%%%%%%%%%%%%%%%%%%%%%%%%%%%%%%%%%%%%%%%%%%%%%%%%%%%%%%%%%%%%%%%%%%%%%%%%%%%%%%%%%%%%%%%%%%%%%
%%%%%%%%%%%%%%%%%%%%%%%%%%%%%%%%%%%%%%%%%%%%%%%%%%%%%%%%%%%%%%%%%%%%%%%%%%%%%%%%%%%%%%%%%%%%%%%%%%%
With the assumption that $|\Mb_1 \Mb_2|\geq \frac{1}{4\pi^2},$ they also obtained the Entropy UP 
\begin{align}
E\left(\left|\mathcal{G}^{ij,\phi}_{M_1,M_2}f\right|^2\right)
\geq\left[\frac{1}{2\pi^2}\left(\ln(2)-\ln(2\pi)-\ln(|\Mb_1 \Mb_2|)\right)-\ln\left(\|f\|^2_{L^2_{\mathbb{H}}(\mathbb{R}^2)}\|\phi\|^2_{L^2_{\mathbb{H}}(\mathbb{R}^2)}\right)\right]\|f\|^2_{L^2_{\mathbb{H}}(\mathbb{R}^2)}\|\phi\|^2_{L^2_{\mathbb{H}}(\mathbb{R}^2)},
\end{align}
for $f\in L^2_{\mathbb{H}}(\mathbb{R}^2).$

\section{Quaternion Linear Canonical Wavelet Transform (QLCWT)}

A function $\psi\in L^2_{\mathbb{H}}(\mathbb{R}^2)$ is said to be an admissible quaternion linear canonical wavelet if it satisfies the admissibility condition $0<c_\psi:=\int_{SO(2)}\int_{\mathbb{R}^+}\left|\left(\mathcal{L}_{M_1,M_2}\psi^{M_1,M_2}_{\A,\theta,\bldb}(\bldx)\right)(\bldomega)\right|^2\frac{d\A d\theta}{\A}<\infty$ a.e. $\bldomega\in \mathbb{R}^2.$ Shah et al. \cite{shah2021linear}, defined the unitary linear operator $U_{\A,\theta,\bldb}:L^2_{\mathbb{H}}(\mathbb{R}^2)\rightarrow L^2_{\mathbb{H}}(\mathcal{G})$ as $U_{\A,\theta,\bldb}(\psi)=\psi^{M_1,M_2}_{\A,\theta,\bldb}(\bldx)=\frac{1}{\A}e^{-i\frac{\Ma_1}{2\Mb_1}(\xx_1^2-b_1^2)}\psi\left(r_{-\theta}\left(\frac{\bldx-\bldb}{\A}\right)\right)e^{-j\frac{\Ma_2}{2\Mb_2}(\xx_2^2-b_2^2)},~\A\in\mathbb{R}^+,\bldb\in\mathbb{R}^2$ and $r_{-\theta}\in SO(2)$ and called the family $\psi^{M_1,M_2}_{\A,\theta,\bldb}$ as the quaternion daughter linear canonical wavelets. They defined the QLCWT of $f\in L^2_{\mathbb{H}}(\mathbb{R}^2)$ with respect to the quaternion admissible linear canonical wavelet $\psi$ as 
\begin{align}\label{R2021_Shah-Teali-Tantary_AACA_QLCWT_Definition}
(\mathcal{W}^{M_1,M_2}_\psi f)(\A,\theta,\bldb)=(f,\psi_{\A,\theta,\bldb}^{M_1,M_2})=\int_{\mathbb{R}^2}f(\bldx)\frac{1}{\A}\overline{e^{-i\frac{\Ma_1}{2\Mb_1}(\xx_1^2-b_1^2)}\psi\left(r_{-\theta}\left(\frac{\bldx-\bldb}{\A}\right)\right)e^{-j\frac{\Ma_2}{2\Mb_2}(\xx_2^2-b_2^2)}}d\bldx.
\end{align} 
The transform is $\mathbb{H}-$linear and anti-linear in $f$ and $\psi$ respectively and have the following translation and scaling formula
$$(\mathcal{W}^{M_1,M_2}_\psi f(\bldx-\boldsymbol k))(\A,\theta,\bldb)=e^{i\frac{\Ma_1 k_1^2}{\Mb_1}}\left(\mathcal{W}^{M_1,M_2}_\psi \left(e^{i\frac{\Ma_1 k_1\xx_1}{\Mb_1}}f(\bldx)e^{j\frac{\Ma_2 k_2\xx_2}{\Mb_2}}\right)\right)(\A,\theta,\bldb-\boldsymbol k)e^{j\frac{\Ma_2 k_2^2}{\Mb_2}}$$ 
and 
$$(\mathcal{W}^{M_1,M_2}_\psi f)(\A,\theta,\bldb)=\frac{1}{\lambda}\left(\mathcal{W}^{M_1',M_2'}_\psi f(\lambda\bldx)\right)(\A\lambda,\theta,\lambda\bldb),~\lambda>0$$
where $M_l'=(\Ma_l',\Mb_l',\Mc_l',\Md_l)$ such that $\frac{\Ma_l'}{\Mb_l'}=\frac{\Ma_l}{\Mb_l\lambda^2},~l=1,2.$
Apart from these basic properties authors in \cite{shah2021linear} obtained the inner product relation 
\begin{align}
\int_{SO(2)}\int_{\mathbb{R}^+}\left(\int_{\mathbb{R}^2}(\mathcal{W}^{M_1,M_2}_\psi f)(\A,\theta,\bldb)\overline{(\mathcal{W}^{M_1,M_2}_\psi g)(\A,\theta,\bldb)}d\bldb\right)\frac{d\A d\theta}{\A^3}=c_\psi(f,g)
\end{align}
for $f,g\in L^2_{\mathbb{H}}(\mathbb{R}^2).$ In particular if $f=g,$ then $\|\mathcal{W}^{M_1,M_2}_\psi f\|^2_{L^2_{\mathbb{H}}(\mathcal{G})}=c_\psi\|f\|^2_{L^2_{\mathbb{H}}(\mathbb{R}^2)},$ as a result of which it follows that, expect for the factor $c_\psi,$ the transform $\mathcal{W}^{M_1,M_2}_\psi$ is an isometry from $L^2_{\mathbb{H}}(\mathbb{R}^2)$ to $L^2_{\mathbb{H}}(\mathcal{G}).$ The reconstruction formula for the transform \eqref{R2021_Shah-Teali-Tantary_AACA_QLCWT_Definition} is given by $f(\bldx)=\frac{1}{c_\psi}\int_{\mathcal{G}}(\mathcal{W}^{M_1,M_2}_\psi f)(\A,\theta,\bldb)\psi^{M_1,M_2}_{\A,\theta,\bldb}(\bldx)d\lambda(\A,\theta,\bldb),$ a.e. Moreover, the range of the transform is a RKHS with the kernel $K_{\psi}(\A,\theta,\bldb;\A',\theta',\bldb')=\frac{1}{c_\psi}(\psi^{M_1,M_2}_{\A,\theta,\bldb},\psi^{M_1,M_2}_{\A',\theta',\bldb'}).$\\

\noindent\textbf{Uncertainty Principles:} Using the Heisenberg's and logarithmic UP for the QFT, the authors obtained the Heisenberg's UP for the  transform \eqref{R2021_Shah-Teali-Tantary_AACA_QLCWT_Definition} for $f\in L^2_{\mathbb{H}}(\mathbb{R}^2)$ as
$$\left(\int_{\mathcal{G}}\left|b^2_l (\mathcal{W}^{M_1,M_2}_\psi f)(\A,\theta,\bldb)\right|^2d\bldb\right)\left(\int_{\mathbb{R}^2}\left|\omega_l(\mathcal{F}_qf)(\bldomega)\right|d\bldomega\right)\geq \frac{\Mb_l^2}{4}\sqrt{c_\psi}\|f\|^4_{L^2_{\mathbb{H}}(\mathbb{R}^2)},~l=1,2$$
and also the logarithmic UP for the QWLCT \eqref{R2021_Shah-Teali-Tantary_AACA_QLCWT_Definition} 
$$\int_{\mathcal{G}}\ln |\bldb|\left| (\mathcal{W}^{M_1,M_2}_\psi f)(\A,\theta,\bldb)\right|^2d\bldb\frac{d\A d\theta}{\A^3}+c_\psi\int_{\mathbb{R}^2}\ln|\bldomega|\left|(\mathcal{F}_qf)(\bldomega)\right|^2d\bldomega\geq c_\psi\left(D+\ln\sqrt{\Mb_1^2+\Mb_2^2}\right)\|f\|^2_{L^2_{\mathbb{H}}(\mathbb{R}^2)}$$
for $f\in \mathcal{S}(\mathbb{R}^2,\mathbb{H})$ followed by its local UP.
%\textbf{Local UP}
%$$\|f\|^2_{L^2_{\mathbb{H}}(\mathbb{R}^2)}\leq \frac{1}{\sqrt{1-|E|}}\|\mathcal{W}^{M_1,M_2}_\psi f\|_{L^2(E^c,\mathbb{H})}$$
%Moreover, for every $\alpha>0,$ there exist $C(\alpha)>0,$ such that 
%$$\|f\|^2_{L^2_{\mathbb{H}}(\mathbb{R}^2)}\leq C(\alpha)\left\{\int_{\mathbb{R}^2}\int_{\mathbb{R}^2}\right\}$$
\section{Quaternion Linear Canonical Stockwell Transform (QLCST)}
Motivated by the definition of the ST \eqref{R1996_Stockwell-Mansinha-Lowe_Definition_STransform} given by Stockwell et al. in \cite{stockwell1996localization}. Bhat and Dar \cite{younus2023quaternion} introduced a LCST for a $\mathbb{H}-$valued functions. For $f\in L^2_{\mathbb{H}}(\mathbb{R}^2)$ and $\psi\in L^1_{\mathbb{H}}(\mathbb{R}^2)\cap L^2_{\mathbb{H}}(\mathbb{R}^2)$ is a non-zero windowed function, the QLCST is defined as 
\begin{align}\label{R2023_Bhat-Dhar_IJWMIP_QLCST_Definition}
\left(\mathcal{S}^{ij,\psi}_{M_1,M_2}f\right)(\bldomega,\bldu)=\int_{\mathbb{R}^2}K^i_{M_1}(\xx_1,\oomega_1)f(\bldx)\overline{\psi(\bldx-\bldu,\bldomega)}K^j_{M_2}(\xx_2,\oomega_2)d\bldx
\end{align} 
where the kernels $K^i_{M_1}(\xx_1,\oomega_1)$ and $K^j_{M_1}(\xx_2,\oomega_2)$ are given respectively by \eqref{R2016Kou-Ou-MoraisDefinitionTsQLCT-KernelI} and \eqref{R2016Kou-Ou-MoraisDefinitionTsQLCT-KernelII} with
$M_l=(\Ma_l,\Mb_l,\Mc_l,\Md_l)\in\mathbb{R}^{4}$ are parameters such that $\Ma_l\Md_l-\Mc_l\Md_l=1,$ for $l=1,2$ and $\int_{\mathbb{R}^2}\psi(\bldx,\bldomega)=1,~\forall~\bldomega\in\mathbb{R}^2.$ The proposed QLCST gives rise to a novel  QFrST. The authors obtained their fundamental properties, like linearity, parity, translation, and modulation, followed by the orthogonality relation and the reconstruction formula. They also obtained the Heisenberg and logarithmic UPs associated with the QLCST \eqref{R2023_Bhat-Dhar_IJWMIP_QLCST_Definition}, respectively, as 
\begin{align*}
\left(\int_{\mathbb{R}^2}\int_{\mathbb{R}^2}\oomega_l^2\left|\left(\mathcal{S}^{ij,\psi}_{M_1,M_2}f\right)(\bldomega,\bldu)\right|^2d\bldomega d\bldu\right)^{\frac{1}{2}}\left(\int_{\mathbb{R}^2}\xx_l^2|f(\bldx)|^2d\bldx\right)^{\frac{1}{2}}\geq \frac{\Mb_l C_\psi}{2}\|f\|^2_{L^2_{\mathbb{H}}(\mathbb{R}^2)},~l=1,2.
\end{align*}
and
\begin{align*}
\int_{\mathbb{R}^2}\int_{\mathbb{R}^2}\ln|\bldomega|\left|\left(\mathcal{S}^{ij,\psi}_{M_1,M_2}f\right)(\bldomega,\bldu)\right|^2d\bldomega d\bldu+C_\psi\int_{\mathbb{R}^2}\ln|\bldx||f(\bldx)|^2d\bldx\geq \left(\frac{\Gamma'(\frac{1}{2})}{\Gamma(\frac{1}{2})}-\ln 2\right)C_\psi\|f\|^2_{L^2_{\mathbb{H}}(\mathbb{R}^2)},
\end{align*}
where $C_\psi$ is some constant depending on $\psi.$

\section{Quaternion Quadratic Phase Windowed Fourier Transform (QQPWFT)}

Let $\wedge_l=(\LmA_l,\LmB_l,\LmC_l,\LmD_l,\LmE_l),\LmA_l,\LmB_l,\LmC_l,\LmD_l,\LmE_l\in\mathbb{R} ~\mbox{and}~B_l\neq 0~\mbox{for}~l=1,2$. The short time quaternion quadratic phase Fourier transform (STQQPFT) of a function $f\in L^2_\mathbb{H}(\mathbb{R}^2)$ with respect to a quaternion window function (QWF) $g\in L^2_\mathbb{H}(\mathbb{R}^2)\cap L^\infty_\mathbb{H}(\mathbb{R}^2)$ is defined by 
\begin{align}\label{R2022_Gupta-Verma_arXiv_QQPFT_Definition}
\left(\mathcal{S}^{\wedge_1,\wedge_2}_{\mathbb{H},g}f\right)(\bldx,\bldxi)=\int_{\mathbb{R}^2}\mathcal{K}^i_{\wedge_1}(\tee_1,\xxi_1)f(\bldt)\overline{g(\bldt-\bldx)}\mathcal{K}^j_{\wedge_2}(\tee_2,\xxi_2)d\bldt,~(\bldx,\bldxi)\in\mathbb{R}^2\times\mathbb{R}^2,
\end{align}
where $\mathcal{K}^i_{\wedge_1}(\tee_1,\xxi_1)$ and $\mathcal{K}^j_{\wedge_2}(\tee_2,\xxi_2)$ are given by equations \eqref{RGupta_VermaKernelofQQPFTinI} and \eqref{RGupta_VermaKernelofQQPFTinJ}, respectively. The transform \eqref{R2022_Gupta-Verma_arXiv_QQPFT_Definition} is linear in $f$ with respect to the scalars in $\{\xx+i\yy:\xx,\yy\in\mathbb{R}\}$ and anti-linear in $\psi$ with respect to the scalars  $\{\xx+j\yy:\xx,\yy\in\mathbb{R}\}.$ It is bounded and satisfies $\left\|\mathcal{S}^{\wedge_1,\wedge_2}_{\mathbb{H},g}f\right\|_{L^\infty_\mathbb{H}(\mathbb{R}^2)}\leq\frac{1}{2\pi}\|g\|_{L^2_\mathbb{H}(\mathbb{R}^2)}\|f\|_{L^2_\mathbb{H}(\mathbb{R}^2)}$ followed by the translation  formula 
$$\left(\mathcal{S}^{\wedge_1,\wedge_2}_{\mathbb{H},g}(\tau_{\boldsymbol k}f)\right)(\bldx,\bldxi)=\phi^i_{\wedge_1}(k_1,\xxi_1)\left(\mathcal{S}^{\wedge_1,\wedge_2}_{\mathbb{H},g}f\right)(\bldx-\boldsymbol k,\bldxi'_{\bldx})\phi^j_{\wedge_2}(k_2,\xxi_2),$$
where $(\tau_{\boldsymbol k}f)(\bldt)=f(\bldt-\boldsymbol k),$ $\bldxi'_{\bldx}=\left(\xxi_1+\frac{2\LmA_1\xx_1}{\LmB_1},\xxi_2+\frac{2\LmA_2\xx_2}{\LmB_2}\right),$  
$\phi^i_{\wedge_1}(k_1,\xxi_1)=e^{-i\left(\LmA_1k_1^2+\LmD_1k_1+\LmB_1k_1\xxi_1-\frac{4\LmA_1^2\LmC_1k_1^2}{\LmB_1^2}-\frac{4\LmA_1\LmC_1k_1\xxi_1}{\LmB_1}-\frac{2\LmA_1k_1}{\LmB_1}\right)}$
and 
$\phi^j_{\wedge_2}(k_2,\xxi_2)=e^{-j\left(\LmA_2k_2^2+\LmD_2k_2+\LmB_2k_2\xxi_2-\frac{4\LmA_2^2\LmC_2k_2^2}{\LmB_2^2}-\frac{4\LmA_2\LmC_2k_2\xxi_2}{\LmB_2}-\frac{2\LmA_2k_2}{\LmB_2}\right)}.$
Along with the scaling formula $\left(\mathcal{S}^{\wedge_1,\wedge_2}_{\mathbb{H},g_\lambda}f_{\lambda}\right)(\bldx,\bldxi)=\left(\mathcal{S}^{\wedge_1',\wedge_2'}_{\mathbb{H},g}f\right)\left(\frac{1}{\lambda}\bldx,\bldxi\right),$ where $(f_{\lambda})(\bldt)=\frac{1}{\lambda} f\left(\frac{1}{\lambda}\bldt\right),$ $\wedge_l'=\left(\lambda^2\LmA_l,\lambda \LmB_l,\LmC_l,\lambda \LmD_l,\LmE_l\right),~l=1,2.$ Authors in \cite{gupta2022short}, obtained the Lieb's inequality for the transform \eqref{R2022_Gupta-Verma_arXiv_QQPFT_Definition} given by 
\begin{align}\label{R2022_Gupta-Verma_arXiv_QQPFT_Lien-Inequality}
\left\|\mathcal{S}^{\wedge_1,\wedge_2}_{\mathbb{H},g}f \right\|_{L^q_\mathbb{H}(\mathbb{R}^2\times\mathbb{R}^2)}\leq\frac{(2\pi)^{\frac{1}{q}-\frac{1}{p}}}{|\LmB_1\LmB_2|^{\frac{1}{q}}}\left(\frac{2}{q}\right)^\frac{2}{q}\|g\|_{L^2_\mathbb{H}(\mathbb{R}^2)}\|f\|_{L^2_\mathbb{H}(\mathbb{R}^2)}, \mbox{for}~f\in L^2_\mathbb{H}(\mathbb{R}^2).
\end{align}
Apart from this, authors also proved that for quaternion window functions $g_1,g_2$  and $f_1,f_2\in L^2_\mathbb{H}(\mathbb{R}^2),$ then $\mathcal{S}^{\wedge_1,\wedge_2}_{\mathbb{H},g_1}f_1,$ $\mathcal{S}^{\wedge_1,\wedge_2}_{\mathbb{H},g_2}f_2\in L^2_\mathbb{H}(\mathbb{R}^2\times\mathbb{R}^2)$ and satisfies the following inner product relation
$\left\langle\mathcal{S}^{\wedge_1,\wedge_2}_{\mathbb{H},g_1}f_1,\mathcal{S}^{\wedge_1,\wedge_2}_{\mathbb{H},g_2}f_2\right\rangle=\frac{1}{|\LmB_1\LmB_2|}\langle f_1(\overline{g_1},\overline{g_2}),f_2\rangle,$ followed by the inversion formula 
$$f(\bldt)=\frac{|\LmB_1\LmB_2|}{\|g\|^2_{L^2_\mathbb{H}(\mathbb{R}^2)}}\int_{\mathbb{R}^2}\int_{\mathbb{R}^2}\overline{\mathcal{K}^i_{\wedge_1}(\tee_1,\xxi_1)}\left(\mathcal{S}^{\wedge_1,\wedge_2}_{\mathbb{H},g}f\right)(\bldx,\bldxi)\overline{\mathcal{K}^j_{\wedge_2}(\tee_2,\xxi_2)}g(\bldt-\bldx)d\bldx d\bldxi.$$
Bhat et al. \cite{bhat2022uncertainty} considered the two-sided quaternion QPFT \eqref{R2022_Gupta-Verma_arXiv_QQPFT_Definition} where the factor $\frac{1}{\sqrt{2\pi}}$ in the kernels $\mathcal{K}^i_{\wedge_1}(\tee_1,\xxi_1)$ and $\mathcal{K}^j_{\wedge_2}(\tee_2,\xxi_2)$ are replaced respectively with $\sqrt{\frac{\LmB_1i}{2\pi}}$ and $\sqrt{\frac{\LmB_2j}{2\pi}}.$ Let us denote the transform in this case by $^{(1)}\mathcal{S}^{\wedge_1,\wedge_2}_{\mathbb{H},g}.$ They also obtained the basic properties along with the inversion formula and Plancherel identity for the transform $^{(1)}\mathcal{S}^{\wedge_1,\wedge_2}_{\mathbb{H},g}$ as in \cite{gupta2022short}. With the aid of the Pitt's inequality for the two-sided QFT and the relation of the QQPWFT with the two-sided QFT, they also obtained the Pitt's inequality  
$$\int_{\mathbb{R}^2}\int_{\mathbb{R}^2}|\bldxi|^{-\alpha}\left|\left(^{(1)}\mathcal{S}^{\wedge_1,\wedge_2}_{\mathbb{H},g}f\right)(\bldx,\bldxi)\right|^2d\bldx d\bldxi\leq \frac{|\boldsymbol{B}|^\alpha}{4\pi^2}\Delta_\alpha\|g\|^2_{L^2_{\mathbb{H}}(\mathbb{R}^2)}\int_{\mathbb{R}^2}|\bldx|^\alpha|f(\bldx)|^2d\bldx,$$
$\Delta_\alpha=\pi^\alpha\left[\Gamma\left(\frac{2-\alpha}{4}\right)\Gamma\left(\frac{2+\alpha}{4}\right)\right],$ $0\leq\alpha\leq 2$ and $f\in \mathcal{S}(\mathbb{R}^2,\mathbb{H}).$\\

\noindent\textbf{Uncertainty Principles:} Using the logarithmic UP for the QQPFT $^{(1)}\mathcal{Q}^{\wedge_1,\wedge_2}_{\mathbb{H}}$ given by equations \eqref{R2023_Bhat-Dar_MMAS_QQPFT_Logarithmic_UP}, authors in \cite{bhat2022uncertainty} obtained the logarithmic UP for $^{(1)}\mathcal{S}^{\wedge_1,\wedge_2}_{\mathbb{H},g}$ as 
$$\frac{\|g\|^2_{L^2_{\mathbb{H}}(\mathbb{R}^2)}}{(2\pi)^2}\int_{\mathbb{R}^2}\ln|\bldx||f(\bldx)|^2d\bldx+\int_{\mathbb{R}^2}\int_{\mathbb{R}^2}\ln|\bldxi|\left|\left(^{(1)}\mathcal{S}^{\wedge_1,\wedge_2}_{\mathbb{H},g}f\right)(\bldx,\bldxi)\right|^2d\bldx d\bldxi\geq \frac{1}{(2\pi)^2}(D-\ln|\boldsymbol{B}|)\|g\|^2_{L^2_{\mathbb{H}}(\mathbb{R}^2)}\|f\|^2_{L^2_{\mathbb{H}}(\mathbb{R}^2)}$$
for $f\in \mathcal{S}(\mathbb{R}^2,\mathbb{H}).$ Also with the assumption $\|g\|_{L^2_\mathbb{H}(\mathbb{R}^2)}=\|f\|_{L^2_\mathbb{H}(\mathbb{R}^2)}=1$ and $\int\int_{\Omega}\left|\left(^{(1)}\mathcal{S}^{\wedge_1,\wedge_2}_{\mathbb{H},g}f\right)(\bldx,\bldxi)\right|^2d\bldx d\bldxi\geq 1-\epsilon,$ for $\Omega$ being a measurable subset of $\mathbb{R}^2\times\mathbb{R}^2$ and $\epsilon\geq 0$, they proved that the measure $|\Omega|$ of $\Omega$ is bounded below by $\frac{2\pi(1-\epsilon)}{\sqrt{|\LmB_1\LmB_2|}}.$ In addition, if $0<|\Omega|<1,$ then the transform $^{(1)}\mathcal{S}^{\wedge_1,\wedge_2}_{\mathbb{H},g}$ satisfies the following local uncertainty inequality $1\leq \frac{1}{\sqrt{1-|\Omega|}}\|^{(1)}\mathcal{S}^{\wedge_1,\wedge_2}_{\mathbb{H},g}f\|_{L^2_{\mathbb{H}}(\Omega^c)}.$ As a consequence of the Lieb's inequality \eqref{R2022_Gupta-Verma_arXiv_QQPFT_Lien-Inequality}, author in \cite{gupta2022short} proved the Lieb's UP, i.e., if $\mathcal{S}^{\wedge_1,\wedge_2}_{\mathbb{H},g}f,$ is $\epsilon-$concentrated on a measurable subset $\Omega$ of $\mathbb{R}^2\times\mathbb{R}^2,$ then for every $q>2,$
$|\Omega|\geq \frac{(2\pi)^2}{|\LmB_1\LmB_2|}(1-\epsilon^2)^{\frac{q}{q-2}}\left(\frac{q}{2}\right)^{\frac{4}{q-2}}.$ With the assumption that $\|g\|_{L^2_\mathbb{H}(\mathbb{R}^2)}\|f\|_{L^2_\mathbb{H}(\mathbb{R}^2)}=1,$ authors also proved that the QQPWFT \eqref{R2022_Gupta-Verma_arXiv_QQPFT_Definition}  satisfies 
$$\mathcal{E}_{S}(f,g,\wedge_1,\wedge_2)\geq\frac{2}{|\LmB_1\LmB_2|},$$
where $\displaystyle\mathcal{E}_{S}(f,g,\wedge_1,\wedge_2)=-\int_{\mathbb{R}^2}\int_{\mathbb{R}^2}\left|\left(\mathcal{S}_{\mathbb{H},g}^{\wedge_1,\wedge_2}f\right)(\bldx,\bldxi)\right|^2\log\left(\left|\left(\mathcal{S}_{\mathbb{H},g}^{\wedge_1,\wedge_2}f\right)(\bldx,\bldxi)\right|^2\right)d\bldx d\bldxi.$

%\section{Quaternion QPWT}
%star symbol using MnSymbol package
\section{Conclusion}
This review is a comprehensive survey of a large class of Fourier-type integral transforms. Looking at the recent interests of researchers, it becomes extremely important for young researchers to have a collective study in this field. Not only the Fourier transform to FrFT, LCT, or QPFT, but also the popular joint time-frequency analyzing tools such as windowed Fourier transform, wavelet transforms, and the Stockwell transforms are generalized in the fractional domain, linear canonical domain, and the quadratic phase domain. This survey paper contains the mathematical properties of all these transforms and their variants. Despite many applications for some transforms, this survey mainly focuses on their mathematical properties. Section-wise presentation for each transform makes it convenient for readers to grasp the direction and find research gaps in the relevant area. We have tried our level best to include almost all the variants of each transform reported in the literature. The pictorial representation given at the beginning serves as the best tool for understanding the direction of the generalization of the transform. Despite several properties and inequalities, like sharp Hausdorff Youngs inequality, Lieb's inequality, Pitt's inequality, etc., of the transformation, this review article also presents several qualitative as well as the quantitative UPs like Heisenberg, local, logarithmic, R\'{e}yni (in particular, Shannon), Lieb's UP, etc. Apart from just giving the result obtained in the literature, we have also tried to include the ideas of proofs in short.

Quaternion is the generalization of complex numbers. Recently, researchers have been actively involved in defining and exploring new integral transforms for quaternion-valued functions. Generalization of the transform for functions of complex-valued to quaternion-valued has fascinated researchers for the last two decades. This paper also discusses a non-trivial generalization of the convolutions in the Fourier domain, fractional Fourier domain, linear canonical domain, and quadratic phase Fourier domain to quaternion-valued functions defined on $\mathbb{R}.$ This survey paper also discusses the ideas of defining some convolution-based wavelet, and Stockwell transforms in quaternion Fourier, fractional Fourier, linear canonical, and quadratic phase Fourier domains. For the function of quaternion-valued defined on $\mathbb{R}^2$, we observe that several transforms can be defined depending upon the positions of the kernels. This article can serve better for the signal processing engineers in selecting an appropriate tool for the quaternion value signal analysis.

While going through the literature, we noticed that the Fourier analysis and its offshoots have also been explored for the function of Clifford algebra valued, a generalization of the quaternion algebra. Scholars researching the Fourier-type integral transform and its variants for Clifford algebra-valued functions will find this survey to be of great use.  
%%%%%%%%%%%%%%%%%%%%%%%%%%%%%%%%%%%%%%%%%%%%%%%%%%%%%%%%%%%%%%%%%%%%%%%%%%%%%%%%%%%%%
%%%%%%%%%%%%%%%%%%%%%%%%%%%%%%%%%%%%%%%%%%%%%%%%%%%%%%%%%%%%%%%%%%%%%%%%%%%%%%%%%%%%%
%%%%%%%%%%%%%%%%%%%%%%%%%%%%%%%%%%%%%%%%%%%%%%%%%%%%%%%%%%%%%%%%%%%%%%%%%%%%%%%%%%%%%
%%%%%%%%%%%%%%%%%%%%%%%%%%%%%%%%%%%%%%%%%%%%%%%%%%%%%%%%%%%%%%%%%%%%%%%%%%%%%%%%%%%%%
%%%%%%%%%%%%%%%%%%%%%%%%%%%%%%%%%%%%%%%%%%%%%%%%%%%%%%%%%%%%%%%%%%%%%%%%%%%%%%%%%%%%%
%\section{Some variants of stars, that can be used throughout}
%   $\pentagram,\;\bigotimes,\;\bigoast$\\[6pt]
%   % Star symbol in text mode using bbding package
%   \FiveStar, \FiveStarCenterOpen, \FiveStarConvex \\[6pt]
%   \EightStar, \EightStarTaper, \SixStar \\[6pt]
%   %Star sumbol using pifont package(use \ding{65 to 107})
%   \ding{65}, \ding{80}, \ding{81}, \ding{107}
%%%%%%%%%%%%%%%%%%%%%%%%%%%%%%%%%%%%%%%%%%%%%%%%%%%%%%%%%%%%%%%%%%%%%%%%%%%%%%%%%%%%%
%%%%%%%%%%%%%%%%%%%%%%%%%%%%%%%%%%%%%%%%%%%%%%%%%%%%%%%%%%%%%%%%%%%%%%%%%%%%%%%%%%%%%
%%%%%%%%%%%%%%%%%%%%%%%%%%%%%%%%%%%%%%%%%%%%%%%%%%%%%%%%%%%%%%%%%%%%%%%%%%%%%%%%%%%%%
%%%%%%%%%%%%%%%%%%%%%%%%%%%%%%%%%%%%%%%%%%%%%%%%%%%%%%%%%%%%%%%%%%%%%%%%%%%%%%%%%%%%%
%%%%%%%%%%%%%%%%%%%%%%%%%%%%%%%%%%%%%%%%%%%%%%%%%%%%%%%%%%%%%%%%%%%%%%%%%%%%%%%%%%%%%

\bibliography{Master_Bivek_Intro}
\bibliographystyle{plain}
\end{document}